\documentclass[letterpaper,11pt]{amsart}
\usepackage{amssymb}
\usepackage{amsthm}
\usepackage{amsmath}
\usepackage{amsfonts}
\usepackage{stmaryrd}
\usepackage{xcolor}
\usepackage{xy}
\xyoption{all}
\usepackage{mathrsfs}
\usepackage[utf8]{inputenc}
\usepackage[margin=1.0in]{geometry}
\usepackage{pst-node}
\usepackage{multirow}
\usepackage{tikz-cd} 
\usepackage{hyperref}
\usepackage[T1]{fontenc}
\usepackage{appendix}
\usepackage[english]{babel}
\usepackage[backend=biber,style=alphabetic,maxalphanames=5]{biblatex}
\usepackage{csquotes}
\usepackage{comment}
\newtheorem{theorem}{Theorem}[section]
\newtheorem{lemma}[theorem]{Lemma}
\newtheorem{proposition}[theorem]{Proposition}
\newtheorem{question}[theorem]{Question}

\newtheorem{corollary}[theorem]{Corollary}
\newtheorem{conjecture}[theorem]{Conjecture}
 
\newtheorem{condition/definition}[theorem]{Condition/Definition} 
\newtheorem{definition/proposition}[theorem]{Definition/Proposition} 

\newtheorem{assumption}[theorem]{Assumption}
\newtheorem{setup}[theorem]{Setup}

\theoremstyle{definition}
\newtheorem{definition}[theorem]{Definition}
\theoremstyle{remark}
\newtheorem{remark}[theorem]{Remark}
\theoremstyle{definition}

\newtheorem{construction}[theorem]{Construction}

\def\ul{\underline}

\def\ol{\overline}

\def\SL{\mathrm{SL}}
\def\lis{\mathrm{lis}}

\def\mf{\mathfrak}
\def\ra{\rightarrow}
\def\la{\leftarrow}
\def\lim{\mathop{\rm lim}\nolimits}
\def\colim{\mathop{\rm colim}\nolimits}
\def\Spec{\mathop{\rm Spec}}
\def\Spa{\mathop{\rm Spa}}
\def\Spd{\mathop{\rm Spd}}
\def\Jbul{\underline{J_{b}(\bb{Q}_{p})}}

\def\Spf{\operatorname{Spf}}
\def\Hom{\mathop{\rm Hom}\nolimits}

\def\Sh{\mathop{\rm Sh}\nolimits}

\def\Bun{\mathrm{Bun}}
\def\Div{\mathrm{Div}}
\def\IC{\mathrm{IC}}
\def\Hck{\mathrm{Hck}}

\def\pre{\mathrm{pre}}
\def\Perf{\mathrm{Perf}}
\def\oc{\mathrm{o.c}}

\def\dim{\mathrm{dim}}
\def\Mant{\mathrm{Mant}}

\def\Sht{\mathrm{Sht}}
\def\Rep{\mathrm{Rep}}
\def\D{\mathrm{D}}
\def\Dconstf{\mathrm{D}_{\mathrm{cons,tf}}}
\def\Dbc{\mathrm{D}^\mathrm{b}_{\mathrm{c}}}
\def\PerfAdm{\mathrm{PfAdm}}
\def\disc{\mathrm{disc}}

\def\RHom{\mathrm{RHom}}

\def\Corr{\mathrm{Corr}}

\def\Convex{\mathrm{Convex}}
\def\GL{\mathrm{GL}}
\def\Shstar{\mathrm{Sh}^{*}}
\def\Perv{\mathrm{Perv}}
\def\HT{\mathrm{HT}}
\def\ULA{\mathrm{ULA}}
\def\Dlis{\mathrm{D}_{\mathrm{lis}}}
\def\cDlis{\mathcal{D}_{\mathrm{lis}}}

\def\LinCat{\mathrm{LinCat}}
\def\perf{\mathrm{perf}}
\def\alg{\mathrm{alg}}
\def\BB{\mathrm{BB}}
\def\bb{\mathbb}
\def\CO{\mathcal{O}}

\def\et{\mathrm{\acute{e}t}}
\def\Detale{\mathrm{D}_{\mathrm{\acute{e}t}}}
\def\cDetale{\mc{D}_\mathrm{\Acute{e}t}}

\def\pcDetale{\phantom{}^{p}\mc{D}}
\def\pH{\phantom{}^{p}\mathcal{H}}

\def\oc{\mathrm{o.c}}
\def\RHomint{R\mathcal{H}\mathrm{om}}
\def\hs{\mathrm{hs}}
\newcommand{\ParG}{\mathfrak{X}_{\widehat{G}}}
\newcommand{\tor}{\mathrm{tor}}
\newcommand{\PP}{{[\mathsf{P}]}}
\newcommand{\Q}{\mathbb{Q}}

\newcommand{\mc}{\mathcal}
\newcommand{\coloneq}{:=}

\newcommand{\gx}{(\mathsf{G},\mathsf{X})}
\newcommand{\gxno}{\mathsf{G},\mathsf{X}}
\newcommand{\A}{\mathbb{A}}
\newcommand{\Ig}{\mathrm{Ig}}
\newcommand{\Igs}{\mathrm{Igs}}

\newcommand{\qp}{\mathbb{Q}_p}
\newcommand{\Fl}{{\mathcal{F}\ell_{G,\mu^{-1}}}}

\newcommand{\Adm}{\mathrm{Adm}}

\newcommand{\bs}{\mbox{\tiny $\blacksquare$}}

\newcommand{\Dcons}{\mathrm{D}_\mathrm{cons}}

\usepackage{graphicx} 

\title{Intersection Cohomology of Igusa Stacks}

\author{Ana Caraiani}
\address{Department of Mathematics, Imperial College London, Huxley 668,
180 Queen's Gate,
London SW7 2AZ}
\email{a.caraiani@imperial.ac.uk}
\thanks{AC was supported in part by ERC Starting Grant 804176, by a Royal Society University Research Fellowship and by a Leverhulme Prize.}

\author{Linus Hamann}
\address{Department of Mathematics, Harvard University, Science Center Room 239,
1 Oxford Street,
Cambridge MA, 02138, USA}
\email{hamann@math.harvard.edu}

\author{Mingjia Zhang}
\address{Department of Mathematics, Princeton university, Fine Hall, Washington Road,
Princeton, NJ, 08544-1000, USA}
\email{mz9413@princeton.edu}
\thanks{MZ is supported by the NSF through the Institute for Advanced Study, under Grant No. DMS-2424441.}

\begin{document}
\begin{abstract} We study the intersection cohomology of minimally compactified Shimura varieties of PEL type AC using Igusa stacks and the work of Fargues--Scholze. More precisely, we construct a sheaf on the moduli stack of $G$-bundles on the Fargues--Fontaine curve, which recovers this intersection cohomology after applying a Hecke operator in the sense of geometric Langlands. We show that this sheaf has several desirable properties; for example, it is Verdier self-dual and perverse. This leads to several applications to intersection cohomology, including a version of the Mantovan product formula, as well as torsion-vanishing and Eichler--Shimura relations. Along the way, we investigate the interaction between Baily--Borel and Newton stratifications on minimally compactified Igusa stacks, and we study perverse t-structures on stratified $v$-stacks.
\end{abstract}
\maketitle
\tableofcontents
\section{Introduction}

A large amount of literature has been devoted to studying the intersection cohomology of Shimura varieties in the context of the Langlands program. In the more classical setup, the starting point is Borel--Casselman's generalization of Matsushima's formula \cite{BorelCasselman}, which describes a clean decomposition of the $L^2$-cohomology of a noncompact locally symmetric space in terms of spectral data. In contrast to the results for the usual or compactly supported Betti cohomology, only the discrete automorphic representations contribute. See Section 2 of \cite{ArthurClay96} for a beautiful account of this picture. 

In the case when the locally symmetric space arises from a Shimura variety, Zucker's conjecture (proven by Looijenga, Saper--Stern) implies that the $L^2$-cohomology identifies canonically with the intersection cohomology of its Baily--Borel compactification, compatibly with the Hecke action and Hodge--Lefschetz structures, see \cite{Looijenga88}, \cite{SaperStern}, \cite{Looijenga25}. This brings the intersection cohomology into play. When compared with the $L^2$-cohomology, the intersection cohomology has several advantages. For example, intersection cohomology can be defined in the setting of $\ell$-adic sheaves, and consequently it carries a Galois action. Moreover, unlike $L^{2}$-cohomology, it can also be defined in the torsion setting. 

In recent years, a number of new tools have been introduced to study the cohomology of Shimura varieties. 
Igusa stacks, introduced in~\cite{zhang2023, DvHKZ,Kim, DvHKZ2}, are one such tool. They connect Shimura varieties with the categorical local Langlands program over $p$-adic fields, as developed by Fargues--Scholze~\cite{FSGeomLLC} and Zhu~\cite{Zhu}. Until now, novel applications such as the torsion-vanishing results of~\cite{KoTorsionPaper, HamannLeeTorsion, YangZhuTorsion} and the Eichler--Shimura relations of~\cite{KoshikawaES,DvHKZ, DvHKZ2}, have concerned the usual or compactly supported cohomology of Shimura varieties. Our goal in this paper is to provide a framework for also understanding the intersection cohomology of Shimura varieties using techniques from the categorical local Langlands program. 

To achieve this goal, we construct a version of the intersection complex in the setting of Igusa stacks\footnote{For technical reasons, we restrict ourselves to the PEL setting of~\cite{zhang2023} in this paper, but we expect that our geometric constructions will carry over to much more general Shimura varieties. See Remark~\ref{rem: generalizations to other Shimura varieties} for more details.}. We study the resulting relative intersection cohomology complex over $\mathrm{Bun}_G$, the moduli stack of $G$-bundles on the Fargues--Fontaine curve, and show that it has a number of desirable properties. See Theorem~\ref{thm: main} for the precise statement. In particular, we show that this relative intersection cohomology recovers the $\ell$-adic intersection cohomology of the minimal compactification of Shimura varieties after applying a Hecke operator in the sense of geometric Langlands. Furthermore, by computing the stalks of this complex, we prove a version of the Mantovan product formula, cf.~\cite{mantovan-thesis, mantovan-PEL}, for the intersection cohomology of Shimura varieties, see Theorem~\ref{thm: product formula intro}. 

We also show that the relative intersection cohomology of the Igusa stack gives rise to a perverse sheaf for the natural perverse $t$-structure on $\mathrm{Bun}_G$.  As additional applications of this framework, we establish torsion-vanishing (Theorem~\ref{thm: genericLocalization}) and Eichler--Shimura relations (Theorem~\ref{thm: ES intro}) for the intersection cohomology of Shimura varieties.  

Motivated by torsion-vanishing, the recent paper of Koshikawa and Shin \cite{KoshikawaShin} formulates precisely a number of conjectures about the degrees in which the intersection cohomology of Shimura varieties with both rational and torsion coefficients can be concentrated. For example, Koshikawa and Shin conjecture that certain local conditions can determine the range of degrees to which an automorphic representation can contribute to the spectral decomposition of intersection cohomology. There is a similar expectation for modulo $\ell$ systems of Hecke eigenvalues, which points towards a torsion analog of Arthur's conjectures. The main result in our current paper makes it possible to apply results from \cite{FSGeomLLC,Zhu} to study local aspects of the intersection cohomology of Shimura varieties, which, we believe, will eventually be helpful in proving results along these lines.

\subsection{The main theorem}

We let $(\mathsf{G}, \mathsf{X})$ be a Shimura datum with reflex field $\mathsf{E}\subset \overline{\Q}$. For a neat compact open subgroup $K\subset \mathsf{G}(\A_f)$, we denote by $\Sh_K:=\Sh\gx_K$ the attached Shimura variety over $\mathsf{E}$, a smooth quasi-projective variety of dimension $d$. We denote by $\mathrm{Sh}^*_K$ its minimal (Baily--Borel) compactification. Fix a prime $\ell$ and let $\Lambda$ be a finite self-injective ring over $\mathbb{Z}_{\ell}$\footnote{For most of our results, we also allow for coefficients in an extension of $\mathbb{Q}_{\ell}$ or its ring of integers, see Remark~\ref{rem: GeneralCoefficientSystems}, as well as coefficients in a local system obtained from an algebraic representation of $\mathsf{G}$.}. On each $\mathrm{Sh}^*_K$,  we have an intersection complex, denoted $\IC_{\mathrm{Sh}_K^*}$, obtained by taking the middle extension of $\Lambda[d]$ towards the minimal compactification. We are interested in studying the \'etale cohomology groups 
\begin{equation}\label{eq: classical IC}
H^i\left(\mathrm{Sh}^*_{K, \ol{\mathsf{E}}}, 
\IC_{\mathrm{Sh}^*_{K}}\right).
\end{equation} 

Fix a prime $p\not = \ell$ and set $G:=\mathsf{G}_{\Q_p}$ to be the local group at $p$. 
We choose an isomorphism $\iota: \mathbb{C}\simeq \ol{\Q}_p$, which induces a $p$-adic place of $\mathsf{E}$. Let $E$ be the completion of $\mathsf{E}$ at this place and $C$ be its completed algebraic closure. We can analytify the previous construction over $C$. For a neat $K\subset \mathsf{G}(\A_f)$, we write $\mathcal{S}_{K}$ for the diamond over $\Spd C$ attached to $\Sh_K$, and $\mathcal{S}^*_{K}$ for the minimal compactification. 
Fixing the tame level $K^p\subset \mathsf{G}(\A_f^p)$, we also have the inverse limits
\[
\mathcal{S}_{K^p}:= \varprojlim_{K_p\ra \{1\}} \mathcal{S}_{K^pK_p}\ \mathrm{and}\  \mathcal{S}^*_{K^p}:= \varprojlim_{K_p\ra \{1\}} \mathcal{S}^*_{K^pK_p},
\]
where $K_p\subset G(\Q_p)$ runs over compact open subgroups. Over each $\mc{S}^*_{K}$, we have a complex $\IC_{\mathcal{S}_{K}^*}$, obtained by analytifying $\IC_{\mathrm{Sh}_K^*}$ (see \S \ref{sec: analytification}). Pulling these back to $ \mathcal{S}_{K^p}^*$ and passing to the colimit, we get the infinite-level intersection complex, which we denote by $\IC_{\mc{S}_{K^p}^*}$.
The intersection cohomology complex
\[
\IC(\mathsf{G}, \mathsf{X})_{K^p,\Lambda}:= \varinjlim_{K_p\ra \{1\}} R\Gamma\left(\mathcal{S}_{K^pK_p}^*, \IC_{\mathcal{S}^*_{K^pK_p}} \right)
\]
is equipped with a prime-to-$p$ Hecke action, a smooth $G(\Q_p)$-action, and a continuous action by the Weil group $W_E$. Via a series of comparison theorems, its cohomology groups agree with the colimit of the intersection cohomology groups in~\eqref{eq: classical IC} after restricting the Galois action on the latter to $W_E$. 

We aim to reconstruct the complex $\IC(\mathsf{G}, \mathsf{X})_{K^p,\Lambda}$ with its additional structures using Igusa stacks and the action of Hecke operators on $\Bun_{G}$. In order to appeal to the results of~\cite{zhang2023} on the existence of the minimal compactification 
of the Igusa stack, we impose the following assumption throughout 
this paper.

\begin{assumption}\label{assumption:codimension}
    The Shimura datum $\gx$ is of PEL type, and the absolute root system of $\mathsf{G}$ is of type AC. The group $G :=\mathsf{G}_{\qp}$  is unramified. The boundary of the minimal compactification of the attached Shimura variety has codimension at least $2$ in the minimal compactification.
\end{assumption}

\noindent Under this assumption, we have a cartesian diagram of $v$-stacks over $\Spd \ol{\mathbb{F}}_{p}$ (\cite[Theorem 9.40]{zhang2023})
\begin{equation}{\label{eqn: CartesiandiagramforClosedShimuraVariety}}
 \begin{tikzcd}
 \mc{S}_{K^p}^\ast \arrow[r,"\pi_{\mathrm{HT}}"] \arrow[d,"\tilde{h}"] & \mathcal{F}\ell \arrow[d,"h"] & \\ 
\Igs^\ast \arrow[r,"\overline{\pi}_{\mathrm{HT}}"] & \Bun_{G}. &
\end{tikzcd}
\end{equation}
Here $\mathrm{Bun}_G$ is the moduli stack parameterizing $G$-bundles on the Fargues--Fontaine curve, 
$\mathcal{F}\ell$ is the analytification of a flag variety attached to the Hodge cocharacter of $\gx$, 
and $\Igs^\ast=\Igs_{K^p}^\ast$ is the minimally compactified Igusa stack attached to $\gx$ at level $K^p$. The top horizontal arrow is the Hodge--Tate period morphism introduced by Scholze. The right vertical arrow is the Beauville--Laszlo map introduced in \cite{FSGeomLLC}. 

To state our main result, we describe the geometry of the cartesian diagram~\eqref{eqn: CartesiandiagramforClosedShimuraVariety} in a little more detail, following~\cite{FSGeomLLC, zhang2023}. The $v$-stack $\Bun_G$ admits a Newton / Harder--Narasimhan stratification, by locally closed strata 
\[i_b:\Bun^b_G\hookrightarrow \Bun_G\] 
indexed by elements $b$ of the Kottwitz set $B(G)$. This induces a semi-orthogonal decomposition on the category of sheaves $\Detale(\Bun_G,\Lambda)$ by excision with respect to this stratification. Furthermore, there is an equivalence of categories $\D(\Bun_G^b,\Lambda)\simeq \D(J_b(\mathbb{Q}_p),\Lambda)$, where $J_b$ is an inner form of a Levi subgroup of $G$ and $\D(J_b(\mathbb{Q}_p),\Lambda)$ is the derived category of smooth representations of $J_b(\mathbb{Q}_p)$ with $\Lambda$-coefficients. For example, if $b=1$ is the neutral element of $B(G)$, then $J_b$ is $G$.  
 
Both $\ol{\pi}_{\HT}$ and $h$ factor through the open substack $i_{\mu^{-1}}: \Bun_{G,\mu^{-1}} \ra \Bun_{G}$, corresponding to the $\mu$-admissible locus $B(G,\mu^{-1}) \subset B(G)$, where $\mu$ denotes the Hodge cocharacter of the Shimura datum transported under the fixed isomorphism $\ol{\bb{Q}}_{p} \simeq \bb{C}$.  For $b\in B(G,\mu^{-1})$, the pullback of $\Igs^*$ to $\Bun^b_G$ can be identified (up to canonical compactification) with the stack quotient $[\mathrm{Ig}^{b,*, \diamond}/\mathcal{J}_b]$, for a partially minimally compactified perfect Igusa variety $\mathrm{Ig}^{b,*}$ appropriately analytified (see \S \ref{ss: analytificationofperfectschemes}), and a certain unipotent thickening $\mathcal{J}_b$ of $J_b(\mathbb{Q}_p)$ that acts on $\mathrm{Ig}^{b,*}$. 

Our main theorem shows that the intersection cohomology complex $\IC(\mathsf{G}, \mathsf{X})_{K^p}$ interacts with Diagram~\eqref{eqn: CartesiandiagramforClosedShimuraVariety} in the best possible way. 

\begin{theorem}{\label{thm: main}}{(Theorem \ref{thm: perversityofICIgs}, Theorem \ref{thm: Hecke operator torsion})}
Let $\mu$ be the Hodge cocharacter attached to $\gx$. The complex $\mathcal{F}_{\IC,\Lambda}\in \Detale(\Bun_{G,\mu^{-1}}, \Lambda)$, defined as the pushforward under the proper morphism $\ol{\pi}_{\HT}$ of an intersection complex $\IC_{\Igs^*}$ on $\Igs^*$, satisfies the following properties: 
\begin{enumerate}
    \item (Verdier self-duality) It is self-dual with respect to Verdier duality on $\Bun_{G,\mu^{-1}}$. If $\Lambda$ is moreover a regular ring, then $\mathcal{F}_{\IC,\Lambda}$ is universally locally acyclic over $\Spd \ol{\mathbb{F}}_{p}$ in the sense of \cite[Section IV.2]{FSGeomLLC}. 

\medskip 
    
    \item (Perversity) It is a perverse sheaf with respect to the natural perverse $t$-structure on $\Bun_{G,\mu^{-1}}$ (see Definition~\ref{defn: existenceofperversetstructure}).

\medskip 

    \item (Hecke operator) Let $T_{\mu}$ be the Hecke operator attached to $\mu$ (see \S \ref{s: HeckeOperatorsandDlisse}) and $i_1$ the inclusion of the neutral stratum into $\Bun_G$, we have an isomorphism
\begin{equation}\label{eq:geom Hecke operator formula}
i_{1}^{*}T_{\mu}\mathcal{F}_{\mathrm{IC},\Lambda} \simeq 
\IC(\mathsf{G},\mathsf{X})_{K^p,\Lambda}    
\end{equation}
of $G(\Q_p) \times W_{E}$-representations, which is also Hecke-equivariant away from $p$. 

\medskip 

\item (Stalks) There is a $J_b(\mathbb{Q}_p)$- and away-from-$p$ Hecke-equivariant isomorphism 
    \[
    i_b^*\mathcal{F}_{\IC,\Lambda} \simeq V_{\IC,\Lambda,b}[-d_b], 
    \]
    where $V_{\IC,\Lambda,b}$ is the intersection cohomology of the partially minimally compactified Igusa variety $\mathrm{Ig}^{b,*}$ with $\Lambda$-coefficients and $d_b = \dim\mathrm{Ig}^{b,*}$, i.e. the colimit of the intersection cohomologies of the partially minimally compactified finite level Igusa varieties, see Equation~\eqref{eqn: defofIC}.
\end{enumerate}
\end{theorem}

\begin{remark}{\label{rem: GeneralCoefficientSystems}}
In fact, we can also work with \'etale local systems coming from algebraic representations of $\mathsf{G}$ and allow for coefficients in an algebraic extension $F/\mathbb{Q}_{\ell}$ or its ring of integers, see Assumption \ref{assumption: coefficients} and Theorem \ref{thm: Hecke operator}. However, in the case of non-self-injective coefficients, we do not have Verdier self-duality of $\mc{F}_{\IC,\Lambda}$, since, like in the case of schemes, the perverse $t$-structure is not preserved under Verdier duality (See \cite[Section~2.1]{JuteauDecompositionNumbers}). Moreover, we also establish analogous results for various boundary cohomology groups of the minimal compactification, see \S \ref{sec: variants}.
\end{remark}

\begin{remark}\label{rem: generalizations to other Shimura varieties} We establish Theorem~\ref{thm: main} under Assumption~\ref{assumption:codimension} primarily because this is the setting in which the minimal compactification $\Igs^*$ of the Igusa stack has been constructed. We also rely on the identification of the fibers of $\pi_{\mathrm{HT}}$ with partial minimal compactifications of Igusa varieties, a result known in the unramified, PEL type AC setting by~\cite{CS2} and~\cite{santos}. Once these geometric ingredients are in place for more general Shimura varieties, one should be able to construct $\mathcal{F}_{\IC,\Lambda}$ more generally and prove that it has analogous properties, following the outline of the present paper. One key ingredient for both the construction of $\Igs^*$ and the computation of the fibers of $\pi_{\mathrm{HT}}$ is the construction and affineness of the partial minimal compactifications $\mathrm{Ig}^{b,*}$. This has been established recently for Shimura varieties of Hodge type in~\cite{shengkai-mao}. 
\end{remark}

In the rest of this introduction, we first state some applications of Theorem~\ref{thm: main} and then explain the proof ideas for our main geometric results. We conclude with a brief discussion of future directions for the study of $\mathcal{F}_{\mathrm{IC},\Lambda}$.  

\subsection{Applications}

\subsubsection{The Mantovan product formula} As our first application, we establish a version of the Mantovan product formula for intersection cohomology. The first version of the Mantovan product formula was implicit in~\cite{HT} for Harris--Taylor Shimura varieties and established in~\cite{mantovan-thesis, mantovan-PEL} for more general compact Shimura varieties of PEL type. This was subsequently extended to the non-compact case and beyond the PEL setting, for both usual and compactly supported cohomology. See~\cite[Theorem 7.1]{KoTorsionPaper}, ~\cite[Theorem 1.15]{HamannLeeTorsion} and \cite[Theorem 8.5.7]{DvHKZ} for state of the art results. 

Roughly speaking, the cohomological version of the Mantovan product formula gives a $G(\bb{Q}_{p}) \times W_{E}$-equivariant filtration (in the derived sense) on the cohomology of a global Shimura variety, indexed by Newton strata and with graded pieces isomorphic to a product of the cohomology of a local Shimura variety, or Rapoport--Zink space, and an Igusa variety. This product formula has played a fundamental role in computing the cohomology of Shimura varieties over the past few decades, providing a viable alternative to the Langlands--Kottwitz method, often better suited than the latter to understanding places of bad reduction. For example, the product formula was used in landmark results such as the proof of local Langlands for $\GL_n$ over $p$-adic fields in~\cite{HT} and the construction of Galois representations in the conjugate self-dual, regular algebraic setting in~\cite{shin-galois}. 

Regarding the intersection cohomology of Shimura varieties, the work of Morel led to the computation of the alternating sum of these cohomology groups in the Grothendieck group in a number of important cases, for example in~\cite{morel-cohomology}. As she follows the Langlands--Kottwitz method, this approach works primarily at places of good reduction. However, for applications to the conjectures of Arthur and Kottwitz, it is desirable to also understand the places of bad reduction. This is where the Mantovan product formula comes in. 

To state our result precisely, we let $R\Gamma_c(G,b,\mu,\Lambda(d_{b}))$ denote the cohomology with $\Lambda$-coefficients of the local Shimura variety determined by the triple $(G, b, \mu)$, in the sense of~\cite{Scholze-Weinstein} twisted by a suitable character if $b$ is non-basic (see~\eqref{eqn: CohofLocalShimuraVariety}). This is a module over the smooth Hecke algebra with $\Lambda$-coefficients, denoted $\mathcal{H}(J_b)$. The intersection cohomology $V_{\IC,\Lambda,b}$ of the partially minimally compactified Igusa variety is also a module over $\mathcal{H}(J_b)$ and the away from $p$ Hecke algebra. 

\begin{theorem}\label{thm: product formula intro}{(Corollary~\ref{cor: product formula})} Assume in addition that $\Lambda$ is a $\mathbb{Z}_{\ell}[\sqrt{p}]$-algebra. We have an away from $p$ Hecke-equivariant and $G(\mathbb{Q}_p)\times W_E$-equivariant filtration on $\IC(\mathsf{G},\mathsf{X})_{K^p,\Lambda}$ indexed by $b\in B(G,\mu^{-1})$, with graded pieces isomorphic to 
\[
R\Gamma_{c}(G,b,\mu,\Lambda(d_{b}))[d]\otimes |\cdot|^{d/2} \otimes_{\mathcal{H}(J_{b})} V_{\IC,\Lambda,b}[d_{b}],
\]
where $|\cdot|: W_{E} \ra \Lambda^{*}$ denotes the norm character.
\end{theorem}

Theorem~\ref{thm: product formula intro} is deduced by combining parts (3) and (4) of Theorem~\ref{thm: main},  via the standard method used in the proof of~\cite[Corollary 3.22]{HamannLeeTorsion} and~\cite[Theorem 8.5.7]{DvHKZ} for usual and compactly supported cohomology. The idea is to apply excision with respect to the Newton stratification on $\Bun_G$ to $i_1^*T_\mu\mathcal{F}_{\IC,\Lambda}$ and to rewrite the graded pieces using the cohomology with compact support of local Shimura varieties and the computation of the stalks $i_b^*\mathcal{F}_{\IC,\Lambda}$.

\begin{remark}\label{rem:kottwitz conjecture}
We emphasize that the cohomology complex $V_{\IC,\Lambda,b}$ has an explicit description as a colimit of intersection cohomology groups of partial minimal compactifications of Igusa varieties that are schemes of finite type over $\ol{\mathbb{F}}_p$. Therefore, setting $\Lambda=\overline{\mathbb{Q}}_{\ell}$, the alternating sum of the cohomology groups of $V_{\IC,\overline{\mathbb{Q}}_{\ell},b}$ in the Grothendieck group is, at least in theory, computable using the Grothendieck--Lefschetz trace formula. Indeed, one should be able to use Morel's weight truncation~\cite{morel-weight} to describe the intersection complexes on these partial compactifications.

Together with progress (some on-going) on the local Kottwitz conjecture~\cite{HKW} and the Harris--Viehmann conjecture~\cite{HHSII}, which concern the purely local contribution $R\Gamma_c(G,b,\mu,\overline{\mathbb{Q}}_{\ell}(d_{b}))$ to the product formula, Theorem~\ref{thm: product formula intro} suggests a possible route to the global Kottwitz conjecture for the intersection cohomology of Shimura varieties. By this, we mean computing the Galois action on intersection cohomology, at least on the level of the Grothendieck group and then verifying that it agrees with the general recipe given in~\cite{kottwitz-galois}.\footnote{Scholze and Shin introduced another method to compute the Galois action in the case of bad reduction in~\cite{scholze-shin}, based on a generalization of the Langlands--Kottwitz method and building on Scholze's approach to local Langlands for $\GL_n$. This has not been implemented for intersection cohomology at present.} Because intersection cohomology satisfies purity, one should even be able to distinguish individual cohomological degrees. 
\end{remark}

\subsubsection{Torsion-vanishing for intersection cohomology}
For simplicity, set $\Lambda:=\overline{\mathbb{F}}_{\ell}$ in this section. The isomorphism in Equation~\eqref{eq:geom Hecke operator formula} decomposes the intersection cohomology of Shimura varieties into a global part, namely $\mathcal{F}_{\mathrm{IC},\Lambda}$, and a purely local part, given by the Hecke operator $T_\mu$. Together 
with the perversity of $\mathcal{F}_{\mathrm{IC},\Lambda}$ from part (2) of Theorem~\ref{thm: main}, this means that the cohomological amplitude of $\IC\gx_{K^p,\Lambda}$ is constrained by the perverse $t$-exactness properties of the Hecke operator $T_\mu$. Therefore, we can apply purely local methods to study the generic part of the intersection cohomology with $\overline{\mathbb{F}}_{\ell}$-coefficients, as was done in~\cite{HamannLeeTorsion, YangZhuTorsion} for compactly supported and usual cohomology. 

For example, if $K = K^pK_p$ with $K_p\subset G(\Q_p)$ hyperspecial, and $\mf{m} \subset \mathcal{H}_{K_{p}}$ is a maximal ideal of the spherical Hecke algebra with $\ol{\mathbb{F}}_\ell$ coefficients, then, by arguing as \cite[Theorem~5.2]{HamannLeeTorsion}, cf. \cite[Theorem 10.1.6]{DvHKZ}, we may deduce the following.

\begin{theorem}[Corollary \ref{cor: HeckeAlgebraResult}]{\label{thm: genericLocalization}}
For $\mf{m}$ generic in the sense of \cite[Definition~1.1]{HamannLeeTorsion}, 
$\ell$ banal for $G(\Q_p)$, $p > 2$, and $G$ a product of groups in \cite[Table (1)]{HamannLeeTorsion} or of those in \cite[Theorem~8.2.1]{peng2025farguesscholzeparameterstorsionvanishing}, the localization
\[ 
R\Gamma\left(\mathrm{Sh}^*_{K, \ol{\mathsf{E}}}, 
\IC_{\mathrm{Sh}^*_{K}}\right)_{\mathfrak{m}}
\]
is concentrated in degree $0$.
\end{theorem}

\begin{remark}\label{rem: Koshikawa-Shin} In their recent work, Koshikawa and Shin study how certain local conditions can give bounds to the cohomological amplitude of the intersection cohomology of locally symmetric spaces and formulate precise conjectures with rational coefficients, namely~\cite[Conjecture~1.2.1]{KoshikawaShin}, and with torsion coefficients, namely~\cite[Conjecture~6.2.1]{KoshikawaShin}. Theorem~\ref{thm: genericLocalization} provides evidence towards the conjecture for torsion coefficients. We discuss this in more detail in \S \ref{ss: comparisonwithKoshikawaShin}.  
\end{remark}

\begin{remark}{\label{rem: EisensteinvsLSType}} If we localize the formula~\eqref{eq:geom Hecke operator formula} at a more general semi-simple local Langlands parameter $\phi$, using the excursion action of Fargues--Scholze, we can establish a similar result, whenever the localized Hecke operator $i_{1}^{*}T_{\mu,\phi}$ is perverse $t$-exact. This is the content of Theorem~\ref{thm: axiomatizedvanishingstatement}.

We expect that the best possible condition for $\phi$ that guarantees this perverse $t$-exactness is it being of weakly Langlands--Shahidi type, cf.~\cite[Definition~6.2, Conjecture 6.4]{HamannLeeTorsion}. Note that when further localized at a non-Eisenstein maximal ideal in the global prime-to-$p$ Hecke algebra, the intersection cohomology is expected (and known in many cases) to agree with both usual and compactly supported cohomology, see \cite[Remark 1.6]{CS2}. However, there are many Eisenstein maximal ideals whose local constituent at $p$ is generic, and our result is novel in this case.
\end{remark}


\begin{remark} Yang--Zhu~\cite{YangZhuTorsion} perform a similar analysis to~\cite{HamannLeeTorsion} by replacing $\mathrm{Bun}_G$ with the stack $\mathrm{Isoc}_G$ on perfect schemes, which parametrizes isocrystals with $G$-structure. They appeal to the unipotent categorical local Langlands equivalence established in~\cite{Zhu} to deduce the analog of perverse $t$-exactness in their setting, whereas~\cite{HamannLeeTorsion} relied on the more restricted results of~\cite{HamGeomES} on geometric Eisenstein series. This leads Yang--Zhu to prove more general torsion-vanishing results for the compactly supported and usual cohomology of Shimura varieties. In the paper~\cite{GHILZIsocComparison}, the authors compare the \'etale sheaf theories on $\mathrm{Isoc}_G$ and on $\mathrm{Bun}_G$, giving an equivalence of categories
\[
\mathrm{Shv}^{!}(\mathrm{Isoc}_{G},\ol{\mathbb{F}}_{\ell}) \xrightarrow{\substack{\Psi \\ \simeq}} \Detale(\Bun_{G},\ol{\mathbb{F}}_{\ell}).
\] 
Hence, one should be able to combine Theorem~\ref{thm: main} with the result of Yang--Zhu on perverse $t$-exactness of the Hecke operators to deduce torsion-vanishing for intersection cohomology in the generality of~\cite{YangZhuTorsion}.
\end{remark}

\subsubsection{Eichler--Shimura congruence relations}
The work of Eichler--Shimura reveals a congruence relation between the Hecke operators and Frobenius action on the cohomology of the modular curve at a prime of good reduction. Blasius and Rogawski~\cite{Blasius-Rogawski} conjectured that such relations hold true more generally for higher dimensional Shimura varieties. 

More precisely, let the notation be as before, except that we set $\Lambda$ to be an algebraic extension of $\mathbb{Q}_{\ell}$. In particular, $\mu$ is the Hodge cocharacter of $\gx$, viewed as a geometric cocharacter of $G$ via $\iota:\mathbb{C}\simeq \ol{\mathbb{Q}}_p$. Let $(r_\mu, V)$ be the representation of the $L$-group of $G$ determined by $\mu$ (see \S\ref{sec:Eichler-Shimura}). We assume the residue field of $E$ has cardinality $q$ and $\operatorname{Frob}_q\in \operatorname{Gal}(\ol{\mathbb{F}}_q/\mathbb{F}_q)$ is the arithmetic Frobenius. Then for each $i=0,\cdots,\operatorname{rk} r_\mu$, consider the $\Lambda$-valued function on the space of unramified Langlands parameters for $G$ 
\[T_i: \phi\mapsto \operatorname{tr}(r_\mu \circ\phi(\mathrm{Frob}_q)\mid \wedge^i V).\]
This can be viewed as an element in a spherical Hecke algebra $\mathcal{H}_{K_p,\Lambda}$ for some hyperspecial subgroup $K_p\subset G(\qp)$ via the Satake isomorphism\footnote{For simplicity, we assume $G$ is unramified. One can also make a conjecture for ramified groups using spherical $L$-parameters, as in \cite{vdHove}.}. For an Iwahori subgroup $I\subset K_p$, we identify $\mc{H}_{K_p}$ with the center of $\mc{H}_I$ (which acts on $R\Gamma(\mathrm{Sh}^*_{K^pI, \ol{E}}, 
\IC_{\mathrm{Sh}^*_{K}})$ from the right) via the Bernstein isomorphism. We define the Hecke polynomial 
\[H_\mu(X):=\sum_{i=0}^{r=\mathrm{dim}V} (-1)^i T_iX^{r-i}\in \mc{H}_{K_p}[X].\] 
Then the conjectural congruence relation can be formulated as follows:

\begin{conjecture}[Blasius--Rogawski]\label{conj: ES}
Let $I\subset G(\qp)$ be an Iwahori subgroup, and $K=IK^p$. Then the action of the inertia subgroup $I_E\subset W_E$ on $R\Gamma(\mathrm{Sh}^*_{K, \ol{E}}, 
\IC_{\mathrm{Sh}^*_{K}})$ is unipotent. For any $\sigma_E\in W_E$ lifting arithmetic Frobenius $\operatorname{Frob}_q$, one has the relation $H_\mu(\sigma_E)=0$ as endomorphism of this cohomology complex.
\end{conjecture}

\begin{remark}
Similar conjectures can be made for usual or compactly supported cohomology, with more general coefficients, including in the case when $\Lambda$ is integral or torsion and our results also address this.    
\end{remark}

As a third application of Theorem~\ref{thm: main}, we obtain the following:
\begin{theorem}[Theorem~\ref{thm: ES}]\label{thm: ES intro}
    Under Assumption~\ref{assumption:codimension}, Conjecture~\ref{conj: ES} is true.
\end{theorem}

\begin{remark}
Conjecture~\ref{conj: ES} and its variants have been extensively studied by \cite{faltings-chai}, \cite{Wedhorn}, \cite{BueltelWedhorn}, \cite{Koskivirta}, \cite{Lee}, \cite{Wu} using algebraic correspondences on the special fiber of Shimura varieties, and more recently by \cite{DvHKZ}, \cite{vdHove}, \cite{DvHKZ2} using Igusa stacks and the Fargues--Scholze local Langlands, following ideas of \cite{XiaoZhu} and \cite{KoshikawaES}. For the usual or compactly supported cohomology, the most general result is obtained in \cite{DvHKZ2}, which deals with abelian type Shimura varieties at Iwahori levels. 

For intersection cohomology, the best known result so far is for Hodge type Shimura varieties at hyperspecial level by \cite{Wu}. Our approach follows \cite{DvHKZ}, and hence applies to Iwahori levels. On the other hand, we do not work with Hodge type Shimura varieties, since minimal compactifications of Igusa stacks are not currently available in this generality. 
\end{remark}

\subsection{Proof ideas}

\subsubsection{The stratified cartesian diagram}
From Diagram (\ref{eqn: CartesiandiagramforClosedShimuraVariety}) and the work of \cite[Theorem~D]{Kim}, one can deduce\footnote{We actually give an independent proof of this result - see Corollary \ref{cor: BeyondGoodRed}.} the existence of an open substack, denoted $j: \Igs \hookrightarrow \Igs^{*}$, whose pullback along $h: \mathcal{F}\ell\to \Bun_G$ recovers the open Shimura variety $\mc{S}_{K^p}$. 

It is therefore natural to attempt to define an intersection complex $\IC_{\Igs^{*}}$ on $\Igs^*$ as an intermediate extension of the constant sheaf $\Lambda$ along $j$\footnote{Recall that, by \cite[Section 8.3]{DvHKZ}, the Igusa stack $\Igs$ is $\ell$-cohomologically smooth of $\ell$-dimension $0$ and has trivial dualizing sheaf.}. To do this, we use a perverse $t$-structure, defined in terms of a boundary stratification of $\Igs^*$ that mimics the boundary stratification on $\mc{S}^*_{K^p}$. 
Indeed, the infinite level Shimura variety $\mc{S}^*_{K^p}$ admits a locally closed, away from $p$ Hecke-equivariant and $G(\Q_p)$-stable stratification (the \textit{Baily--Borel stratification})
\begin{equation}\label{Eq: Shimura Stratification Intro}
\mathcal{S}^{*}_{K^p} = \bigsqcup_{[\mathsf{P}]} \mathcal{S}_{K^p, [\mathsf{P}]}, 
\end{equation}
where $[\mathsf{P}]$ runs over conjugacy classes of admissible rational parabolics $\mathsf{P} \subset \mathsf{G}$ (see \S \ref{Subsub:StratificationShimura}). Our first main geometric result is to prove that this stratification is pulled back from a similar locally closed, away from $p$ Hecke-stable stratification 
\begin{equation}\label{Eq: Igusa Stratification Intro}
\Igs^{*} := \bigsqcup_{[\mathsf{P}]} \Igs_{[\mathsf{P}]}
\end{equation}
that already exists on the level of the minimally compactified Igusa stack. More precisely, we prove the following theorem, which should be of independent interest. 

\begin{theorem}[Theorem \ref{Thm:StratifiedCartesian}]{\label{thm: stratifiedCartesianIntro}}
For each rational conjugacy class of admissible parabolics $\PP$, there exists a locally closed substack $\Igs_{[\mathsf{P}]} \hookrightarrow \Igs^{*}$ such that Diagram (\ref{eqn: CartesiandiagramforClosedShimuraVariety}) restricts to a cartesian diagram
    \begin{equation}{\label{Eq: stratifiedCartesianDiagramIntro}}
    \begin{tikzcd}
        \mathcal{S}_{K^p, \PP} \ar[r,"\pi_{[\mathsf{P}]}"] \ar[d,"\tilde{h}_{[\mathsf{P}]}"] & \mathcal{F}\ell\ar[d,"h"]\\
        \Igs_{\PP} \ar[r,"\overline{\pi}_{[\mathsf{P}]}"] & \Bun_{G}.
    \end{tikzcd}
    \end{equation}
    of $v$-stacks over $\Spd \overline{\mathbb{F}}_{p}$. 
\end{theorem}

It is not too hard to construct the commutative diagram~\eqref{Eq: stratifiedCartesianDiagramIntro}. In fact, one may define $\Igs_\PP$ simply to be the ($v$-sheaf-theoretic) image of $\mc{S}_{K^p,\PP}$ in $\Igs^\ast$. To prove that the diagram is cartesian, we can specialize to individual Harder--Narasimhan strata on $\mathrm{Bun}_{G}$ and use the description of the fibers of $\ol{\pi}_{\HT}$ as $[\mathrm{Ig}^{b,*, \diamond}/\mathcal{J}_b]$ (up to canonical compactification).  It turns out that the key point needed to prove that Diagram~\eqref{Eq: stratifiedCartesianDiagramIntro} is cartesian is to show that the $\mc{J}_b$-action on $\Ig^{b,\ast,\diamond}$ respects the boundary stratification by the conjugacy classes $[\mathsf{P}]$. Since this action is not constructed moduli theoretically, but rather using a form of Hartogs extension principle, this is far from obvious from the construction. However, this result can be proved on the level of partial toroidal compactifications, where a moduli-theoretic description of the action is available after passing to formal completions of the boundary strata.

\subsubsection{Universal local acyclicity, Verdier self-duality and comparison with $\IC(\mathsf{G},\mathsf{X})_{K^p}$} For simplicity, let us explain the ideas in the case $\Lambda$ is regular and torsion. Theorem~\ref{thm: stratifiedCartesianIntro} implies that the boundary strata $\mathrm{Igs}_{\PP}\subset \mathrm{Igs}^*$ are $\ell$-cohomologically smooth of dimension $d_{[\mathsf{P}]}-d$, where $d_{[\mathsf{P}]}$ is the dimension of the boundary stratum labeled by $[\mathsf{P}]$ in the Shimura variety. It is therefore natural to define a perverse $t$-structure on $\mathrm{Igs}^*$ with respect to the stratification~\eqref{Eq: Igusa Stratification Intro}, by gluing the standard $t$-structures on the boundary strata $\mathrm{Igs}_{\PP}$ with a degree shift by $d_{[\mathsf{P}]}-d$. Using this ad hoc perverse $t$-structure, we may now define an intersection complex $\mathrm{IC}_{\mathrm{Igs^\ast}}\in \Detale(\Igs^\ast,\Lambda)$ as the middle extension of the constant sheaf $\Lambda$ under $j: \Igs\hookrightarrow\Igs^\ast$

To construct the relative intersection cohomology of the Igusa stack, we can then simply define 
\[
\mathcal{F}_{\IC,\Lambda}:= R\ol{\pi}_{\HT*}(\mathrm{IC}_{\mathrm{Igs^\ast}}).
\] 
While this definition is straightforward, the properties established in Theorem~\ref{thm: main} are more subtle. 
One obstruction is the lack of a well-behaved formalism of perverse sheaves on Artin $v$-stacks, and another obstacle is the difference between the stack quotient $[\mc{S}^\ast_{K^p}/K_p]$ and its coarse moduli space $\mc{S}_{K^pK_p}^\ast$. 

For example, 
the Verdier self-duality of $\mathcal{F}_{\IC,\Lambda}$ follows from that of $\IC_{\Igs^\ast}$. This can be checked by a formal argument from the $\ell$-cohomological smoothness of the boundary strata, but only after we know a strong finiteness condition, which we call \emph{stratified ULA-ness} (Definition~\ref{defn: StratifiedULA}) or a slightly weaker condition which we call \emph{stratified reflexivity} (Definition \ref{defn: stratifiedreflexive}). This is a consequence of the fact that we work within $p$-adic geometry rather than in the more familiar algebraic geometric set-up of constructible sheaves on schemes of finite type and our sheaves behave more like a colimit of constructible sheaves with smooth group actions. In particular, constructible sheaves in this setting do not interact well with Verdier duality and this notion of universal local acyclicity is designed to fix this. 

To establish the stratified ULA property, we give a criterion to verify universal local acyclicity in a certain geometric setting (see Proposition~\ref{prop: ULA-Criterion}), inspired by the correspondence between ULA-ness for sheaves on the classifying stack of a $p$-adic group and admissibility for smooth representations of the $p$-adic group. To verify this criterion in our case, we compare with the analytification of the usual intersection complexes, after pulling back $\IC_{\Igs^*}$ under $\tilde{h}$ and pushing it forward to the coarse quotients $\mathcal{S}^*_{K^pK_p}$. This comparison is also key to establishing part (3) of Theorem~\ref{thm: main}, i.e. showing that $\mathcal{F}_{\IC,\Lambda}$ recovers the complex  $\IC(\mathsf{G},\mathsf{X})_{K^p}$ after applying the Hecke operator $i_1^*T_\mu$. 


\subsubsection{The computation of stalks and perversity} 
For the computation of the stalks of $\mathcal{F}_{\IC,\Lambda}$, we use three ingredients. The first ingredient is smooth base change, together with the fact that we have an open immersion $[\mathrm{Ig}^{b,*,\diamond}/\mathcal{J}_b]\hookrightarrow \Igs^*$ for every $b\in B(G,\mu^{-1})$ (Proposition~\ref{prop: openstrata}). Indeed, recall that the locally closed stratum of $\Igs^{*}$ given by pulling back the stratum $i_{b}: \Bun_{G}^{b} \hookrightarrow \Bun_{G}$ differs from $[\mathrm{Ig}^{b,*,\diamond}/\mathcal{J}_b]$ by some higher rank points in the analytic space, and it is precisely the oddities caused by these higher rank points that allow $[\mathrm{Ig}^{b,*,\diamond}/\mathcal{J}_b]\hookrightarrow \Igs^*$ to always be an open immersion. In particular, the proof of this fact relies on the observation that there are, a priori, two ways to define a Newton stratification on $\mathcal{S}^*_{K^p}$: one purely on the generic fiber
\[
\mathcal{S}^*_{K^p} = \bigsqcup_{b\in B(G,\mu^{-1})} \mathcal{S}_{K^p}^{*,b,\eta},
\]
by pulling back from the Harder--Narasimhan stratification on $\Bun_G$, and one by pulling back from the Newton stratification on the special fiber
\[
\mathcal{S}^*_{K^p} = \bigsqcup_{b\in B(G,\mu^{-1})} \mathcal{S}_{K^p}^{*,b,s}.
\] 
The closure relations between the generic fiber and the special fiber are reversed precisely due to the presence of higher rank points and, therefore, the intersection $\mathcal{S}_{K^p}^{*,b,\eta}\cap \mathcal{S}_{K^p}^{*,b,s}$ can be shown to be an open subspace of $\mathcal{S}^*_{K^p}$. 

The second ingredient is a \emph{transversality} result between Newton and boundary stratifications on $\Igs^*$ (Proposition~\ref{prop: dimensionStrataIg}). If we denote by $d_{b,\PP}$ the dimension of the boundary stratum in $\mathrm{Ig}^{b,*}$ labeled by the conjugacy class of admissible parabolics $\PP$, then we have the equality 
\[
d_{b,\PP} - d_b = d_{\PP} - d. 
\]
This formula follows from Theorem~\ref{thm: stratifiedCartesianIntro}. Heuristically, the underlying geometric reason for this equality is that the degenerate part of a semi-abelian scheme is a torus, with $p$-divisible group of multiplicative type. This strongly
constrains how the Newton and boundary strata intersect. Together with a version of Deligne's formula in the setting of $v$-stacks and smooth base change, this gives a concrete comparison between the ad hoc intersection complexes on $\Igs^*$ and on $[\mathrm{Ig}^{b,*,\diamond}/\mathcal{J}_b]$. 

The third ingredient comes from~\cite{GleasonTubularNeighborhoods, GHILZIsocComparison}. These allow us to compare the colimit of intersection cohomology complexes of finite-level versions of $\Ig^{b,*}$, which are algebraic varieties over $\ol{\mathbb{F}}_p$, with the ad hoc intersection cohomology complex of the $v$-stack $[\Ig^{b,*,\diamond}/\mathcal{J}_b]$. Here, we emphasize that the functor $(-)^{\diamond}$ is an analytification operation that interacts with the theory of $\ell$-adic \'etale sheaves on schemes and analytic spaces in a non-obvious way. The references \cite{GleasonTubularNeighborhoods, GHILZIsocComparison} allow us to control the subtle behavior of this analytification.

Finally, we discuss the perversity of $\mathcal{F}_{\IC,\Lambda}$. Recall that the perverse $t$-structure on $\mathrm{Bun}_G$ is constructed by gluing the standard $t$-structures on the Harder--Narasimhan strata 
$\mathrm{Bun}^b_G$ with a degree shift by $d_b$, which can be defined as $\langle 2\rho_{G},\nu_{b} \rangle$ and which, when $b\in B(G,\mu^{-1})$, agrees with the dimension of the Igusa variety $\mathrm{Ig}^{b,*}$. As $\mathcal{F}_{\IC,\Lambda}$ is Verdier self-dual on $\mathrm{Bun}_{G,\mu^{-1}}$, to prove that it is perverse it suffices to establish upper semi-perversity. This follows from the computation of stalks in part (4) of Theorem~\ref{thm: main} together with the affineness of $\mathrm{Ig}^{b,*}$ and Artin vanishing. 

\subsection{Concluding remarks}
One of the most important properties of the intersection cohomology (of minimal compactifications) of Shimura varieties with $\ol{\Q}_\ell$-coefficients is that, unlike compactly supported or usual cohomology, it admits a direct sum decomposition over discrete automorphic representations, via its relation to $L^2$-cohomology. On the other hand, besides the Hecke and Galois actions, it carries a Lefschetz structure and a Hodge structure. 

It is particularly interesting to study the interaction between all these structures. This is related to Arthur's (conjectural) classification of automorphic representations in terms of the Arthur parameters. Compared with the Weil--Deligne form of Langlands parameters, Arthur parameters have an extra copy of $\SL_2$, which should encode the Lefschetz structure on intersection cohomology. It is expected that the range of cohomological degrees to which a given discrete automorphic representation $\pi$ contributes is controlled by the Arthur $\mathrm{SL}_2$ of $\pi$, an explicit algebraic representation attached to $\pi$ under Arthur's conjectures, see \cite{ArthurUnipotentConjectures} and \cite{KoshikawaShin}. 

We expect that some of these structures exist already on the level of the perverse sheaf $\mathcal{F}_{\mathrm{IC},\Lambda}$. More precisely, we expect that $\mathcal{F}_{\mathrm{IC},\Lambda}$ admits a direct sum decomposition over discrete automorphic representations and, within each summand, a further decomposition in terms of direct summands of generalized eigensheaves on $\mathrm{Bun}_G$ restricted to $\Bun_{G,\mu^{-1}}$, which occur with globally determined (automorphic) multiplicities. This is inspired by the local-global compatibility part of Fargues' conjecture~\cite[Section 7]{Fargues-Overview}, and we believe that this is the correct interpretation of Fargues' vision.

More precisely, we believe that one will see the direct summands of the generalized eigensheaves attached to (local) Arthur parameters, that are being constructed by Bertoloni-Meli and Koshikawa~\cite{CuspidalSheavesBertiKoshikawa} on the spectral side of the categorical conjecture. These generalized eigensheaves have a shearing property that leads to a spreading out of cohomological degrees after applying Hecke operators. This is determined by a local Arthur $\mathrm{SL}_2$ in a way that precisely mirrors the way that the global Arthur $\mathrm{SL}_2$ determines the Lefschetz structure on intersection cohomology. 

See \S \ref{ss: EigensheafConjecture} for a more detailed discussion of this topic. In joint work in progress with Bertoloni-Meli and Koshikawa, we formulate a precise conjecture that relates $\mathcal{F}_{\mathrm{IC},\Lambda}$ to generalized eigensheaves, and verify that it is consistent with various other conjectures in the literature concerning the cohomology of Shimura varieties and Igusa varieties.

\subsection{Organization of the paper} In \S\ref{sec: prelim}, we discuss various geometric and sheaf-theoretic preliminaries, including analytification functors, a construction of intersection complexes on $v$-stacks equipped with nice stratifications, and our criterion for verifying ULA-ness. 

In \S\ref{sec: stratifications}, we review the structure of the minimal compactifications $\mathcal{S}^*_{K}$, following~\cite{PinkHigherDirectImages}, and we use this to establish the $\ell$-cohomological smoothness of boundary strata in the stacky quotients $[\mathcal{S}^*_{K^p}/\underline{K_p}]$.

In \S\ref{sec: pc Igusa varieties}, we review properties of partial minimal and toroidal compactifications of Igusa varieties. We construct a moduli-theoretic action of (a truncated form of) $\mathcal{J}_b$ on the partial toroidal compactifications and show that it preserves the boundary strata labeled by conjugacy classes $\PP$.  

In \S\ref{sec: IgusaStack}, we introduce the minimally compactified Igusa stack and establish the stratified cartesian diagram of Theorem~\ref{thm: stratifiedCartesianIntro}. We construct the intersection complex $\IC_{\Igs^*}$ and prove part (1) of Theorem~\ref{thm: main}.  

In \S\ref{sec: Stalks}, we prove parts (2) and (4) of Theorem~\ref{thm: main}. In \S\ref{sec: applications}, we prove part (3) of Theorem~\ref{thm: main}. We deduce our main applications, namely Theorems~\ref{thm: product formula intro}, \ref{thm: genericLocalization} and~\ref{thm: ES intro}. Finally, we conclude with a speculative discussion on how our work connects Arthur's conjectures on intersection cohomology with the categorical local Langlands program in the setting of Fargues--Scholze.

\section*{Acknowledgments}
We would like to thank the Simons Collaboration on Perfection in Algebra, Geometry and Topology, during whose annual meeting the project started. We would like to thank Alexander Bertoloni Meli, Bharghav Bhatt, George Boxer, Patrick Daniels, Mark Goresky, Ian Gleason, Tom Haines, David Hansen, Pol van Hoften,  Dongryul Kim, Teruhisa Koshikawa, Si Ying Lee, Lucas Mann, Peter Scholze, Matteo Tamiozzo and Xinwen Zhu for various helpful discussions surrounding this project. We thank Ofer Gabber, Kai-Wen Lan, Si Ying Lee,  Teruhisa Koshikawa, Shengkai Mao, Matteo Tamiozzo and Peihang Wu for providing helpful feedback on a first draft of this paper. The third named author thanks the University of California, Berkeley and IHES for hospitality. Part of this project was carried out during her visit to these institutes.

\section*{Notation and conventions} \label{Notation}
\begin{itemize}
\item Fix rational primes $\ell\neq p$. 
\item We denote the category of perfectoid spaces in characteristic $p$ by $\Perf$, equipped with the v-topology as defined in \cite{Ecod}.
\item Let $X$ be a space or a v-sheaf, with an action by some group $\Gamma$. We denote by $[X/\Gamma]$ the stack quotient of $X$ by $\Gamma$, and denote by $X\sslash\Gamma$ its coarse moduli, as a space or a v-sheaf.
\item We call a PEL type Shimura datum $\gx$ \textit{of type AC}, if the absolute root datum of each simple factor of $\mathsf{G}$ is either of type $A$ or $C$. 
\item We set $\Breve{\bb{Q}}_{p} := W(\ol{\bb{F}}_{p})[\frac{1}{p}]$ to be the completion of the maximal unramified extension of $\bb{Q}_{p}$. We write $\sigma$ for its Frobenius endomorphism.
For $E/\bb{Q}_{p}$ a finite extension, we let $\Breve{E} := E\Breve{\bb{Q}}_{p}$ be the compositum with the maximal unramified extension.  
\item For a linear algebraic group $H/\bb{Q}_{p}$, we let $B(H) := H(\Breve{\bb{Q}}_{p})/(b \sim hb\sigma(h)^{-1})$ denote the set of $\sigma$-conjugacy classes (the Kottwitz set). 
\item For $G/\bb{Q}_{p}$ a connected reductive group with a fixed choice of maximal torus $T \subset G$, we recall \cite{KottwitzIsocrystalsI,RR} that the set $B(G)$ is equipped with two maps 
\begin{enumerate}
\item The slope homomorphism
    \[ \nu: B(G) \rightarrow \bb{X}_*(T_{\overline{\mathbb{Q}}_{p}})^{+,\Gamma}_{\mathbb{Q}} \]
    \[ b \mapsto \nu_{b}, \]
    where $\Gamma := \mathrm{Gal}(\overline{\mathbb{Q}}_{p}/\mathbb{Q}_{p})$ and 
    $\bb{X}_*(T_{\overline{\mathbb{Q}}_{p}})_{\mathbb{Q}}^{+}$ is the set of rational dominant cocharacters of $G$. 
    \item The Kottwitz invariant
    \[ \kappa: B(G) \rightarrow \pi_{1}(G)_{\Gamma} \]
    where $\pi_1(G) = \bb{X}_*(T_{\bar{\mathbb{Q}}_{p}})/X_*(T_{\bar{\mathbb{Q}}_{p},sc})$ for $T_{\bar{\mathbb{Q}}_{p},sc}$ the maximal simply connected subtorus of $T_{\bar{\mathbb{Q}}_{p}}$. 
\end{enumerate}
Moreover, the product map $\nu \times \kappa$ is injective. We can endow $B(G)$ with a partial ordering where $b \geq b'$ if $\nu_{b} - \nu_{b'}$ (in this order!) is a $\bb{Q}_{\geq 0}$-linear combination of positive coroots and $\kappa(b) = \kappa(b')$. 
\item For $G/\qp$ and $b\in B(G)$ as above, we denote by $J_b$ the $\sigma$-centralizer of $b$. This is again a reductive group over $\qp$.
\item For $\mu \in \bb{X}_{*}(T_{\ol{\mathbb{Q}}_{p}})^{+}$ a geometric dominant cocharacter, we write $B(G,\mu) \subset B(G)$ for the $\mu$-admissible locus (\cite[Definition~2.3]{RV}).
\item Throughout, we will work with various derived categories that admit $\infty$-categorical enhancements, which will occasionally be essential to work with (e.g if we wish to take inverse limits of the derived categories). In order to emphasize this, we will usually use the notation $\mathcal{D}$ to denote this $\infty$-categorical enhancement and the notation $\mathrm{D}$ to denote the underlying homotopy category. E.g. $\mathcal{D}(X,\Lambda)$ is the $\infty$-category of $\Lambda$-sheaves on $X$ with underlying homotopy category $\mathrm{D}(X,\Lambda)$.
\item All our functors will be derived unless otherwise stated. Due to the force of habit, we will still write $R\Gamma_{c}(-,\Lambda)$ and $R\Gamma(-,\Lambda)$ for compactly supported and regular cohomology, and $R\mathcal{H}om(-,-)$ for internal Hom in derived categories of sheaves. Similarly, we write $\RHom(-,-)$ for external Hom. However, for the usual pushforward and pullback as well as limits, colimits, and tensor products we will suppress the derived notation.

\item As coefficients, we will take the following: 

\begin{setup}{\label{assumption: coefficientsystemsingeneral}}
We let $\Lambda/\bb{Z}_{\ell}$ be one of the following.
\begin{enumerate}
\item $\Lambda$ is a finite ring that is $\ell$-power torsion.
\item  $\Lambda/\bb{Z}_{\ell}$ is an $\ell$-adically complete $\ell$-torsion-free ring, whose mod $\ell$ reduction is finite.
\end{enumerate}
\end{setup}
For a diamond or $v$-stack, in (1) we take $\Detale(X,\Lambda)$ to be the category of \'etale sheaves constructed in \cite{Ecod}. In (2), we will work with $\Detale(X,\Lambda)$ the category of $\Lambda$-adic sheaves.  By~\cite[Proposition~26.2]{Ecod}, this can be described as the underlying homotopy category of the inverse limit $\cDetale(X,\Lambda) := \lim_{n \geq 1} \cDetale(X,\Lambda/\ell^{n})$ for $n \geq 1$, where $\cDetale(X,\Lambda/\ell^{n})$ is the $\infty$-categorical enhancement of the category of sheaves described in (1). This upgrades to a full six functor formalism equipped with excision (since mod $\ell$-reduction is conservative and excision holds in the torsion case). We will often cite results from \cite{Ecod} that concern the torsion case; however, as remarked at the end of \cite[Section~26]{Ecod}, it is easy to see that they carry over to this context. We will also cite results from \cite{Mann2022NuclearSheaves} which uses the nuclear six functor formalism; however, we will only invoke these results in the torsion case, and implicitly pass between the comparison of nuclear and the \'etale category in this case (\cite[Proposition~3.20]{Mann2022NuclearSheaves}), implicitly deducing results from the torsion case for the adic category by taking inverse limits.

We will also make use of ULA objects with respect to the structure map $X \ra \ast$ of an Artin $v$-stack in all of these contexts. If $\Lambda$ is as in (1) and $X$ is an Artin $v$-stack, this is as in \cite[Definition IV.2.31]{FSGeomLLC}. If $\Lambda$ is as in (2) then we define $\Detale^{\ULA}(X,\Lambda)$ to be the homotopy category of the inverse limit of $\cDetale^{\ULA}(X,\Lambda/\ell^{n})$ for all $n \geq 1$.

For our ultimate results, we will use the following setup.  
\begin{setup}{\label{assumption: coefficients}}
We assume that $\Lambda/\mathbb{Z}_{\ell}$ is of the following form: 
\begin{enumerate}
\item $\Lambda$ is a finite ring that is self-injective and $\ell$-power torsion.
\item $\Lambda = \mathcal{O}_{F}$ for $F/\bb{Q}_{\ell}$ a finite extension.
\item $\Lambda = F$ for $F/\bb{Q}_{\ell}$ a finite extension. 
\item $\Lambda = \mathcal{O}_{F},F$ for $F/\bb{Q}_{\ell}$ a (possibly infinite) algebraic extension.
\item $\Lambda$ is a (possibly infinite) algebraic extension of $\bb{F}_{\ell}$.
\end{enumerate}
\end{setup}

\noindent Self-injectivity in case (1) will allow us to define a perverse sheaf theory that behaves well with respect to Verdier duality.

\item For $\Lambda$ a coefficient system as in Setup~\ref{assumption: coefficientsystemsingeneral} (1) and $X$ a $v$-stack, we define $\Detale^{\oc}(X,\Lambda) \subset \D(X,\Lambda)$ to be the full subcategory of overconvergent sheaves. If $X = \Spa(C,C^{+})$ then this is the full subcategory of sheaves such that the adjunction map $j_{*}j^{*} \ra \mathrm{id}$ is an isomorphism, where $j: \Spa(C,\mathcal{O}_{C}) \ra \Spa(C,C^{+})$ is the inclusion of the rank one geometric point. In general, it is the full subcategory of sheaves whose pullback to each geometric point is overconvergent. For $\Lambda$ a coefficient system as in Setup~\ref{assumption: coefficientsystemsingeneral} (2), we define $\Detale^{\oc}(X,\Lambda)$ by taking inverse limits.
\item For a locally profinite group $K$ and $\Lambda/\bb{Z}[\frac{1}{p}]$ a ring, we write $\D(K,\Lambda)$ for the left-completed derived category of smooth representations of $K$ on $\Lambda$-modules. If $\Lambda$ is as in Setup~\ref{assumption: coefficientsystemsingeneral} (2), we write $\widehat{\D}(K,\Lambda)$ for the derived category of $\ell$-adically-complete smooth $\Lambda$-representations, cf. \cite[Section~2.2]{HansenSupercuspidalCohomology}. Formally, it is the underlying homotopy category of the inverse limit $\lim_{n \geq 1} \mathcal{D}(K,\Lambda/\ell^{n})$. We also use this notation when $\Lambda$ is as in Setup~\ref{assumption: coefficientsystemsingeneral} (1); however, we note that in this case $\widehat{\D}(K,\Lambda) = \D(K,\Lambda)$. We omit $K$ when $K=\{1\}$ is trivial. Then $\D(K,\Lambda)=\D(\Lambda)$ is the usual derived category of (discrete) $\Lambda$ modules, and $\widehat{\D}(\Lambda)$ is the subcategory of $\ell$-complete objects.
\end{itemize}

\section{Preliminaries}\label{sec: prelim}
\subsection{Variants on Analytification}{\label{ss: Analytification}}
In this section, we discuss analytifications of perfect schemes in characteristic $p$ and analytification of finite type schemes over mixed characteristic fields. 

\subsubsection{Perfect schemes over $k$}{\label{ss: analytificationofperfectschemes}}
We will need some variants of the analytification functor in the context of perfect schemes, as described extensively in \cite[Section~4.2]{GHILZIsocComparison}. We let $k$ be a perfect field of characteristic $p$, equipped with the discrete topology. Let $X$ be a groupoid valued fpqc sheaf on perfect schemes over $k$, e.g the Yoneda sheaf attached to a perfect scheme. We may consider the functors on $\Perf_{k}$
	\begin{align*}
		\label{some functors pre}
		X^{\diamond_\pre}:(R,R^{+})& \mapsto X(\Spec{R^{+}}) \\
		X^{\dagger_\pre}:(R,R^{+})& \mapsto X(\Spec{R^\circ}) \\
		X^{\Diamond_\pre}:(R,R^{+})& \mapsto X(\Spec{R}).
	\end{align*}
    We denote their associated $v$-sheafifications by $X^{\diamond}$, $X^{\dagger}$, and $X^{\Diamond}$. There are natural morphisms 
    \begin{equation}{\label{eqn: NaturalTransformationsofAnalytifications}}
    X^{\diamond} \xrightarrow{a_{X}} X^{\dagger} \xrightarrow{b_{X}} X^{\Diamond}. 
    \end{equation}
    induced by the natural maps 
    of rings $R^{+} \subset R^{\circ} \subset R$. We write $t_{X}$ for the composition $b_{X} \circ a_{X}$. 
    \begin{remark}{\label{rem: representablevsheafifyingautomatic}}
    We note that if $X$ is representable by an affine perfect scheme then it follows by \cite[Theorem 8.7]{Ecod} that $v$-sheafifying is unnecessary. 
    \end{remark}
    The following claim describes the relation between $(-)^{\diamond}$ and $(-)^{\dagger}$.
    \begin{lemma}\label{lem: DaggerisCC}
        Let $X$ be a perfect scheme separated over $k$. The map $a_X$ identifies $X^\dagger$ with the canonical compactification of $X^\diamond$ in the sense of \cite[Section 18]{Ecod}.
    \end{lemma}
    \begin{proof}
        This is clear from the definition of the canonical compactification.
    \end{proof}
    
    An important point for us will be that the transformations $(-)^{\diamond} \ra  (-)^{\dagger} \ra (-)^{\Diamond}$ also interact well with closed immersions of perfect schemes $Z \hookrightarrow X$, as essentially follows from the valuative criterion for properness and the fact that finitely presented morphisms of schemes interact well with colimits of points.
    \begin{lemma}{\cite[Lemma~4.2.3]{GHILZIsocComparison}}{\label{lemma: CartesianGivesProper}}
		Let $f \colon X \to Y$ be a perfectly finitely presented proper morphism of perfectly finitely presented schemes over $k$. Then the natural transformations (\ref{eqn: NaturalTransformationsofAnalytifications}) fit into a commutative diagram
		\[ 
			\begin{tikzcd}
				X^{\diamond} \arrow[r] \arrow[d,"f^{\diamond}"] & X^{\dagger} \arrow[r] \arrow[d,"f^{\dagger}"] & X^{\Diamond} \arrow{d}{f^{\Diamond}} \\
				Y^\diamond \arrow{r} & Y^\dagger \arrow{r} & Y^{\Diamond}
			\end{tikzcd}
		\]
		with cartesian squares.
	\end{lemma}

We will later apply these lemmas to (partial minimally compactified) Igusa varieties to compute the stalks of the sheaf $\mc{F}_\mathrm{IC}$ in terms of the intersection cohomology of Igusa varieties, see \S \ref{sec: Stalks}.
    
\subsubsection{Analytification in Mixed Characteristic}{\label{ss: AnalytificationinMixedCharacteristic}}
We now discuss analytification in a slightly more general setup following \cite[Lecture 18]{SW}. Namely, to a pre-adic space $X$ over $\bb{Z}_p$, we can associate a v-sheaf $X^{\Diamond}/\Spd \bb{Z}_p$ by taking the functor on $\Perf$
\[\Spa(R,R^+)\mapsto\{(R^\sharp, \iota, f)\},\]
where $(R^\sharp,\iota)$ is an untilt of $R$ over $\bb{Z}_p$, and $f$ is a map $\Spa(R^\sharp, R^{\sharp+})\to X$. 

\begin{remark}
We note that if $X=\Spec A$ is an affine perfect scheme over $\bb{F}_p$, one can attach to it a pre-adic space via $\Spa(A,A)$, where $A$ is equipped with the discrete topology; or $\Spa(A,A^+)$, where $A^+$ is the integral closure of $\bb{F}_p$ in $A$. Then it has two attached v-sheaves: the data $(R^\sharp,\iota)$ is uninteresting, i.e. $R^\sharp=R$ and $\iota$ is the identity. In the first case, $f$ is equivalent to a map $A\to R^{\sharp+}$. This gives a $v$-sheaf, by Remark~\ref{rem: representablevsheafifyingautomatic} and recovers the definition $X^\diamond$ in \S~\ref{ss: analytificationofperfectschemes}. In the second case, $f$ amounts to a map $A\to R$, and hence the attached v-sheaf recovers $X^\Diamond$ in \S~\ref{ss: analytificationofperfectschemes}.
\end{remark}

When $X$ is a finite type scheme over a non-archimedean field over $\qp$, one can attach to it an adic space $X^\mathrm{an}$ by taking analytification, see for example \cite[Definition 10.3]{hubner}. It then has an associated v-sheaf $X^\Diamond=(X^\mathrm{an})^\Diamond$. When $X=\Spec A$ is affine (sheafification is needed in general), this can also be directly constructed from the scheme $X$ as in \cite[Definition 2.10(2)]{AGLRLocalModels}, cf. Lemma 2.11 in \textit{loc. cit.}. 

Another observation is that the functor $(-)^\Diamond$ preserves open and closed decompositions. 
\begin{lemma}\label{lem: ExcisionBigDiamond}
    Suppose $j: U\hookrightarrow X$ is an open immersion of schemes over $k$, and $i:Z\hookrightarrow X$ is the complementary closed immersion.
    The attached map of v-sheaves $j^\Diamond$ (resp. $i^\Diamond$) is an open (resp. closed) immersion. Moreover, there is a set theoretic decomposition of the underlying topological space of the small v-sheaf $X^\Diamond$
    \[|X^\Diamond|=|U^\Diamond|\coprod |Z^\Diamond|.\]
\end{lemma}
\begin{proof}
    The statement is local on $X$, so we may assume $X=\Spec A$ and $Z=\Spec A/I$ to be affine. Then, for any affinoid perfectoid space $S=\Spa(R,R^+)$ mapping to $X^\Diamond$ with associated untilt $S^{\sharp} = \Spa(R^{\sharp},R^{\sharp+})$, we have that the pullback of $i^\Diamond:Z^\Diamond\to X^\Diamond$ to $S$ is represented by the Zariski closed immersion 
    \[\Spa(R^{\sharp} \otimes_A A/I,\overline{R^{\sharp+}})\hookrightarrow S^{\sharp},\]
    under the isomorphism $|S^{\sharp}| \simeq |S|$ induced by (un)tilting. Here $\overline{R^{\sharp+}}$ is the integral closure of $R^{\sharp+}$ in $R\otimes_A A/I$, hence $i^\Diamond$ is a closed immersion. Similarly, for $j^\Diamond$, we may reduce to the case where $U$ is a standard open $U=\Spec A[1/f]$. Then for $S\to X^\Diamond$ as above, the pullback of $j^\Diamond$ to $S$ is representable by the open immersion
    \[S^{\sharp}\backslash V(\alpha(f))=\bigcup_n \Spa(R^{\sharp}\langle\varpi^n/\alpha(f)\rangle, R^{\sharp+}\langle\varpi^n/\alpha(f)\rangle)\hookrightarrow S^{\sharp},\]
    under the isomorphism $|S^{\sharp}| \simeq |S|$ induced by $i$, where $\alpha:A\to R^{\sharp}$ is the map corresponding to $S\to X^\Diamond$. This shows that $j^\Diamond$ is an open immersion. The statement on topological spaces is clear, by testing on geometric points (cf. \cite[Lemma~3.3.2]{GHILZIsocComparison}).
\end{proof}

\begin{remark}
    This result allows us to use excision of \'etale sheaves on $X^\Diamond$ for the open-closed decomposition given by $U^\Diamond$ and $Z^\Diamond$. Note that for the functors $(-)^\diamond$, $(-)^\dagger$, it is not the case that open-closed decompositions are preserved by them. For example, let $X=\Spec A$ be a perfect affine scheme in characteristic $p$, $Z = V(f)$ a divisor cut out by a function $f$, and $U=\Spec A[1/f]$ its open complement. Let $C$ be an algebraically closed perfectoid field in characteristic $p$. A $\Spa C$-point $x$ of $X^\diamond$ amounts to a map $\alpha: A\to \mc{O}_C$. Then $x$ factors through $Z^\diamond$ if and only if $\alpha(f)$ is zero, while $x$ factors through $U^\diamond$ means $\alpha(f)$ is invertible. Clearly these do not exhaust all possibilities, since it could be the case that $0\neq \alpha(f)$ lies in the maximal ideal of $\mc{O}_C$.
\end{remark}

We have the following relation that describes how the $(-)^{\diamond}$ functor on perfect schemes interacts with the $(-)^{\Diamond}$ functor for schemes or adic spaces over a non-archimedean base in characteristic zero via deformation theory. 
Let $k$ be a perfect field of characteristic $p$ and $X$ be a perfect $k$-scheme. Let $C/W(k)$ be a perfectoid field over the Witt ring. Let $W(X)_{\CO_C}$ be the formal scheme obtained by taking the Witt vector lift of $X$ and then base-changing to $\Spf \CO_C$. Denote by $W(X)_C$ its adic generic fiber, i.e.
\[W(X)_C:=W(X)_{\CO_{C}} \times_{\Spa(\mathcal{O}_{C},\mathcal{O}_{C})} \Spa(C,\mathcal{O}_{C}).\] 
This is a perfectoid space over $C$. We now have the following.

\begin{lemma}\label{lemma: PerfectDeformation}
With notation as above, we have an isomorphism
\[ X^\diamond \times_{\Spd k} \Spd C \simeq W(X)_{{C}}^{\Diamond} \] 
of diamonds over $\Spd C$.
\end{lemma}
\begin{proof}
Without loss of generality, we may assume $X$ is affine. Then perfectoidness of $W(X)_{C}$ is clear from the definition and the fact that $X$ is perfect. Let $S=\Spa(R,R^+)$ in $\Perf_k$ be an affinoid perfectoid space. Giving an $S$-point of $X^\diamond \times_{\Spd k}\Spd C$ is the same as giving a $\CO_{C^\flat}$-algebra structure on $R^+$, together with a $\Spec R^+$-point of $X$ (using Remark \ref{rem: representablevsheafifyingautomatic}). Since $X$ is perfect, a $\Spec R^+$-point is the same as a  $\Spec R^+/\varpi$-point, for any pseudo-uniformizer $\varpi\in R^+$ in light of the isomorphism $\lim_{x \mapsto x^{p}} R^{+}/\varpi \simeq R^{+}$. Hence, the datum $(\CO_{C^\flat}\to R^+, \Spec R^+\to X)$ amounts to a $\Spec R^+/\varpi$-point of the base-change $X_{\CO_{C^\flat}/\varpi}$, where we choose $\varpi$ to be in $\CO_{C^\flat}$. Let $R^{\sharp+}$ be the unique untilt of $R^+$ over $\CO_C$. Then a $\Spec R^+/\varpi$-point of $X_{\CO_{C^\flat}/\varpi}$ uniquely deforms to a $\Spf R^{\sharp +}$-point of $W(X)_{\CO_C}$, or equivalently, a $\Spa (R^\sharp, R^{\sharp+})$-point of its generic fiber.
\end{proof}

\subsection{Analytification and six functors}\label{sec: analytification} 

\subsubsection{Various sheaf categories}
\label{ss: VariousSheafCategories}
Let $k$ be a fixed field, either of characteristic $p$ with discrete topology, or over $\mathbb{Q}_{p}$ with non-archimedean topology. Let $X$ be a separated scheme over $k$.

We let $\Lambda$ be a coefficient system as in Setup~\ref{assumption: coefficientsystemsingeneral} and consider $\Detale(X,\Lambda)$, which is the left completed derived category of \'etale $\Lambda$-modules on $X$ in case (1) and the homotopy category of the inverse limit $\varprojlim_n\mc{D}_{\et}(X,\Lambda/\ell^n)$ in case (2). The functor $(-)^{\Diamond}$ from \S~\ref{ss: Analytification} defines a natural map $X_{v}^{\Diamond} \ra X_{\text{proet}}$ of sites, and in turn gives rise to a natural analytification functor 
    \[ c_{X}^\ast: \Detale(X,\Lambda) \ra \Detale(X^{\Diamond},\Lambda), \]
where we recall that the target is defined similarly, see our convention in \S~\ref{Notation}. 

We moreover consider the full subcategories 
\[\Dconstf(X,\Lambda)\subset \Dcons(X,\Lambda) \subset \Detale(X,\Lambda)\]
of perfect constructible and constructible sheaves. The former subcategory is spanned by objects $A\in \Detale(X,\Lambda)$ such that $A\otimes_{\Lambda} \Lambda/\ell$ is \'etale locally, after passing to a constructible stratification, a locally constant perfect complex. The latter one is similar, dropping the perfectness condition and replacing it with finitely generated.

We denote by $\Detale^{\oc}(X^{\Diamond},\Lambda)\subset \Detale(X^{\Diamond},\Lambda)$ the full subcategory of overconvergent sheaves, i.e. the subcategory spanned by objects that are constant on geometric points, see our conventions in \S~\ref{Notation}.

Finally, we write $\Detale^{\ULA}(X,\Lambda) \subset \Detale(X,\Lambda)$ (respectively $\Detale^{\ULA}(X^{\Diamond},\Lambda) \subset \Detale(X^\Diamond,\Lambda)$) for the full subcategory of complexes that are universally locally acyclic (ULA) over $\Spec k$ (resp. $\Spd k$). This is defined using adjunction in the $2$-category of kernels (which is a notion native to any six functor formalism, see \cite[Section~4]{HeyerMann}, where they use the terminology ``suave'' for ULA).  In particular, see \cite[Section~3]{HansenScholzeRelativePerversity} for a discussion of ULA sheaves in the algebraic setting, and \cite[Section~IV.2]{FSGeomLLC} for the analogous discussion in the analytic setting. 
\begin{remark}{\label{rem: PerfectConstructibleSameasULA}}
When $X/k$ is of finite presentation (or the perfection thereof if $\operatorname{char}k=p$) and $\Lambda$ is as in Setup~\ref{assumption: coefficientsystemsingeneral} (1), then we have the relation
\[\Detale^{\ULA}(X,\Lambda) \simeq \Dconstf(X,\Lambda)\]
by \cite[Proposition 3.4(iii), Theorem 4.4]{HansenScholzeRelativePerversity}. Indeed, Proposition 3.4(iii) implies that being ULA in $\Detale(X,\Lambda)$ is equivalent to being ULA in $\Dcons(X,\Lambda)$. Then the relation follows from the equivalence between part (i) and (ii) in \cite[Theorem~4.4]{HansenScholzeRelativePerversity}, where we note that the condition on the specialization map being an isomorphism in part (ii) is vacuous over a point. 

In Setup~\ref{assumption: coefficientsystemsingeneral} (2), using \cite[Lemma 4.4.5]{HeyerMann} to characterize ULA objects, then the same proof as in \cite[Proposition 7.14]{Mann2022NuclearSheaves} shows that $A\in \Detale(X,\Lambda)$ is ULA over $k$ if and only if $A/\ell \in \Detale^{\ULA}(X,\Lambda/\ell)$. This implies that the relation
\[\Detale^{\ULA}(X,\Lambda) \simeq \Dconstf(X,\Lambda)\]
still holds true.
\end{remark}
\subsubsection{Properties of Analytification}
Let $\Lambda$ be as in Setup \ref{assumption: coefficientsystemsingeneral}. Now we study some basic properties of the functor $c_X^\ast$. 

\begin{proposition}{\label{prop: propertiesofthealgebraizationfunctor}}
We assume that all schemes are separated over $k$. The following is true.
\begin{enumerate}
\item The functor $c_{X}^\ast$ commutes with the completed (derived) tensor product $-\hat{\otimes}_\Lambda-$ and admits a right adjoint $c_{X,\ast}$. 
\item For any map $f: X \ra Y$ with associated map of v-sheaves $f^{\Diamond}: X^{\Diamond} \ra Y^{\Diamond}$, we have a natural identification:
\[ c_{X}^\ast f^\ast \simeq f^{\Diamond,\ast}c_{Y}^\ast.\]
\item If $k$ is of characteristic $p$, then $c_X^\ast$ is fully faithful; if $k$ is a non-archimedean field over $\qp$, and $X$ is locally of finite type over $k$, then $c_X^\ast$ is fully faithful on $\Dcons(X,\Lambda)$.
\item If $f:X \ra Y$ is of finite type (or the perfection thereof if $\operatorname{char}k=p$), then we have natural isomorphisms
\[ f^{\Diamond}_{!}c_{X}^\ast \simeq c_{Y}^\ast f_{!}\]
\[ c_{X\ast} f^{\Diamond!}  \simeq f^{!}c_{Y\ast}.\]
\item If $f$ is as in (4) and $Y$ is of finite type (or the perfection thereof if $\operatorname{char}k = p$), then on $\Dcons(X,\Lambda)$ we have
\[ f^{\Diamond}_\ast c_{X}^\ast \simeq c_{Y}^\ast f_\ast\]
\[ c_{Y}^{*}f^{!}\simeq f^{\Diamond!}c_{X}^{*}, \]
and
\[ \bb{D}_{X^{\Diamond}/Y^\Diamond}c_{X}^{*} \simeq c_{X}^{*}\bb{D}_{X/Y}. \]
\item The functor $c_{X}^\ast$ factors through the full subcategory $\Detale^{\oc}(X^{\Diamond},\Lambda)$ of overconvergent sheaves. 
\end{enumerate}
\end{proposition}
\begin{proof}
Part (1) and (2) follow from \cite[Proposition 27.1, 27.2]{Ecod} (the case $k/\qp$ is not explicitly written, but the argument applies, see also Page 165 of \textit{loc. cit.}). 

Part (3) follows from \cite[Proposition 27.2, Proposition 27.7]{Ecod}, when the coefficient ring is finite. In the case $k/\qp$ and $\Lambda$ is as in Setup~\ref{assumption: coefficientsystemsingeneral} (2), we need to take an inverse limit of \cite[Proposition 27.7]{Ecod}.

For Part (4), the statements are local on $Y$, so we may assume $Y$ to be qcqs. Then it follows from \cite[Proposition 27.4, 27.5]{Ecod}. 

For part (5), if $k$ is over $\qp$, then the first isomorphism follows from \cite[Proposition 27.6]{Ecod} and an inverse limit of the cases of finite coefficient rings. If $k$ has characteristic $p$, one can argue similarly. Namely, by factoring $f$ into an open immersion followed by a proper map, one reduces to the case $f$ is an open immersion. By a further d\'evissage, one can reduce to checking that the relation 
\[f^{\Diamond}_\ast c_{X}^\ast \bb{F}_\ell\simeq c_{Y}^\ast f_\ast\bb{F}_\ell\]
holds on the constant sheaf $\bb{F}_\ell$ on $X$. In this case, we can conclude by using an excision sequence and applying Part (4) to the complement closed immersion (see also \cite[Lemma~3.3.2]{GHILZIsocComparison}).

For the second isomorphism in Part (5), we use the following argument as suggested on \cite[Page~128,Footnote~2]{FSGeomLLC}: We have a natural transformation 
\begin{equation}{\label{eqn: BaseChangeMapofUpper!}}
 c_{Y}^{*}f^{!} \ra f^{\Diamond!}c_{X}^{*},
\end{equation}
obtained by adjunction from the composition 
\[ f_{!}^{\Diamond}c_{X}^{*}f^{!} \simeq c_{Y}^{*}f_{!}f^{!} \ra c^{*}_{Y},\]
where the first map is the base-change map in (4). Localizing on $Y$ and using Part (2), we may assume $f$ factorizes as $X \ra \bb{A}^{n}_{Y} \ra Y$, where $X \ra \bb{A}^{n}_{Y}$ is a closed embedding. That Equation~\eqref{eqn: BaseChangeMapofUpper!} is an isomorphism for the second map is clear. This reduces us to the case where $f=i$ is a closed immersion. Since $i_\ast: \Detale(X,\Lambda) \ra \Detale(Y,\Lambda)$ is fully faithful and hence conservative, it suffices to check that (\ref{eqn: BaseChangeMapofUpper!}) is an isomorphism  after applying $i_{*}$. Now we write $j$ for the complementary open to $i$ and consider the excision sequence
\[ i_{*}i^{!} \ra \mathrm{id} \ra j_{*}j^{*}. \]
Using part (2) and (4), this reduces us to check that $j_{*}$ commutes with analytification, which follows from the first isomorphism in Part (5) we just established. 

The last claim about Verdier duality follows from the computation that 
\begin{align*}
    \bb{D}_{X^{\Diamond}/Y^\Diamond}c_{X}^\ast= \RHomint(c_X^\ast(-), f^{\Diamond!}\omega_{Y^\Diamond})
    &\simeq c_X^\ast\RHomint(-, c_{X,\ast}f^{\Diamond!}\omega_{Y^\Diamond})\\
    &\simeq c_X^\ast\RHomint(-, f^{!}c_{Y\ast}\omega_{Y^\Diamond})\simeq c_{X}^{*}\bb{D}_{X/Y},
\end{align*}
where we have used the adjunction $c_{X,\ast}\vdash c_X^\ast$ for the first isomorphism, Part (4) for the second isomorphism, and, for the last isomorphism, the identification $c_{Y\ast}\omega_{Y^\Diamond}=c_{Y\ast}c_Y^\ast\omega_{Y}\simeq\omega_Y$ from Part (3).

Part (6) is straightforward using partial properness. See for example \cite[Proposition~4.4.6]{GHILZIsocComparison} for the torsion case and then take inverse limits.
\end{proof}

We have the following results about how analytification interacts with the ULA subcategory. Below for a morphism $X\to S$, we write $\mathbb{D}_{X/S}$ (resp. $\mathbb{D}_{X^{\Diamond}/S^\Diamond}$) for Verdier duality on $\Detale(X,\Lambda)$ (resp. $\Detale(X^{\Diamond},\Lambda)$). 

\begin{proposition}{\label{prop: constructiblealgebrizestoULA}}
Let $S$ be a quasi-compact separated scheme over $k$. Let $f: X\to S$ be a separated scheme of finite presentation (when $\operatorname{char}k=p$, we assume both $S$ and $X$ are perfect and $f$ is perfectly of finite presentation), then 
\begin{enumerate}
\item The functor $c_{X}^{*}$ preserves the universally locally acyclic (over $S$) objects and hence induces a functor
\[ c_{X}^\ast:\Detale^{\mathrm{ULA}}(X,\Lambda) \ra \Detale^{\mathrm{ULA}}(X^{\Diamond},\Lambda).\] 
\item There is a natural identification of functors
\[c_{X}^{*}\mathbb{D}_{X/S} \simeq \mathbb{D}_{X^{\Diamond}/S^\Diamond}c_{X}^{*}: \Detale^{\mathrm{ULA}}(X,\Lambda) \ra \Detale^{\mathrm{ULA}}(X^{\Diamond},\Lambda)^{\mathrm{op}}.\]
\end{enumerate}
\end{proposition}

\begin{proof}
Part (2) is a special case of Proposition~\ref{prop: propertiesofthealgebraizationfunctor} (5), but in the case of ULA objects, there is an alternative argument. Namely, for both (1) and (2), by $\ell$-adic completeness, it reduces to checking the analogous statements for $\Lambda/\ell^n$ for each $n$, so we may assume $\Lambda$ to be finite. If $\operatorname{char}k=p$, the statement is \cite[Proposition~A.1]{AGLRLocalModels}. In general, we can argue exactly as in \textit{loc. cit.}. For convenience, we recall the argument here: Let $\mathcal{C}_{S}$ denote the 2-category of correspondences over $S$, whose objects are given by schemes as in the hypothesis and whose morphisms are given by objects in $\mathrm{Hom}_{\mathcal{C}_{S}}(X,Y) := \Detale(X \times_S Y,\Lambda)$. The composition law of morphisms is given by convolution, namely
\[ \Hom_{\mathcal{C}_{S}}(X,Y) \times \mathrm{Hom}_{\mathcal{C}_{S}}(Y,Z) \ra \Hom_{\mathcal{C}_{S}}(X,Z)\] 
\[ (A,B) \mapsto \pi_{XZ!}(\pi_{XY}^\ast(A) \otimes \pi_{YZ}^\ast(B)), \]
where the maps $\pi_{XZ}$ etc. are the obvious projection maps. Let $\mathcal{C}_{S^\Diamond}$ denote the analogous 2-category for diamonds which are fdcs over $S^\Diamond$, as in \cite[Section~IV.2.23]{FSGeomLLC}, see \cite[Definition 5.4]{Mann2022NuclearSheaves} for the definition of fdcs. Then Proposition~\ref{prop: propertiesofthealgebraizationfunctor} (2) and (4) above tells us that $c_{X}^{*}$ upgrades to a functor of 2-categories $c_{X}^{*}: \mathcal{C}_{S} \to \mathcal{C}_{S^\Diamond}$.

Now $A \in \Detale(X,\Lambda) =:\mathrm{Hom}_{\mc{C}_{S}}(X,S)$ is ULA with respect to the structure map $X \ra S$ if and only if $A$ admits a right adjoint in $\mc{C}_{S}$, which must necessarily be given by $\bb{D}_{X/S}(A)$, by \cite[Definition~3.2, Proposition~3.4(ii)]{HansenScholzeRelativePerversity}. Similarly, $c_X^\ast A \in \Detale(X^{\Diamond},\Lambda)$ is ULA with respect to $X^{\Diamond} \ra S^\Diamond$ if and only if $c_{X}^{*}A \in \Detale(X^{\Diamond},\Lambda) = \mathrm{Hom}_{\mc{C}_{S^\Diamond}}(X^{\Diamond},S^\Diamond)$ admits a right adjoint, given necessarily by $\bb{D}_{X^\Diamond/S^\Diamond}(c^\ast_XA)$, by \cite[Theorem~IV.2.23]{FSGeomLLC}. Since functors of 2-categories preserve adjunctions between the 1-morphisms, we deduce the preservation of the ULA property and Verdier duals.
\end{proof}
\begin{remark}
We note that in light of the identification $\Dconstf(X,\Lambda) \simeq \Detale^{\ULA}(X,\Lambda)$ discussed in Remark \ref{rem: PerfectConstructibleSameasULA}. Part (2) of Proposition \ref{prop: constructiblealgebrizestoULA} gives an alternative proof of Proposition \ref{prop: constructiblealgebrizestoULA} (4) on the perfect constructible subcategory.
\end{remark}

\begin{corollary}{\label{prop: Compatabilitieswithupper!s}}
If $f: X \ra Y$ is a smooth morphism of separated $k$-schemes of pure dimension $d$ (or the perfection thereof when $\operatorname{char}k=p$), then $f^{\Diamond}: X^{\Diamond} \ra Y^{\Diamond}$ is $\ell$-cohomologically smooth of pure $\ell$-dimension $d$.  
\end{corollary}
\begin{proof}
The statement is local on $Y$, so we can assume it is quasi-compact. It follows by standard results in \'etale cohomology of schemes that the map $f$ is $\ell$-cohomologically smooth of pure dimension $d$ in the six functor formalism $\mc{D}_\et(-,\Lambda)$. In particular, the constant sheaf $\Lambda$ on $X$ is ULA over $Y$ and the dualizing sheaf $f^!\Lambda$ is invertible, locally concentrated in degree $-2d$. By Proposition~\ref{prop: constructiblealgebrizestoULA} and Proposition~\ref{prop: propertiesofthealgebraizationfunctor}(5), this implies the same for $X^{\Diamond}$, which in turn implies that it is $\ell$-cohomologically smooth of pure dimension $d$ (cf. \cite[Proposition~ IV.2.33]{FSGeomLLC} and \cite[Theorem~3.3.12 (7)]{GHILZIsocComparison}).
\end{proof}
\subsubsection{Miscellany}
We have the following claim, which will be essential for computing the stalks of our intersection cohomology sheaf.
\begin{proposition}\label{prop: SmallDiamondvsPerfectScheme}
        Let $k$ be a perfect field of characteristic $p$, and $X$ be a perfect scheme separated over $k$. Then the functor $t_X^\ast c_X^\ast$ is fully faithful, where $t_X$ is the morphism defined below Equation~\eqref{eqn: NaturalTransformationsofAnalytifications}. In particular, for all $A \in \Detale(X,\Lambda)$ we have an isomorphism 
        \[ R\Gamma(X^{\diamond},t_{X}^{*}c_{X}^{*}(A)) \simeq R\Gamma(X,A) \]
        of $\Lambda$-modules.
    \end{proposition}
    \begin{proof}
        This is {\cite[Proposition 4.2]{GleasonTubularNeighborhoods}}. The second claim follows by expressing cohomology as external Hom from the constant sheaf.
    \end{proof}

\subsection{\'Etale sheaves attached to representations}\label{sec: EtaleSheavesRep}
In this section, we recall the construction of certain \'etale sheaves attached to smooth representations of certain topological groups considered in \cite[Section~1.10]{PinkHigherDirectImages}. We then compare it with an alternative construction involving classifying stacks under the analytification functor. This will allow us to compare the intersection complex of a local system attached to some algebraic representation of $\mathsf{G}$ on the Igusa stack to that on the Shimura variety considered by \cite{PinkHigherDirectImages}. The latter is constructible on the boundary stratification by the main results of \cite{PinkHigherDirectImages}, which gives rise to important finiteness properties for our intersection complex on the Igusa stack. Since this subsection is only of rather subtle technical importance, it can be skipped on a first reading. 

\subsubsection{Pink's construction}{\label{ss: PinksConstruction}}
Let $K$ be a profinite group, and $X$ be a scheme, together with a pro-\'etale Galois covering with covering group $K$, denoted $\phi: \widetilde{X} \ra X$. Given a ring $\Lambda$ as in Setup \ref{assumption: coefficientsystemsingeneral}(1), Pink~\cite{PinkHigherDirectImages} defines a functor on the bounded below derived category of smooth $K$-representations on $\Lambda$-modules
\[ \mu_{K}: \D^+(K,\Lambda) \ra \Detale^+(X,\Lambda). \]

We briefly recall his construction: Writing $K := \varprojlim_{i \in I} K/K_{i}$ as an inverse limit of finite groups, for a system of open normal subgroups $\{K_i\}_{i\in I}$, one has the corresponding \'etale $K/K_{i}$-covering 
\[ \phi_{i}: X_{i}:=\widetilde{X}\times^K (K/K_i) \ra X.\]
This allows us to define the \'etale sheaf
\[ \lambda_{\phi,\Lambda}(M) := \underline{M} \otimes \varinjlim_{i \in I} \phi_{i*}(\Lambda) \simeq \underline{M} \otimes \phi_{*}(\Lambda) \in \Detale^+(X,\Lambda),  \]
where the isomorphism follows from \cite[Tag~09Z1]{stacks-project}. As explained in \cite[Section~1.8]{PinkHigherDirectImages}, the object $\lambda_{\phi,\Lambda}(M)$ carries a continuous $K$-action, which is diagonal on the two tensor factors. Composed with the continuous group cohomology functor $R\Gamma(K,-)$ on sheaves with continuous $K$-action, one obtains the desired functor
\[\mu_K:= R\Gamma(K,-)\circ \lambda_{\phi,\Lambda}.\]
By taking left completions and $\ell$-adic completions, this induces a functor
\[ \mu_{K}: \widehat{\D}(K,\Lambda) \ra \Detale(X,\Lambda) \]
for $\Lambda$ any coefficient system in Setup \ref{assumption: coefficientsystemsingeneral}.

\subsubsection{Classifying stacks}
We now give an alternative description of the functor $\mu_{K}$ using classifying stacks, which will more easily adapt to the context in this paper. Let $\Lambda$ be as in Setup~\ref{assumption: coefficientsystemsingeneral}. 

Let us assume from now on that $X$ is a finitely presented scheme over an algebraically closed field $k$, which is either an extension of $\mathbb{Q}_{p}$ equipped with the non-archimedean topology or an algebraically closed field of characteristic $p$ equipped with the discrete topology. We use $\ast$ to denote $\Spd k$. Consider as before a pro-\'etale Galois covering of $k$-schemes $\phi: \widetilde{X} \ra X$ whose covering group is a profinite group $K$. Passing to the associated diamonds, the induced pro-\'etale $\underline{K}$-torsor $\phi^{\Diamond}: \widetilde{X}^{\Diamond} \ra X^{\Diamond}$ defines a map $\eta_{K}: X^\Diamond \ra [\ast /\underline{K}].$ We have a cartesian diagram 
\begin{equation}{\label{eqn: CartesianDiagramofTorsors}}
\begin{tikzcd}
\widetilde{X}^{\Diamond} \arrow[d] \arrow[r, "\phi^{\Diamond}"] & X^{\Diamond} \arrow[d,"\eta_{K}"] & \\
\ast \arrow[r] & \left[\ast /\underline{K}\right]. & 
\end{tikzcd}
\end{equation}

Note that by \cite[Lemma~10.5]{Mann2022NuclearSheaves} and \cite[Proposition 3.20]{Mann2022NuclearSheaves} and taking inverse limits in the setting of \ref{assumption: coefficients} (2), one has a natural equivalence of categories $\Detale([\ast/\underline{K}],\Lambda) \simeq \widehat{\D}(K,\Lambda)$ (in fact, for any locally profinite $K$). As a consequence, pullback defines a functor 
\[ \eta_{K}^{\ast}: \widehat{\D}(K,\Lambda) \ra \Detale(X^{\Diamond},\Lambda). \]
We would like to compare this with the analytification $\mu_{K}^{\Diamond} := c_{X}^{*}\mu_{K}$ of Pink's functor. 

\begin{lemma}\label{lemma: ComparisonbetweenPinksFunctorandClassifyingStack}
Let $K$, $\Lambda$, and $\phi: \widetilde{X} \ra X$ be as above. There is a natural equivalence $\mu_{K}^{\Diamond} \simeq \eta_{K}^{*}$, between functors from $\widehat{\D}(K,\Lambda)$ to $\Detale(X^\Diamond,\Lambda)$.
\end{lemma}
\begin{proof}
By left and $\ell$-adic completeness, we can reduce to the case $\Lambda$ is finite and just need to compare the two functors after restricting to the bounded below subcategories.

Observe that there is a commutative diagram of v-stacks
\begin{equation}\label{eq:pink commutative diagram} 
\begin{tikzcd}
X^{\Diamond} \times \left[\ast/\ul{K}\right] \ar[d,"\ol{\phi}",swap] & X^{\Diamond} \arrow[d,"{(\mathrm{id}_X,\eta_K)}"] \ar[l,"{(\mathrm{id}_X,\eta_K)}",swap]&  \\
 X^{\Diamond} \times \left[\ast/\ul{K}\right] \times \left[\ast/\ul{K}\right] \arrow[d,"p_3",swap] & X^{\Diamond} \times \left[\ast/\ul{K}\right] \arrow[d,"\gamma"]\ar[l,"{g}"] \\
 \left[\ast/\ul{K}\right] & X^{\Diamond}, 
\end{tikzcd} 
\end{equation}
where the square is cartesian; the map labeled by $\ol{\phi}$ is defined as $(\mathrm{id}_X,\eta_K)\times\mathrm{id}_{[\ast/\ul{K}]}$; the map $g$ is defined as $\mathrm{id}_X\times \Delta_{[\ast/\ul{K}]}$, where $\Delta_{[\ast/\ul{K}]}$ indicates the diagonal morphism on $[\ast/\ul{K}]$; the map $\gamma$ is induced by the structure morphism $[\ast/\ul{K}]\to \ast$; and the composition of the right vertical maps is the identity on $X$. 

We first note that the analytification of Pink's functor $\mu_K^\Diamond$ can be rewritten as
\[\mu_K^\Diamond = \gamma_\ast g^\ast(p_3^{*}(-) \otimes_\Lambda \ol{\phi}_\ast\Lambda). \]
Indeed, this follows from the fact that pushforward along $[\ast/\underline{K}] \ra \ast$ identifies with continuous group cohomology with respect to $K$, Proposition~\ref{prop: propertiesofthealgebraizationfunctor}, the fact that $\phi^{\Diamond}$ and $\phi_{i}^{\Diamond}$ are proper as maps of diamonds since $K$ is profinite, and \cite[Proposition 8.5(ii)]{Ecod}.

Now we can simplify this formula, using a chain of natural isomorphisms
\begin{align}{\label{align: naturaltransformation}}
\gamma_\ast g^\ast(p_3^{*}(-) \otimes_\Lambda \ol{\phi}_\ast\Lambda) \simeq \gamma_\ast g^\ast \ol{\phi}_\ast \ol{\phi}^\ast p_3^{*}(-)  \simeq \gamma_\ast (\mathrm{id}_X,\eta_K)_\ast(\mathrm{id}_X,\eta_K)^\ast \ol{\phi}^\ast p_3^{*}(-)\simeq \mathrm{id}_{X^\Diamond,\ast}\eta_K^\ast=\eta_K^\ast.
\end{align}
Here the first isomorphism uses the projection formula for $\ol{\phi}_\ast$, where, since $K$ is profinite, $(\mathrm{id}_{X} \times \eta_{K})$ is proper and so is $\ol{\phi}$. Hence the projection formula applies to $\ol{\phi}_\ast\simeq \ol{\phi}_!$. The second isomorphism is given by proper base-change along the cartesian square. Finally, the third isomorphism follows by the identification $p_3\circ \ol{\phi}\circ (\mathrm{id}_X,\eta_K)=\eta_K$. 
\end{proof}

\subsection{Perverse $t$-structures on v-stacks}\label{sec: tStructureVstacks}
For the purpose of defining intersection complexes on minimally compactified Igusa stacks and relating them to the intersection complexes on the corresponding Shimura varieties, we need to give some appropriate definition of perverse $t$-structures in the context of $v$-stacks. In fact, even for rigid spaces over an algebraically closed field, making sense of such a thing is a delicate issue and only works properly after restricting to Zariski constructible sheaves, see \cite{BhattHanSixFunctors}.

Nonetheless, for a stratified v-stack with cohomologically smooth boundary strata, one can often define a perverse $t$-structure in an ad hoc way, by gluing together the standard $t$-structures on the boundaries with the appropriate degree shifts. Such constructions already exist in various contexts, cf. \cite[Definition/Proposition~VI.7.1]{FSGeomLLC}. We recall this construction below, in the generality that is convenient for us.

\subsubsection{Perverse \texorpdfstring{$t$-structures}{}}\label{subsub: PerverseSetup}
Let $k$ be an algebraically closed field over $\mathbb{F}_p$, or an algebraically closed non-archimedean field over $\qp$. Let $X$ be a small v-stack over $\Spd k$ with a finite collection of locally closed sub-functors $\{i_\alpha:X_\alpha\hookrightarrow X\}_{\alpha\in I}$, $|I|<\infty$.  We assume they form a set theoretic stratification on the level of underlying topological spaces
\[|X|=\coprod_{\alpha\in I} |X_\alpha|,\]
so that
\[X_\alpha\simeq X\times_{\underline{|X|}}\underline{|X_\alpha|},\]
by combining \cite[Lemma~2.7]{AGLRLocalModels} and \cite[Proposition~12.9]{Ecod}. Here $X\to \underline{|X|}$ is the map of v-stacks as defined in \cite[Page 7, Equation (2.4)]{AGLRLocalModels}. We further assume that $X$ is an Artin stack over $\Spd k$. It follows that each $X_\alpha$ is an Artin v-stack (\cite[Example~IV.1.9(iii)]{FSGeomLLC}).  We assume that the strata $X_{\alpha}$ are $\ell$-cohomologically smooth of pure $\ell$-dimension $d_{\alpha}\in \bb{Z}$.

Let $\Lambda$ be any coefficient system as in Setup~\ref{assumption: coefficientsystemsingeneral}, and consider the sheaf category $\cDetale(-,\Lambda)$ and its homotopy category $\Detale(-,\Lambda)$ as introduced after Setup~\ref{assumption: coefficientsystemsingeneral}. We endow $I$ with the partial ordering topology such that the natural map $|X| \ra |I|$ sending each stratum $|X_{\alpha}|$ to $\alpha$ is continuous.
\begin{lemma}\label{lem: excision}
    Via excision, there is a semi-orthogonal decomposition of $\mc{D}_\et(X,\Lambda)$ in the sense of \cite[Definition~6.2.1]{GHILZIsocComparison} with respect to the poset $I$. The graded pieces of this semi-orthogonal decomposition are given by $\mc{D}_\et(X_\alpha,\Lambda)$s, for $\alpha\in I$.
\end{lemma}
\begin{proof}
    The v-stack $X$ is weakly-stratified in the sense of \cite[\S 6.3]{GHILZIsocComparison}. When $\Lambda$ is finite, this is {\cite[Proposition 6.3.1(1)]{GHILZIsocComparison}}. The case for $\ell$-adically complete $\Lambda$ follows from the finite case by using that the mod $\ell$-reduction functor is conservative and the projection formula to verify that the excision triangle in \cite[Definition~6.2.1 (3)]{GHILZIsocComparison} holds. 
\end{proof}

\begin{definition/proposition}{\label{defn: tstructuresonIgusaStacks}}
There exists a unique $t$-structure $(\pcDetale^{\leq 0},\pcDetale^{\geq 0})$
on $\cDetale(X,\Lambda)$ defined as follows:
\[\pcDetale^{\leq 0}\coloneq \{A \in \cDetale(X,\Lambda)\mid i_\alpha^\ast A \in \cDetale^{\leq -d_{\alpha}}(X_\alpha,\Lambda)\text{ for all $\alpha\in I$}\},\] 
\[\pcDetale^{\geq 0}\coloneq \{A \in \cDetale(X,\Lambda)\mid i_\alpha^!A \in \cDetale^{\geq -d_{\alpha}}(X_\alpha,\Lambda)\text{ for all $\alpha\in I$}\}.\] 
We write $\phantom{}^{p}\Detale^{\leq 0}(X,\Lambda)$ and $\phantom{}^{p}\Detale^{\geq 0}(X,\Lambda)$ for the resulting halves of the $t$-structure on the homotopy category $\Detale(X,\Lambda)$.
\end{definition/proposition}
\begin{proof}
This is similar to \cite[VI.7.1]{FSGeomLLC}. We first need to check that $\cDetale(X,\Lambda)$ is a presentable stable $\infty$-category, and the full subcategory $\pcDetale^{\leq 0}$ is presentable, closed under small colimits and extensions. The statement about $\cDetale(X,\Lambda)$ is \cite[Lemma 17.1]{Ecod}. That $\pcDetale^{\leq 0}$ is closed under small colimits and extensions is clear from the definition. It is presentable since, via excision (Lemma~\ref{lem: excision} above), it is generated by the categories $i_{\alpha,!}\cDetale^{\leq d_\alpha}(X_\alpha,\Lambda)$, each of which is presentable (since up to shift, it is part of the standard $t$-structure). It follows from \cite[Proposition 1.4.4.11(1)]{LurieHA}\footnote{Note that, in \textit{loc. cit.}, Lurie uses the homological indexing convention instead of the cohomological one, while we use the latter.} that there is a unique $t$-structure $(\pcDetale^{\leq 0},\pcDetale_1^{\geq 0})$ on $\cDetale(X,\Lambda)$. 

It remains to identify $\pcDetale_1^{\geq 0}$ with the more explicitly defined subcategory $\pcDetale^{\geq 0}$ above. For this we only need to check that $A\in \pcDetale^{\geq 0}$ if and only if $\Hom_{\cDetale(X,\Lambda)}(B,A)=0$ for all $B\in \pcDetale^{\leq -1}$, cf. \cite[Remark 1.2.1.3]{LurieHA}. But this follows immediately from excision and induction.
\end{proof}

\begin{remark}{\label{rem: dualityinterchangingtstructures}}
The assumption on cohomological smoothness of the strata ensures that 
\[ \bb{D}_{X}(\phantom{}^{p}\mc{D}^{\leq 0}) \subset \phantom{}^{p}\mc{D}^{\geq 0}. \]
Indeed, the naive duality on each stratum $X_\alpha$ sends the subcategory $\mc{D}^{\leq 0}$ to $\mc{D}^{\geq 0}$ in the standard $t$-structure, as can be checked stalkwise (The other direction is incorrect unless $\Lambda$ is self-injective!). Cohomological smoothness of $X_\alpha$ implies that the dualizing complex of $X_\alpha$ is invertible and is concentrated in degree $-2d_\alpha$. Therefore, by combining these facts, we see that Verdier duality on $X_\alpha$ sends the subcategory $\mc{D}^{\leq -d_\alpha}(X_\alpha,\Lambda)$ to $\mc{D}^{\geq -d_\alpha}(X_\alpha,\Lambda)$ (and vice-versa if $\Lambda$ is self-injective). The desired statement is then a direct consequence of the relation  $\bb{D}_{X_{\alpha}}i_{\alpha}^{*} = i_{\alpha}^{!}\bb{D}_{X}$. 

However, even if $\Lambda$ is self-injective, since the identity $\bb{D}_{X_{\alpha}}i_{\alpha}^{!} = i_{\alpha}^{*}\bb{D}_{X}$ does not always hold in general, we cannot generally conclude that
\[ \bb{D}_{X}(\phantom{}^{p}\mc{D}^{\geq 0}) \subset \phantom{}^{p}\mc{D}^{\leq 0}.\]
We will discuss this matter in more detail below.
\end{remark}

The resulting $t$-structure only depends on the homotopy category $\Detale(X,\Lambda)$. We will work mostly with this triangulated category instead of the infinity category. We write $\Perv(X,\Lambda)$ for the heart of the perverse $t$-structure on $\Detale(X,\Lambda)$. For all $n \in \mathbb{Z}$, this gives us perverse truncation functors denoted $(\phantom{}^{p}\tau^{\leq n},\phantom{}^{p}\tau^{\geq n})$, where the functor
\[ \phantom{}^{p}\tau^{\leq n}: \Detale(X,\Lambda) \ra \phantom{}^{p}\Detale^{\leq n}(X,\Lambda) \]
is defined to be the right adjoint to the natural inclusion $\phantom{}^{p}\Detale^{\leq n}(X,\Lambda) \hookrightarrow \Detale(X,\Lambda)$. Similarly, the functor $\phantom{}^{p}\tau^{\geq n}$ will be the left adjoint to the inclusion $\phantom{}^{p}\Detale^{\geq n}(X,\Lambda) \hookrightarrow \Detale(X,\Lambda)$.
As usual, by composing these functors, we obtain perverse cohomology functors. 
\[ \pH^{n}(-): \Detale(X,\Lambda) \ra \Perv(X,\Lambda).\]

\subsubsection{Basic properties} We now establish some standard lemmas on these perverse $t$-structures in the required level of generality. Assume we are given a decomposition of $X$ into a closed immersion and its open complement 
\[Z\xrightarrow{i} X \xleftarrow{j} U,\] 
where we assume that both $Z$ (closed) and $U$ are unions of strata. This allows us to analogously define perverse $t$-structures on $Z$ and $U$ in the same way as above. By postcomposing the adjoint functors $(i^{*},i_{*},i^{!})$ with the 0-th perverse cohomology, we obtain functors $\phantom{}^{p}i^{*}$, $\phantom{}^{p}i_\ast$, and $\phantom{}^{p}i^{!}$. Similarly, we obtain $\phantom{}^{p}j_{!}$, $\phantom{}^{p}j^\ast$, and $\phantom{}^{p}j_\ast$ from the adjoint functors $(j_!,j^\ast,j_\ast)$. They have the following basic interactions with the perverse $t$-structure, as in the classical setting. 
\begin{lemma}{\label{lemma: tstructureLemmasClosedImmersions}}
\begin{enumerate}
\item The functor $i_\ast$ is perverse $t$-exact, so that $\phantom{}^{p}i_\ast = i_\ast$. 
\item The functor $i^\ast$ is right perverse $t$-exact; in particular, $\phantom{}^{p}i^\ast =  \phantom{}^{p}\tau_{Z,\geq 0}i^\ast$.
\item The functor $i^{!}$ is left perverse $t$-exact; in particular, $\phantom{}^{p}i^{!} = \phantom{}^{p}\tau_{Z, \leq 0}i^{!}$
\end{enumerate}
\end{lemma}
\begin{proof}
Claim (1) follows from the fact that, since $i$ is proper, we have that $i_{!} \simeq i_{*}$ so $i_{!}$ commutes with both $*$ and $!$-pullback to the strata by proper base-change. Claim (2) (resp. Claim (3)) follows from the fact that $i^{*}$ (resp. $i^{!}$) commutes with $*$-pullback (resp. $!$-pullback).
\end{proof}
We also have analogous claim for the open immersion.
\begin{lemma}{\label{lemma: tstructureLemmasOpenImmersions}}
\begin{enumerate}
\item The functor $j^\ast$ is perverse $t$-exact, so that $\phantom{}^{p}j^\ast = j^\ast$. 
\item The functor $j_{!}$ is right perverse $t$-exact; in particular, $\phantom{}^{p}j_{!} =  \phantom{}^{p}\tau_{\geq 0}j_{!}$.
\item The functor $j_\ast$ is left perverse $t$-exact; in particular, $\phantom{}^{p}j_\ast = \phantom{}^{p}\tau_{\leq 0}j_\ast$.
\end{enumerate}
\end{lemma}
\begin{proof}
Claim (1) follows from the fact that $j^{!} \simeq j^\ast$, since $j$ is an open immersion. Claim (2) (resp. Claim (3)) follows from the fact that $j_{!}$ (resp. $j_\ast$) commutes with $*$-pullback (resp. $!$-pullback) to the stratum by proper base-change.
\end{proof}
This gives us the following consequence. 
\begin{corollary}{\label{cor: adjoints}}
There exists a sequence of adjunctions 
\[ \phantom{}^{p}i_\ast \dashv \phantom{}^{p}i_\ast = i_\ast \dashv \phantom{}^{p}i^{!} \]
and 
\[ \phantom{}^{p}j_{!} \dashv \phantom{}^{p}j^\ast = j^\ast \dashv \phantom{}^{p}j_\ast. \]
\end{corollary}
\begin{proof}
This follows from Lemmas \ref{lemma: tstructureLemmasClosedImmersions} and \ref{lemma: tstructureLemmasOpenImmersions} and the fact that the composite of adjoint pairs is again an adjoint pair (Recall that  $\phantom{}^{p}\tau_{\leq n}$ is a right adjoint and $\phantom{}^{p}\tau_{\geq n}$ is a left adjoint.). 
\end{proof}
This leads to the following properties of the perverse six functors:
\begin{lemma}{\label{lemma: KeyPropertiesofOpenandClosedTstructures}}
\begin{enumerate}
    \item $\phantom{}^pi^\ast: \mathrm{Perv}(X)\to \mathrm{Perv}(Z)$ is right exact;
    \item $\phantom{}^pi^!: \mathrm{Perv}(X)\to \mathrm{Perv}(Z)$ is left exact;
    \item $\phantom{}^pi^\ast \phantom{}^pj_!=0$;
    \item $\phantom{}^pi^! \phantom{}^pj_\ast=0$.
\end{enumerate}
\end{lemma}
\begin{proof}
(1) and (2) follow from the fact that $\phantom{}^pi^\ast$ is a left adjoint and that $\phantom{}^pi^!$ is a right adjoint, by Corollary \ref{cor: adjoints}. (3) and (4) follow from the fact that $j^\ast i_\ast =0$ and adjunction.
\end{proof}

\subsubsection{Intersection complex}\label{subsub: IC}
Let $X$, $\{X_\alpha\}_{\alpha\in I}$ be as in \S \ref{subsub: PerverseSetup} before. We assume in the stratification of $X$, there is a single open topologically dense stratum $j: U\hookrightarrow X$ of pure $\ell$-dimension $d_U$. They both have perverse $t$-structures defined as above. We can now define the intersection complexes for \'etale local systems on $U$, using the same construction as in the case of schemes. 

\begin{definition}{\label{defn: intermediateextension}}
The intermediate extension functor is 
    \[j_{!\ast}\coloneq \operatorname{Im}({}^pj_!\to {}^pj_\ast): \operatorname{Perv}(U,\Lambda)\to \operatorname{Perv}(X,\Lambda).\]
\end{definition}
Since $U$ is assumed to be $\ell$-cohomologically of pure $\ell$-dimension $d_U$, $\operatorname{Perv}(U,\Lambda)$ contains the abelian category of local systems on $U$, i.e. \'etale sheaves of $\Lambda$-modules on $U$ that become constant on an \'etale cover when $\Lambda$ is torsion, up to a shift by $d_U$. 

\begin{definition}\label{defn: IntersectionComplex}
For each $\mathcal{L} \in \Perv(U,\Lambda)$, we define its intersection complex to be 
\[\IC^{X_\ast}_X(\mc{L})\coloneq j_{!\ast} (\mc{L})\in \operatorname{Perv}(X,\Lambda).\]
\end{definition}
\begin{remark}
    The more usual convention is to define the intersection complex of $\mc{L}\in \Detale^{=0}(U,\Lambda)$ and insert a shift in the right-hand side of the above formula. We use the above normalization for convenience. We drop $X_\ast$ from notation unless we would like to stress on the stratification involved. 
\end{remark}

We have the following characterization of the intersection complex. 
\begin{lemma}\label{lemma: uniqueIC}
    For $\mc{L} \in \Perv(U,\Lambda)$, there exists a unique $\mathcal{F}\in \Perv(X,\Lambda)$ that extends $\mc{L}$ and such that, for each $\alpha \in I$ such that $X_{\alpha}\subset X\setminus U$, we have $i_{\alpha}^\ast\mathcal{F}\in \Detale^{\leq -d_{\alpha}-1}(X_{\alpha}, \Lambda)$ and $i_{\alpha}^!\mathcal{F}\in \Detale^{\geq -d_{\alpha}+1}(X_\alpha, \Lambda)$.
\end{lemma}
\begin{proof}
    We argue by induction on the number of strata. The statement is empty when there is a single stratum, so let us consider the case where there are at least two strata. 
    
    Take a stratum $Z$ that is minimal under the closure relations. Then $i: Z\hookrightarrow X$ is a closed immersion, whose open complement we denote by $j: U'\hookrightarrow X$. Then, by the induction hypothesis, there is a unique extension $\mathcal{F}_{U'}$ of $\mc{L}$ from $U$ to $U'$, such that, for all strata $X_\alpha$ in $U'$, the conditions on $i_{\alpha}^\ast \mathcal{F}_{U'}$ and $i_{\alpha}^! \mathcal{F}_{U'}$ are satisfied. Now let $\mathcal{F}$ be any complex extending $\mc{L}$ that satisfies the condition in the lemma. 
    Then we must have $j^\ast \mathcal{F}\cong \mathcal{F}_{U'}$ by uniqueness. Consider the exact triangle
    \[i_\ast i^!\mathcal{F}\to \mathcal{F}\to j_\ast \mathcal{F}_{U'}.\]
    Apply $i^\ast$ to it and shift by one, we get another exact triangle
    \[ i^\ast \mathcal{F}\to i^\ast j_\ast \mathcal{F}_{U'}\to i^! \mathcal{F}[1].\]
    Now applying the truncation functor $\tau_{\leq -d_Z-1}$ to this triangle, the condition that $i^\ast \mathcal{F}\in \Detale^{\leq -d_Z-1}(Z,\Lambda)$ and $i^! \mathcal{F}\in \Detale^{\geq -d_Z+1}(Z,\Lambda)$ implies that $i^\ast \mathcal{F}\simeq \tau_{\leq -d_Z-1}i^\ast j_\ast\mathcal{F}_{U'}$. The extension $\mathcal{F}$ is then uniquely determined by the gluing triple 
    \[(i^\ast \mathcal{F}\cong \tau_{\leq -d_Z-1}i^\ast j_\ast\mathcal{F}_{U'}, \mathcal{F}_{U'}, \tau_{\leq -d_Z-1}i^\ast j_\ast\mathcal{F}_{U'}\to i^\ast j_\ast\mathcal{F}_{U'}),\]
    where the last map is induced by the natural map $\tau_{\leq -d_Z-1}\to \mathrm{id}$. Using this gluing triple as a construction of $\mc F$ from $\mc{F}_{U'}$, we have shown the unique existence of $\mathcal{F}$. 
\end{proof}

Write $X$ as an increasing union of open sub-functors
\[X=U_n\supset U_{n-1}\supset \cdots \supset U_1\supset U_0 = U,\]
with inclusion maps $j_i: U_{i-1}\hookrightarrow U_i$. We assume that for each $i$, $U_i\backslash U_{i-1}$ is a union of strata of pure $\ell$-codimension $i$. Let $\mc{L} \in \Perv(U,\Lambda)$. 
\begin{corollary}[Deligne's formula]\label{cor: DeligneFormula}
When $\mathcal{L}$ is a local system (up to shift) on $U$, then we have:
\[\IC_X(\mc L)\simeq \tau_{\leq -d_U+n-1}j_{n,\ast}\cdots \tau_{\leq -d_U+1} j_{2,\ast}\tau_{\leq -d_U} j_{1,\ast}\mc{L}\]
\end{corollary}
\begin{proof}
    It suffices to check both sides satisfy the characterization of the unique extension as in Lemma~\ref{lemma: uniqueIC} above, which can be done inductively on the number of strata. This is clear for the right-hand side (called Deligne's complex), using that $\mathcal{L}$ is locally constant on $U$. The condition on $\ast$-pullbacks is clear and the condition on $!$-pullbacks may be checked by using the truncation triangles $\tau_{\leq d} \ra \mathrm{id} \ra \tau_{\geq d + 1}$ in the standard $t$-structure together with proper base-change. We explain the proof for $\IC_{X}(\mathcal{L})$.
    
    The case when there is a single stratum is clear. 
    Consider the inductive step. We assume that the claim holds on an open $j: U' \hookrightarrow X$ with complementary closed $i:Z\hookrightarrow X$, which are both assumed to be unions of strata. Set $\mathcal{L}':=\mathrm{IC}_{U'}(\mathcal{L})$. Then we have a surjection $\phantom{}^pj_{!}\mathcal{L}' \twoheadrightarrow \IC_X(\mathcal{L}')$, to which we can apply $\phantom{}^pi^\ast$. By Lemma \ref{lemma: KeyPropertiesofOpenandClosedTstructures} (1) and (3), this gives a surjection $0=\phantom{}^pi^\ast\phantom{}^pj_{!}\mathcal{L}' \twoheadrightarrow \phantom{}^pi^\ast\IC_X(\mathcal{L}')$, which implies $\phantom{}^pi^\ast\IC_{X}(\mathcal{L}') =0$. But by definition 
    \[\phantom{}^pi^\ast\IC_{X}(\mathcal{L}') = \phantom{}^p\mathcal{H}^0(i^\ast\IC_X(\mc{L}'))=\tau_{\geq -d_Z}\tau_{\leq -d_Z}i^\ast\IC_X(\mc{L}'). \]
    Since $\IC_X(\mc{L}')$ is perverse,  $i^\ast\IC_X(\mc{L}') \in \Detale^{\leq -d_Z}(Z,\Lambda)$. Hence, the sheaf $\phantom{}^pi^\ast\IC_{X}(\mc{L}')$ being zero forces 
    \[i^\ast\IC_X(\mc{L}') =\tau_{< -d_Z}i^\ast\IC_X(\mc{L}') \in  \Detale^{\leq -d_Z-1}(Z,\Lambda).\]
    For the claim about $i^!$, one can work dually, and apply $\phantom{}^pi^!$ to the injection $\IC_X(\mc{L}') \hookrightarrow \phantom{}^pj_{\ast}\mc{L}'$ and use the fact that $\phantom{}^pi^! \phantom{}^pj_\ast=0$. Here we use Lemma \ref{lemma: KeyPropertiesofOpenandClosedTstructures} (2) and (4).
\end{proof}

\subsubsection{Verdier duality}
Let $\Lambda$ be as in Setup~\ref{assumption: coefficientsystemsingeneral} as before. We now discuss the behavior of the intersection complex with respect to Verdier duality, which will depend on the reflexivity of the IC complex and its restriction to strata. In the context of constructible sheaves on finite type schemes, reflexivity of the intersection complex follows from the fact that Verdier duality is a biduality on the whole category of bounded constructible sheaves, see \cite[Th\'eor\`eme 6.1.1]{TravauxdeGabberXVII}. As we are interested in studying intersection complexes where the sheaves can have infinite dimensional stalks, we will need to identify some analogous finiteness condition in order to get self-duality. 

Recall that, if $f:X \ra\Spd k$ is an Artin $v$-stack, there is a natural (Verdier) biduality map 
\begin{equation}{\label{eqn: DoubleDual}}
 A \ra \bb{D}_{X}\bb{D}_{X}(A), 
\end{equation}
which is given by applying Hom-tensor adjunction to the evaluation pairing $A \otimes \bb{D}_{X}(A) \ra \omega_X:=Rf^{!}\Lambda$. We say an object $A\in \Detale(X,\Lambda)$ is reflexive (over $\Spd k$) if the map~\eqref{eqn: DoubleDual} is an isomorphism. 

\begin{definition}{\label{defn: stratifiedreflexive}}
    An object $A\in \Detale(X,\Lambda)$ is \emph{stratified reflexive} with respect to $X\to \Spd k$ and the stratification $\{X_\alpha\}_{\alpha\in I}$, if it is reflexive over $\Spd k$, and, for all $\alpha\in I$, the sheaf $i_\alpha^\ast A \in \Detale(X_\alpha,\Lambda)$ is reflexive with respect to $X_\alpha\to \Spd k$.
\end{definition}
This is guaranteed by the following stronger condition.
\begin{definition}\label{defn: StratifiedULA}
    An object $A\in \Detale(X,\Lambda)$ is \emph{stratified universally locally acyclic} (stratified ULA) with respect to $X\to \Spd k$ and the stratification $\{X_\alpha\}_{\alpha\in I}$, if it is itself ULA over $\Spd k$, and, for any $\alpha\in I$, the sheaf $i_\alpha^\ast A \in \Detale(X_\alpha,\Lambda)$ is ULA with respect to $X_\alpha\to \Spd k$, in the sense of \cite[Definition IV.2.31]{FSGeomLLC}.
\end{definition}
\begin{remark}\label{remark: ULAReflexive}
    By \cite[Corollary IV.2.25]{FSGeomLLC} and \cite[Proposition IV.2.15]{FSGeomLLC}\footnote{We need \cite[Proposition IV.2.15]{FSGeomLLC} to descend \cite[Corollary IV.2.25]{FSGeomLLC} to a version that works for relative Verdier duality with respect to maps between Artin v-stacks, not necessarily representable. Also, although the statement in \textit{loc. cit.} is given for torsion $\Lambda$, it works for a general six-functor formalism.}, if $A\in \Detale(X,\Lambda)$ is stratified ULA, then $A$ (resp. for all $\alpha\in I$, $i_\alpha^\ast A$) is reflexive with respect to the Verdier duality $\mathbb{D}_X=\mathbb{D}_{X/k}$ (resp. $\mathbb{D}_{X_\alpha}$).
\end{remark}

\begin{construction}\label{cons: DualityMap}
There is a natural transformation between functors from the subcategory of $\Perv(U,\Lambda)$ consisting of objects whose Verdier dual are again perverse, to $\Detale(X,\Lambda)^\mathrm{op}$
\begin{equation}\label{eq: ComparisonIC}
j_{!\ast}\bb{D}_{U}(-)\to \mathbb{D}_{X}j_{!\ast}(-)    
\end{equation}
induced by the following commutative diagram of natural transformations
\[
\begin{tikzcd}
    {}^pj_!\bb{D}_U \ar[r, two heads]\ar[dr] & j_{!\ast}\bb{D}_U \ar[d, dashed]\ar[r,hook] & {}^pj_\ast \bb{D}_U \ar[d]\\
    & \phantom{}^p\mathcal{H}^0(\bb{D}_Xj_{!\ast}) \ar[r,hook] & {}^p\mc{H}^0(\bb{D}_X{}^pj_!),
\end{tikzcd}
\]
where the bottom horizontal map is obtained by applying Verdier duality to the map
\[\phantom{}^pj_!\twoheadrightarrow  j_{!\ast} \]
and taking the 0-th perverse cohomology. Note that this is a pointwise injective natural transformation, since, if we write $\mc{K}$ for the kernel of the surjection $\phantom{}^pj_!\twoheadrightarrow j_{!\ast}$ (which is perverse) applied to some object $A$, then applying Verdier duality, we have a distinguished triangle
\begin{equation}\label{eq: ExactTriangleDIC}
    \bb{D}_Xj_{!\ast}A \to \bb{D}_X\phantom{}^pj_!A \to \bb{D}_X\mc{K}.
\end{equation}
However, $\bb{D}_X$ sends $\phantom{}^p\D^{\leq 0}$ to $\phantom{}^p\D^{\geq 0}$ (see Remark~\ref{rem: dualityinterchangingtstructures}), so $\bb{D}_X \mc{K}$ has no perverse cohomology in negative degrees. Therefore, the (perverse) cohomology long exact sequence induced by \eqref{eq: ExactTriangleDIC} gives
\[0\to \phantom{}^p\mc{H}^0(\bb{D}_Xj_{!\ast})\to \phantom{}^p\mc{H}^0(\bb{D}_X\phantom{}^pj_!).\]

The left diagonal arrow is obtained by appealing to the adjunction $\phantom{}^pj^\ast=j^\ast\vdash \phantom{}^pj_!$ (Corollary \ref{cor: adjoints}) and the identification
\begin{equation}
  \bb{D}_U=\bb{D}_U j^\ast j_{!\ast}\simeq j^\ast \bb{D}_X j_{!\ast},
\end{equation}
followed by taking perverse cohomology. The rightmost vertical arrow is obtained by considering the composition
\[\phantom{}^pj_\ast\bb{D}_U=\phantom{}^p\tau_{\leq 0}j_\ast\bb{D}_U\simeq \phantom{}^p\tau_{\leq 0}\bb{D}_Xj_!\to \phantom{}^p\tau_{\leq 0}\bb{D}_X\phantom{}^p\tau_{\geq 0}j_!=\phantom{}^p\tau_{\leq 0}\bb{D}_X\phantom{}^pj_!\to \bb{D}_X\phantom{}^pj_!\]
and taking perverse cohomology. The diagram with solid arrows is commutative, since going both ways from $\phantom{}^pj_!\bb{D}_U$ to $\phantom{}^p\mc{H}^0(\bb{D}_X \phantom{}^pj_!)$, this is the map adjoint to the identity on $\bb{D}_U$ via the isomorphism
\[j^\ast\phantom{}^p\mc{H}^0(\bb{D}_X \phantom{}^pj_!)\simeq \phantom{}^p\mc{H}^0(j^\ast\bb{D}_X \phantom{}^pj_!)\simeq \phantom{}^p\mc{H}^0(\bb{D}_Uj^\ast \phantom{}^pj_!)=\bb{D}_U.\]

Now the dashed arrow exists by commutativity of the outer square and the injectivity of the lower horizontal map. Namely, one can consider the following composition of morphisms
\[j_{!\ast}\mathbb{D}_U \to \operatorname{Im}(\phantom{}^pj_!\bb{D}_U\to \phantom{}^p\mc{H}^0(\bb{D}_X\phantom{}^pj_!))\to \phantom{}^p\mathcal{H}^0(\bb{D}_Xj_{!\ast}),\]
where both maps use the universal property of images. This induces the desired comparison map by taking composition with the natural map 
\[\phantom{}^p\mc{H}^0(\bb{D}_Xj_{!\ast})\to \bb{D}_Xj_{!\ast},\]
where we have used again the fact that $\bb{D}_X$ sends $\phantom{}^p\D^{\leq 0}$ to $\phantom{}^p\D^{\geq 0}$, and hence $\bb{D}_Xj_{!\ast}$ has no perverse cohomology in negative degrees.
\end{construction}


We now have the following key proposition. 
\begin{proposition}\label{prop: VerdierDualIC}
Assume $\Lambda$ is as in Setup~\ref{assumption: coefficientsystemsingeneral} and is self-injective. Let $\mathcal{L} \in \Perv(U,\Lambda)$. Suppose  
\begin{enumerate}
    \item Both the intersection complex $\IC_X(\mc{L})$ and its Verdier dual $\mathbb{D}_{X}\IC_X(\mc{L})$ are stratified reflexive;
    \item $\mc{L}$ has a finite filtration whose graded pieces are Schur irreducible objects (in the sense that the endomorphism ring is just $\Lambda$), which all satisfy (1).
\end{enumerate}
Then the comparison map \eqref{eq: ComparisonIC} induces an isomorphism
\[\IC_{X}(\bb{D}_{U}\mathcal{L})\simeq \mathbb{D}_{X}\IC_X(\mc{L}),\] 
where the right-hand side is well-defined, since $\bb{D}_U\mc{L}$ is perverse by the cohomological smoothness of $U$ and Remark \ref{rem: dualityinterchangingtstructures}.
\end{proposition}

\begin{remark}
    The condition for the complex $\mathbb{D}_X\IC_X(\mc{L})$ to be stratified reflexive is equivalent to asking that each $i_\alpha^! \IC_X(\mc{L}) \in \Detale(X_\alpha,\Lambda)$ is reflexive, for all $\alpha\in I$.
\end{remark}

\begin{lemma}\label{lemma: relationVerdierDuality}
Let $\mc{L} \in \Perv(U,\Lambda)$ such that both $\IC_X(\mc{L})$ and $\mathbb{D}_X\IC_X(\mc{L})$ are stratified reflexive. Then, for each $\alpha\in I$, there is an isomorphism
    \[i_\alpha^!\mathbb{D}_X\IC_X(\mathcal{L}) \simeq \mathbb{D}_{X_\alpha} i_\alpha^\ast\IC_X(\mathcal{L}), \,i_\alpha^\ast\mathbb{D}_X\IC_{X}(\mathcal{L}) \simeq \mathbb{D}_{X_\alpha}i_\alpha^!\IC_{X}(\mathcal{L}).\]
\end{lemma}
\begin{proof}
    The first isomorphism is an easy application of Yoneda's lemma and projection formula. 
    For the second relation, we use the first relation and the reflexivity of $\IC_X(\mathcal{L})$ and $i_\alpha^\ast\mathbb{D}_X\IC_X(\mathcal{L})$. Namely, we have     
    \[\mathbb{D}_{X_\alpha}i_\alpha^! \IC_X(\mathcal{L}) \simeq\mathbb{D}_{X_\alpha}i_\alpha^!\mathbb{D}_X \mathbb{D}_X \IC_X(\mathcal{L}) \simeq\mathbb{D}_{X_\alpha} \mathbb{D}_{X_\alpha} i_\alpha^\ast\mathbb{D}_X  \IC_X(\mathcal{L}) \simeq i_\alpha^\ast \mathbb{D}_X\IC_X(\mathcal{L}).\]
\end{proof}

\begin{proof}[Proof of Proposition~\ref{prop: VerdierDualIC}]
We first show that $\mathbb{D}_{X}\IC_{X}(\mathcal{L})$ is isomorphic to $\IC_{X}(\bb{D}_{U}\mathcal{L})$ as a complex (not necessarily induced by the prescribed map). For this, it suffices to check that both $\IC_{X}(\bb{D}_{U}\mathcal{L})$ and the Verdier dual $\mathbb{D}_{X}\IC_{X}(\mathcal{L})$ satisfy the unique characterization in Lemma~\ref{lemma: uniqueIC}. The former follows from the proof of Corollary~\ref{cor: DeligneFormula}. The latter follows from the result for $\IC_X(\mathcal{L})$ and Lemma~\ref{lemma: relationVerdierDuality}. Here we have used the $\ell$-cohomological smoothness of the strata $X_{\alpha}$ and the self-injectivity of the coefficient system $\Lambda$, see Remark \ref{rem: dualityinterchangingtstructures}. We fix an isomorphism $\alpha$ between the two complexes.

Now in order to show the map in Equation~\eqref{eq: ComparisonIC} is an isomorphism when applied to $\mc{L}$, we first assume $\mc{L}$ is Schur irreducible. Let us compute the endomorphisms of $\IC_{X}(\bb{D}_{U}\mathcal{L})$. We have a short exact sequence
\[0\to \IC_{X}(\bb{D}_{U}\mathcal{L})\to {}^pj_\ast \bb{D}_U\mc{L}\to \mc{Q}\to 0,\]
for some perverse sheaf $\mc{Q}$. Applying $\operatorname{RHom}(\IC_{X}(\bb{D}_{U}\mathcal{L}),-)$ to it, we obtain a long exact sequence
\[0\to \operatorname{End}_{X}(\IC_{X}(\bb{D}_{U}\mathcal{L}))\to \operatorname{Hom}_{X}(\IC_{X}(\bb{D}_{U}\mathcal{L}),{}^pj_\ast\bb{D}_U\mc{L})\to \operatorname{Hom}_{X}(\IC_{X}(\bb{D}_{U}\mathcal{L}),\mc{Q}) \ra.\]
By adjunction, we have 
\[\operatorname{Hom}_X(\IC_{X}(\bb{D}_{U}\mathcal{L}),{}^pj_\ast\bb{D}_U\mc{L})\simeq \operatorname{Hom}_U(j^\ast\IC_{X}(\bb{D}_{U}\mathcal{L}),\bb{D}_U\mc{L})\simeq \operatorname{End}_U(\bb{D}_{U}\mathcal{L})\simeq \Lambda,\]
where the last isomorphism uses our assumption that $\mc{L}$ (and hence $\bb{D}_U \mc{L}$) is Schur irreducible. Since $\operatorname{End}(\IC_{X}(\bb{D}_{U}\mathcal{L}))$ is a sub-$\Lambda$-module of the above group, this implies that endomorphisms of $\IC_{X}(\bb{D}_{U}\mathcal{L})$ are constants. By composing with the isomorphism $\alpha$, we see that morphisms from $\IC_{X}(\bb{D}_{U}\mathcal{L})$ to $\bb{D}_X\IC_{X}(\mathcal{L})$ are also constants. In particular, one can check whether a map between them is an isomorphism by restricting to $U$. The desired result follows, therefore, from the fact that \eqref{eq: ComparisonIC} restricts to the identity on $U$.

For a general $\mc{L}$, write it as an iterated extension of Schur irreducible objects. We then apply the argument above and use induction together with \cite[Tag~014A]{stacks-project}.
\end{proof}

\subsubsection{Comparison with the algebraic IC-sheaf}{\label{ss: ComparisonWithAlgebraicIC}}

We apply the analytification map from \S \ref{sec: analytification} to compare the perverse $t$-structure described above with algebraic ones, for v-sheaves attached to schemes. 

More precisely, let $k$ be either an algebraically closed field over $\bb{F}_p$, or an algebraically closed non-archimedean field over $\qp$. Suppose $X$ is a separated finite type scheme over $k$ equipped with a stratification $\{X_{\alpha}\}$ such that each $X_{\alpha}$ is smooth of pure dimension $d_{\alpha}$, and that there exists an open dense $U \subset X$ in the stratification. We assume that $\Lambda$ is as in Setup~\ref{assumption: coefficientsystemsingeneral}. 
Consider $\Dbc(X,\Lambda)$ its bounded derived category of constructible \'etale $\Lambda$-sheaves.

We can then define a perverse $t$-structure $(\phantom{}^{p}\D^{X_\ast, \geq 0}, \phantom{}^{p}\D^{X_\ast, \leq 0})$ on $\Dbc(X,\Lambda)$ by insisting that an object $A$ satisfies $i_\alpha^\ast A \in \D^{\leq -d_{\alpha}}(X_\alpha,\Lambda)$ and $i_\alpha^{!}A \in \D^{\geq -d_\alpha}(X_\alpha,\Lambda)$, just as in Definition/Proposition \ref{defn: tstructuresonIgusaStacks}. We write $j^{X_\ast}_{!*}: \Dbc(U,\Lambda) \ra \Dbc(X,\Lambda)$ for intermediate extension operation attached to this $t$-structure. 

We also recall that there is a more intrinsic perverse $t$-structure $(\phantom{}^{p}\D^{\leq 0},\phantom{}^{p}\D^{\geq 0})$ for the middle perversity on $\Dbc(X,\Lambda)$, defined without the need to choose any stratification in $X$. Indeed, for a point $x\in X$, let $i_x:x\hookrightarrow X$ denote the  the inclusion. Then one can define  
\[
\phantom{}^{p}\D^{\leq 0} = \{A\in \Dbc(X,\Lambda)\mid H^j(i_x^*A) = 0\text{ for all }j>-\mathrm{dim}\overline{\{x\}}\text{ and for all }x\in X \},
\]
\[
\phantom{}^{p}\D^{\geq 0} = \{A\in \Dbc(X,\Lambda)\mid H^j(i_x^{!}A) = 0\text{ for all }j<-\mathrm{dim}\overline{\{x\}}\text{ and for all }x\in X \}. 
\]
This is the perverse $t$-structure studied in \cite[\S 2.2]{BBD}.

\begin{remark}{\label{rem: perverseconstructibletstructure}}
If we assume $\Lambda$ is cohomologically regular (in the sense that the category of perfect complexes of $\Lambda$-modules is preserved under standard truncation), then the intrinsic perverse $t$-structure restricts to a perverse $t$-structure $(\phantom{}^{p}\Dconstf^{\leq 0},\phantom{}^{p}\Dconstf^{\geq 0})$ on $\Dconstf(X,\Lambda)$, the full subcategory of perfect  constructible sheaves (e.g by \cite[Theorem~1.9]{HansenScholzeRelativePerversity} and the fact that perfect constructible sheaves coincide with the ULA objects, as explained in \S\ref{ss: VariousSheafCategories}). In the setup~\ref{assumption: coefficientsystemsingeneral}, being cohomologically regular is guaranteed by being regular.
\end{remark}

For the intrinsic perverse $t$-structure and the open immersion $j:U\hookrightarrow X$, we also have the intermediate extension operation $j_{!\ast}$. One can therefore attach an intersection complex to any $\mc{L}\in \mathrm{D}^{\mathrm{b}, =-d_U}_\mathrm{c}(U,\Lambda)$ by
\[ \IC_{X}(\mathcal{L}) := j_{!*}(\mathcal{L}).\]

The lemma below compares it with $\IC^{X_\ast}_X(\mc{L}):= j_{!\ast}^{X_\ast}(\mc{L})$, constructed from the stratification $X_\ast$, as in the previous section.
\begin{lemma}\label{lemma: NaturalandStratifiedIC}
Suppose that $\mathcal{L}$ is lisse (i.e a shift by $d_{U}$ of an \'etale local system in Setup \ref{assumption: coefficientsystemsingeneral} (1) or an inverse system of such sheaves in Setup \ref{assumption: coefficientsystemsingeneral} (2)). If $j_{*}\mathcal{L}$ is constructible with respect to the stratification $X_{*}$, then we have
    \[\IC_{X}(\mathcal{L})\simeq \IC^{X_\ast}_X(\mc{L}).\]
    In particular, $\IC_{X}(\mathcal{L})$ is constructible with respect to $X_{*}$.
\end{lemma} 
\begin{proof}
    By \cite[Proposition~2.2.2]{BBD}, the two perverse $t$-structures $(\phantom{}^{p}\D^{\leq 0,X_{*}},  \phantom{}^{p}\D^{\geq 0,X_{*}})$ and $(\phantom{}^{p}\D^{\leq 0},  \phantom{}^{p}\D^{\geq 0})$ coincide after restricting to the full subcategory of $ \Dbc(X,\Lambda)$ with objects whose cohomology sheaves are constructible with respect to the stratification $X_{*}$. In particular, if we have an  $\mathcal{L}$ such that $j_{*}\mathcal{L}$ is constructible with respect to $X_{*}$, then it follows that the perverse truncations of $j_{!}\mathcal{L}$ and $j_*\mathcal{L}$ under these two $t$-structures agree.     
\end{proof}

Now, applying the functor $X\mapsto X^\Diamond$, we have an associated stratification $X_\ast^{\Diamond} := \{X^{\Diamond}_{\alpha}\}_{\alpha\in I}$ of $X^{\Diamond}$, and each $X_\alpha^\Diamond$ is $\ell$-cohomologically smooth of pure $\ell$-dimension $d_\alpha$ by Corollary~\ref{prop: Compatabilitieswithupper!s}. Then applying Definition/Proposition \ref{defn: tstructuresonIgusaStacks} gives a $t$-structure $(\phantom{}^{p}\Detale^{X^\Diamond_\ast,\leq 0},\phantom{}^{p}\Detale^{X^\Diamond_\ast,\geq 0})$ on $\Detale(X^{\Diamond},\Lambda)$, where we have implicitly invoked Lemma \ref{lem: ExcisionBigDiamond} to see the analytified stratification satisfies the assumptions of the previous section. We write $j_{!*}^{X_\ast,\Diamond}$ for the intermediate extension operation with respect to this $t$-structure, along the map $U^\Diamond\to X^\Diamond$. 
\begin{remark}{\label{rem: stratificationsanalytify}}
Note that the stratification $X_\ast^{\Diamond}$ is a Zariski stratification, which we warn the reader is very different from being constructible with respect to the analytic topology on $X^{\Diamond}$, as defined in \cite[Section~20]{Ecod}, due to the fact that analytification functor $(-)^{\Diamond}$ does not preserve quasi-compactness. The analytification functors $(-)^{\diamond}$ and $(-)^{\dagger}$ of \S \ref{ss: analytificationofperfectschemes} do actually preserve quasi-compactness \cite[Lemma~4.4.12]{GHILZIsocComparison}. However, the $(-)^{\Diamond}$ functor has the important property that, given an open-closed decomposition $U=X \setminus Z$, we obtain $|U^{\Diamond}| = |X^{\Diamond}| \setminus |Z^{\Diamond}|$ by Lemma \ref{lem: ExcisionBigDiamond}. Therefore, a locally closed stratification analytifies to a locally closed stratification. 
\end{remark}

We have the following claim. 
\begin{lemma}{\label{lemma: intermediateextensioncompatiblewithanalytification}}
For the perverse $t$-structures defined by the stratifications $X_\ast$ and $X_\ast^\Diamond$ above, the analytification functor 
\[ c_{X}^\ast: \Dbc(X,\Lambda) \ra \Detale(X^\Diamond,\Lambda) \]
is perverse $t$-exact. Moreover, we have a natural isomorphism  
\[ c_{X}^{*}j_{!*}^{X_{*}}(-) \simeq j_{!*}^{X_{*},\Diamond}c_{U}^{*}(-) \]
of functors 
\[ \Perv(U,\Lambda) \ra \Detale(X^\Diamond,\Lambda). \]
\end{lemma}

\begin{proof}
The first part of the claim easily follows from combining Corollary~\ref{prop: Compatabilitieswithupper!s}, the natural isomorphisms given by Proposition \ref{prop: propertiesofthealgebraizationfunctor}
\[ i_\alpha^{\Diamond*}c_{X}^{*} \simeq c_{X_\alpha}^{*}i_\alpha^{*} \]
\[ i_\alpha^{\Diamond!}c_{X}^{*} \simeq c_{X_\alpha}^{*}i_\alpha^{!} \]
between functors 
\[ \Dcons(X,\Lambda) \ra \Detale(X_\alpha^{\Diamond},\Lambda), \]
for all $\alpha\in I$, and the fact that analytification is exact for the standard $t$-structure, since it comes from a pullback of sites (see \S \ref{ss: AnalytificationinMixedCharacteristic}). This in particular implies that $c_{X}^{*}$ commutes with the perverse truncation functors and therefore reduces the second part of the claim to the natural isomorphisms of functors on $\Dcons(U,\Lambda)$
\[ c_{X}^{*}j_{!} \simeq j_{!}^{\Diamond}c_{U}^{*} \]
\[ c_{X}^{*}j_{*} \simeq j_{*}^{\Diamond}c_{U}^{*},\]
which again follows from Proposition \ref{prop: propertiesofthealgebraizationfunctor}.
\end{proof}

\begin{corollary}{\label{cor: diamonddescriptionoftheICsheaf}}
Let $\mathcal{L}\in \Perv(U,\Lambda)$ be lisse and assume that $j_{*}\mathcal{L}$ is constructible with respect to $X_{*}$, then the following statements are true.
\begin{enumerate}
\item There is an isomorphism in $\Detale(X^{\Diamond},\Lambda)$
\[ c_{X}^{*}\IC_{X}(\mathcal{L}) \simeq j_{!*}^{X_{*},\Diamond}c_{U}^{*}(\mathcal{L})=: \IC_{X^\Diamond}^{X_\ast}(c_U^\ast \mc{L}). \]
\item If $\Lambda$ is self-injective, then the complexes $c_{X}^{*}\IC_{X}(\mathcal{L})$ and $\bb{D}_{X^{\Diamond}}c_{X}^{*}\IC_{X}(\mathcal{L})$ are stratified reflexive with respect to $X^\Diamond\to \Spd k$ and $X_\ast^{\Diamond}$, in the sense of Definition~\ref{defn: stratifiedreflexive}. 
\item If $\Lambda$ is regular, then the complexes $c_{X}^{*}\IC_{X}(\mathcal{L})$ and $\bb{D}_{X^{\Diamond}}c_{X}^{*}\IC_{X}(\mathcal{L})$ are stratified ULA with respect to $X^\Diamond\to \Spd k$ and $X_\ast^{\Diamond}$, in the sense of Definition~\ref{defn: StratifiedULA}. 
\end{enumerate}
\end{corollary}
\begin{proof}
Part (1) follows from Lemma~\ref{lemma: NaturalandStratifiedIC} and Lemma~\ref{lemma: intermediateextensioncompatiblewithanalytification}. 

Part (2) then follows from (1) and Proposition \ref{prop: propertiesofthealgebraizationfunctor} (5), together with the fact that Verdier duality is an anti-involution on $\Dbc(X,\Lambda)$, see \cite[Th\'eor\`eme XVII.6.1.1]{TravauxdeGabberXVII}. 

Part (3) follows from (1), where we note that, by the regularity assumption on $\Lambda$ the perverse $t$-structure restricts to a perverse $t$-structure on the perverse constructible subcategory, as explained in Remark \ref{rem: perverseconstructibletstructure}. In particular, both $\IC_{X}(\mathcal{L})$ and its Verdier dual are perfect constructible in this case. Therefore, the claim follows from Proposition~\ref{prop: propertiesofthealgebraizationfunctor} (3), Proposition \ref{prop: constructiblealgebrizestoULA} and the identification of $\Dconstf(X,\Lambda)$ with $\Detale^\mathrm{ULA}(X,\Lambda)$ as in Remark \ref{rem: PerfectConstructibleSameasULA}.

\end{proof}

\subsection{A criterion for universal local acyclicity}
In this section, we establish a criterion for universal local acyclicity on stacky quotients of diamonds by locally pro-$p$ groups, which will ultimately be applied to the pullback of our intersection cohomology sheaf on the Igusa stack to a cohomologically smooth atlas given by a stacky quotient of the minimal compactification of the Shimura variety. We take $\Lambda$ to be a coefficient system as in Setup~\ref{assumption: coefficientsystemsingeneral}. We start with some lemmas. 
\begin{lemma}{\label{lemma: basicCohomologylemma}}
If $K$ is a profinite group whose pro-order is coprime to $\ell$ and $\Lambda$ is torsion, then the functor
\[ R\Gamma(K,-): \D(K,\Lambda) \ra \D(\Lambda) \]
defined by taking continuous $K$-cohomology has trivial higher cohomology.
\end{lemma}
\begin{proof}
    This is \cite[Corollary~III.3.3.7]{CohomologyofNumberFields}.
\end{proof}

This translates into the following lemma in terms of six functors on classifying stacks. We consider the following setup. Let $X \ra \Spd k$ be a fdcs map as defined in \cite[Definition~5.4]{Mann2022NuclearSheaves} for $k$ an algebraically closed field equipped with the discrete topology or over $\bb{Q}_{p}$  with the non-archimedean topology. Let $H$ be a locally pro-$p$ group of locally finite $\ell$-cohomological dimension. For a compact open subgroup $K \subset H$, we form the diagram 
\[ 
\begin{tikzcd}
    \left[X/H\right] & \left[X/K\right]\ar[l,"{f_{K}}",swap] \ar[r,"{q_{K}}", left] & X.
\end{tikzcd} \]
of $v$-stacks, where the quotients are formed with respect to the trivial action, and the maps are the natural ones. 
\begin{lemma}{\label{lemma: classifiyingstackmapproperties}}
If $K$ has finite $\ell$-cohomological dimension (e.g it is pro-$p$, by Lemma \ref{lemma: basicCohomologylemma}) then the map $q_{K}$ is $\ell$-cohomologically proper in the sense of \cite[Definition~9.5]{Mann2022NuclearSheaves}. In particular, by \cite[Proposition~9.7]{Mann2022NuclearSheaves}, there is a natural identification $q_{K!} \simeq q_{K*}$. Moreover, its cohomological amplitude is $0$ if $K$ has pro-order coprime to $\ell$ and $\Lambda$ is torsion. 
\end{lemma}
\begin{proof}
We first explain the case where $X = \Spd k'$, for an algebraically closed field $k'/k$. The $\ell$-cohomological properness follows from  \cite[Proposition~10.9]{Mann2022NuclearSheaves}. The claim on the cohomological amplitude in the torsion case follows from combining Lemma~\ref{lemma: basicCohomologylemma} with the fact that, under the identification $\Detale([\Spd k'/\ul{K}],\Lambda) \simeq \D(K,\Lambda)$ (\cite[Theorem~V.1.1]{FSGeomLLC}) via $\Detale(\Spd k',\Lambda) \simeq \D(\Lambda)$,  the functor $q_{K*}$ identifies with continuous group cohomology, by \cite[Example~4.2.4]{HKW}.

In general, we note that $[X/\ul{K}] \ra X$ is $\ell$-cohomologically proper, since $\ell$-cohomologically proper maps are stable under base-change by \cite[Lemma~9.8 (ii)]{Mann2022NuclearSheaves} and $[\Spd k/\underline{K}] \ra \Spd k$ is $\ell$-cohomologically proper, as explained above. To see the claim on the cohomological amplitude if $K$ is pro-$p$ and $\Lambda$ is torsion, we first note by \cite[Theorem~1.13 (ii)]{Ecod} and proper base-change applied to $[X/\ul{K}] \ra X$ that we may assume that $X$ is fdcs over a rank one geometric point $\Spa(C,\mathcal{O}_{C})$. Moreover, using \cite[Proposition~14.3]{Ecod},
it suffices to show that, for $A \in \Detale([X/\ul{K}],\Lambda)$ in degree $0$, the stalk of $q_{K*}A$ at each geometric point of $X$ is concentrated in degree $0$. By proper base-change, this reduces us to the case of a (possibly higher rank) geometric point $X = \Spa(C',C'^{+})$. To reduce to the previous case, we consider the cartesian diagram 
\begin{equation}{\label{eqn: StupidCartesianDiagram}}
\begin{tikzcd}
\Spa(C',C'^{+}) \arrow[r] \arrow[d,"g"] & \left[\Spa(C',C'^{+})/\underline{K}\right] \arrow[d,"g_{K}"] \\
\Spa(C,\mathcal{O}_{C}) \arrow[r] & \left[\Spa(C,\mathcal{O}_{C})/\underline{K}\right],
\end{tikzcd}
\end{equation}
where $g$ is the structure map, and the horizontal arrows are the obvious ones. Denote by $|X|$ the underlying topological space of $X$, and write $\D(|X|,\Lambda)$ for the left-completion of the category of sheaves of $\Lambda$-modules on the topological space $|X|$. Recall that, since \'etale covers of $\Spa(C',C'^{+})$ and of $\Spa C$ all split, we have identifications 
\[\Detale(\Spa(C',C'^{+}),\Lambda) \simeq \D(|\Spa(C^{'},C'^{+})|,\Lambda),\,\Detale(\Spa(C,\mathcal{O}_{C}),\Lambda) \simeq \D(|\Spa(C,\mathcal{O}_{C})|,\Lambda) \simeq \D(\Lambda),\] 
in such a way that $|g|^{*}$ is identified with $g^{*}$. By passing to right adjoints, we get an identification $R\Gamma(|\Spa(C',C'^{+})|,-) \simeq g_{\ast}$. However, by \cite[Lemma~7.2]{Ecod}, $g_{*}$ has cohomological amplitude zero. Hence, using qcqs base-change \cite[Theorem~1.9 (i)]{Ecod} applied to (\ref{eqn: StupidCartesianDiagram}), the functor $g_{K,\ast}$ also has cohomological amplitude zero. We may therefore replace $A$ by $g_{K,\ast}A$ and thus reduce to the already established case of $X = \Spa(C,\mathcal{O}_{C})$, as desired.
\end{proof}
We have the following very useful conservativity criterion. 
\begin{lemma}{\label{lemma: conservativitylemma}}
With notation as above, the functor 
\[\varinjlim_Kq_{K*}f_{K}^{*}: \Detale([X/H], \Lambda)\to \Detale(X,\Lambda)\]
is conservative, for $K$ varying over a cofinal system of pro-$p$ compact open subgroups of $H$. 
\end{lemma}
\begin{proof}
Since the functor given by mod $\ell$-reduction is conservative on $\Detale([X/H], \Lambda)$, we may assume without loss of generality that $\Lambda$ is torsion. Let $A\in \Detale([X/H], \Lambda)$ such that $\varinjlim_Kq_{K*}f_{K}^{*}A=0$. We need to show $A=0$. 

First, when $X = \Spd k'$ is an algebraically closed field then $\Detale([X/H],\Lambda)$ identifies with $\D(H,\Lambda)$ by \cite[Theorem~V.1.1]{FSGeomLLC}, the left-completed derived category of smooth representations of $H$ on $\Lambda$-modules. The claim then reduces to the observation that $A \simeq \varinjlim_{K} A^{K}$, where $K$ runs over a cofinal collection of open compact pro-$p$ subgroups of $H$.

Now suppose that $X = \Spa(C,C^{+})$ is a geometric point of higher rank. By \cite[Theorem 19.5]{Ecod} and Lemma \ref{lemma: classifiyingstackmapproperties}, we may perform a base-change and assume $X$ lives over a rank one geometric point $\Spa F$, for some algebraically closed perfectoid field $F/k$. Write $\ast$ for $\Spa F$ and consider the following commutative diagram
\[
\begin{tikzcd}
     \left[X/H\right]\ar[d,"g_H"] & \left[X/K\right]\ar[l,"{f_{K}}",swap] \ar[r,"{q_{K}}", left] \ar[d,"g_K"]& X\ar[d,"g"]\\
     \left[\ast/H\right] & \left[\ast/K\right]\ar[l,"{a_K}",swap] \ar[r,"{b_K}", left] & \ast,
\end{tikzcd}
\]
of cartesian squares. Here the vertical arrows are induced by the structure map $X\to \ast$. 

We have the following.
\begin{align*}
    0 &= g_\ast\varinjlim_Kq_{K*}f_{K}^{*}A\simeq \varinjlim_Kb_{K*}g_{K,\ast}f_{K}^{*}A\\
    &\simeq \varinjlim_Kb_{K*}a_{K}^{*} g_{H,\ast}A.
\end{align*}
Here the isomorphism in the first line follows from the identification $R\Gamma(|\Spa(C,C^{+})|,-) \simeq g_{*}(-)$ (see the proof of Lemma~\ref{lemma: classifiyingstackmapproperties} above) and the fact that $|\Spa(C,C^{+})|$ is a qcqs topological space, so taking cohomology on it commutes with colimits, see \cite[Tag~0739]{stacks-project}. The second isomorphism follows from the fact that $f_{K}$ is \'etale and smooth base-change. Now, by the case $X=\Spd k'$, applied to $k'=F$, the functor $\varinjlim_K b_{K*}a_{K}^{*}$ is conservative; hence, $g_{H,\ast}A=0$.

We note, however, by virtue of the fact that $\Detale(\Spa(C,C^{+}),\Lambda)$ identifies with the category of sheaves on the topological space $|\Spa(C,C^{+})|$ as in the previous proof, that the functor $g_\ast$ is conservative on sheaves which have non-zero stalk at the unique closed point of $|\Spa(C,C^{+})|$. In particular, $g_{H,\ast}A=0$ implies that $A$ is not supported at the closed point of $X$. It is therefore of the form $j_!B$, where $j: \Spa(C,C^{+}_{1}) \hookrightarrow \Spa(C,C^{+})$ is the open immersion corresponding to some larger valuation ring $C^{+} \subset C^{+}_{1}$. Therefore, by proper base change, we may replace $\Spa(C,C^{+})$ with $\Spa(C,C^{+}_{1})$, and the desired result then follows from induction on
the rank of the valuation ring $C^+$ (where the base case is the already established rank one case). The rank of this valuation is finite by our assumption that $X \ra \ast$ is fdcs and \cite[Lemma 21.6]{Ecod}.

For general $X$, consider the quotient map $f: X \ra [X/H]$. Since $H$ is locally pro-$p$ it is a pro-\'etale cover, and therefore pullback is conservative (using that $\cDetale(-,\Lambda)$ is a $v$-sheaf). Therefore, if $A$ is non-zero, it follows that $f^{*}A$ is also non-zero, which by \cite[Proposition~14.3]{Ecod}, implies that the stalk of $f^{*}A$ at some geometric point $\Spa(C,C^{+}) \xrightarrow{x} X$ is non-zero. This, in turn, implies that the pullback of $A$ along the induced map $x_{H}: [\Spa(C,C^{+})/H] \ra [X/H]$ is non-zero. Let $x_K$ be the induced map $[\Spa(C,C^{+})/K]\to [X/K]$. Then, by cohomological properness of $q_K$ (Lemma \ref{lemma: classifiyingstackmapproperties}) and the fact that $f_{K}$ is \'etale, we have
\[ x^\ast\varinjlim_Kq_{K*}f_{K}^{*}A\simeq \varinjlim_K(q_{x,K})_\ast(f_{x,K})^\ast x_H^\ast A.\] 
Hence the last term is non-zero by the result for geometric points above. This shows that $\varinjlim_Kq_{K*}f_{K}^{*}A$ is non-zero, as desired.
\end{proof}

We use this to prove the following criterion for universal local acyclicity.
\begin{proposition}{\label{prop: ULA-Criterion}}
Let $X \ra \Spd k$ be a fdcs map and $H$ be a locally pro-$p$ group, that acts trivially on $X$. A sheaf $A \in \Detale([X/H],\Lambda)$ is universally locally acyclic over $\Spd k$ if and only if $q_{K*}f_{K}^{*}A$ is universally locally acyclic over $\Spd k$ for $K$ varying over a cofinal system of open pro-$p$-subgroups of $H$.
\end{proposition}
\begin{proof}
To simplify the notation, we write $X_{K} := [X/K]$ and $X_{H} := [X/H]$. We write $p_{i,X}: X \times X \ra X$ for $i = 1,2$ for the projection maps. Similarly, we write $p_{i,X_{K}}$ and $p_{i,X_{H}}$ for the analogous projection maps for $X_{K}$ and $X_{H}$, respectively.

Using the criterion \cite[Theorem IV.2.23(ii)]{FSGeomLLC} for universal local acyclicity (cf. \cite[Lemma 4.4.5]{HeyerMann}), to show $A$ is ULA over $\Spd k$, it is equvalent to show that the natural map
\[\alpha_{X_H/k,A}:p_{1,X_{H}}^\ast \mathbb{D}_{X_H/k}(A)\otimes^\mathbb{L}_\Lambda p_{2,X_{H}}^\ast A\rightarrow R\mathcal{H}om(p_{1,X_{H}}^\ast A, p_{2,X_{H}}^! A)\]
is an isomorphism. For this, we consider, for each $K$ in the cofinal system, the diagram
\begin{align*}
    X_H\times X_H \xleftarrow{f_K \times f_K} X_K \times X_K \xrightarrow{q_K \times q_K} X \times X.
\end{align*}
The ULAness of $q_{K*}f_{K}^{*}(A)$ is equivalent to the condition that the natural map
\[\alpha_{X/k,q_{K*}f_{K}^{*}(A)}:p_{1,X}^\ast \mathbb{D}_{X/k}q_{K*}f_{K}^{*}(A) \otimes^\mathbb{L}_\Lambda p_{2,X}^\ast A\rightarrow R\mathcal{H}om(p_{1,X}^*q_{K*}f_{K}^{*}(A), p_{2,X}^!q_{K*}f_{K}^{*}(A))\]
be an isomorphism. The ``only if'' direction of the proposition follows from Lemma~\ref{lem:ULA intertwine} below. The ``if'' direction of the proposition also follows from Lemma~\ref{lem:ULA intertwine} below, together with the conservativity result of Lemma~\ref{lemma: conservativitylemma}. 
\end{proof}

\begin{lemma}\label{lem:ULA intertwine}
The functor $(q_K \times q_K)_\ast (f_K \times f_K)^\ast$ intertwines the maps $\alpha_{X_H/k,A}$ and $\alpha_{X/k,q_{K*}f_{K}^{*}(A)}$.
\end{lemma}
\begin{proof}
We apply $(q_{K} \times q_{K})_{\ast} (f_{K} \times f_{K})^{*}$ to both sides of $\alpha_{X_{H}/k,A}$. We can rewrite the RHS as 
\begin{align*} 
(q_{K} \times q_{K})_{*}(f_{K} \times f_{K})^{*}R\mathcal{H}om(p_{1,X_{H}}^\ast A, p_{2,X_{H}}^! A) &\simeq \\
 (q_{K} \times q_{K})_{*}R\mathcal{H}om(p_{1,X_{K}}^{\ast}f_{K}^{*}A, p_{2,X_{K}}^{!}f_{K}^{!}A) & \simeq \\
R\mathcal{H}om(q_{K,*}p_{1,X_{K}}^{\ast}f_{K}^{*}A, q_{K,*}p_{2,X_{K}}^{!}f_{K}^{!}A) &\simeq \\
R\mathcal{H}om(p_{1,X}^*q_{K*}f_{K}^{*}A, p_{2,X}^!q_{K*}f_{K}^{*}A).
\end{align*}
Here we used the projection formula \cite[Proposition~23.3 (2)]{Ecod} for the first isomorphism, several applications of the projection formula \cite[Proposition~23.3 (1)]{Ecod} for the second isomorphism, and proper base change for the third isomorphism, where we note that we have an identification $(q_{K} \times q_{K})_{!} \simeq (q_{K} \times q_{K})_{*}$ by Lemma~\ref{lemma: classifiyingstackmapproperties}. For the LHS, we can similarly rewrite this as 
\begin{align*} 
(q_{K} \times q_{K})_{*}(f_{K} \times f_{K})^{*}(p_{1,X_{H}}^\ast \mathbb{D}_{X_H/k}(A)\otimes^\mathbb{L}_\Lambda p_{2,X_{H}}^\ast A) & \simeq \\
(q_{K} \times q_{K})_{*}(p_{1,X_{K}}^{*}f_{K}^{*}\mathbb{D}_{X_H/k}(A)\otimes^\mathbb{L}_\Lambda p_{2,X_{K}}^*f_{K}^{*}A) & \simeq \\
q_{K*}p_{1,X_{K}}^{*}f_{K}^{*}\mathbb{D}_{X_H/k}(A)\otimes^\mathbb{L}_\Lambda q_{K*}p_{2,X_{K}}^*f_{K}^{*}A)
& \simeq  \\
p_{1,X}^\ast \mathbb{D}_{X/k}q_{K*}f_{K}^{*}(A) \otimes^\mathbb{L}_\Lambda p_{2,X}^\ast A.
\end{align*}
Here we have used the commutation of upper $*$ with tensor products for the first isomorphism, commutation of lower $!$ with external tensor product for the second isomorphism (which follows from several applications of projection formula in its usual form), and proper base-change for the third isomorphism, where we note that $f_{K}$ is \'etale so that we have $f_{K}^{!} \simeq f_{K}^{*}$. The desired claim follows (cf. \cite[Proposition~IV.2.28]{FSGeomLLC}).
\end{proof}

\section{Minimal compactifications of Shimura varieties}{\label{sec: stratifications}}

Let $\gx$ be a Shimura datum of PEL type with reflex field $\mathsf E$. In this section, we briefly recall the construction of the minimal compactification $\Shstar\gx_K$, for a neat compact open subgroup $K\subset\mathsf{G}(\mathbb{A}_f)$, and its stratification by boundary strata, following \cite{PinkHigherDirectImages}. We then move to the $p$-adic situation and analyze the boundary strata there. 

\subsection{Stratification by boundary components}\label{Subsub:StratificationShimura}
We let $\mathsf{X}^\ast$ denote the union of the rational boundary components of $\mathsf{X}$, as in \cite[Section~3.5]{PinkHigherDirectImages} endowed with the Satake topology \cite[Section~3.1]{PinkHigherDirectImages}. The complex analytification of the $\mathbb{C}$-points of $\Shstar\gx_{K}$, the minimal compactification of level $K$, is given by the double cosets 
\begin{equation}{\label{eqn: adelicquotientdefiningminimalcompactification}}
\mathsf{G}(\bb{Q}) \backslash \mathsf{X}^{*} \times \mathsf{G}(\bb{A}_{f})/K. 
\end{equation}
The resulting space is stratified according to the type of rational boundary components. To describe this more precisely, recall the following notion.
\begin{definition}{\label{defn: admissibleparabolic}}
    A parabolic subgroup $\mathsf{P}$ of $\mathsf{G}$ (defined over $\mathbb{Q}$) is called \emph{admissible} if its projection to every simple factor in $\mathsf{G}^\mathrm{ad}$ is either a maximal parabolic or equal to that simple factor.
\end{definition}

The action of the group $\mathsf{G}(\mathbb{R})$ on $\mathsf X$ extends to an action on $\mathsf X^\ast$. The normalizer of the stabilizer of each boundary component is given by the $\mathbb{R}$-points of some admissible parabolic $\mathsf{P}$. We write $\mathsf{X}_{\mathsf{P}}$ for the union of the rational boundary components stabilized by $\mathsf{P}(\mathbb{Q})$. Then the boundary strata $\Sh\gx_{K,\PP}$ of 
$\Shstar\gx_{K}$ have complex points given by 
\[\Sh\gx_{K,\PP}(\bb{C}) := \mathsf{G}(\bb{Q})\backslash \bigsqcup_{\mathsf{P}\in [\mathsf{P}]} \mathsf{X}_{\mathsf{P}} \times \mathsf{G}(\bb{A}_{f})/K, \]
where $\PP$ runs through $\mathsf{G}(\mathbb{Q})$-conjugacy classes of admissible rational parabolic subgroups of $\mathsf{G}$. We will denote the closure of the stratum $\Sh(\mathsf{G},X)_{K,[\mathsf{P}]}$ by $\Sh(\mathsf{G},X)_{K,\leq [\mathsf{P}]}$. 

Consider a choice of admissible parabolic $\mathsf{P} \in [\mathsf{P}]$. Then there exists a connected normal subgroup $\mathsf{G}_{1}$ of the Levi quotient of $\mathsf{P}$, which acts transitively on an open and closed subspace $\mathsf{X}_{1}\subset \mathsf{X}_{\mathsf{P}}$. Moreover, the pair $(\mathsf{G}_{1},\mathsf{X}_{1})$ can be upgraded to a Shimura datum, as in \cite[Section~3.6]{PinkHigherDirectImages}. We write $\mathsf{P}_1$ for the preimage of $\mathsf{G}_{1}$ in $\mathsf{P}$. Given also a choice of $g\in G(\bb{A}_f)$, we let $K_{\mathsf{P}_1}$ denote the intersection $gKg^{-1} \cap \mathsf{P}_1(\bb{A}_{f})$. We set $\ol{K}_{\mathsf{P}_1}$ to be the image of $K_{\mathsf{P}_1}$ in the Levi quotient of $\mathsf{P}_1(\bb{A}_f)$. Then projection to the Levi quotient induces an isomorphism
\[ \mathsf{G}_{1}(\bb{Q})\backslash \mathsf{X}_{1} \times \mathsf{G}_{1}(\bb{A}_{f})/\ol{K}_{\mathsf{P}_1} \simeq \mathrm{Stab}_{\mathsf{P}(\bb{Q})}(\mathsf{X}_{1})\backslash \mathsf{X}_{1} \times \mathrm{Stab}_{\mathsf{P}(\bb{Q})}(\mathsf{X}_{1})\mathsf{P}_1(\bb{A}_{f})/K_{\mathsf{P}_1}, \]
as in \cite[Equation~3.7.1]{PinkHigherDirectImages}, where the RHS admits a natural map to the adelic quotient in (\ref{eqn: adelicquotientdefiningminimalcompactification}), which also depends on $g$. By \cite[Theorem~3.7.2]{PinkHigherDirectImages}, this descends to a map 
\begin{equation}{\label{eqn: finitemap}}
\mathrm{Sh}(\mathsf{G}_{1},\mathsf{X}_{1})_{\ol{K}_{\mathsf{P}_1}} \ra \Shstar\gx_{K} 
\end{equation}
of schemes over $\mathsf E$ that is finite over its image. The images of these maps exhaust, up to varying $\mathsf{P} \in [\mathsf{P}]$, $X_1\subset X_{\mathsf{P}}$ and $g\in G(\A_f)$, the components of the boundary strata $\Sh\gx_{K,[\mathsf{P}]}$. 

\subsection{Structure of the strata}\label{subsub: SVStrata} 
We will now need a  more explicit description of these strata. Maintain the choices of $\mathsf{P}\in [\mathsf{P}]$, $X_1\subset X_{\mathsf{P}}$ and $g\in G(\bb{A}_f)$. Set $\mathsf{U}$ (resp. $\mathsf{U}_1$) to be the unipotent radical of $\mathsf{P}$ (resp. $\mathsf{P}_1$). We define the following subgroups (where, by an abuse of notation, we suppress $g$):
\[H_{\mathsf{P}} := gKg^{-1} \cap \mathrm{Stab}_{\mathsf{P}(\bb{Q})}(\mathsf{X}_{1})\mathsf{P}_1(\bb{A}_{f}) \]
\[ H_{\mathsf{U}} := gKg^{-1} \cap \mathrm{Cent}_{\mathsf{P}(\bb{Q})}(\mathsf{X}_{1})\mathsf{U}_1(\bb{A}_{f}) \]
\[ K_{\mathsf{P}_{1}} := gKg^{-1} \cap \mathsf{P}_{1}(\bb{A}_{f})  \]
\[ K_{\mathsf{U}_1} := gKg^{-1} \cap \mathsf{U}_1(\bb{A}_{f}).  \]
The natural map (\ref{eqn: finitemap}) factors through a locally closed embedding
\begin{equation}{\label{eqn: locallyclosedembedding}}
\mathrm{Sh}(\mathsf{G}_{1},\mathsf{X}_{1})_{\ol{K}_{\mathsf{P}_1}}\sslash H_{\mathsf{P}} \hookrightarrow \Shstar\gx_{K},
\end{equation}
where the LHS denotes the scheme theoretic quotient of $H_{\mathsf{P}}$ acting on $\mathrm{Sh}(\mathsf{G}_{1},\mathsf{X}_{1})_{\ol{K}_{\mathsf{P}_1}}$; the action factors through the quotient of $H_{\mathsf{P}}$ by $K_{\mathsf{P}_1}H_{\mathsf{U}}$. 

\begin{remark}\label{rem:trivial quotient}
Because we are working with a Shimura datum of PEL type, when the neat level $K$ is principal, the action of the finite group $H_{\mathsf{P}}/K_{\mathsf{P}_1}H_{\mathsf{U}}$ on $\mathrm{Sh}(\mathsf{G}_{1},\mathsf{X}_{1})_{\ol{K}_{\mathsf{P}_1}}$ is trivial. See~\cite[Remark 4.3.3]{LanStroh} and also the proof of Theorem 4.3.10 of \emph{loc. cit.}. Therefore, the morphism~\eqref{eqn: finitemap} is already a locally closed embedding. 
\end{remark}

The groups $K_{\mathsf{U}_1}$  and $K_{\mathsf{P}_1}$ are compact subgroups of $\mathsf{G}(\bb{A}_{f})$ and are therefore profinite with the induced topology. However, the subgroups $H_{\mathsf{U}}$ and $H_{\mathsf{P}}$ of $\mathsf{G}(\bb{A}_{f})$ are not profinite. Rather, they are extensions 
\[ 1 \ra K_{\mathsf{P}_1} \ra H_{\mathsf{P}} \ra \Gamma_\mathsf{P} \ra 1 \]
and 
\[ 1 \ra K_{\mathsf{U}_1} \ra H_{\mathsf{U}} \ra \Gamma_{\mathsf{U}} \ra 1, \]
where $\Gamma_\mathsf{P}$ and $\Gamma_\mathsf{U}$ are arithmetic groups. Nevertheless, the quotient $H_{\mathsf{P}}/H_{\mathsf{U}}$ with induced topology will be profinite, see the explanation on \cite[Page 222]{PinkHigherDirectImages}. We write $\hat{H}_\mathsf{P}$ and $\hat{H}_\mathsf{U}$ for the closure of these groups in $\mathsf{G}(\bb{A}_{f})$ (equivalently in $gKg^{-1}$), which will be profinite. We note that the injection $H_{\mathsf{U}} \hookrightarrow H_{\mathsf{P}}$ extends to an injection $\hat{H}_\mathsf{U} \hookrightarrow \hat{H}_\mathsf{P}$ and we have an isomorphism $H_{\mathsf{P}}/H_{\mathsf{U}} \simeq \hat{H}_\mathsf{P}/\hat{H}_\mathsf{U}$. 

For an inclusion of compact open subgroups $K' \subset K$, we have a cartesian diagram of the form
\[
\begin{tikzcd}
\mathrm{Sh}(\mathsf{G}_{1},\mathsf{X}_{1})_{\ol{K}'_{\mathsf{P}_1}}\sslash H'_{P} \arrow[r] \arrow[d] & \Shstar\gx_{K'} \arrow[d] \\
\Sh(\mathsf{G}_{1},\mathsf{X}_{1})_{\ol{K}_{\mathsf{P}_1}}\sslash H_{\mathsf{P}} \arrow[r] & \Shstar\gx_{K}, 
\end{tikzcd}
\]
where we have used the superscript $(-)^{'}$ to denote analogous notations for $K'$, and the horizontal maps are as in (\ref{eqn: locallyclosedembedding}). We record the following lemma.

\begin{lemma}\label{lemma: proetaletorsorboundary}
Taking the limit over open normal subgroups $K' \subset K$, we obtain the map
\begin{equation}
\varprojlim_{K' \subset K} \Sh(\mathsf{G}_{1},\mathsf{X}_{1})_{\ol{K}'_{\mathsf{P}_1}}\sslash H'_{P} \ra \Sh(\mathsf{G}_{1},\mathsf{X}_{1})_{\ol{K}_{\mathsf{P}_1}}\sslash H_{\mathsf{P}},
\end{equation}
which is a pro-\'etale $H_{\mathsf{P}}/H_{\mathsf{U}} \simeq \hat{H}_\mathsf{P}/\hat{H}_\mathsf{U}$-torsor.    
\end{lemma}
\begin{proof}
This follows from \cite[Proposition~3.7.5]{PinkHigherDirectImages}. 
\end{proof}

\subsection{The minimal compactification as a diamond}\label{sec: SVasDiamond}
Let us move to the $p$-adic situation. Choose a $p$-adic place $v$ of $\mathsf{E}$ and denote by $E$ the completion of $\mathsf{E}$ at $v$. We fix a completed algebraic closure $C$ of $E$ and base-change the Shimura variety $\Sh\gx_{K}$ to $C$.  We will replace $\Sh$ with the symbol $\mathcal{S}$ when discussing the diamonds over $C$ attached to the Shimura variety at various levels by applying $(-)^{\Diamond}$ as defined in \S \ref{ss: AnalytificationinMixedCharacteristic}, often omitting the Shimura datum when the context is clear. We will fix an away-from-$p$ level $K^p$ and denote by $\mc{S}_{K^p}$ the inverse limit $\varprojlim_{K_p}\mc{S}_{K^pK_p}$. Similarly, for the minimal compactifications, we use the parallel notations decorated by the superscript ${}^\ast$. 

Applying the functor $(-)^\Diamond$ to the stratification in \S \ref{Subsub:StratificationShimura} above, we obtain a set-theoretic locally closed stratification 
\[ \mathcal{S}^{*}_{K} = \bigsqcup_{[\mathsf{P}]} \mc{S}_{K,\PP}, \]
where $[\mathsf{P}]$ runs over conjugacy classes of admissible rational parabolics $P \subset \mathsf{G}$, and $\mc{S}_{K,[\mathsf{G}]} = \mc{S}_K$ is the open stratum. We write $\mc{S}_{K,\leq [\mathsf{P}]}$ for the $(-)^{\Diamond}$ functor applied to the closure of the strata.

Similarly, since the $K_p$-action on $\mc{S}_{K^p}^\ast$ preserves each stratum labeled by $\PP$, we can define a set-theoretic locally closed stratification for the stacky quotient $[\mc{S}_{K^p}^\ast/\ul{K_{p}}]$, by pullback along the natural map $q_{K_p}: [\mc{S}_{K^p}^\ast/\ul{K_{p}}]\to \mc{S}^\ast_{K}$ from the stacky quotient at level $K_{p}$ to the coarse quotient. More precisely, the map $q_{K_{p}}$ is defined by the factorization
\begin{equation}{\label{eqn: factorizationofqKp}}
q_{K_p}: [\mc{S}_{K^p}^{*}/\ul{K_{p}}] \xrightarrow{q^{1}_{K_{p}}} [\mc{S}^\ast_{K}/\ul{K_{p}}] \xrightarrow{q_{K_{p}}^{2}} \mc{S}^\ast_{K},  
\end{equation}
where the first map is induced by the natural $\ul{K_p}$-equivariant (with trivial action on the target) map $\mc{S}_{K^p}^{*} \ra \mc{S}^\ast_{K}$, and the second map is the base-change of the obvious map $[\ast/\underline{K_{p}}] \ra \ast$. The pullback of $\mc{S}_{K,\PP}$ along $q_{K_p}$ can be written as the quotient of 
\[ \mathcal{S}_{K^p,[\mathsf{P}]} := \varprojlim_{K_{p} \ra \{1\}} \mc{S}_{K,\PP}, \]
by the natural action of $K_{p}$. In particular, we have a natural locally closed stratification
\[ [\mc{S}_{K^p}^\ast/\ul{K_{p}}] = \bigsqcup_{[\mathsf{P}]} [\mathcal{S}_{K^p,[\mathsf{P}]}/\ul{K_{p}}]. \]
Analogously, we write $[\mc{S}^{*}_{K^p,\leq [\mathsf{P}]}/\ul{K_{p}}]$ for the pullback of the closure of the strata corresponding to $[\mathsf{P}]$ along $q_{K_{p}}$.

We let $q_{K_{p},[\mathsf{P}]}: [\mathcal{S}_{K^p,[\mathsf{P}]}/\ul{K_{p}}] \ra \mc{S}_{K,\PP}$ be the base-change of $q_{K_{p}}$, and let $q_{K_{p},[\mathsf{P}]}^{1}$ (resp. $q_{K_{p},[\mathsf{P}]}^{2}$) be the base-change of $q_{K_{p}}^{1}$ (resp. $q_{K_{p}}^{2}$) appearing in the factorization (\ref{eqn: factorizationofqKp}). We use the analogous notation for the closed strata.  

We record two propositions on the structure of $q_{K_{p}}$ for later use. Let $\Lambda$ be as in Setup \ref{assumption: coefficientsystemsingeneral}.
\begin{proposition}\label{prop: qKpProper}
For a compact open subgroup $K_p\subset G(\qp)$ of finite $\ell$-cohomological dimension, the map $q_{K_p}$ is $\ell$-cohomologically proper and $q_{K_p,\ast} \simeq q_{K_p,!}$ has cohomological amplitude $0$ when $\Lambda$ is as in Setup \ref{assumption: coefficientsystemsingeneral} (1) and the pro-order of $K_p$ is coprime to $\ell$.
\end{proposition}
\begin{proof}
Consider the factorization in Equation~(\ref{eqn: factorizationofqKp}). It suffices to show the claims for $q^{1}_{K_{p}}$ and $q^2_{K_{p}}$ separately, using that cohomologically proper maps are stable under composition by \cite[Lemma~9.8 (ii)]{Mann2022NuclearSheaves}. For the map $q^2_{K_P}$, this follows immediately from Lemma~\ref{lemma: classifiyingstackmapproperties}. The desired claim for $q^1_{K_p}$ follows from the next lemma using \cite[Lemma~9.8 (i)]{Mann2022NuclearSheaves}, and \cite[Remark~21.14]{Ecod} for the claim on the cohomological amplitude. 
\end{proof}

\begin{lemma}{\label{lemma: q1isquasiproetaleandproper}}
For $K_{p} \subset G(\mathbb{Q}_{p})$ any compact open subgroup, the map $q^1_{K_p}: [\mathcal{S}_{K^p}^{*}/\ul{K_p}] \ra [\mathcal{S}^{*}_{K}/\ul{K_p}]$ is quasi-pro-\'etale, fdcs, and proper.
\end{lemma}

\begin{proof}
We note that, by definition of $q^{1}_{K_{p}}$, we have a natural cartesian diagram 
\[
\begin{tikzcd}
\mathcal{S}_{K^{p}}^{*} \arrow[r,"f_{K_{p}}"] \arrow[d] & \mathcal{S}_{K}^{*} \arrow[d] & \\
\left[\mathcal{S}_{K^{p}}^{*}/\underline{K_{p}}\right] \arrow[r,"q^{1}_{K_{p}}"] & \left[\mathcal{S}_{K}^{*}/\underline{K_{p}}\right]. & 
\end{tikzcd}
\]
In particular, it is clear that the map $f_{K_{p}}$ is fdcs and that the same is true for $q^{1}_{K_{p}}$, since being fdcs may be checked after base-changing along a $v$-surjective quasi-pro-\'etale map (see \cite[Proposition~3.2.7]{GHILZIsocComparison}). For the remaining claims on $q^{1}_{K_{p}}$, we note that, using \cite[Proposition~10.11]{Ecod}, it suffices to show the analogous claim for the map $f_{K_{p}}$, since $\mathcal{S}^{*}_K \ra [\mathcal{S}^{*}_K/\ul{K_{p}}]$ is a v-cover. The properness of $f_{K_p}$ follows from the fact that both $\mathcal{S}^{*}_{K^p}$ and $\mathcal{S}^{*}_K$ are proper so that the map is qcqs. Then, we conclude by invoking the valuative criterion of properness \cite[Proposition~18.3]{Ecod}, recalling that we may write $f_{K_{p}}$ as an inverse limit over the analytification of finite maps, which are in particular proper. 

As already discussed, the map $f_{K_{p}}$ is qcqs, in order to check that the map is quasi-pro-\'etale, using \cite[Lemma~7.19]{Ecod}, this reduces to the fact that the fiber of $f_{K_{p}}: \mathcal{S}^{*}_{K^p} \ra \mathcal{S}^{*}_K$ over a rank one geometric point is a profinite set, which is clear since it is an  inverse limit of the analytification of finite morphisms.
\end{proof}

For the next proposition, we fix any conjugacy class $[\mathsf{P}]$ of admissible rational parabolics of $\mathsf{G}$. For a choice of $\mathsf{P}\in [\mathsf{P}]$ $X_1\subset X_{\mathsf{P}}$ and $g\in G(\A_f)$, we adopt the notation from \S\ref{subsub: SVStrata}. We write $g=(g^p,g_p)$ with $g^p\in G(\bb{A}_f^p)$ and $g_p\in G(\bb{Q}_p)$. For the various groups introduced  there, we use the decoration $(-)^p$, resp. $(-)_p$ to denote their variants obtained by intersecting with $g^pK^p(g^p)^{-1}$, resp. $g_pK_p(g_p)^{-1}$. For example, we have 
\[H_{\mathsf{P}}^p= H_{\mathsf{P}}\cap g^pK^p(g^p)^{-1}\]
and so on. Note that this is still an extension of a discrete group by a profinite group. Moreover, we will use $\hat{(-)}$ on such a group to denote its closure in the corresponding level subgroup, for example as before $\hat{H}_\mathsf{P}$ is the closure of $H_\mathsf{P}$ in $gKg^{-1}$, $\hat{H}^p_\mathsf{P}$ is the closure of $H^p_\mathsf{P}$ in $g^pK^p(g^p)^{-1}$, $\hat{H}_{\mathsf{P},p}$ is that of $H_{\mathsf{P},p}$ in $g_pK_p(g_p)^{-1}$ and so on.

We first need a lemma, where for a choice $\mathsf{P} \in [\mathsf{P}]$ and $\mathsf{X}_1\subset \mathsf{X}_P$, we write $\mathcal{S}(\mathsf{G}_{1},\mathsf{X}_{1})_{\ol{K}_{\mathsf{P}_1}}$ for the diamond over $C$ attached to the variety $\mathrm{Sh}(\mathsf{G}_{1},\mathsf{X}_{1})_{\ol{K}_{\mathsf{P}_1}}$ appearing in the target of (\ref{eqn: finitemap}). To define the map~\eqref{eqn: finitemap}, we also need a choice of $g\in G(\bb{A}_f)$. 

\begin{lemma}\label{lemma: BoundaryInfinitePLevel} Assume that both levels $K^p\subset G(\bb{A}_f^p)$ and $K_p\subset G(\Q_p)$ are principal. 
The maps (\ref{eqn: finitemap}) induce a locally closed embedding
\[\mathcal{S}(\mathsf{G}_{1},\mathsf{X}_{1})_{\ol{K}^{p}_{\mathsf{P}_1}}:= \varprojlim_{K_{p} \ra \{1\}} \mathcal{S}(\mathsf{G}_{1},\mathsf{X}_1)_{\overline{K}_{\mathsf{P}_1}}\hookrightarrow \mc{S}_{K^p}^\ast,\]  
which is $\hat{H}_{\mathsf{P},p}$-equivariant, with the action on the RHS induced by the composition $\hat{H}_{\mathsf{P},p}\hookrightarrow g_pK_p(g_p)^{-1}\stackrel{\sim}{\to} K_p$, where the first map is the natural inclusion and the second map is conjugation. Hence, we obtain a map of the stack quotients
\[ [\mathcal{S}(\mathsf{G}_{1},\mathsf{X}_{1})_{\ol{K}^{p}_{\mathsf{P}_1}}/\ul{\hat{H}_{\mathsf{P},p}}] \ra [\mc{S}_{K^p}^\ast/\ul{K_{p}}],\]
which, upon passing to coarse moduli spaces, recovers the analytification of (\ref{eqn: locallyclosedembedding}). Moreover, 
\[\mathcal{S}(\mathsf{G}_{1},\mathsf{X}_{1})_{\ol{K}^{p}_{\mathsf{P}_1}}\to \mathcal{S}(\mathsf{G}_{1},\mathsf{X}_{1})_{\ol{K}_{\mathsf{P}_1}}\sslash H_p\] 
is a $H_{\mathsf{P},p}/H_{\mathsf{U},p}\simeq \hat{H}_{\mathsf{P},p}/\hat{H}_{\mathsf{U},p}$-torsor.
\end{lemma}
\begin{proof}
The first part follows by passing to the inverse limit over the locally closed embeddings in~\eqref{eqn: locallyclosedembedding} as $K'_p\to \{1\}$ runs over principal levels. Note that at each finite level $K^pK'_p$, which is principal by assumption, Remark~\ref{rem:trivial quotient} applies to tell us that there is no need to take an additional finite quotient in~\eqref{eqn: locallyclosedembedding}. 
    The $\hat{H}_{\mathsf{P},p}$-equivariance is clear. The last statement follows from Lemma~\ref{lemma: proetaletorsorboundary} by taking the $p$-part and analytifying.
\end{proof}

\begin{proposition}\label{lemma: stackystratasmooth}
For any neat level $K = K^pK_p$, the diamond $\mc{S}_{K,\PP}$ is $\ell$-cohomologically smooth of $\ell$-dimension $d_{\PP}$ over $\Spd C$. Moreover, if the level subgroup $K_{p}$ is pro-$p$, then the map $q_{K_{p},[\mathsf{P}]}$ is $\ell$-cohomologically smooth of $\ell$-dimension $0$. In particular, in this case, $[\mathcal{S}_{K^p,\PP}/\ul{K_{p}}]$ is also  $\ell$-cohomologically smooth of $\ell$-dimension $d_{[\mathsf{P}]}$ over $\Spd C$.
\end{proposition}

\begin{proof}
The first part follows from smoothness of $\Sh\gx_{K,\PP}$, cf. \cite[Proposition~3.7.5]{PinkHigherDirectImages} and Corollary \ref{prop: Compatabilitieswithupper!s}. 

For the second part, assume first that the level $K=K^pK_p$ is principal. Choose $\mathsf{P} \in [\mathsf{P}]$, $\mathsf{X}_{1}\subset \mathsf{X}_{\mathsf{P}}$ and $g\in G(\A_f)$, and consider the map
\[ [\mathcal{S}(\mathsf{G}_{1},\mathsf{X}_{1})_{\ol{K}^{p}_{\mathsf{P}_1}}/\ul{\hat{H}_{\mathsf{P},p}}] \ra [\mc{S}_{K^p}^\ast/\ul{K_{p}}]\]
as in Lemma~\ref{lemma: BoundaryInfinitePLevel} above. It follows from the last part of Lemma~\ref{lemma: BoundaryInfinitePLevel} that we have an isomorphism 
\[ [\mathcal{S}(\mathsf{G}_1,\mathsf{X}_1)_{\ol{K}^{p}_{\mathsf{P}_1}}/\ul{\hat{H}_{\mathsf{P},p}}] \simeq [(\mathcal{S}(\mathsf{G}_1,\mathsf{X}_1)_{\ol{K}_{\mathsf{P}_1}}\sslash H_{\mathsf{P}})/\ul{\hat{H}_{\mathsf{U},p}}],  \]
of $v$-stacks, with $\ul{\hat{H}_{\mathsf{U},p}}$ acting trivially on the RHS. This gives a cartesian diagram 
\begin{equation}{\label{eqn: ClassifyingStackDescriptionofBoundary}}
 \begin{tikzcd}
\left[(\mathcal{S}(\mathsf{G}_1,\mathsf{X}_1)_{\ol{K}_{\mathsf{P}_1}}\sslash H_{\mathsf{P}})/\ul{\hat{H}_{\mathsf{U},p}}\right] \arrow[d] \arrow[r] & \mathcal{S}(\mathsf{G}_1,\mathsf{X}_1)_{\ol{K}_{\mathsf{P}_1}}\sslash H_{\mathsf{P}} \arrow[d] & \\
\left[\mc{S}_{K^p}^\ast/\ul{K}_{p}\right] \arrow[r,"q_{K_p}"] & \mc{S}^\ast_{K}, & 
\end{tikzcd} 
\end{equation}
and every connected component of the map $q_{K_{p},[\mathsf{P}]}$ is represented by one of the top horizontal maps as one varies over $\mathsf{P} \in [\mathsf{P}]$, $\mathsf{X}_1\subset \mathsf{X}_P$ and $g\in G(\A_f)$. Now, to conclude, note that the map
\[ [(\mathcal{S}(\mathsf{G}_1,\mathsf{X}_1)_{\ol{K}_{\mathsf{P}_1}}\sslash H_{\mathsf{P}})/\ul{\hat{H}_{\mathsf{U},p}}] \ra \mathcal{S}(\mathsf{G}_1,\mathsf{X}_1)_{\ol{K}_{\mathsf{P}_1}}\sslash H_{\mathsf{P}} \]
is $\ell$-cohomologically smooth of $\ell$-dimension $0$. Indeed, it is the base-change of $[\ast/\ul{\hat{H}_{\mathsf{U},p}}] \ra \ast$, since $\hat{H}_{\mathsf{U},p}$ acts trivially. Hence, we conclude by \cite[Theorem~10.13, Example~10.11 (b)]{Mann2022NuclearSheaves}. 

For a general neat level $K = K^pK_p$ with $K_p$ pro-$p$, we can find a principal level $K'\triangleleft K$. Writing $K' = (K^p)'K'_p$, we have that $K'$ is neat and $K'_p$ is pro-$p$. Therefore, $[\mathcal{S}_{(K^p)',[\mathsf{P}]}/K'_p]$ is $\ell$-cohomologically smooth of dimension $d_{[\mathsf{P}]}$ over $\mathrm{Spd}\ C$. The change-of-level map on the stacky quotients
\[
[\mathcal{S}_{(K^p)',[\mathsf{P}]}/K'_p]
\to 
[\mathcal{S}_{K^p,[\mathsf{P}]}/K_p]
\]
is finite \'etale. Therefore, we deduce that $[\mathcal{S}_{K^p,[\mathsf{P}]}/K_p]$ is also $\ell$-cohomologically smooth of dimension $d_{[\mathsf{P}]}$ over $\mathrm{Spd}\ C$. Similarly, the $\ell$-cohomological smoothness of the map $q_{K^p, [\mathsf{P}]}$ follows from the $\ell$-cohomological smoothness of $q_{(K^p)', [\mathsf{P}]}$ and the finite \'etaleness of the change-of-level maps on the stacky quotients as well as on the finite-level boundary strata labeled by $[\mathsf{P}]$.  
\end{proof}

\section{Partial compactifications of Igusa varieties}\label{sec: pc Igusa varieties}

Igusa varieties and their partial minimal compactifications are relevant to us because they define $v$-sheaves which can be used to describe the fibers of the morphism from the Igusa stack to $\Bun_{G}$. The relationship between Igusa varieties and the Igusa stack is recalled in detail in \S\ref{sec: Stratified}. 

In this section, we recall the definition and properties of Igusa varieties and their partial minimal and toroidal compactifications in the PEL case, following~\cite{CS17, CS2, santos}. We then describe a natural and explicit extension of the moduli-theoretic action of $\widetilde{J}_b$ on open Igusa varieties to their partial compactifications by lifting to toroidal compactifications. This will play a key role in the construction of the Baily--Borel stratification on Igusa stacks which we carry out in \S 5.

\subsection{Igusa varieties and their partial compactifications}\label{Sec:Igusa background}

We let $\gx$ be a PEL type AC Shimura datum in Kottwitz's classification and $p$ be a rational prime, for which $G:=\mathsf{G}_{\qp}$ is unramified, as in Assumption \ref{assumption:codimension}. Let  $\mathsf{E}$ be the reflex field. We fix an isomorphism $\ol{\mathbb{Q}}_{p} \simeq \bb{C}$ and let $v\mid p$ denote resulting place of $\mathsf{E}$; set $E:=\mathsf{E}_v$. Denote by $\mathcal{O}$ the ring of integers of $E$, and by $\mathbb{F}$ its residue field. We choose an algebraic closure $\overline{\mathbb{F}}_p$ of $\mathbb{F}$.

Fix a neat compact open subgroup $K\subset \mathrm{G}(\mathbb{A}_f)$ of the form $K = K^pK_p$, with $K_p\subset G(\mathbb{Q}_p)$ a hyperspecial maximal compact subgroup. We let $(\mc{O}_B,\ast,\Lambda,(\cdot,\cdot),h)$ be an integral PEL datum that gives rise to $\gx$ and consider the canonical integral model over $\mathcal{O}$ of the Shimura variety attached to it, see \cite{Kottwitz}, cf. \cite[Definition 5.9]{zhang2023}.\footnote{The Hasse principle can fail for $\mathsf{G}$'s that are of type A, in which case the generic fiber of the moduli problem consists of several copies of the Shimura variety, see \cite[Remark 5.7]{zhang2023}. We ignore this difference below when it does not affect the arguments.} We let $\mathscr{S}:=\mathscr{S}_{K, \overline{\mathbb{F}}_p}$be its special fiber, base-changed to $\overline{\mathbb{F}}_p$. We denote by $\mathcal{A}$ the universal abelian scheme over $\mathscr{S}$, equipped with its PEL structures and by $\mathscr{G}:=\mathcal{A}[p^\infty]$ the universal $p$-divisible group, equipped with $G$-structures. 

We let $\mu$ denote the (conjugacy class of) Hodge cocharacters of our Shimura datum viewed as a geometric dominant cocharacter of $G$ via the fixed choice of isomorphism $\ol{\mathbb{Q}}_{p} \simeq \bb{C}$. We fix $b\in B(G,\mu^{-1})$ and a completely slope divisible $p$-divisible group with $G$-structures $\mathbb{X}/\overline{\mathbb{F}}_p$, in the isogeny class defined by $b$ and compatible with $\mu$ (in the sense of satisfying the Kottwitz determinant condition). 

\subsubsection{Igusa varieties}{\label{ss: IgusaVarieties}} We consider the Oort central leaf $\mathscr{C}:=\mathscr{C}^{\mathbb{X}}\subset \mathscr{S}$, which is the smooth subscheme of $\mathscr{S}$ where the universal $p$-divisible group $\mathscr{G}$ is geometrically fiberwise isomorphic to $\mathbb{X}$, compatibly with the $G$-structures. Over $\mathscr{C}$, we have a perfect scheme $\mathrm{Ig}^b$ (the \emph{perfect Igusa variety}), which parameterizes trivializations 
\[
\gamma_\infty:
\mathbb{X}\times_{\overline{\mathbb{F}}_p} \mathrm{Ig}^b
\stackrel{\sim}{\to}
\mathscr{G}\times_{\mathscr{C}}\mathrm{Ig}^b 
\]
compatible with the $G$-structures on both sides. 
Explicitly, such a trivialization is required to be compatible with the $\mathcal{O}_B\otimes_{\mathbb{Z}}\mathbb{Z}_p$-action and with the polarizations up to an element of $\underline{\mathbb{Z}_p^\times} (\mathrm{Ig}^b)$, the global sections of the sheaf $\underline{\bb{Z}_{p}^\times}$ on the perfect Igusa variety $\Ig^{b}$. There is an action of the $\overline{\mathbb{F}}_p$-group scheme $\underline{\mathrm{Aut}}_{G}(\mathbb{X})$ of automorphisms of $\mathbb{X}$ that respect the $G$-structure on $\mathrm{Ig}^b$ given by pre-composing $\gamma_\infty$. The morphism $\mathrm{Ig}^b\to \mathscr{C}$ is an fpqc torsor under $\underline{\mathrm{Aut}}_{G}(\mathbb{X})$, by~\cite[Corollary 2.3.2]{CS2}.  

We also introduce a variant of finite-level Igusa varieties, denoted $\Ig^{b}_{m}$, where we trivialize only the $p^m$-torsion in the universal $p$-divisible group $\mathscr{G}$. For each $m\in \mathbb{Z}_{\geq 1}$, we consider the moduli problem on $\mathscr{C}$-schemes $S$ that parameterizes trivializations 
 \begin{equation}{\label{eqn: trivializations}}
 \gamma_m: \mathbb{X}[p^m]\times_{\overline{\mathbb{F}}_p} S \stackrel{\sim}{\to}
  \mathscr{G}[p^m]\times_{\mathscr{C}}S
 \end{equation}
 that (fppf locally on $S$) lift to arbitrary $m'\geq m$, that commute with the $\mathcal{O}_B\otimes_{\mathbb{Z}}\mathbb{Z}_p$-action, and that commute with the polarization up to an element of $\underline{(\mathbb{Z}/p^m\mathbb{Z})^\times}(S)$. 

\begin{lemma}\label{lem:automorphisms at level m} For any $m\in \mathbb{Z}_{\geq 1}$, there exists a finite group scheme $H_m$ over $\ol{\mathbb{F}}_p$, characterized by the following properties: 
\begin{itemize}
\item $H_m\subset \underline{\mathrm{Aut}}_G(\mathbb{X}[p^m])$ is an affine closed subgroup;
\item the natural map $ \underline{\mathrm{Aut}}_G(\mathbb{X})\to  \underline{\mathrm{Aut}}_G(\mathbb{X}[p^m])$ factors through $H_m$;
\item the map $ \underline{\mathrm{Aut}}_G(\mathbb{X})\to H_m$ is faithfully flat. 
\end{itemize}
\end{lemma}

\begin{proof}
    This is proved in the same way as~\cite[Lemma 5.5.3]{DAVH}, which is the case $m=2$. 
\end{proof}

\begin{lemma}\label{lem:representability}\leavevmode
\begin{enumerate}
\item For any $m\in \mathbb{Z}_{\geq 1}$, the moduli problem of trivializations $\gamma_m$ as in (\ref{eqn: trivializations}) is representable by an Igusa variety $\mathrm{Ig}^b_m$, which is finite and faithfully flat over $\mathscr{C}$. 
\item For $m\geq 2$, the Igusa variety $\mathrm{Ig}^b_m$ is smooth over $\ol{\mathbb{F}}_p$. 
\end{enumerate}
\end{lemma}

\begin{proof} Using the fact that $\mathrm{Ig}^b\to \mathscr{C}$ is an $\underline{\mathrm{Aut}}_G(\mathbb{X})$-torsor, we define $\mathrm{Ig}^b_m$ as the pushout
\[
\mathrm{Ig}^b_m:=\left(\mathrm{Ig}^b\times H_m\right)/\underline{\mathrm{Aut}}_G(\mathbb{X}),
\]
where the action of $\underline{\mathrm{Aut}}_G(\mathbb{X})$ is via the diagonal. As $\mathrm{Ig}^b\to \mathscr{C}$ is an fpqc cover, we can check finiteness and faithful flatness after pullback to $\mathrm{Ig}^b$, where we can appeal to Lemma~\ref{lem:automorphisms at level m}. 

For the second part, as $\mathrm{Ig}^b_m$ is of finite type over $\ol{\mathbb{F}}_p$, it is enough to prove that its completed local rings at closed geometric points are formally smooth. When $m\geq 2$, the argument in~\cite[\S 5.5]{DAVH}, particularly Construction 5.5.6, shows that the completed local ring $(\mathrm{Ig}^b_m)^{\hat{}}_{x}$ of $\Ig^{b}_{m}$ at a closed geometric point $x$ can be identified with $Z(p^m\mathfrak{b}^+)$, in the notation of \emph{loc. cit.}. Here, $\mathfrak{b}^+$ is a completely slope divisible and nilpotent Dieudonn\'e--Lie $\breve{\mathbb{Z}}_p$-algebra in the sense of~\cite[\S 4.2]{DAVH}, constructed as in Example 5.5.8 of \emph{loc. cit.}. Whenever $m\geq 2$, the argument in Lemma 4.2.10 of \emph{loc. cit.} shows that  $p^m\mathfrak{b}^+$ is furthermore integrable, and then $Z(p^m\mathfrak{b}^+) \simeq \tilde{\Pi}(\mathfrak{b})/\Pi(p^m\mathfrak{b}^+)$ is the formal Lie variety assigned to $p^m\mathfrak{b}^+$ by Corollary 4.3.9. As formal Lie varieties are, in particular, formally smooth, we conclude. 
 \end{proof}

  The perfect Igusa variety $\mathrm{Ig}^b$ can be identified with the inverse limit of the system $(\mathrm{Ig}^b_m)_{m\in \mathbb{Z}_{\geq 1}}$ of finite-level Igusa varieties.

  \begin{remark}\label{rem: comparison with Mantovan} In the literature on Igusa varieties, for example in~\cite{mantovan-PEL} and~\cite{CS17}, it is common to consider a variant of finite-level Igusa varieties that are finite \'etale covers of the leaf $\mathscr{C}$. The moduli problem for these finite \'etale covers is given in terms of trivializations of the $p^m$-torsion on each graded piece of the slope filtration on $\mathscr{G}$, compatibly with the extra structures. For our purposes, it is more convenient to trivialize the full $\mathscr{G}[p^m]$. We denote these latter finite level Igusa varieties by $\Ig^{b}_{m,\mathrm{Mant}}$. The difference between the two moduli problems lies in splitting the slope filtration, which amounts to taking a purely inseparable cover of the finite \'etale cover of $\mathscr{C}$. In particular, we note that the inverse limit $\lim_{m \geq 1} \Ig^{b}_{m,\mathrm{Mant}}$ gives rise to $\Ig^{b}$ after taking perfections (see \cite[Remark~2.3.7.]{CS2} and \cite[Proposition~4.3.8]{CS17}) and we have an equality on perfections $(\Ig^{b}_{m,\mathrm{Mant}})^{\perf} = (\Ig^{b}_{m})^{\perf}$.
  \end{remark}

We denote by $\widetilde{\mathbb{X}}$ the universal cover of $\mathbb{X}$, in the sense of~\cite[\S 3.1]{Scholze-Weinstein}. We recall the following construction from~\cite[Definition 9.15]{zhang2023}, which gives the special fiber of the formal group scheme considered in~\cite[Definition 4.2.9]{CS17}. 

\begin{definition}{\label{defn: automorphismsfpqcsheaf}}
 We define $\widetilde{J}_{b, \ol{\mathbb{F}}_p}=\underline{\mathrm{Aut}}_{G}(\widetilde{\mathbb{X}})$ as the fpqc sheaf on the opposite category of $\ol{\mathbb{F}}_p$-algebras defined by 
 \[
 R\mapsto \underline{\mathrm{Aut}}_G(\widetilde{\mathbb{X}}\times_{\ol{\mathbb{F}}_p} R),
 \]
 where the latter is the group of automorphisms of $\widetilde{\mathbb{X}}\times_{\ol{\mathbb{F}}_p} R$ that preserve the $G$-structures, in the sense that they commute with the $\mc{O}_B$-endomorphisms, and preserve the polarization up to an element of $\underline{\mathrm{Aut}}(\tilde{\mu}_{p^\infty})(R)= \underline{\mathbb{Q}^\times_p}(R)$. 
\end{definition}

It follows from~\cite[Lemma 4.2.10]{CS17} that $\widetilde{J}_{b,\ol{\mathbb{F}}_p}$ is representable by a formal group scheme over $\ol{\mathbb{F}}_p$. Furthermore, for an $\ol{\mathbb{F}}_p$-scheme $S$, we can identify $\widetilde{J}_{b, \ol{\mathbb{F}}_p}(S)$ with the group of quasi-self isogenies of $\mathbb{X}\times_{\ol{\mathbb{F}}_p} S$ that preserve the $G$-structures, in the sense that they commute with the endomorphisms, and preserve the polarization up to an element of $\underline{\mathbb{Q}_p^\times}$. 

There is a natural \emph{right} action of $\widetilde{J}_{b, \ol{\mathbb{F}}_p}$ on $\mathrm{Ig}^b$ that extends the action of $\underline{\mathrm{Aut}}_G(\mathbb{X})$. This can be constructed as follows. One proves, as in~\cite[Lemma 4.3.4]{CS17} and~\cite[Lemma 4.2.2]{Caraiani-Tamiozzo} (which provides a number of details of the argument), that $\mathrm{Ig}^b$ also represents a different moduli problem, where the $p$-divisible group is trivialized via a quasi-isogeny rather than via an isomorphism, and where the isomorphism classes of objects are given by $p$-power quasi-isogenies. 

For the sake of precision, we make this statement explicit. One can construct a functorial bijection between $\mathrm{Ig}^b(S)$, for an $\overline{\mathbb{F}}_p$-scheme $S$, and the set of isomorphism classes of pairs $(A,\gamma)$, where $A/S$ is an abelian scheme equipped with $\mathsf{G}$-structures and with a $K^p$-level structure\footnote{Note, however, that in this alternative moduli problem, the abelian scheme $A/S$ is no longer required to satisfy the Kottwitz determinant condition.}, and $\gamma$ is a quasi-isogeny 
\[
\gamma: \mathbb{X}\times_{\overline{\mathbb{F}}_p} S \to A[p^\infty]
\]
that respects the $\mathsf{G}$-structures, in the sense that the polarization is respected up to an element of $\underline{\mathbb{Q}_p^\times}(S)$. Two such pairs $(A, \gamma)$ and $(A', \gamma')$ are considered isomorphic if there exists a $p$-power quasi-isogeny $A\to A'$ that commutes with the level structures, with the quasi-isogenies $\gamma$ and $\gamma'$, and with $\mathsf{G}$-structures, in the sense that the polarizations are respected up to an integral power of $p$ on each connected component of $S$. 
This alternative moduli-theoretic description makes it clear that there is a right action of $\widetilde{J}_{b,\overline{\mathbb{F}}_p}$ on $\mathrm{Ig}^b$, given by pre-composition: $\rho\in \widetilde{J}_{b, \overline{\mathbb{F}}_p}(S)$ acts via $\gamma \mapsto \gamma\circ \rho$.

By chasing through the proof that the Igusa moduli problem via quasi-isogenies is equivalent to the one via isomorphisms, e.g. as in the proof of~\cite[Lemma 4.2.2]{Caraiani-Tamiozzo}, one sees that the action of $\widetilde{J}_{b, \overline{\mathbb{F}}_p}$ on $\mathrm{Ig}^b$ described above is equivalent to the one that can be deduced from~\cite[\S 4]{mantovan-PEL}. Indeed, assume that we have a trivialization 
\[
\gamma: \mathbb{X}\times_{\overline{\mathbb{F}}_p} S \stackrel{\sim}{\to} A[p^\infty],
\]
which is an isomorphism (not just a quasi-isogeny), and that $\rho\in \widetilde{J}_{b, \overline{\mathbb{F}}_p}(S)$ is a genuine isogeny. For some $m\geq 0$, there exists an isogeny $\rho': \mathbb{X}\to \mathbb{X}$, compatible, as usual, with the $G(\mathbb{Q}_p)$-structures and such that $\rho\circ \rho' = \rho'\circ \rho = [p^m]$.  
Then we can define a finite flat subgroup scheme $K_\rho:=\gamma(\mathrm{ker}\rho')$ of $A$ and form the quotient $B:=A/K_\rho$. The new abelian scheme $B$ over $S$ can be shown to inherit $\mathsf{G}$-structures and a $K^p$-level structure from the $\mathsf{G}$-structures and $K^p$-level structure of $A$, and is furthermore equipped with an induced isomorphism 
$\gamma': \mathbb{X}\times_{\overline{\mathbb{F}}_p} S\stackrel{\sim}{\to} B[p^\infty]$
that respects these structures, in the usual sense. The data of the abelian scheme $B$, its additional structures and the trivialization $\gamma'$ determine another $S$-valued point of $\mathrm{Ig}^b$. This is the same construction as the one in~\cite[\S 4]{mantovan-PEL}, up to the fact that Mantovan restricts to those $\rho$ that are contained in a certain submonoid of $J_b(\mathbb{Q}_p)$. Furthermore, we have a commutative diagram 
\begin{equation}
\begin{tikzcd}
    \mathbb{X}\times_{\overline{\mathbb{F}}_p}S \arrow[r, "\gamma'"] \arrow[d, "\rho"] & B[p^\infty]\arrow[d] \\ 
      \mathbb{X}\times_{\overline{\mathbb{F}}_p}S \arrow[r, "\gamma"] & A[p^\infty],
\end{tikzcd}
\end{equation}
where the horizontal maps are isomorphisms and the vertical maps are $p$-power quasi-isogenies. This shows that the new $S$-valued point of $\mathrm{Ig}^b$ is equivalent to the one determined by $\gamma\circ \rho$ when we consider the moduli description of $\mathrm{Ig}^b(S)$ using quasi-isogenies. 

\begin{remark}\label{rem: comparison with Mantovan Jb} The reason for the restriction to a submonoid of $J_b(\mathbb{Q}_p)$ in Mantovan's work is that she does not pass to perfection on the level of the Igusa varieties, so the slope filtration on the universal $p$-divisible group is not split (cf. Remark \ref{rem: comparison with Mantovan}). However, the formation of $K_{\rho}$, when $\rho$ is an element of a certain submonoid of $J_b(\mathbb{Q}_p)$, is still unambiguous, even when the slope filtration is not split. 
\end{remark}

\subsubsection{Partial toroidal compactifications}

In~\cite[\S 3.2]{CS2} and~\cite[\S 3.2]{santos}, partial toroidal compactifications of the Igusa varieties considered above are constructed. This relies on the construction of toroidal compactifications of canonical integral models of Shimura varieties of PEL type in~\cite{lan-thesis} and, in particular, on the auxiliary choice of a compatible system of cone decompositions. 

We assume, for simplicity, that the level $K^p$ is principal. Recall that $K_p$ was chosen to be hyperspecial. We have a notion of equivalence classes of cusp labels at level $K^pK_p$, cf.~\cite[\S 5.4]{lan-thesis}. We choose a compatible system of cone decompositions $\Sigma = \{\Sigma_Z\}_{Z}$ indexed by these equivalence classes of cusp labels at level $K^pK_p$. Assume that $\Sigma$ is a \emph{good} choice of compatible system, in the sense that it satisfies the conditions in~\cite[Remark 2.5.6]{CS2}. 

Given a good choice of $\Sigma$, we have an open dense embedding
\[
\mathscr{S}\hookrightarrow \mathscr{S}_{\Sigma}^{\mathrm{tor}},
\]
with the target a smooth and projective scheme over $\overline{\mathbb{F}}_p$. The universal abelian scheme $\mathcal{A}$ over $\mathscr{S}$ extends to a semi-abelian scheme over $\mathscr{S}_{\Sigma}^{\mathrm{tor}}$, which we denote by $\mathcal{A}_{\Sigma}$. The toroidal compactification $\mathscr{S}_{\Sigma}^{\mathrm{tor}}$ admits a stratification in terms of (equivalence classes of) cusp labels, where the strata essentially measure the degeneracy of the semi-abelian scheme $\mathcal{A}_{\Sigma}$. This refines the decomposition 
\begin{equation}{\label{eqn: decompositionoftoroidalintoparabolics}}
\mathscr{S}^{\mathrm{tor}}_{\Sigma} = \bigsqcup_{\PP} \mathscr{S}^{\mathrm{tor}}_{\Sigma, \PP},
\end{equation}
where $\PP$ runs through $\mathsf{G}(\mathbb{Q})$-conjugacy classes of admissible rational parabolic subgroups of $\mathsf{G}$, as in Definition \ref{defn: admissibleparabolic}. In fact, a cusp label $Z$ determines a split symplectic admissible filtration (\cite[Section~5.2.1]{lan-thesis}), which gives rise to a conjugacy class of parabolics $\PP$ by looking at the stabilizer of the filtration.\footnote{We note that this filtration has a rational structure, as follows from~\cite[Prop. A.5.9]{Lan-ANT} for the generic fiber of the Shimura variety, together with~\cite[Lemma A.4.7]{Lan-ANT}, which reduces to the generic fiber.}

The formal completion of $\mathscr{S}_{\Sigma}^{\mathrm{tor}}$ along a boundary stratum determined by a cusp label $Z$ can be described via a toroidal boundary chart. The properties of toroidal boundary charts are stated precisely, for example, in~\cite[Theorem 2.5.9]{CS2} and \cite[Theorem 6.4.1.1]{lan-thesis}. 

Given the choice of $b\in B(G,\mu^{-1})$ and a completely slope divisible $p$-divisible group $\mathbb{X}/\ol{\mathbb{F}}_p$ in the isogeny class determined by $b$, Lan and Stroh define a partial toroidal compactification
\[
\mathscr{C}\hookrightarrow \mathscr{C}_{\Sigma}^{\mathrm{tor}}
\]
of the corresponding Oort central leaf in~\cite[\S 3.4]{LanStroh}, by taking, in a precise sense, its closure towards the boundary of $\mathscr{S}_{\Sigma}^{\mathrm{tor}}$. The partial toroidal compactification inherits the stratification in terms of cusp labels and its formal completion along a boundary stratum can also be described via a toroidal boundary chart. This is simply the pullback of the toroidal boundary chart from the special fiber of the Shimura variety to the Oort central leaf. 

\begin{remark} The fact that the structure of toroidal compactifications for Oort central leaves is analogous to the structure of toroidal compactification for the Shimura variety is a consequence of the fact that Oort central leaves are \emph{well-positioned} in the sense of~\cite[\S 2.2]{LanStroh} (see also~\cite[\S 3.4]{Boxer}). 
\end{remark}

We now explain the result about toroidal compactifications of \emph{perfect} Igusa varieties, and then comment on the finite-level variants. For a good choice of $\Sigma$, the partial toroidal compactification $\mathrm{Ig}^{b,\mathrm{tor}}_{\Sigma}$ is a perfect scheme over $\overline{\mathbb{F}}_p$ that fits in a cartesian diagram 
\[
\begin{tikzcd}
\mathrm{Ig}^b \arrow[d] \arrow[r, hook] & \mathrm{Ig}^{b,\mathrm{tor}}_{\Sigma}\arrow[d] \\ 
\mathscr{C} \arrow[r, hook] & \mathscr{C}^{\mathrm{tor}}_{\Sigma}.
\end{tikzcd}
\]
Pulling back the stratification~\eqref{eqn: decompositionoftoroidalintoparabolics} from $\mathscr{S}_\Sigma^{\tor}$, the perfect scheme $\mathrm{Ig}^{b,\mathrm{tor}}_{\Sigma}$ has a stratification
\begin{equation}{\label{eqn: BoudnaryStrataofToroidal}}
\mathrm{Ig}^{b,\mathrm{tor}}_{\Sigma} = \bigsqcup_{\PP}
\mathrm{Ig}^{b,\mathrm{tor}}_{\Sigma, \PP},
\end{equation}
where $\PP$ runs through $\mathsf{G}(\mathbb{Q})$-conjugacy classes of admissible rational parabolic subgroups of $\mathsf{G}$.  

It also inherits a stratification in terms of cusp labels from $\mathscr{C}^{\mathrm{tor}}_{\Sigma}$, and this can be further refined in terms of \emph{Igusa cusp labels}. We review this notion, since it will be used in the next section, see also~\cite[Definition 3.2.15]{santos}, cf.~\cite[Definition 3.3.10]{CS2} for the special case when the polarizations that are part of our PEL moduli problem are principal.

Note that, because we are assuming $p$ is an unramified prime for the Shimura datum, the $p$-divisible group $\mathbb{X}$ is equipped with a principal polarization. 

\begin{definition}\label{defn: IgusaCuspLabel} An Igusa cusp label is an equivalence class of tuples 
\[
\widetilde{Z} = (Z^p, Z_b, (X,Y,\phi, \varphi_{-2},\varphi_{0}), \psi, \delta^p, \delta_b),
\] 
where
\begin{enumerate}
\item the tuple $Z:=(Z^p, (X,Y,\phi, \varphi_{-2},\varphi_{0}),\delta^p)$ is a (representative of a) cusp label in the sense of~\cite[\S 5.4.1]{lan-thesis}\footnote{Note that, unlike \emph{loc. cit.}, we are only taking the prime-to-$p$ part here. However, this agrees with equivalence classes of cusp labels with hyperspecial level at $p$ in the sense of~\cite{lan-thesis}, via the Iwasawa decomposition - see~\cite[Lemma A.4.7]{Lan-ANT}}. 
\item $Z_b$ is an $\mathcal{O}_{B}$-stable filtration on $\mathbb{X}$, of the form 
\[
Z_{b,-2}\subset Z_{b,-1}\subset Z_{b,0}:=\mathbb{X},
\] 
where $\mathrm{Gr}_{-2}^{\mathbb{X}} = Z_{b,-2}$ is multiplicative and $\mathrm{Gr}_{0}^{\mathbb{X}} = Z_{b,0}/Z_{b,-1}$ is identified with its Cartier dual under the principal polarization on $\mathbb{X}$; in particular, $\mathrm{Gr}_0^{\mathbb{X}}$ is \'etale. 
\item $\psi$ is an $\mathcal{O}_{B}$-equivariant isomorphism 
$\mathrm{Gr}_0^{\mathbb{X}}\stackrel{\sim}{\to} Y\otimes \mathbb{Q}_p/\mathbb{Z}_p$. 
\item $\delta_b$ is an $\mathcal{O}_{B}$-linear splitting of the filtration corresponding to $Z_b$, i.e. an $\mathcal{O}_B$-linear isomorphism  
\[
\delta_b:\mathbb{X}\simeq \bigoplus_{i=-2}^0 \mathrm{Gr}_i^{\mathbb{X}}.
\]
We furthermore assume that $\delta_b$ is compatible with the principal polarization $\lambda$ on $\mathbb{X}$; see Lemma~\ref{lem:splitting anti-diagonal} below for a proof that we can always construct such a splitting. 
\end{enumerate}
Two such tuples $\widetilde{Z}$ and $\widetilde{Z}'$ are equivalent if the filtrations $Z^p$ and $(Z')^{p}$, respectively $Z_b$ and $Z'_b$, agree, and if there exists a pair of $\mathcal{O}_B$-linear isomorphisms $f_X:X'\stackrel{\sim}{\to} X$ and $f_Y:Y\stackrel{\sim}{\to}Y'$ such that $\phi = f_X\circ \phi'\circ f_Y$, $\varphi'_{-2} =f_X^{\vee}\circ \varphi_{-2} $, $\varphi'_{0} =f_Y\circ \varphi_0$, and $\psi' = f_Y\circ \psi$. There is no additional constraint relating $\delta_b$ and $\delta'_b$. 
\end{definition}

\begin{lemma}\label{lem:splitting anti-diagonal} There exists a choice of splitting 
\[
\delta_b:\mathbb{X}\simeq \bigoplus_{i=-2}^0 \mathrm{Gr}_i^{\mathbb{X}}
\]
which is $\mathcal{O}_B$-linear and compatible with the anti-symmetric principal polarization $\lambda$ on $\mathbb{X}$. Via $\delta_b$, $\lambda$ can be written as an anti-diagonal matrix. 
\end{lemma}

\begin{proof} 
As the slope decomposition on $\mathbb{X} / \overline{\mathbb{F}}_p$ respects both the $\mathcal{O}_B$-action and the polarization $\lambda$, it is enough to construct an $\mathcal{O}_B$-linear splitting of the surjection 
\[
\mathbb{X}^{\et} \twoheadrightarrow \mathrm{Gr}_0^{\mathbb{X}} \stackrel{\psi}{\simeq} Y\otimes \mathbb{Q}_p/\mathbb{Z}_p. 
\]
The polarization $\lambda$ will then induce, by Cartier duality, a corresponding $\mathcal{O}_B$-linear splitting of the injection 
$\mathrm{Gr}^{\mathbb{X}}_{-2} \hookrightarrow \mathbb{X}^{\mu}$
of multiplicative $p$-divisible groups. 

Recall that $p$ was chosen to be an unramified prime for the Shimura datum $\gx$. Therefore, $\mathcal{O}_B$ is a maximal order, hence hereditary; the $\mathcal{O}_B$ lattice $Y\otimes \mathbb{Z}_p$ is a projective $\mathcal{O}_B$-module. This means that the surjection
\[
T_p\mathbb{X}^{\et}\twoheadrightarrow T_p\mathrm{Gr}_0^{\mathbb{X}} \simeq Y\otimes \mathbb{Z}_p
\]
on the level of Tate modules admits an $\mathcal{O}_B$-linear splitting. This induces the desired $\mathcal{O}_B$-linear splitting on the level of \'etale $p$-divisible groups.     
\end{proof}

There is an action of $G(\A_f^p)\times J_b(\mathbb{Q}_p)$ on the set of Igusa cusp labels. This is a variant of the Hecke action at good primes on usual cusp labels described explicitly in~\cite[\S 5.4.3]{lan-thesis}. Indeed, one can formulate a variant of Definition~\ref{defn: IgusaCuspLabel} after extending scalars to $\mathbb{Q}$, and there the group action is obvious. One then checks that Igusa cusp labels defined integrally are in bijection with Igusa cusp labels defined rationally. 
This group action allows us to define the notion of \emph{Igusa cusp label at level $K$} for a compact open subgroup $K\subset G(\mathbb{A}_f^p)\times J_b(\mathbb{Q}_p)$, or even at level $H$ for a closed subgroup $H\subset G(\mathbb{A}_f^p)\times J_b(\mathbb{Q}_p)$, by taking orbits under these groups. Note that the group action preserves the rational conjugacy class of admissible parabolic subgroup $\PP$ determined by the Igusa cusp label. 

For each good choice of $\Sigma$ and $m\in \mathbb{Z}_{\geq 1}$, one can also construct a finite-level toroidal compactification $\mathrm{Ig}^{b,\mathrm{tor}}_{\Sigma, m}$ of $\mathrm{Ig}^{b}_m$, for example by taking the relative normalization of $\mathscr{C}^{\mathrm{tor}}_{\Sigma}$ in the Igusa variety $\mathrm{Ig}^b_m$, as defined in \S~\ref{ss: IgusaVarieties}. We have a natural isomorphism 
\[
\mathrm{Ig}^{b,\mathrm{tor}}_{\Sigma} \stackrel{\sim}{\to}
\varprojlim_m\ \mathrm{Ig}^{b,\mathrm{tor}}_{\Sigma,m}.
\]

\begin{remark}
The Igusa variety $\mathrm{Ig}^{b,\mathrm{tor}}_{\Sigma, m}$ is a variant of the finite-level toroidal compactification constructed in~\cite{CS2, santos}. The difference comes from the fact that the finite-level Igusa varieties considered in \emph{loc. cit.} only trivialize the graded pieces of the slope filtration on $\mathscr{G}[p^m]$, whereas $\mathrm{Ig}^b_m$ trivializes the entire $\mathscr{G}[p^m]$, see Remark~\ref{rem: comparison with Mantovan}. In particular, their perfections agree, and they give rise to the same underlying topological space. 
\end{remark}

 Furthermore, each $\mathrm{Ig}^{b,\mathrm{tor}}_{\Sigma,m}$ has a locally closed stratification with boundary strata indexed by Igusa cusp labels at level $K^pK_b(p^m)$, where $K_b(p^m)$ is the level $m$ congruence subgroup in $\underline{\mathrm{Aut}}_G(\mathbb{X})(\overline{\mathbb{F}}_p)$, and therefore a compact open subgroup of $J_b(\mathbb{Q}_p)$. For each such Igusa cusp label, there is also a description via toroidal boundary charts analogous to the ones that arise in the setting of Shimura varieties. See~\cite[Theorem 3.3.12]{CS2}, respectively~\cite[Theorem 3.2.18]{santos}, for the precise description of each toroidal boundary chart in the case of the standard moduli problem. The arguments there can be adapted in our case to give Theorem~\ref{thm:IgusaStratumDiagram} below, which applies specifically to the boundary strata of $\mathrm{Ig}^{b,\mathrm{tor}}_{\Sigma, m}$. As in the work of Faltings--Chai and Lan, one key input to constructing the isomorphism~\eqref{eq:boundary chart iso} of Theorem~\ref{thm:IgusaStratumDiagram} is the smoothness, hence normality, of $\mathrm{Ig}^b_{m}$ and of the corresponding boundary charts, which follows from Lemma~\ref{lem:representability}.   

The choice of Igusa cusp label $\widetilde{Z}$ induces a decomposition 
\begin{equation}\label{eq:IgusaCuspLabelDecomposition}
\mathbb{X}\simeq \mathbb{X}_Y \oplus \mathbb{X}_{\widetilde{Z}} \oplus \mathbb{X}_X
\end{equation}
over $\overline{\mathbb{F}}_p$, where $\mathbb{X}_Y:=Y\otimes (\mathbb{Q}_p/\mathbb{Z}_p)$ is \'etale, $\mathbb{X}_X:=\mathrm{Hom}(X,\mu_{p^\infty})$ is multiplicative, and $\mathbb{X}_{\widetilde{Z}}$ is isomorphic to the (constant) geometric fibers of the $p$-divisible group of the universal abelian scheme over a smaller Igusa variety $\mathrm{Ig}^b_{\widetilde{Z}}$. In the notation of Definition~\ref{defn: IgusaCuspLabel}, we identify $\mathbb{X}_Y\simeq \mathrm{Gr}_0^{\mathbb{X}}$ using $\psi$ and hence $\mathbb{X}_X\simeq \mathrm{Gr}^{\mathbb{X}}_{-2}$ using Cartier duality and the fact that $\phi\otimes \mathbb{Z}_p$ is an isomorphism.

\begin{theorem}\label{thm:IgusaStratumDiagram}
   Fix a tame principal level $K^p=K(N)$, $p\nmid N$. For each Igusa cusp label $\widetilde{Z}$ at level $K^pK_b(p^m)$, we have a diagram 
   \begin{equation}\label{eq:IgusaStratumDiagram}
   \begin{tikzcd}
       \Xi_{\widetilde{Z}}\arrow[d]\arrow[r] & \Xi_{\widetilde{Z}, \Sigma_{\widetilde{Z}}}\arrow[dl] \\
       C_{\widetilde{Z}}\arrow[d] &  \\
       \mathrm{Ig}^b_{\widetilde{Z}} & ,
   \end{tikzcd}
   \end{equation}
   where $\mathrm{Ig}^b_{\widetilde{Z}}$ denotes a level $m$ smaller Igusa variety determined by the Igusa cusp label $\widetilde{Z}$ (in the sense of our moduli problem, i.e. which trivializes the full $p^m$-level structure), 
    $C_{\widetilde{Z}}\to \mathrm{Ig}^b_{\widetilde{Z}}$ is an abelian scheme, $\Xi_{\widetilde{Z}}\to C_{\widetilde{Z}}$ is a torus torsor, and $\Xi_{\widetilde{Z}} \hookrightarrow \Xi_{\widetilde{Z}, \Sigma_{\widetilde{Z}}}$ is a relative torus embedding. There is an action of $\Gamma_{\widetilde{Z}}$ on this diagram, and, letting $\mathfrak{X}_{\widetilde{Z}, \Sigma_{\widetilde{Z}}}$ denote the formal completion of $\Xi_{\widetilde{Z}, \Sigma_{\widetilde{Z}}}$ along its toroidal boundary, there is an isomorphism 
   \begin{equation}\label{eq:boundary chart iso}
   \mathfrak{X}_{\widetilde{Z}, \Sigma_{\widetilde{Z}}} / \Gamma_{\widetilde{Z}} \simeq \widehat{\mathrm{Ig}}^{b,\mathrm{tor}}_{m,\widetilde{Z}}.
   \end{equation}
\end{theorem}

 Given an Igusa cusp label $\widetilde{Z}$ at level $K^pK_b(p^m)$, we let $Z$ denote the underlying cusp label for the Shimura variety obtained from $\widetilde{Z}$ by forgetting the part at $p$. We let $\mathscr{C}_Z$ denote the corresponding boundary stratum in the partial minimal compactification of the underlying leaf. Over $\mathscr{C}_{Z}$, we have a similar diagram 
 \begin{equation}\label{eq:LeafStratumDiagram}
      \begin{tikzcd}
       \Xi_{Z}\arrow[d]\arrow[r] & \Xi_{Z, \Sigma_Z}\arrow[dl] \\
       C_{Z}\arrow[d] & \\
       \mathscr{C}_{Z} &, 
   \end{tikzcd}
 \end{equation}  
   in which the objects are pulled back from the corresponding boundary stratum in the Shimura variety (see \cite[Theorem~6.4.1.1]{lan-thesis}).
 
 We now give a precise description of each term appearing in Diagram~\eqref{eq:IgusaStratumDiagram}. The moduli-theoretic data below is compatible with the PEL structures, in the usual sense. 

 \begin{enumerate}
     \item Over $\mathscr{C}_{Z}$, there is a universal abelian scheme $\mathcal{B}_{Z}$. Its pullback $\mathcal{B}_{\widetilde{Z}}$ to $\mathrm{Ig}^b_{\widetilde{Z}}$ is equipped with a trivialization 
     \[
     \alpha_m: \mathbb{X}_{\widetilde{Z}}[p^m]\times_{\overline{\mathbb{F}}_p} \mathrm{Ig}^b_{\widetilde{Z}} \stackrel{\sim}{\to} \mathcal{B}_{\widetilde{Z}}[p^m]. 
     \]
     Furthermore, the morphism $\mathrm{Ig}^b_{\widetilde{Z}}\to \mathscr{C}_Z$ is universal for the existence of such trivializations. 
     \item Over $C_{Z}$, there is a Raynaud extension 
     \[
     0\to T_Z\to \mathcal{G}_Z \to \mathcal{B}_Z\to  0,
     \]
     where $T_Z:=T\times_{\overline{\mathbb{F}}_p} C_Z$ and $T$ is the split torus with character group $X$. The Raynaud extension is parameterized by an element $f_0\in \mathrm{Hom}(\underline{X}, \mathcal{B}^\vee_{Z})$, where the homomorphisms are taken as \'etale sheaves of groups over $\mathscr{C}_Z$. The Raynaud extension induces a short exact sequence 
     \[
     0\to \mathbb{X}_{X}[p^m]\times_{\overline{\mathbb{F}}_p}C_Z\to \mathcal{G}_Z[p^m]\to \mathcal{B}_{Z}[p^m]\to 0.
     \]
     The pullback of this short exact sequence to $C_{\widetilde{Z}}$ is equipped with a splitting 
     \[
     \beta_m: \mathcal{B}_{\widetilde{Z}}[p^m]\to \mathcal{G}_{\widetilde{Z}}[p^m]
     \] 
     that is parameterized by an extension $f_{0,m}: \frac{1}{p^m} X\to \mathcal{B}_{\widetilde{Z}}$ of $f_0$. The set of possible such extensions is a torsor under $\mathrm{Hom}\left(\underline{X/p^mX}, \mathcal{B}^\vee_{\widetilde{Z}}[p^m]\right)$, where the homomorphisms are taken as \'etale sheaves of groups over $C_Z\times_{\mathscr{C}_Z}\mathrm{Ig}^b_{\widetilde{Z}}$. This space of homomorphisms can be identified with $\mathrm{Hom}\left(\mathcal{B}_{\widetilde{Z}}[p^m], \mathbb{X}_{X}[p^m]\times_{\overline{\mathbb{F}}_p}\left(C_Z\times_{\mathscr{C}_Z}\mathrm{Ig}^b_{\widetilde{Z}}\right)\right)$. 
     
There is also dual data, parameterized by an element $g_0\in \mathrm{Hom}(\underline{Y}, \mathcal{B}_Z)$ and by an extension $g_{0,m}:\underline{\frac{1}{p^m} Y}\to \mathcal{B}_{\widetilde{Z}}$. The space of such extensions is a torsor under $\mathrm{Hom}\left(\mathbb{X}_Y[p^m]\times_{\overline{\mathbb{F}}_p} \left(C_Z\times_{\mathscr{C}_Z}\mathrm{Ig}^b_{\widetilde{Z}}\right),\mathcal{B}_{\widetilde{Z}}[p^m]\right)$. The choices of $f_{0,m}$ and $g_{0,m}$ are equivalent under $\phi: Y\hookrightarrow X$ and the prime-to-$p$ polarization on $\mathcal{B}_{\widetilde{Z}}$. Explicitly, we have $\lambda_{\mathcal{B}_{\widetilde{Z}}}\circ g_{0,m} = f_{0,m}\circ (\frac{1}{p^m}\phi)$.    

The morphism $C_{\widetilde{Z}}\to C_Z\times_{\mathscr{C}_Z}\mathrm{Ig}^b_{\widetilde{Z}}$ is universal for the existence of splittings $\beta_m$ or equivalently extensions $g_{0,m}$ of $g_0$.  

\item Over $\Xi_Z$, there is a lift $g:\underline{Y}\to \mathcal{G}_{Z}$ of $g_0$ that satisfies a certain symmetry condition. Its pullback to $\Xi_{\widetilde{Z}}$ has an extension to a lift $g_m:\underline{\frac{1}{p^m} Y}\to \mathcal{G}_{\widetilde{Z}}$ of $g_{0,m}$. The space of possible extensions $g_m$ of $g$ is a torsor under 
$\mathrm{Hom}\left(\mathbb{X}_Y[p^m]\times_{\overline{\mathbb{F}}_p}\left(\Xi_Z\times_{C_Z} C_{\widetilde{Z}}\right),\mathcal{G}_{\widetilde{Z}}[p^m]\right)$ over $\Xi_Z\times_{C_Z} C_{\widetilde{Z}}$. If we further impose the condition that this extension should lift $g_{0,m}$ then one obtains a torsor under $\mathrm{Hom}\left(\mathbb{X}_Y[p^m],\mathbb{X}_X[p^m]\right)(\Xi_Z\times_{C_Z}C_{\widetilde{Z}})$, and one further cuts out a closed subscheme by imposing the symmetry condition and the compatibility with endomorphisms. 

The morphism $\Xi_{\widetilde{Z}}\to \Xi_{Z}\times_{C_Z}C_{\widetilde{Z}}$ is universal for the existence of extensions $g_m$ of $g$ that lift a fixed $g_{0,m}$ and satisfy the desired compatibilities. 

\item The last piece of data, $\Xi_{\widetilde{Z}}\hookrightarrow \Xi_{\widetilde{Z},\Sigma_{\widetilde{Z}}}$ is a relative torus embedding. 
 \end{enumerate}

 \begin{remark}\label{rem:comparison with Lan}
     The reader might compare the above moduli description with the moduli description of the boundary charts for principal level structures constructed in~\cite[\S 6.2.3]{lan-thesis} in the PEL case. The relationship between $f_{0,m}$ and $g_{0,m}$ and the symmetry condition satisfied by $g_m$ are simplified in our setting. This is because we are working at an unramified prime $p$ and we are assuming that the splitting $\delta_b$ used in the definition of an Igusa cusp label is compatible with the polarization; by Lemma~\ref{lem:splitting anti-diagonal}, such splittings exist in the unramified case. In Lan's case, the construction of the spaces $C_Z$ and $\Xi_Z$ is more subtle; this is because, depending on the choice of splitting $\delta^p$, the polarization could induce additional non-trivial pairings between the graded pieces of the filtration induced by the degeneration data. 
 \end{remark}

 As a consequence of Theorem~\ref{thm:IgusaStratumDiagram} and Lemma \ref{lem:representability} (2), $\mathrm{Ig}^{b,\mathrm{tor}}_{\Sigma,m}$ is smooth and hence normal for $m \geq 2$. Furthermore, one can show that it satisfies the analog of the universal property in Part 6 of~\cite[Theorem 6.4.1.1]{lan-thesis}. This is done in exactly the same way as in \emph{loc. cit.}, using the toroidal boundary charts for $\mathrm{Ig}^{b,\mathrm{tor}}_{\Sigma', m'}$ of Theorem~\ref{thm:IgusaStratumDiagram} and using the universal property of the toroidal boundary chart established in~\cite[Prop. 6.2.5.11]{lan-thesis}. 

\subsubsection{Partial minimal compactifications}{\label{ss: partialminimalcompactifications}}

Fix a principal tame level $K^p=K(N)$, $p\nmid N$ and a hyperspecial level $K_p$ at $p$. Let $K=K_pK^p$ and consider the corresponding Shimura variety (resp. Igusa variety). In~\cite[\S 3.3]{CS2} and~\cite[\S 3.3]{santos}, partial minimal compactifications of the Igusa varieties $\Ig^{b}$ introduced in \S~\ref{ss: IgusaVarieties} are constructed. This relies, once again, on the construction of  minimal compactifications of canonical integral models of Shimura varieties due to~\cite{lan-thesis}. On the generic fiber, Lan's construction recovers the construction in Section~\ref{sec: stratifications}. On the special fiber (base-changed to $\ol{\mathbb{F}}_p$), we have an open dense embedding of $\overline{\mathbb{F}}_p$-schemes
\[
\mathscr{S}\hookrightarrow\mathscr{S}^*:= \mathscr{S}^*_{K, \overline{\mathbb{F}}_p},
\] 
with $\mathscr{S}^*$ projective. Furthermore, $\mathscr{S}^*$ admits a stratification into locally closed strata indexed by cusp labels, and, in particular, a set-theoretic decomposition
\[
\mathscr{S}^* = \bigsqcup_{\PP}\mathscr{S}_{\PP}, 
\]
with $\PP$ running over $\mathsf{G}(\mathbb{Q})$-conjugacy classes of admissible rational parabolic subgroups of $\mathsf{G}$. Given a good choice of compatible system $\Sigma$ of cone decompositions, there is a proper morphism $\mathscr{S}^{\mathrm{tor}}_{\Sigma}\to \mathscr{S}^*$, whose restriction to the open dense subscheme $\mathscr{S}$ is an isomorphism and which is compatible with the decomposition into cusp labels.  

Given a choice of $b\in B(G,\mu^{-1})$ and $\mathbb{X}$ in the isogeny class determined by $b$, Lan--Stroh also define a partial minimal compactification $\mathscr{C}\hookrightarrow\mathscr{C}^*$ of the corresponding Oort central leaf in~\cite[\S 3.4]{LanStroh}. This partial minimal compactification inherits the stratification in terms of cusp labels and in turn there is a stratification 
\[
\mathscr{C}^* = \bigsqcup_{\PP} \mathscr{C}_{\PP},
\]
by conjugacy classes of admissible parabolics, as in the discussion surrounding (\ref{eqn: decompositionoftoroidalintoparabolics}). 
For any good choice of $\Sigma$, there is a proper morphism $\mathscr{C}^{\mathrm{tor}}_{\Sigma}\to \mathscr{C}^*$, whose restriction to $\mathscr{C}$ is an isomorphism and which is compatible with the decomposition into cusp labels.

On the level of Igusa varieties, we have a perfect scheme $\mathrm{Ig}^{b,*}$ over $\overline{\mathbb{F}}_p$ that fits in a cartesian diagram 
\begin{equation}{\label{eqn: BoundaryDecompositionofPerfectIgusaVariety}}
\begin{tikzcd}
\mathrm{Ig}^b \arrow[d] \arrow[r, hook] & \mathrm{Ig}^{b,*}\arrow[d] \\ 
\mathscr{C} \arrow[r, hook] & \mathscr{C}^{*}.
\end{tikzcd}
\end{equation}
The perfect scheme $\mathrm{Ig}^{b,*}$ inherits a stratification in terms of cusp labels from $\mathscr{C}^*$, and this can be further refined to a stratification in terms of Igusa cusp labels. It in particular also has a stratification 
\begin{equation}{\label{eqn: BailyBorelStratificationonIgusavariety}}
\mathrm{Ig}^{b,*} = \bigsqcup_{\PP}\mathrm{Ig}^b_{\PP},
\end{equation}
and, as before, we denote the closure of each stratum by $\mathrm{Ig}^b_{\leq \PP}$. For any good choice of $\Sigma$, there is a proper, surjective morphism 
\begin{equation}{\label{eqn: propermorphismofIgusavarieties}}
\mathrm{Ig}^{b,\mathrm{tor}}_{\Sigma}\to \mathrm{Ig}^{b,*},
\end{equation}
whose restriction to $\mathrm{Ig}^b$ is an isomorphism and which respects the stratification by Igusa cusp labels. Furthermore, we have the following result. 
\begin{proposition}\label{prop: globalsections}
The partial minimal compactification
$\mathrm{Ig}^{b,*}$ is affine and its underlying ring agrees with the ring of global sections $H^0\left(\mathrm{Ig}^{b,\mathrm{tor}}_{\Sigma}, \mathcal{O}_{\mathrm{Ig}^{b,\mathrm{tor}}_{\Sigma}}\right)$.     
\end{proposition}
\begin{proof}
    This is \cite[Proposition 3.3.4]{CS2} and~\cite[Proposition 3.3.5]{santos}. 
\end{proof}

\begin{remark}{\label{rem: WhatarethepartialminimalCompactifications}}
There are also finite level partial minimal compactifications $\mathrm{Ig}^{b,*}_{m}$ of the $\Ig^{b}_{m}$, which can once again be defined by taking normalizations with respect to the map $\Ig^{b}_{m} \ra \mathscr{C}$ and the inclusion $\mathscr{C} \hookrightarrow \mathscr{C}^{*}$, as in (\ref{eqn: BoundaryDecompositionofPerfectIgusaVariety}). These satisfy the analog of Proposition~\ref{prop: globalsections} with respect to $\mathrm{Ig}^{b,\mathrm{tor}}_{\Sigma,m} \ra \Ig^{b,*}_{m}$, as all of these claims are deduced from the analogous claim for $\mathscr{C}_{\Sigma}^{\tor} \ra \mathscr{C}^{*}$.
\end{remark}
\begin{remark}
    When $K^p$ is a general neat, but not necessarily principal, tame level, we can still construct a partial minimal compactification of the Igusa variety $\Ig_{K^p}^b$. Indeed, we can find a neat tame level $(K^p)'\triangleleft K^p$ that is principal and then we construct $\Ig^{*,b}_{K^p}$ by taking the coarse quotient of $\Ig^{b,*}_{(K^p)'}$ under the action of the finite group $K^p/(K^p)'$. This is affine, contains  $\mathrm{Ig}^b_{K^p}$ as an open dense subscheme, and can be shown to agree with the relative normalization of $\mathscr{C}^*_{K^p}$ in $\mathrm{Ig}^b_{K^p}$.   
\end{remark}  

\subsection{The action of $\widetilde{J}_b$}\label{Sec:Action}

We maintain the assumptions on the PEL datum from the previous subsection and we impose, in addition, Assumption~\ref{assumption:codimension}. This implies, by~\cite[Corollary 9.24]{zhang2023}, that global sections over the open Igusa variety $\mathrm{Ig}^b$ extend uniquely to the affine scheme $\mathrm{Ig}^{b,*}$. This leads to the following result. 
\begin{lemma}\label{lemma: UniqueExistenceJbAction}
Under Assumption~\ref{assumption:codimension}, there is a unique extension of the action of the fpqc sheaf $\widetilde{J}_{b, \overline{\mathbb{F}}_p}$ on the perfect scheme $\mathrm{Ig}^b$ to the partial minimal compactification $\mathrm{Ig}^{b,*}$.     
\end{lemma}
\begin{proof}
The action of $\widetilde{J}_{b, \overline{\mathbb{F}}_p}$ on $\mathrm{Ig}^b$ is given by a morphism 
\[\widetilde{J}_{b, \overline{\mathbb{F}}_p}\times_{\ol{\mathbb{F}}_p} \mathrm{Ig}^b\to \mathrm{Ig}^b,\]
taking global section of which gives a ring homomorphism
\[\mc{O}(\mathrm{Ig}^{b,\ast})\simeq \mc{O}(\mathrm{Ig}^b)\to \mc{O}(\widetilde{J}_{b, \overline{\mathbb{F}}_p}\times_{\ol{\mathbb{F}}_p} \mathrm{Ig}^b)\simeq \mc{O}(\widetilde{J}_{b,\overline{\mathbb{F}}_p}\times_{\ol{\mathbb{F}}_p} \mathrm{Ig}^{b,\ast}), \]
where the isomorphisms follow from Assumption~\ref{assumption:codimension} and \cite[Proposition 1.9]{zhang2023}. This in turn gives the desired action map \[\widetilde{J}_{b, \overline{\mathbb{F}}_p}\times_{\ol{\mathbb{F}}_p} \mathrm{Ig}^{b,\ast}\to \mathrm{Ig}^{b,\ast}\]
by affineness of $\mathrm{Ig}^{b,\ast}$. The uniqueness is clear.
\end{proof}
We recall that we have a decomposition 
 \begin{equation}\label{eq:Jb decomposition}
 \widetilde{J}_{b,\overline{\mathbb{F}}_p} = \underline{J_b(\mathbb{Q}_p)}\ltimes \widetilde{\mathcal{U}}_b,   
 \end{equation}
where $J_b(\mathbb{Q}_p)$ can be identified with the group of self-quasi-isogenies of $\bb{X}$ over $\ol{\bb{F}}_p$ that respect the $G$-structure endowed with its profinite topology, $\ul{J_b(\mathbb{Q}_p)}$ is the associated constant group scheme, and $\widetilde{\mc{U}}_b$ is a unipotent formal group scheme, see \cite[Proposition 4.2.11]{CS17}.

For the rest of this subsection, we assume the prime-to-$p$ level $K^p=K(N)$, $p\nmid N$ is principal, in order to simplify the geometric structure of the toroidal boundary strata, in particular, so that we can refer to Theorem~\ref{thm:IgusaStratumDiagram}. The main result is the following theorem. 
\begin{theorem}\label{thm: JbPreservesStrata}
There exists a submonoid $J_{b}(\bb{Q}_{p})^{+} \subset J_{b}(\bb{Q}_{p})$, which generates $J_{b}(\bb{Q}_{p})$ as a group  and a subgroup $T_{p}\mathcal{U}_{b} \subset \tilde{\mathcal{U}}_{b}$ satisfying $\tilde{\mathcal{U}}_{b} = \bigcup_{g \in \underline{J_{b}(\bb{Q}_{p})}} g(T_{p}\mathcal{U}_{b})g^{-1}$ such that the action of 
$\underline{J_{b}(\bb{Q}_{p})^{+}}$ and $T_{p}\mathcal{U}_{b}$ on $\mathrm{Ig}^b$ extend to the projective system 
\[ \left(\mathrm{Ig}^{b,\mathrm{tor}}_{\Sigma}\right)_{\Sigma},\] where $\Sigma$ runs over good compatible choices of cone decomposition. This extended action preserves the projective system $\left(\mathrm{Ig}^{b,\mathrm{tor}}_{\Sigma,\PP}\right)_{\Sigma}$ consisting of boundary strata labeled by $\PP$, for each conjugacy class of admissible parabolics $\PP$. 
\end{theorem}

\begin{proof}
As stated, it is enough to construct the extension of the action and prove the desired preservation of the boundary strata for:
\begin{enumerate}
\item the constant scheme attached to a submonoid $J_b(\Q_p)^+$ that generates $J_b(\Q_p)$ as a group;
\item the subgroup $T_p\mathcal{U}_b:=\underline{\mathrm{Aut}}_{G}(\mathbb{X})\cap \widetilde{\mathcal{U}}_b$, which has the claimed property that 
\[
\widetilde{\mathcal{U}}_b = \bigcup_{g\in \underline{J_b(\Q_p)}} g\left(T_p\mathcal{U}_b\right)g^{-1}.
\]
\end{enumerate}
For (1), this is done in \S \ref{subsub: J(qp)}, in particular in Proposition~\ref{prop: HeckeAction}. 
For (2), this is done in \S \ref{subsub: Ub}, in particular in Corollary~\ref{cor: TpUbAction}.
\end{proof}

\subsubsection{The action of \texorpdfstring{$J_b(\qp)$}{}}\label{subsub: J(qp)}
We first consider the case of $J_b(\Q_p)$, which is identified with the group of self-quasi-isogenies of $\mathbb{X}$ over $\overline{\mathbb{F}}_p$ that respect the $G$-structure. We define the submonoid $J_b(\Q_p)^+$ as 
\[\{\rho\in J_b(\Q_p)\mid \rho^{-1}\ \mathrm{is}\ \mathrm{an}\ \mathrm{isogeny}\}.
\]
This submonoid clearly generates $J_b(\Q_p)$ as a group. For each $\rho\in J_b(\mathbb{Q}_p)^+$, we have a finite flat subgroup scheme $K_{\rho}\coloneq \operatorname{ker}(\rho^{-1})\subset \mathbb{X}$ that is compatible with the $G$-structures. We choose $m_0\gg 0$ such that $K_\rho\subset \mathbb{X}[p^{m_0}]$. For any $m\geq m_0$ and we can define a finite flat subgroup scheme 
\[
\mathcal{K}_{\rho}\subset \left(\mathscr{G}\times_{\mathscr{C}}\mathrm{Ig}^b_m\right)[p^m]
\]
by taking the image of $K_{\rho}\times_{\overline{\mathbb{F}}_p}\mathrm{Ig}^b_m$ under the trivialization $\gamma_m: \mathbb{X}\times_{\overline{\mathbb{F}}_p}\mathrm{Ig}^b_m \stackrel{\sim}{\to} \mathcal{A}[p^m]$. The quotient $\mathcal{A}/\mathcal{K}_{\rho}$ inherits the PEL structures of $\mathcal{A}$ and induces a morphism $[\rho]:\mathrm{Ig}^b_m\to \mathrm{Ig}^b_{m-m_0}$. Varying $m$ and taking an inverse limit, we obtain a morphism 
$[\rho]: \mathrm{Ig}^b\to \mathrm{Ig}^b$ on the perfect Igusa variety. Varying $\rho$, we obtain an action of the monoid $J_b(\mathbb{Q}_p)^+$ on $\mathrm{Ig}^b$. By extending from the case of connected test objects, this construction gives a scheme-theoretic action of the constant scheme $\underline{J_b(\mathbb{Q}_p)^+}$ on $\mathrm{Ig}^b$. This action agrees with the restriction of the action of $\widetilde{J}_{b, \overline{\mathbb{F}}_p}$ described in \S\ref{Sec:Igusa background}. We will show that each morphism $[\rho]$ extends to partial toroidal compactifications and preserves the boundary strata labeled by each $\PP$. This will extend to an explicit action of $\underline{J_b(\mathbb{Q}_p)^+}$ on these partial compactifications that preserves the boundary strata labeled by each $\PP$.

Given a good compatible choice of cone decompositions $\Sigma = \{\Sigma_{Z}\}$ and an Igusa cusp label $\widetilde{Z}$, we define $\Sigma_{\widetilde{Z}}:= \Sigma_Z$. On the level of partial toroidal compactifications of Igusa varieties, it is not true in general that the Hecke action by elements of $J_b(\mathbb{Q}_p)^+$ will preserve a single compatible choice of cone decompositions. This phenomenon already occurs for the Hecke action of $G(\A_f^p)$ on the level of the Shimura variety. However, modeled on this case, we formulate below, for each $\rho\in J_b(\mathbb{Q}_p)^+$, a precise relationship between two such choices $\Sigma$ and $\Sigma'$ such that, when the relationship is satisfied, we will be able to extend $[\rho]$ to a morphism 
\[
[\rho]^{\mathrm{tor}}:\mathrm{Ig}^{b,\mathrm{tor}}_{\Sigma}\to \mathrm{Ig}^{b,\mathrm{tor}}_{\Sigma'}.
\]

\begin{definition}\label{defn:rho-refinement} If $\Sigma$ and $\Sigma'$ are two good compatible choices of cone decompositions and $\rho\in J_b(\Q_p)$, we say that $\Sigma$ is a \emph{$\rho$-refinement} of $\Sigma'$ if the following condition is satisfied:
\begin{itemize}
\item For any Igusa cusp label $\widetilde{Z}$, which is mapped by $\rho$ to an Igusa cusp label $\widetilde{Z}'$, we obtain isomorphisms 
\[
f_X:X\otimes_{\mathbb{Z}} \mathbb{Q}\stackrel{\sim}{\to} X'\otimes_{\mathbb{Z}} \mathbb{Q}\  \mathrm{and} \ 
f_Y:Y\otimes_{\mathbb{Z}} \mathbb{Q}\stackrel{\sim}{\to}  Y'\otimes_{\mathbb{Z}} \mathbb{Q}
\]
The pair $(f_X, f_Y)$ induces an identification $\mathbf{P}_{{\widetilde{Z}}}\simeq \mathbf{P}_{{\widetilde{Z}'}}$ 
and we assume that the cone decomposition $\Sigma_{\widetilde{Z}}$ is a refinement of the cone decomposition $\Sigma'_{\widetilde{Z}'}$ under this identification.   
\end{itemize}
We note that this condition is independent of the choices of isomorphisms $(f_X,f_Y)$. 
\end{definition}

\begin{proposition}\label{prop: HeckeAction}
    Let $\rho\in J_b(\Q_p)^+$. The morphism $[\rho]:\mathrm{Ig}^b\to \mathrm{Ig}^b$ induced by the isogeny $\mathcal{A}\to \mathcal{A}/\mathcal{K}_{\rho}$ extends to a morphism 
    \[
    [\rho]^{\mathrm{tor}}: \mathrm{Ig}^{b,\mathrm{tor}}_{\Sigma} \to \mathrm{Ig}^{b,\mathrm{tor}}_{\Sigma'} 
    \] 
    under the assumption that the compatible choice of cone decompositions $\Sigma$ is a $\rho$-refinement of the compatible choice of cone decompositions $\Sigma'$. The extended morphism respects the boundary strata labeled by any conjugacy class $\PP$. 
\end{proposition}

\begin{proof} The main subtlety is showing that the isogeny $\mathcal{A}\to \mathcal{A}/\mathcal{K}_{\rho}$ extends to an isogeny of semi-abelian schemes $\mathcal{A}_{\Sigma}\to \mathcal{A}_{\Sigma}/ \mathcal{K}_{\rho, \Sigma}$. In general, this is not obvious -- see~\cite[\S 6]{deJong-Oort}. However, in our case, we can exploit the fact that the families of $p$-divisible groups living over Igusa varieties are constant, as they are equipped with trivializations of their full $p$-divisible group. 

We choose $m_0\gg 0$ such that $\mathcal{K}_{\rho}\subseteq \mathcal{A}[p^{m_0}]$. For any $m\geq m_0$, we consider the finite-level Igusa variety $\mathrm{Ig}^b_{m}$ and its partial toroidal compactification $\mathrm{Ig}^{b,\mathrm{tor}}_{\Sigma, m}$. These are, in particular, normal noetherian schemes. We have a semi-abelian scheme $\mathcal{A}_{\Sigma}$ over $\mathrm{Ig}^{b,\mathrm{tor}}_{\Sigma, m}$, which restricts to an abelian scheme $\mathcal{A}$ over $\mathrm{Ig}^b_m$. Over $\mathrm{Ig}^b_{m}$, we have a trivialization 
\[
\gamma_m: \mathbb{X}[p^m]\times_{\ol{\mathbb{F}}_p} \mathrm{Ig}^b_m\stackrel{\sim}{\to} \mathcal{A}[p^m].
\]
Furthermore, the connected part $\mathcal{A}[p^m]^\circ$ extends to a finite flat subgroup scheme $\mathcal{A}_{\Sigma}[p^m]^\circ$ over $\mathrm{Ig}^{b,\mathrm{tor}}_{\Sigma, m}$ by~\cite[Prop. 3.2.1]{CS2}. The latter, which is defined as the connected part of the quasi-finite flat group scheme $\mathcal{A}_{\Sigma}[p^m]$, is equipped with a trivialization 
\[
\bar{\gamma}_m^\circ : \mathbb{X}[p^m]^\circ \times_{\ol{\mathbb{F}}_p} \mathrm{Ig}^{b,\mathrm{tor}}_{\Sigma, m}\stackrel{\sim}{\to} 
\mathcal{A}_{\Sigma}[p^m]^\circ,
\]
which is compatible with $\gamma_m$. The trivialization $\bar{\gamma}_m^\circ$ can be seen to exist over completions of strict local rings of $\mathrm{Ig}^{b,\mathrm{tor}}_{\Sigma, m}$, using the theory of degeneration, and it is glued to a trivialization over the whole of $\mathrm{Ig}^{b,\mathrm{tor}}_{\Sigma, m}$ in the usual way, e.g. via Lemma~\ref{lem:extendability is a local property}, using the compatibility with $\gamma_m$ over the open, schematically dense complement of the boundary. Over a completion of a strict local ring $\mathrm{Spec}\ R^{\wedge}$ at a point in the boundary stratum cut out by a pair $(\widetilde{Z}, \sigma)$, we have a canonical identification between the pullback of $\mathcal{A}_{\Sigma}[p^m]^\circ$ and $\mathcal{G}[p^m]^\circ$ coming from the theory of degeneration, where $\mathcal{G}$ is the corresponding Raynaud extension, as in the proof of~\cite[Prop. 3.2.1]{CS2}. Let $\Xi_{\widetilde{Z}}(\sigma)$ be the restriction to the cone $\sigma\in \Sigma_{\widetilde{Z}}^+$\footnote{Here, $\Sigma_Z$ is a cone decomposition of some space of symmetric, positive semi-definite pairings with rational radicals, and $\Sigma^+_{Z}\subset \Sigma_{Z}$ is the subset consisting of those cones contained in the interior of this space. We can always assume that the pair $(Z,\sigma)$ satisfies $\sigma\in \Sigma^+_{Z}$, e.g. by~\cite[Prop. 2.1.2 (4)]{LanStroh}.} of the boundary chart of Theorem~\ref{thm:IgusaStratumDiagram}. We have a universal trivialization of $\mathcal{G}_{\widetilde{Z}}[p^m]^\circ$ over $\Xi_{\widetilde{Z}}(\sigma)$ coming from the moduli-theoretic description of the latter - see point (2) below Theorem~\ref{thm:IgusaStratumDiagram}. We obtain a trivialization of $\mathcal{G}[p^m]^\circ$ by pullback along the canonical morphism $\mathrm{Spec}\ R^\wedge\to \Xi_{\widetilde{Z}}(\sigma)$ and this induces the desired trivialization on the pullback of $\mathcal{A}_{\Sigma}[p^m]^\circ$ to $\mathrm{Spec}\ R^{\wedge}$.   

We can therefore apply Lemma~\ref{lem:Extensibility of Isogenies} below, with $\overline{U}:=\mathrm{Ig}^{b,\mathrm{tor}}_{\Sigma, m}$, $U:=\mathrm{Ig}^b_m$ and $\mathcal{Z}$ the pullback of $\mathcal{A}$ to $U$. We deduce that we can extend the isogeny  $\mathcal{A}\to \mathcal{A}':=\mathcal{A}/\mathcal{K}_{\rho}$ to an isogeny of semi-abelian schemes over $\mathrm{Ig}^{b, \mathrm{tor}}_{\Sigma, m}$. The PEL-structures inherited by the quotient $\mathcal{A}'_{\Sigma}:=\mathcal{A}_{\Sigma}/\mathcal{K}_{\rho, \Sigma}$ over the open dense subscheme $\mathrm{Ig}^b_m$ extend uniquely to  PEL-structures over the semi-abelian scheme, by~\cite[Chapter I, Proposition 2.7]{faltings-chai}. We also get an induced Igusa level structure of level $m-m_0$ over $\mathrm{Ig}^{b}_{m}$, i.e. a trivialization 
\[
\gamma'_{m-m_0}:\mathbb{X}[p^{m-m_0}]\times_{\ol{\mathbb{F}}_p} \mathrm{Ig}^{b}_{m} \stackrel{\sim}{\to} \mathcal{A}'[p^{m-m_0}]
\]
compatible with the PEL structures in the usual sense. 

By Part 6 of~\cite[Theorem 6.4.1.1]{lan-thesis} and the assumption that $\Sigma$ is a $\rho$-refinement of $\Sigma'$, we obtain a morphism $\mathrm{Ig}^{b, \mathrm{tor}}_{\Sigma, m}\to \mathscr{S}^{\mathrm{tor}}_{\Sigma'}$. This factors through the leaf $\mathscr{C}^{b,\mathrm{tor}}_{\Sigma'}\subset \mathscr{S}^{\mathrm{tor}}_{\Sigma'}$. We wish to show that this morphism lifts to an Igusa cover $\mathrm{Ig}^{b,\mathrm{tor}}_{\Sigma', m'}$ of the leaf, at least when $m-m_0\gg m'$. 

When $m-m_0\gg m'$, the Igusa level structure $\gamma'_{m-m_0}$ on the complement of the boundary induces a lift
\[
\mathrm{Ig}^b_m\to \mathrm{Ig}^b_{m'}.
\]
To prove that this lift extends to the boundary of the partial toroidal compactification, we use the normality of $\mathrm{Ig}^{b,\mathrm{tor}}_{\Sigma, m}$ and the fact that $\Sigma$ is a $\rho$-refinement of $\Sigma'$. The extendability now follows from the analog of the universal property in Part 6 of~\cite[Theorem 6.4.1.1]{lan-thesis} for $\mathrm{Ig}^{b,\mathrm{tor}}_{\Sigma, m}$. We conclude the construction of $[\rho]^{\mathrm{tor}}$ by taking a projective limit over $m$ and $m'$. 

It remains to see that the morphism $[\rho]^{\mathrm{tor}}$ preserves the stratification labeled by $\PP$ on $\mathrm{Ig}^{b,\mathrm{tor}}_{\Sigma}$. In fact, more is true: the action of $[\rho]^{\mathrm{tor}}$ is compatible with the action of $\rho\in J_b(\mathbb{Q}_p)^+$ on Igusa cusp labels. This can be proved in the same way as~\cite[Proposition 6.4.3.4]{lan-thesis}, which concerns the compatibility between the Hecke action away from $p$ and usual cusp labels, and which follows from the constructions of \S~5.4.3 of \emph{loc. cit.}. 

Specifically, in our case, the compatibility between $[\rho]^{\mathrm{tor}}$ and the stratification labeled by $\PP$ follows from the compatibility of filtrations induced from the theory of degeneration under an isogeny of semi-abelian schemes, cf.~\cite[Construction 5.4.3.3]{lan-thesis}.
(Note that Construction 5.4.3.3 of \emph{loc. cit.} works with a general isogeny between semi-abelian schemes, even though it is only applied in Proposition 6.4.3.4 to the case when  the isogeny has \'etale kernel.) 
The  filtration on the rational Tate module away from $p$ is encoded in the usual cusp label via the filtration $Z^p$,  and the conjugacy class $\PP$ of admissible parabolics can be recovered from this, cf.~\cite[\S A.4, A.5]{Lan-ANT}. 


\end{proof}

\begin{lemma}\label{lem:Extensibility of Isogenies} 
Let $\overline{U}$ be a normal noetherian scheme over $\overline{\mathbb{F}}_p$ containing an open dense subscheme $U$. Assume the following conditions:
\begin{enumerate}
\item  We have an isogeny of abelian schemes $f:\mathcal{Z}\to \mathcal{Z}'$ over $U$, with kernel being a finite flat group scheme $\mathcal{K}\subseteq \mathcal{Z}[p^{m_0}]$ for some $m_0\gg0$.
\item  The finite flat group scheme $\mathcal{H}:=\mathcal{Z}[p^{m_0}]$ is constant over $U$, equipped with a trivialization 
\[
\gamma: \mathbb{X}[p^{m_0}]\times_{\overline{\mathbb{F}}_p} U
\stackrel{\sim}{\to} \mathcal{H},
\]  
such that $\mathcal{K}$ is the image under $\gamma$ of $K\times_{\overline{\mathbb{F}}_p} U$, for a finite flat group scheme $K/\ol{\bb{F}}_{p}$. 
\item The abelian scheme $\mathcal{Z}$ extends to a semi-abelian scheme $\overline{\mathcal{Z}}$ over $\overline{U}$.
\item Consider the quasi-finite flat group scheme  $\overline{\mathcal{H}}:=\overline{\mathcal{Z}}[p^{m_0}]$ over $\overline{U}$. Its connected part $\overline{\mathcal{H}}^\circ$ is a finite flat group scheme, which is equipped with a trivialization 
\[
\bar{\gamma}^\circ: \mathbb{X}[p^{m_0}]^\circ\times_{\overline{\mathbb{F}}_p}\overline{U}
\stackrel{\sim}{\to} \overline{\mathcal{H}}^\circ, 
\] 
which is compatible with $\gamma$. 
\end{enumerate}
Then the abelian scheme $\mathcal{Z}'$ extends to a semi-abelian scheme $\overline{\mathcal{Z}}'$ over $\overline{U}$, and is equipped with an isogeny $\bar{f}:\overline{\mathcal{Z}}\to \overline{\mathcal{Z}}'$ that extends $f$.  
\end{lemma}

\begin{proof} Let $\overline{\mathcal{K}}$ be the schematic closure of $\mathcal{K}$ in $\overline{\mathcal{Z}}$. This is a quasi-finite subgroup scheme of the quasi-finite flat group scheme $\overline{\mathcal{H}}$. As in the proof of~\cite[Lemma 3.1.3.2]{lan-ordinary} (which considers the special case where $b$ is ordinary), the extendability of $f:\mathcal{Z}\to \mathcal{Z}'$ to $\overline{U}$ is equivalent to the flatness of $\overline{\mathcal{K}}$ over $\overline{U}$, so we will prove the latter. As in \emph{loc. cit.}, we may assume that $\overline{U}$ is its strict localization at a point $u$, that $U = \overline{U}\setminus u$, and that $u$ has codimension at least $2$ in $\overline{U}$. 

We have a filtration 
\[
0\subseteq \overline{\mathcal{H}}^\circ \subseteq \overline{\mathcal{H}}^f \subseteq \overline{\mathcal{H}},
\]
where $\overline{\mathcal{H}}^f$ is the maximal finite flat subgroup scheme of $\overline{\mathcal{H}}$ and $\overline{\mathcal{H}}^\circ$ is its connected part. The existence of the trivialization $\bar{\gamma}^\circ$ and its compatibility with $\gamma$ imply that we can identify $\overline{\mathcal{H}}^\circ$ with the schematic closure $\overline{\mathcal{H}^\circ}$ of the connected part $\mathcal{H}^\circ$ of $\mathcal{H}$ in $\overline{\mathcal{H}}$. The quotient $\overline{\mathcal{H}}/\overline{\mathcal{H}}^f$ is quasi-finite \'etale and its special fiber over $u$  consists only of the identity section; the quotient $\overline{\mathcal{H}}^{f,et}:=\overline{\mathcal{H}}^f/\overline{\mathcal{H}}^\circ$ is finite \'etale. By restricting to $U$ and intersecting with $\mathcal{K}$, we get an induced filtration 
\[
0\subseteq \mathcal{K}^\circ\subseteq \mathcal{K}^f\subseteq \mathcal{K}, 
\]
where $\mathcal{K}^f$ is finite flat over $U$. As in \emph{loc. cit.}, we will prove that $\overline{\mathcal{K}}^f$, the schematic closure of $\mathcal{K}^f$ in $\overline{\mathcal{H}}^f$, is finite flat. 

The schematic closure $\overline{\mathcal{K}^\circ}$ of $\mathcal{K}^\circ$ in $\overline{\mathcal{H}^\circ} = \ol{\mathcal{H}}^\circ$ is finite flat, as it can be identified via $\bar{\gamma}^\circ$ with $K^\circ\times_{\overline{\mathbb{F}}_p} \overline{U}$, where $K^{\circ}$ denotes the neutral component of $K$. We claim that $\overline{\mathcal{K}^\circ} = \overline{\mathcal{K}}^\circ$. This follows from the identifications 
\[
\overline{\mathcal{K}^\circ} = \overline{\mathcal{K}}\cap \overline{\mathcal{H}^\circ} = \overline{\mathcal{K}}\cap \overline{\mathcal{H}}^\circ = \overline{\mathcal{K}}^\circ.
\]
The quotient $\mathcal{K}^{f,et}:=\mathcal{K}^f/\mathcal{K}^\circ$ is a finite \'etale group scheme and so is its schematic closure $\overline{\mathcal{K}}^{f,et} = \ol{\mathcal{K}}/\ol{\mathcal{K}}^\circ$ in $\overline{\mathcal{H}}^{f,et}$, since $\overline{U}$ is strict local and normal, making $\overline{\mathcal{H}}^{f,et}$ a constant group scheme. 

Let $\overline{\mathcal{K}}^+$ be the pre-image of $\overline{\mathcal{K}}^{f,et}$ under the projection $\overline{\mathcal{H}}^f\twoheadrightarrow \overline{\mathcal{H}}^{f,et}$. Over $U$, we can realize $\mathcal{K}^f$ as the kernel of an isogeny of finite flat group schemes 
\[
\mathcal{K}^+\to \mathcal{H}^\circ/\mathcal{K}^\circ. 
\]
By Hartogs, using the normality of $\overline{U}$ and the assumption on the codimension of $u$, we can extend this uniquely to a homomorphism of finite flat group schemes 
\[
\Phi: \overline{\mathcal{K}}^+\to \overline{\mathcal{H}}^\circ/\overline{\mathcal{K}}^\circ. 
\]
This extended homomorphism is surjective for the fppf topology, as its restriction to $\overline{\mathcal{H}}^\circ$ is uniquely determined to be the isogeny $\overline{\mathcal{H}}^\circ\to \overline{\mathcal{H}}^\circ/ \overline{\mathcal{K}}^\circ$. 
The kernel of $\Phi$,  which agrees with $\overline{\mathcal{K}}^f$, fits into a short exact sequence of sheaves for the fppf topology 
\[
0\to \overline{\mathcal{K}}^\circ \to \overline{\mathcal{K}}^f\to \overline{\mathcal{K}}^{f,et}\to 0. 
\]
This shows that $\overline{\mathcal{K}}^f$ is a finite flat group scheme. 
Recall that $\overline{\mathcal{K}} = \overline{\mathcal{K}}^f\sqcup \overline{\mathcal{K}}'$, where by definition $\overline{\mathcal{K}}'$ has empty fiber over the closed point $u$ and agrees with $\mathcal{K}\setminus \mathcal{K}^f$ over its open complement $U = \overline{U}\setminus u$. Therefore, we deduce that $\overline{\mathcal{K}}$ is flat as well, as desired. 
\end{proof}



\subsubsection{The action of \texorpdfstring{$\widetilde{\mc{U}}_b$}{}}\label{subsub: Ub} We now consider the subgroup scheme $T_p\mathcal{U}_b = \underline{\mathrm{Aut}}_{G}(\mathbb{X})\cap \widetilde{\mathcal{U}}_b$ of $\widetilde{J}_{b,\overline{\mathbb{F}}_p}$. We note that $T_p\mathcal{U}_b\subset \underline{\mathrm{Aut}}_{G}(\mathbb{X})^\circ$. For any $m\in \mathbb{Z}_{\geq 2}$, we have the smooth finite-level Igusa variety $\mathrm{Ig}^b_m$ that parametrizes trivializations 
\[
\gamma_m: \mathbb{X}[p^m]\times_{\overline{\mathbb{F}}_p}\mathrm{Ig}^b_m \stackrel{\sim}{\to} \mathcal{A}[p^m].
\]
compatible with the $G$-structures.
The group scheme $T_p\mathcal{U}_b$ already acts on $\mathrm{Ig}^b_m$. Indeed, we describe this from the perspective of the functor of points. Let $S\to \mathrm{Ig}^b_m$ be a $S$-point for a scheme $S/\ol{\bb{F}}_{p}$, and consider the pullbacks $\mathcal{A}_S$ and $\gamma_{m,S}$ of the universal abelian scheme $\mathcal{A}$ with its $\mathsf{G}$-structures and, respectively, of the Igusa level structure $\gamma_m$. For any $\rho\in T_p\mathcal{U}_b(S)$, we modify the Igusa level structure via the pre-composition  $\gamma_{m,S} \mapsto \gamma_{m,S}\circ \left(\rho\mid_{\mathbb{X}[p^m]\times_{\ol{\mathbb{F}}_p}S}\right)$ and, by the universal property of $\mathrm{Ig}^b_m$ (i.e. the fact that $\mathrm{Ig}^b_m\to \mathscr{C}$ represents the moduli problem of trivializations of $\mathcal{A}[p^m]$ compatible with the $G$-structures), this gives a new morphism $S\to \mathrm{Ig}^b_m$. 

When we take a projective limit over $m$, we recover the restriction to $T_p\mathcal{U}_b$ of the action of $\widetilde{J}_{b,\overline{\mathbb{F}}_p}$ described in \S\ref{Sec:Igusa background}. 

Fix a good compatible
choice of cone decompositions $\Sigma = (\Sigma_{\widetilde{Z}})_{\widetilde{Z}}$ indexed by Igusa cusp labels at level $K^pK_{b}(p^m)$. 
Fix a (representative of an) Igusa cusp label $\widetilde{Z}$ at this level, which determines a boundary stratum $\mathrm{Ig}^b_{\widetilde{Z},m} \subset \mathrm{Ig}^{b,*}_m$ that is identified with a smaller Igusa variety, as in~\cite[Theorem 3.3.15]{CS2}. We will show that $T_p\mathcal{U}_b$ has an explicit, moduli-theoretic action on the formal completion $\widehat{\mathrm{Ig}}^{b,\mathrm{tor}}_{m,\widetilde{Z}}$ of $\mathrm{Ig}^{b,\mathrm{tor}}_{\Sigma, m}$ along the boundary stratum labeled by $\widetilde{Z}$ in the partial toroidal compactification. This explicit action will visibly extend the moduli-theoretic action on $\mathrm{Ig}^b_m$ described above. 

 Using the decomposition~\eqref{eq:IgusaCuspLabelDecomposition}, we can visualize the automorphisms of the universal cover $\underline{\mathrm{Aut}}(\widetilde{\mathbb{X}})$ as a $3\times 3$ matrix (which acts by left multiplication on $\widetilde{\bb{X}}$ when it is viewed as a column vector)
\begin{equation}\label{eq:ShapeAutX}
    \begin{pmatrix}
    \underline{\mathrm{Aut}}(\widetilde{\mathbb{X}}_Y) & \ & \ \\
    \widetilde{\mathcal{H}}_{\mathbb{X}_{Y},\mathbb{X}_{\widetilde{Z}}} & \underline{\mathrm{Aut}}(\widetilde{\mathbb{X}}_{\widetilde{Z}}) & \ \\
    \widetilde{\mathcal{H}}_{\mathbb{X}_{Y},\mathbb{X}_X} & \widetilde{\mathcal{H}}_{\mathbb{X}_{\widetilde{Z}},\mathbb{X}_X} & \underline{\mathrm{Aut}}(\widetilde{\mathbb{X}}_X)     
    \end{pmatrix}.
\end{equation}
Here, $\widetilde{\mathcal{H}}_{\mathbb{X}_{Y},\mathbb{X}_{\widetilde{Z}}}$ and so on denote universal covers of the internal Hom $p$-divisible groups in the sense of~\cite[\S 4.1.6]{CS17}. The principal polarization on $\mathbb{X}$ induces identifications 
\[
\underline{\mathrm{Aut}}(\widetilde{\mathbb{X}}_Y)\simeq \underline{\mathrm{Aut}}(\widetilde{\mathbb{X}}_X)\ \mathrm{and}\ \widetilde{\mathcal{H}}_{\mathbb{X}_{Y},\mathbb{X}_{\widetilde{Z}}}\simeq \widetilde{\mathcal{H}}_{\mathbb{X}_{\widetilde{Z}},\mathbb{X}_X}.
\]



\noindent Inside $\mathrm{Aut}(\widetilde{\mathbb{X}})$, we can visualize $\rho\in T_p\mathcal{U}_b$ as having the form
\begin{equation}\label{eq:TpUb}
 \rho =   \begin{pmatrix}
        1& \ &\ \\
        \rho_{Y, \widetilde{Z}} & \rho^\circ_{\widetilde{Z}}& \ \\
        \rho_{Y,X} & \rho_{\widetilde{Z},X} & 1
    \end{pmatrix} \in 
    \begin{pmatrix}
        1 & \ & \ \\
        T_p\mathcal{H}_{\mathbb{X}_{Y}, \mathbb{X}_{\widetilde{Z}}} & 
        \underline{\mathrm{Aut}}(\mathbb{X}_{\widetilde{Z}})^\circ & \ \\
        T_p\mathcal{H}_{\mathbb{X}_Y, \mathbb{X}_X} & 
    T_p\mathcal{H}_{\mathbb{X}_{\widetilde{Z}},\mathbb{X}_{X}} & 1 
    \end{pmatrix}.
\end{equation}
\begin{remark}
Here, we are passing to the connected component of the identity $\underline{\mathrm{Aut}}(\mathbb{X}_{\widetilde{Z}})^\circ$, because $T_p\mathcal{U}_b$ is contained in $\underline{\mathrm{Aut}}_G(\mathbb{X})^\circ$, as noted above. Of course, for $\rho$ to be an element of $T_p\mathcal{U}_b$, it has to satisfy additional compatibilities with the $G$-structures of $\mathbb{X}$, i.e. with its endomorphisms and polarization. The endomorphisms of $\mathbb{X}$ respect the decomposition in~\eqref{eq:IgusaCuspLabelDecomposition}, so the compatibility with the endomorphisms can be imposed on each of the entries of the matrix describing $\rho$. The compatibility with the polarization introduces two constraints that relate $\rho^\circ_{\widetilde{Z}}$, $\rho_{Y,\widetilde{Z}}$, $\rho_{\widetilde{Z},X}$ and $\rho_{Y,X}$. One of these constraints is spelled out explicitly in Lemma~\ref{lem: compatibility with polarization}, see in particular the equation~\eqref{eq: compatibility with polarization}. 
\end{remark}

Using the moduli-theoretic description in (1), (2) and (3) below Theorem~\ref{thm:IgusaStratumDiagram}, we first define an action of $T_p\mathcal{U}_b$ on $\Xi_{\widetilde{Z}}$. 

\begin{definition}\label{defn:dd action}
Let $S\to \Xi_{\widetilde{Z}}$ be a scheme and let $\rho\in T_p\mathcal{U}_b(S)$. We pull back the universal tuple $(\alpha_m, \beta_m, g_{0,m}, g_m)$ along the given morphism $S\to \Xi_{\widetilde{Z}}$ to obtain a tuple $(\alpha_{m,S}, \beta_{m,S},g_{0,m,S}, g_{m,S})$. Our chosen isomorphism $\rho$ induces a new tuple as follows:

\begin{enumerate}
\item a trivialization  $\mathbb{X}_{\widetilde{Z}}[p^m]\times_{\ol{\mathbb{F}}_p}S\stackrel{\sim}{\to} \mathcal{B}_{Z}[p^m]\times_{\mathscr{C}_Z} S$, compatible with $G$-structures, given by 
\[
[\rho]_{\widetilde{Z}}(\alpha_{m,S}): = \alpha_{m,S} \circ \left(\rho_{\widetilde{Z}}\mid_{\mathbb{X}_{\widetilde{Z}}[p^m]\times_{\overline{\mathbb{F}}_p}S}\right);
\]

\item a splitting $\mathcal{B}_Z[p^m]\times_{\mathscr{C}_Z} S\to \mathscr{G}_{Z}[p^m]\times_{C_Z}S$ compatible with endomorphisms on both sides, given by 
\[
[\rho]_{\widetilde{Z}}(\beta_{m,S}) := \beta_{m,S} + \left(\rho_{\widetilde{Z},X}\mid_{\mathbb{X}_{\widetilde{Z}}[p^m]\times_{\overline{\mathbb{F}}_p}S}\right)\circ \left(\rho_{\widetilde{Z}}\mid_{\mathbb{X}_{\widetilde{Z}}[p^m]\times_{\overline{\mathbb{F}}_p}S}\right)^{-1}\circ \alpha_{m,S}^{-1};
\]

\item[(\theenumi')] an extension $\underline{\frac{1}{p^m}Y}\to \mathcal{B}_{Z}\times_{\mathscr{C}_Z}S$ of $g_{0,S}: \underline{Y}\to \mathcal{B}_{Z}\times_{\mathscr{C}_Z}S$ compatible with endomorphisms, given by  
\[
[\rho]_{\widetilde{Z}}(g_{0,m,S}) := g_{0,m,S} + \alpha_{m,S}\circ \left(\rho_{Y, \widetilde{Z}}\mid_{\mathbb{X}_Y[p^m]\times_{\overline{\mathbb{F}}_p}S}\right);
\]

\item an extension $\underline{\frac{1}{p^m}Y}\to \mathcal{G}_{Z}\times_{C_Z}S$ of $g_{S}:\underline{Y} \to \mathcal{G}_{Z}\times_{C_Z}S$, compatible with endomorphisms and lifting $[\rho]_{\widetilde{Z}}(g_{0,m,S})$, given by 
\[
[\rho]_{\widetilde{Z}}(g_{m,S}): = g_{m,S} + \beta_{m,S}\circ \alpha_{m,S} \circ \left(\rho_{Y,\widetilde{Z}}\mid_{\mathbb{X}_Y[p^m]\times_{\overline{\mathbb{F}}_p}S}\right) + \rho_{Y, X}\mid_{\mathbb{X}_Y[p^m]\times_{\overline{\mathbb{F}}_p}S}.
\]  
\end{enumerate}
By the universal property of $\Xi_{\widetilde{Z}}$ (i.e. because $\Xi_{\widetilde{Z}}\to \Xi_{Z}$ parametrizes tuples as in (1), (2), (2') and (3) by construction), this new tuple
\[
\left([\rho]_{\widetilde{Z}}(\alpha_{m,S}), [\rho]_{\widetilde{Z}}(\beta_{m,S}), [\rho]_{\widetilde{Z}}(g_{0,m,S}), [\rho]_{\widetilde{Z}}(g_{m,S})
\right)
\]
induces a new $\Xi_Z$-linear morphism $S\to \Xi_{\widetilde{Z}}$. Varying our test object $S$, we obtain an action of the group scheme $T_p\mathcal{U}_b$ on $\Xi_{\widetilde{Z}}$. 
\end{definition}

\begin{lemma}\label{lem: compatibility with polarization} The new structures described in points (2) and (2') above are compatible under duality. 
\end{lemma}

\begin{proof}

    To simplify our notation, we suppress the test object $S$ and the restriction to $p^m$-torsion in what follows. We have 
    \[
    [\rho]_{\widetilde{Z}}(\beta_m) - \beta_m =\rho_{\widetilde{Z},X}\circ \rho_{\widetilde{Z}}^{-1}\circ \alpha_m^{-1} \in \mathrm{Hom}\left(\mathcal{B}_{\widetilde{Z}}[p^m],\mathbb{X}_{X}[p^m]\times_{\overline{\mathbb{F}}_p}\mathrm{Ig}^b_{\widetilde{Z}}\right).
    \]
    The dual morphism is $(\rho_{\widetilde{Z},X}\circ \rho_{\widetilde{Z}}^{-1}\circ \alpha_m^{-1})^\vee \in \mathrm{Hom}(\underline{X/p^m X}, \mathcal{B}^\vee_{\widetilde{Z}}[p^m])$. We claim that this will be identified, up to a scalar in $(\mathbb{Z}/p^m\mathbb{Z})^\times$, with $[\rho]_{\widetilde{Z}}(g_{0,m})-g_{0,m} = \alpha_m\circ \rho_{Y, \widetilde{Z}}$ under the isomorphisms $\phi: Y\otimes \mathbb{Z}_p \stackrel{\sim}{\to}X\otimes \mathbb{Z}_p$ and $\lambda_B: \mathcal{B}_{\widetilde{Z}}[p^m]\stackrel{\sim}{\to}\mathcal{B}_{\widetilde{Z}}^\vee[p^m]$. This is the desired equivalence between points (2) and (2'). 
    
    We prove the claim. By Lemma~\ref{lem:splitting anti-diagonal}, we have a principal, anti-symmetric polarization $\lambda$ on $\mathbb{X}$ given by the formula 
    \[
    \lambda = \begin{pmatrix}
        0 & 0 & -\phi^\vee \\ 0& \lambda_{\widetilde{Z}} & 0 \\ 
        \phi & 0 & 0
    \end{pmatrix},
    \]
    where $\lambda_{\widetilde{Z}}$ is an anti-symmetric polarization on $\mathbb{X}_{\widetilde{Z}}$. The compatibility $\rho^\vee \circ \lambda \circ \rho = \lambda$,  up to a scalar in $(\mathbb{Z}/p^m\mathbb{Z})^\times$, implies that we have
    \begin{equation}\label{eq: compatibility with polarization}
    \rho_{\widetilde{Z}}^\vee\circ \lambda_{\widetilde{Z}}\circ \rho_{Y, \widetilde{Z}} = \rho_{\widetilde{Z},X}^\vee \circ \phi,
    \end{equation}
    once again up to a scalar. Together with the compatibility $\alpha_m^\vee \circ \lambda_{\widetilde{Z}}\circ \alpha_m =\lambda_B$, this implies that we have 
    \[
    (\rho_{\widetilde{Z},X}\circ \rho_{\widetilde{Z}}^{-1}\circ \alpha_m^{-1})^\vee \circ \phi = \lambda_B\circ \alpha_m\circ \rho_{Y,\widetilde{Z}},
    \]
    up to a scalar in $(\mathbb{Z}/p^m\mathbb{Z})^\times$. This proves the desired claim, and implies that the actions described in points (2) and (2') above are equivalent.
\end{proof}

\begin{remark}  Since the action described in point (3) above clearly lifts the action in point (2'), this also
    establishes the compatibility between (3) and (2).
\end{remark}

We now wish to extend the action of $T_p\mathcal{U}_b$ on $\Xi_{\widetilde{Z}}$ constructed in Definition~\ref{defn:dd action} to the relative torus embedding $\Xi_{\widetilde{Z},\Sigma_{\widetilde{Z}}}$. Let us denote by $T\hookrightarrow T_{\Sigma_{\widetilde{Z}}}$ the torus embedding used to construct $\Xi_{\widetilde{Z},\Sigma_{\widetilde{Z}}}$. We can then define 
$\Xi_{\widetilde{Z},\Sigma_{\widetilde{Z}}}$ as the contracted product 
\begin{equation}\label{eq:contracted product}
\left(\Xi_{\widetilde{Z}}\times T_{\Sigma_{\widetilde{Z}}}\right)/ T. 
\end{equation}
In light of (\ref{eq:contracted product}), in order to show that the action of $T_p\mathcal{U}_b$ extends to $\Xi_{\widetilde{Z},\Sigma_{\widetilde{Z}}}$, it is enough to show that it commutes with the action of $T$ on $\Xi_{\widetilde{Z}}$. 

\begin{lemma}\label{lem: commutation with torus} The action of $T_p\mathcal{U}_b$ on $\Xi_{\widetilde{Z}}$ from Definition~\ref{defn:dd action} commutes with the action of $T$ on $\Xi_{\widetilde{Z}}$. 
\end{lemma}

\begin{proof} For a scheme $S\to \Xi_{\widetilde{Z}}$, the action of a $S$-point $t:S\to T$ on the tuple $(\alpha_{m,S}, \beta_{m,S}, g_{0,m,S}, g_{m,S})$ leaves $\alpha_{m,S}$, $\beta_{m,S}$ and $g_{0,m,S}$ fixed and acts on $g_{m,S}$ by translation via the graph $t_S: S\to T_S$ of $t$. This can be easily seen to commute with the action of any $\rho\in T_p\mathcal{U}_b(S)$.  
\end{proof}

 By construction, the extended action of $T_p\mathcal{U}_b$ preserves the boundary of $\Sigma_{\widetilde{Z}, \Sigma_{\widetilde{Z}}}$ scheme-theoretically, so it induces an action on the formal completion $\mathfrak{X}_{\widetilde{Z},\Sigma_{\widetilde{Z}}}$. 
 Furthermore, for fixed $m$, this action depends on the restriction of an automorphism from $\underline{\mathrm{Aut}}_G(\mathbb{X})$ to $\underline{\mathrm{Aut}}_G(\mathbb{X}[p^m])$. On the other hand, the arithmetic group $\Gamma_{\widetilde{Z}}$ is contained inside the congruence subgroup of full level $p^m$, so it acts trivially on $\mathbb{X}_{Y}[p^m]$ and $\mathbb{X}_X[p^m]$. Therefore, the action of $T_p\mathcal{U}_b$ commutes with the  action of $\Gamma_{\widetilde{Z}}$ and we
 get an induced action of $T_p\mathcal{U}_b$ on $\mathfrak{X}_{\widetilde{Z},\Sigma_{\widetilde{Z}}}/\Gamma_{\widetilde{Z}}$. 

\begin{corollary}\label{cor: dd action boundary} For each Igusa cusp label $\widetilde{Z}$, the action in Definition~\ref{defn:dd action} induces an action of 
    $T_p\mathcal{U}_b$ on the formal completion $\widehat{\mathrm{Ig}}^{b,\mathrm{tor}}_{m,\widetilde{Z}}$ of $\mathrm{Ig}^{b,\mathrm{tor}}_{m}$ along the locally closed boundary stratum determined by $\widetilde{Z}$. 
\end{corollary}

\begin{proof} This follows from Theorem~\ref{thm:IgusaStratumDiagram} and the above discussion. 
\end{proof}

 It remains to glue the action of $T_p\mathcal{U}_b$ on the open Igusa variety $\mathrm{Ig}^b_m$ and on the formal completions $\widehat{\mathrm{Ig}}^{b,\mathrm{tor}}_{m, \widetilde{Z}}$ along boundary strata to an action on the whole partial toroidal compactification $\mathrm{Ig}^{b,\mathrm{tor}}_{\Sigma, m}$. The following proposition handles the gluing step.   

\begin{proposition}\label{prop:gluing action}
The action of $T_p\mathcal{U}_b$ on $\mathrm{Ig}^b_m$ described in \S\ref{Sec:Igusa background} and the actions of $T_p\mathcal{U}_b$ on the formal completions $\widehat{\mathrm{Ig}}^{b,\mathrm{tor}}_{m, \widetilde{Z}}$, as $\widetilde{Z}$ runs over (representatives of) equivalence classes of Igusa cusp labels at level $K^pK_{b}(p^m)$, glue uniquely to an action 
\[
T_p\mathcal{U}_b \times \mathrm{Ig}^{b,\mathrm{tor}}_{\Sigma, m}\to \mathrm{Ig}^{b,\mathrm{tor}}_{\Sigma, m}. 
\]
In particular, this action respects the boundary strata labeled by $\PP$. 
\end{proposition}

\begin{proof}
As in Lemma~\ref{lem:automorphisms at level 
m}, there exists a finite subgroup scheme $H'_m \subset \underline{\mathrm{Aut}}_G(\mathbb{X}[p^m])$ 
over $\ol{\mathbb{F}}_p$ characterized by the fact that the natural map $T_p\mathcal{U}_b \to \underline{\mathrm{Aut}}_G(\mathbb{X}[p^m])$ factors through $H'_m$, such that the induced map $T_p\mathcal{U}_b\to H'_m$ is faithfully flat. By construction, all the different actions of $T_p\mathcal{U}_b$ on $\mathrm{Ig}^b_m$ and on its formal completions along boundary strata factor through $H'_m$. We apply Lemma~\ref{lem:extendability is a local property} below to $T:=\mathscr{C}^{\mathrm{tor}}_{\Sigma}$, to the schematically dense quasi-compact open subset  
\[
U:=H'_m\times \mathrm{Ig}^b_m\hookrightarrow H'_m\times \mathrm{Ig}^{b,\mathrm{tor}}_{\Sigma, m}=:X,
\]
to the separated target $Y:=\mathrm{Ig}^{b,\mathrm{tor}}_{\Sigma, m}$, and to the morphism given by the action map described in \S\ref{Sec:Igusa background}. Therefore, it suffices to check the extendability of the morphism 
\begin{equation}\label{eq:action on open Igusa}
H'_m\times \mathrm{Ig}^b_m\to \mathrm{Ig}^{b,\mathrm{tor}}_{\Sigma, m},
\end{equation}
induced by the action map after taking base change to completions of strict local rings of $\mathscr{C}^{\mathrm{tor}}_{\Sigma}$ along boundary strata. 

We let $Z$ denote a representative of an equivalence class of usual cusp labels at level $K^p$ and $\sigma \in \Sigma_{Z}^+$ be a cone. We let $\bar{t}$ be a geometric point in the boundary stratum of $\mathscr{C}^{\mathrm{tor}}_{\Sigma}$ determined by the pair $(Z,\sigma)$. We let $R^{\wedge}$ denote the completion of the strict local ring of $\mathscr{C}^{\mathrm{tor}}_{\Sigma}$ at $\bar{t}$ along the ideal $I$ cut out by the boundary stratum labeled by $(Z,\sigma)$. 
As $\mathscr{C}^{\mathrm{tor}}_{\Sigma}$ is both excellent and normal (being even smooth over $\overline{\mathbb{F}}_p$), the completions of its strict local rings satisfy the hypotheses needed to apply the theory of degenerations of abelian varieties, cf.~\cite[\S 4.1]{lan-thesis}. 

In particular, $R^{\wedge}$ is a normal noetherian domain, complete with respect to an ideal $I$, and we let $\eta = \mathrm{Spec}\ K$ denote the generic point of $\mathrm{Spec}\ R^\wedge$. Over $\eta$, we have an abelian scheme $\mathcal{A}_{\eta}$ equipped with $G$-structures, which extends to a semi-abelian scheme over $\mathrm{Spec}\ R^\wedge$. The theory of degeneration induces, cf.~\cite[Prop. 5.2.2.1]{lan-thesis}, two short exact sequences on the level of the $p^m$-torsion:
    \begin{equation}\label{eq:first ses}
    0\to \mathcal{G}_{\eta}[p^m]\to \mathcal{A}_{\eta}[p^m]\to \mathbb{X}_{Y, \eta}[p^m] \to 0 
    \end{equation}
    and 
    \begin{equation}\label{eq:second ses}
    0\to \mathbb{X}_{X,\eta}[p^m]\to\mathcal{G}_{\eta}[p^m]\to \mathcal{B}_{\eta}[p^m] \to 0.  
    \end{equation}
By construction, all these objects are pulled back from the ``universal objects'' over the toroidal boundary chart $\Xi_{Z}(\sigma)$ along a canonical morphism $f:\mathrm{Spec}\ R^{\wedge}\to \Xi_{Z}(\sigma)$ that makes $\mathrm{Spec}\ R^{\wedge}$ into a ``good formal model'' in the sense of Faltings--Chai and Lan. Note that $R^\wedge$ is abstractly isomorphic to the completion of a strict local ring of $\Xi_{Z}(\sigma)$ at a geometric point $\bar{t}_{Z}$, but the canonical morphism $f$ differs from the morphism induced by this abstract isomorphism by an automorphism of $\mathrm{Spec}\ R^{\wedge}$.     

Fix a representative $\widetilde{Z}$ of an equivalence class of Igusa cusp labels at level $K^pK_{b}(p^m)$ that refines the equivalence class $Z$. Form the cartesian diagram 
\begin{equation}\label{eq: Cartesian diagram chart}
\begin{tikzcd}
    \mathrm{Spec}\ R^\wedge_{\widetilde{Z}}\arrow[r,"f_{\widetilde{Z}}"]\arrow[d] & \Xi_{\widetilde{Z}}(\sigma) \arrow[d] \\ 
    \mathrm{Spec}\ R^\wedge \arrow[r, "f"] & \Xi_{Z}(\sigma).
\end{tikzcd}
\end{equation}
We claim that we have, in addition, a commutative diagram
\begin{equation}\label{eq: Cartesian diagram non-chart}
\begin{tikzcd}
   \mathrm{Spec}\ R^\wedge_{\widetilde{Z}} \arrow[r]\arrow[dr] & \mathrm{Spec} R^\wedge_Z \arrow[r]\arrow[d] & \mathrm{Ig}^{b,\mathrm{tor}}_{\Sigma, m}\arrow[d] \\ \ &\mathrm{Spec}\ R^\wedge \arrow[r] & \mathscr{C}^{\mathrm{tor}}_{\Sigma},
\end{tikzcd}
\end{equation}
where the first horizontal morphism on the top is an open and closed immersion, and where the square is cartesian. The open and closed immersion corresponds to the choice of Igusa cusp label $\widetilde{Z}$ refining $Z$. 

To explain how to construct the diagram~\eqref{eq: Cartesian diagram non-chart}, note first that $R^\wedge_{\widetilde{Z}}$ is isomorphic to a finite product of completions of strict local rings of $\Xi_{\widetilde{Z}}(\sigma)$ at geometric points of $\Xi_{\widetilde{Z}}(\sigma)$ above $\bar{t}_{Z}$. This observation follows from the fact that the morphism $\Xi_{\widetilde{Z}}(\sigma)\to \Xi_{Z}(\sigma)$ is finite and from~\cite[Tag 05WR, Tag 00MA]{stacks-project}. In particular, $\mathrm{Spec}\ R^{\wedge}_{\widetilde{Z}}$ is normal.   
Let $U_{\widetilde{Z}}:=(\mathrm{Spec}\ R^{\wedge}_{\widetilde{Z}})_{\eta}$ be the pullback of $\mathrm{Spec}\ R^\wedge_{\widetilde{Z}}$ to the generic point $\eta$. Over $U_{\widetilde{Z}}$, we have additional data $(\alpha_{m,\eta}, \beta_{m,\eta}, g_{0,m,\eta}, g_{m,\eta})$ induced by pullback along the morphism $f_{\widetilde{Z}}$. Together with the decomposition~\eqref{eq:IgusaCuspLabelDecomposition}, this data can be reassembled into giving a trivialization 
\[
\gamma_{m,\eta}: \mathbb{X}_{\eta}[p^m]
\stackrel{\sim}{\to} \mathcal{A}_{\eta}[p^m]. 
\]
Indeed, $g_{m,\eta}$ induces a splitting of~\eqref{eq:first ses} under the identification $\mathcal{A}_{\eta}\simeq \mathcal{G}_{\eta}/Y_{\eta}$ coming from the theory of degeneration, and $\beta_{m,\eta}$ induces a splitting of~\eqref{eq:second ses}. The trivialization $\gamma_{m,\eta}$ induces a morphism $U_{\widetilde{Z}}\to \mathrm{Ig}^b_{m}$ that fits into a commutative diagram 
\[
\begin{tikzcd}
  U_{\widetilde{Z}} \arrow[r]\arrow[dr] & U_Z \arrow[r]\arrow[d] & \mathrm{Ig}^{b}_{m}\arrow[d] \\ \ & \eta \arrow[r] & \mathscr{C},
\end{tikzcd}
\]
where the first horizontal morphism on the top row is an open and closed immersion, determined by the choice of Igusa cusp label $\widetilde{Z}$, and the square is cartesian. This diagram extends to the desired diagram~\eqref{eq: Cartesian diagram non-chart} using the normality of $R_{\widetilde{Z}}$ and the universal property satisfied by $\mathrm{Ig}^{b,\mathrm{tor}}_{\Sigma, m}$ (in the sense of part 6 of~\cite[Theorem 6.4.1.1]{lan-thesis}). 

The formation of $\mathrm{Spec}\ R^\wedge_{\widetilde{Z}}$ as a fiber product in~\eqref{eq: Cartesian diagram chart} implies that the action of $H'_m$ on $\Xi_{\widetilde{Z}}(\sigma)$ constructed in Definition~\ref{defn:dd action} induces an action of $H'_m$ on $\mathrm{Spec}\ R^\wedge_{\widetilde{Z}}$. To conclude the desired extendability, it is enough to check that this action, when restricted to the complement of the boundary $U_{\widetilde{Z}}$, agrees with the action of $H'_m$ induced by the diagram~\eqref{eq: Cartesian diagram non-chart} from the action on $\mathrm{Ig}^b_m$. 

Choose a test scheme $S$ and an element $\rho\in H'_m(S)$, which determines an automorphism of $\mathbb{X}_{S}[p^m]$ compatible with the $G$-structures. Choose also an element of $U_{\widetilde{Z}}(S)$, which induces by pullback data $(\alpha_{m,S}, \beta_{m,S}, g_{0,m,S}, g_{m,S})$ and therefore also a trivialization $\gamma_{m,S}$. The action of $\rho$ induced from the action on $\mathrm{Ig}^b_m$ amounts to the pre-composition of $\gamma_{m,S}$ by $\rho$. To be able to compare this with the action of $\rho$ induced from Definition~\ref{defn:dd action},  we choose a non-canonical splitting of~\eqref{eq:first ses} and~\eqref{eq:second ses} over $\eta$, which is pulled back to a non-canonical splitting over $S$:
\begin{equation}\label{eq:non-canonical splitting}
\mathcal{A}_{S}[p^m]\simeq \mathbb{X}_{Y,S}[p^m]\oplus \mathcal{B}_{S}[p^m]\oplus \mathbb{X}_{X,S}[p^m].
\end{equation} 
Under the non-canonical splitting~\eqref{eq:non-canonical splitting} and under the decomposition~\eqref{eq:IgusaCuspLabelDecomposition} that is induced by the choice of Igusa cusp label $\widetilde{Z}$, we can write $\gamma_{m,S}$ as the $3\times 3$-matrix
\[
\gamma_{m,S} = \begin{pmatrix} 1 & \ & \ \\
\beta_{m,S, Y} & \alpha_{m,S} & \ \\
\gamma_{m,S, Y,X} & \beta_{m,S, X}& 1
\end{pmatrix} \in 
\begin{pmatrix}
    1 & \ & \ \\ 
    \mathrm{Hom}\left(\mathbb{X}_{Y,S}[p^m], \mathcal{B}_{S}[p^m]\right) & \mathrm{Isom}\left( \mathbb{X}_{\widetilde{Z},S}[p^m], \mathcal{B}_{S}[p^m]\right)& \ \\
    \mathrm{Hom}\left(\mathbb{X}_{Y, S}[p^m], \mathbb{X}_{X, S}[p^m]\right)& \mathrm{Hom}\left(\mathbb{X}_{\widetilde{Z},S}[p^m], \mathbb{X}_{X, S}[p^m]\right) & 1
\end{pmatrix},
\]
where the first column is induced by $g_{m,\bar{\eta}}$ and the second column by $\beta_{m,\bar{\eta}}\circ \alpha_{m,\bar{\eta}}$. Under the decomposition~\eqref{eq:IgusaCuspLabelDecomposition}, we identify $\rho$ with the $3\times3$ matrix
\begin{equation}\label{eq:matrix of rhoZ}
 \begin{pmatrix}
        1& \ &\ \\
        \rho_{Y, \widetilde{Z}} & \rho_{\widetilde{Z}}& \ \\
        \rho_{Y,X} & \rho_{\widetilde{Z},X} & 1
    \end{pmatrix}.    
\end{equation}
We can now check through a direct computation that the effect of $\rho$ on the trivialization $\gamma_{m,S}$, according to Definition~\ref{defn:dd action}, is the same as the effect given by pre-composition 
with the matrix~\eqref{eq:matrix of rhoZ}:
\[
\gamma_{m,S}\circ \rho = 
\begin{pmatrix}
    1 & \ & \ \\
    \beta_{m,S, Y} + \alpha_{m,S}\circ \rho_{Y, \widetilde{Z}} & \alpha_{m,S}\circ \rho_{\widetilde{Z}}& \ \\
    \gamma_{m,S, Y,X} + \beta_{m,S, X}\circ \rho_{Y,\widetilde{Z}}  + \rho_{Y,X}  & \beta_{m,S,X}\circ \rho_{\widetilde{Z}} + \rho_{\widetilde{Z},X} & 1
\end{pmatrix}. 
\]
This concludes the proof of the compatibility of the two actions of $H'_m$ on $U_{\widetilde{Z}}$. We conclude that the morphism~\eqref{eq:action on open Igusa} extends to a morphism 
\[
H'_m \times \mathrm{Ig}^{b,\mathrm{tor}}_{\Sigma, m} \to  \mathrm{Ig}^{b,\mathrm{tor}}_{\Sigma, m} . 
\]

We now claim that the action of $T_p\mathcal{U}_b$ respects the boundary strata in $\mathrm{Ig}^{b,\mathrm{tor}}_{\Sigma, m}$ labeled by the equivalence classes of parabolics $\PP$. Indeed, by construction, this action preserves the formal completion of $\mathrm{Ig}^{b,\mathrm{tor}}_{\Sigma, m}$ along the locally closed stratum labeled by $\widetilde{Z}$. To check that this action preserves the boundary stratum labeled by $\widetilde{Z}$ scheme-theoretically within this formal neighborhood, we use the boundary chart in Theorem~\ref{thm:IgusaStratumDiagram} and the explicit description of the action on $\Xi_{\widetilde{Z}, \Sigma_{\widetilde{Z}}}$. The claim now follows from the fact that the latter scheme is defined as a contracted product, cf.~(\ref{eq:contracted product}), and the action of $T_p\mathcal{U}_b$ is induced from the action on $\Xi_{\widetilde{Z}}$. This, in particular, means that the strata labeled by individual cones in $\Sigma_{\widetilde{Z}}$ are preserved scheme-theoretically. We conclude because the decomposition into cones refines the one labeled by $\PP$. 
\end{proof}

The proof of Proposition~\ref{prop:gluing action} is based on the following two abstract gluing results. 

\begin{lemma}\label{lem:abstract gluing} Let $X$ and $Y$ be qcqs schemes. Let $U\subset X$ be a schematically dense quasi-compact open subset. Assume that $Y$ is separated and that we have a morphism $f: U\to Y$. Assume also that, for an fpqc cover $X'\to X$, the pullback  $f\times_{X}X': U\times_X X'\to Y$ extends to a morphism $\bar{f}':X'\to Y$. Then $f$ extends uniquely to a morphism $\bar{f}:X\to Y$. 
\end{lemma} 

\begin{proof} The uniqueness of the extension follows from the schematic density of $U$ in $X$ and from the fact that the target $Y$ is separated, see~\cite[Tag 01RH]{stacks-project}. For the existence, we use fpqc descent with respect to the cover $X'\to X$. Note that, under the hypotheses, the property of being schematically dense is preserved under flat pullback, so that $U\times_{X}X'$ is schematically dense in $X'$ and $U\times_X X'\times_{X} X'$ is schematically dense in $X'\times_X X'$. To verify the compatibility with the descent datum to $X$, we use uniqueness once more to deduce that it is enough to verify the compatibility with the descent datum after restriction to the pre-image of $U$. On the pre-image of $U$, this compatibility is guaranteed by the existence of the original morphism $f:U\to Y$. 
\end{proof}

\begin{lemma}\label{lem:extendability is a local property} 
Let $T$ be a noetherian scheme and let $Y\to T$ be of finite presentation and separated. 
Let $X\to T$ be of finite presentation and $U\subset X$ a schematically dense quasi-compact open subset. 
Assume that we have a morphism $f:U\to Y$ over $T$ which satisfies the following property:
\begin{itemize}
\item  For every geometric point $\bar{t}\in T$, the base change 
\[
f\times_{T} \mathrm{Spec}\ R^\wedge : U\times_T \mathrm{Spec}\ R^{\wedge}\to Y
\] 
to a completion of the strict local ring of $T$ at $\bar{t}$ extends to a morphism
$\bar{f}\times_{T}\ \mathrm{Spec}\ R^{\wedge}: X\times_T \mathrm{Spec}\ R^\wedge \to Y$ over $T$. 
\end{itemize}
Then $f$ extends uniquely to a morphism $\bar{f}: X\to Y$. 
\end{lemma}

\begin{proof} Consider a geometric point $\bar{t}\in T$, with strict local ring $R$ and with completion $R^{\wedge}$. 
The morphism $\mathrm{Spec}\ R\to X$ is flat, so the schematic density 
of $U$ is preserved under pullback. The morphism $\mathrm{Spec}\ R^{\wedge}\to \mathrm{Spec}\ R$ is an fpqc
cover, because $R$ is a noetherian local ring. By Lemma~\ref{lem:abstract gluing}, we can descend $\bar{f}\times_T\mathrm{Spec}\ R^{\wedge}$ to 
a morphism $\bar{f}\times_{T}\ \mathrm{Spec}\ R$. 

We now write $R$ as a filtered colimit $R \stackrel{\sim}{\to} \varinjlim_{i} R_{i}$, 
where $\mathrm{Spec}\ R_i$ runs over affine \'etale neighbourhoods of $\bar{t}\in T$. Then
\[
X\times_{T}\mathrm{Spec}\ R \stackrel{\sim}{\to} \varprojlim_{i} \left(X\times_{T}\mathrm{Spec}\ R_i\right),
\]
where each $X\times_{T}\mathrm{Spec}\ R_i$ is qcqs (because $X$ over $T$ is) and the 
transition morphisms in the inverse limit are affine. By~\cite[Tag 01ZC]{stacks-project}, using the hypothesis
that $Y\to T$ is locally of finite presentation, we have 
\[
\mathrm{Hom}_{T}\left(X\times_T \mathrm{Spec}\ R, Y\right) \simeq 
\varinjlim_{i}\mathrm{Hom}_{T}\left(X\times_T \mathrm{Spec}\ R_i, Y\right). 
\]
This means that the morphism $\bar{f}\times_{T}\ \mathrm{Spec}\ R$ descends to a morphism 
$\bar{f}\times_{T}\ \mathrm{Spec}\ R_i$ for some $i$, i.e. descends to an \'etale neighbourhood of $\bar{t}\in T$. 
We conclude by another application of Lemma~\ref{lem:abstract gluing}, using \'etale descent instead of fpqc descent. 
\end{proof}
By taking inverse limits over all $m$ in Proposition~\ref{prop:gluing action}, we obtain the following. 
\begin{corollary}\label{cor: TpUbAction} The action of $T_p\mathcal{U}_b$ on $\mathrm{Ig}^b$ described in \S\ref{Sec:Igusa background} extends to an action of $T_p\mathcal{U}_b$ on $\mathrm{Ig}^{b,\mathrm{tor}}_{\Sigma}$ for any good compatible choice of cone decompositions $\Sigma$. In particular, this action preserves the boundary strata labeled by $\PP$, as in (\ref{eqn: BoudnaryStrataofToroidal}).
\end{corollary}

\subsubsection{The action on v-sheaves}\label{sec: JbActionV}
We draw from Theorem~\ref{thm: JbPreservesStrata} some conclusions about the $\widetilde{J}_{b,\ol{\bb{F}}_{p}}$-action on $\Ig^{b,\ast}$ when considered as v-sheaves. This will be crucial for establishing the stratified cartesian diagram~\eqref{Eq: stratifiedCartesianDiagramIntro}. More precisely, fix a neat prime-to-$p$ level $K^p$ (not necessarily principal). We first consider the v-sheaf attached to the perfect scheme $\Ig^{b,\ast}=\Ig^{b,\ast}_{K^p}$ and the formal group scheme $\widetilde{J}_{b,\ol{\bb{F}}_{p}}$ via the small diamond functor of \S \ref{ss: analytificationofperfectschemes}, which we denote by $\Ig^{b,\ast,\diamond}$ and $\mc{J}_b$ respectively. These are v-sheaves over $\Spd \ol{\bb{F}}_{p}$. We further take the canonical compactification of the former towards $\Spd k$ in the sense of \cite[Proposition~18.6]{Ecod} and denote it by $\ol{\Ig^{b,\ast,\diamond}}$. Note that alternatively this identifies with $\Ig^{b,\ast,\dagger}$ in light of Lemma \ref{lem: DaggerisCC}, where $(-)^{\dagger}$ is as in \S  \ref{ss: analytificationofperfectschemes}

Note that the action of $\widetilde{J}_{b,\ol{\bb{F}}_{p}}$ on $\Ig^{b,\ast}$ discussed in \S\ref{Sec:Action} induces an action of $\mc{J}_b$ on $\Ig^{b,\ast,\diamond}$ as follows: By definition, on an affinoid perfectoid test object $S=\Spa(R,R^+)$, we have 
\[\mc{J}_b(S)=\widetilde{J}_b(\Spf R^+)=\varprojlim_n \widetilde{J}_b(\Spec R^+/\varpi^n)\simeq \widetilde{J}_b(\Spec R^+/\varpi),\]
for some pseudo-uniformizer $\varpi$, where the last equivalence follows from Serre--Tate lifting of $p$-divisible groups along $p$-nilpotent rings, see \cite[Theorem 2.4.1]{CS2}. Whereas by definition, $\Ig^{b,\ast,\diamond}(S)=\Ig^{b,\ast}(\Spec R^+)$, but the latter also identifies with the $R^+/\varpi$-points $\Ig^{b,\ast}(\Spec R^+/\varpi)$, since $\Ig^{b,\ast}$ is perfect. Hence, the desired action comes from
\[\widetilde{J}_b(R^+/\varpi)\times \Ig^{b,\ast}(R^+/\varpi)\to \Ig^{b,\ast}(R^+/\varpi),\]
namely the $\Spec R^+/\varpi$-points of the action in Lemma~\ref{lemma: UniqueExistenceJbAction}. Moreover, since $\mc{J}_b$ is partially proper, this action extends to $\Ig^{b,\ast,\dagger}$. 
\begin{remark}[Warning!]
    Since $R^+$ is an $\ol{\bb{F}}_p$-algebra, we can also take the $\Spec R^+$-point of $\widetilde{J}_b$, considered as an fpqc sheaf over $\ol{\bb{F}}_p$. This is genuinely different from its $\Spf R^+$-points in light of the unipotent part $\mc{U}_{b}$ appearing in the decomposition (\ref{eq:Jb decomposition}), and generally contains much fewer points since $R^+$ is reduced.
\end{remark}
We recall that $\tilde{J}_{b,\ol{\bb{F}}_{p}}$ has a subgroup given by $\underline{J_{b}(\bb{Q}_{p})} \subset \tilde{J}_{b,\ol{\bb{F}}_{p}}$, as in (\ref{eq:Jb decomposition}). The $\Spec(R^{+})$ and $\Spf(R^{+})$-points do not differ on this subgroup. In particular, we may view the action of $\ul{J_{b}(\bb{Q}_{p})}^{\diamond} =\ul{J_{b}(\bb{Q}_{p})} \subset \mathcal{J}_{b}$ on $\Ig^{b,*,\diamond}$ as just being induced by applying the functor $(-)^{\diamond}$ to the action of $\ul{J_{b}(\bb{Q}_{p})} \subset \widetilde{J}_{b,\ol{\bb{F}}_{p}}$ on $\Ig^{b,*}$. Similarly, this induces an action of 
\[\ul{J_{b}(\bb{Q}_{p})} = \ul{J_{b}(\bb{Q}_{p})}^{\diamond} = \ul{J_{b}(\bb{Q}_{p})}^{\dagger} = \ul{J_{b}(\bb{Q}_{p})}^{\Diamond}\] 
on $\Ig^{b,*,\dagger} = \ol{\Ig^{b,*,\diamond}}$ and $\Ig^{b,*,\Diamond}$. Here the equalities follow from \cite[Proposition~2.19]{GleasonIvanovMeromorphicVectorBundles} and partial properness. By construction, the natural maps 
\[ \Ig^{b,*,\diamond} \ra \Ig^{b,*,\dagger} \ra \Ig^{b,*,\Diamond} \]
given by evaluating (\ref{eqn: NaturalTransformationsofAnalytifications}) on $X = \Ig^{b,*}$ are $\ul{J_{b}(\bb{Q}_{p})}$-equivariant. It in turn defines a sequence of natural maps
\begin{equation}{\label{eqn: NaturalTransformationofStacks}} 
[\Ig^{b,*,\diamond}/\ul{J_{b}(\bb{Q}_{p})}] \ra [\Ig^{b,*,\dagger}/\ul{J_{b}(\bb{Q}_{p})}] \ra [\Ig^{b,*,\Diamond}/\ul{J_{b}(\bb{Q}_{p})}]  
\end{equation}
of $v$-stacks. This has the following basic properties.
\begin{lemma}\label{lemma: IgCCqcqs}
The following is true. 
\begin{enumerate}
\item The natural map 
\[[\Ig^{b,\ast,\diamond}/\mc{J}_b]\hookrightarrow [\Ig^{b,\ast,\dagger}/\mc{J}_b]\]
is a qcqs open immersion.
\item The natural map 
\[ [\Ig^{b,\ast,\diamond}/\ul{J_{b}(\bb{Q}_{p})}] \hookrightarrow [\Ig^{b,\ast,\Diamond}/\ul{J_{b}(\bb{Q}_{p})}] \]
is a qcqs open immersion.
\end{enumerate}
\end{lemma}
\begin{proof}
We note that, using \cite[Proposition~10.11]{Ecod}, it suffices to show that $\Ig^{b,\ast,\diamond}\hookrightarrow \Ig^{b,\ast,\dagger}$ and $\Ig^{b,*,\diamond} \ra \Ig^{b,*,\Diamond}$ are qcqs open immersions. To see this, we recall that it admits a natural pro-proper map $\Ig^{b,\ast} \ra \mathscr{C}^{b,\ast}$ to the central leaf given by the presentation $\Ig^{b,\ast} = \lim_{m \geq 1} \Ig^{b,\ast}_{m}$. 

If we look at the natural maps $\mathscr{C}^{b,\ast,\diamond} \ra \mathscr{C}^{b,\ast,\dagger}$ and $\mathscr{C}^{b,\ast,\diamond} \ra \mathscr{C}^{b,\ast,\Diamond}$ then these are qcqs open immersions by \cite[Lemma 4.2.2]{GHILZIsocComparison}, and the fact that $\mathscr{C}^{b,\ast}$ is affine and of finite type. Moreover, since each of the maps $\Ig^{b,\ast}_{m} \ra \mathscr{C}^{b,\ast}$ are finitely-presented proper, we have a cartesian diagram 
\[ 
\begin{tikzcd}
\Ig^{b,\ast,\diamond}_{m} \arrow[r] \arrow[d] & \Ig^{b,\ast,\dagger}_{m} \arrow[d] \arrow[r]  & \Ig^{b,\ast,\Diamond}_{m} \arrow[d] \\
\mathscr{C}^{b,\ast,\diamond} \arrow[r] & \mathscr{C}^{b,\ast,\dagger} \arrow[r] & \mathscr{C}^{b,\ast,\Diamond}
\end{tikzcd}
\]
for all $m \geq 1$ in light of Lemma \ref{lemma: CartesianGivesProper}. By virtue of the fact that the analytification functors of \S \ref{ss: analytificationofperfectschemes} all commute with limits, this in turn gives a cartesian diagram, 
\[ 
\begin{tikzcd}
\Ig^{b,\ast,\diamond} \arrow[r] \arrow[d] & \Ig^{b,\ast,\dagger} \arrow[d] \arrow[r]  & \Ig^{b,\ast,\Diamond} \arrow[d] \\
\mathscr{C}^{b,\ast,\diamond} \arrow[r] & \mathscr{C}^{b,\ast,\dagger} \arrow[r] & \mathscr{C}^{b,\ast,\Diamond},
\end{tikzcd}
\]
which implies the desired claim.
\end{proof}

We now have the following.
\begin{proposition}\label{prop: JbPreserveBoundary}
    For $? \in \{\diamond,\dagger\}$, the $\mc{J}_b$-action on $\mathrm{Ig}^{b,\ast,?}$ preserves the closure of the Baily--Borel boundary strata, i.e. for each $\PP$, it restricts to an action on the closed subdiamond $\mathrm{Ig}^{b,?}_{\leq\PP}$, defined by analytifying the closed strata attached to (\ref{eqn: BailyBorelStratificationonIgusavariety}). Moreover, the $\ul{J_{b}(\bb{Q}_{p})}$-action preserves the Baily-Borel strata of $\mathrm{Ig}^{b,\ast,?}$ for $? \in \{\diamond,\dagger,\Diamond\}$.
\end{proposition}
\begin{proof}
For the first part, we explain the case of $? = \diamond$, and the case of $\dagger$ follows by taking canonical compactifications. We need to show that the action map 
    \[\mc{J}_b\times \mathrm{Ig}^{b,\diamond}_{\leq\PP}\to \mathrm{Ig}^{b,\ast,\diamond} \]
    factors through the closed immersion 
    \[ \mathrm{Ig}^{b,\diamond}_{\leq\PP} \hookrightarrow \mathrm{Ig}^{b,\ast,\diamond}. \]
    Since the latter is a quasi-compact injection, it follows by \cite[Proposition 11.20]{Ecod} that it suffices to check on the level of geometric points that the map factors. In particular, the statement can be checked after passing to a deeper level. Indeed, suppose $K'^{p}\subset K^p$ is another prime-to-$p$ level, such that the result is true for $K'^{p}$. Then for any geometric point $x$ of $ \Ig^{b,\diamond}_{K^p,\leq \PP}$, one can lift it to a geometric point $\tilde{x}$ of $ \Ig^{b,\diamond}_{K'^p,\leq \PP}$. Then the desired result for $\Ig^{b,\diamond}_{K^p,\leq \PP}$ follows from the fact that the map $\Ig^{b,\diamond}_{K'^p,\leq \PP}\to \Ig^{b,\diamond}_{K^p,\leq \PP}$ is $\mc{J}_b$-equivariant and respects the Baily--Borel stratification.
    
    We may therefore assume $K^p=K(N)$, $p\nmid N$ is a principal level. In this case, take a $\Spd (C,C^{+})$-point $(x,g)$ of $\mathrm{Ig}^{b,\diamond}_{\leq\PP}\times \mc{J}_b$, for some complete algebraically closed non-archimedean field $C$ of characteristic $p$, with a pseudo-uniformizer $\varpi$. Since $\mc{J}_b(C,C^{+}) = \mc{J}_{b}(C,\mathcal{O}_{C})$ is generated by $J_b(\qp)^+$ and $T_p\mc{U}_b(\mc{O}_C/\varpi)$\footnote{Indeed, the unipotent part $\mc{J}^{>0}_b(C)$ can be identified as the union $\cup_{n} p^{-n}T_p\mc{U}_b(\mc{O}_C/\varpi)$.} as a group, we may assume $g$ lies in either one of them. On the other hand, the point $x$ corresponds to a $\Spec C^{+}$-point of $\Ig^b_{\leq \PP}$. Enlarging $C$ if necessary, we may assume it lifts to a  $\Spec C^{+}$-point $\tilde{x}$ of $\mathrm{Ig}^{b,\tor}_{\Sigma,\leq\PP}$, for some suitable choice of cone decomposition $\Sigma$. This is possible, since the map $\Ig^{b,\tor}_{\Sigma,\leq \PP}\to \Ig^b_{\leq \PP}$ from the corresponding stratum on the partial toroidal compactification is a v-cover (using \cite[Proposition~3.3.4]{CS2}, see \cite[Definition 2.1, Example 2.3(ii)]{BSWitt}). By perfectness of $\mathrm{Ig}^{b,\tor}_{\Sigma,\leq\PP}$, one can also think of $\tilde{x}$ as a $\Spec C^{+}/\varpi$-point. The desired result now follows from Theorem~\ref{thm: JbPreservesStrata} and Proposition \ref{prop: globalsections}.

    For the remaining claims, since we are only concerned with $\ul{J_{b}(\bb{Q}_{p})}$-action, we may (using Remark \ref{rem: representablevsheafifyingautomatic} and the above discussion) reduce to showing that, given a perfect ring $R$ and a $\Spec R$-point of $\Ig^{b,*}$ factoring through the closed subscheme $\Ig^{b,*}_{\leq \PP} \hookrightarrow \Ig^{b,*}$ that the action by an $R$-point of $\underline{J_{b}(\bb{Q}_{p})}$ also lies inside $\Ig^{b}_{\leq \PP}$. However, using \cite[Lemma~3.10]{HamannLeeTorsion}, this reduces to checking the analogous claim on the underlying topological spaces. Now we may argue as above by lifting such a geometric point along the $v$-cover $\mathrm{Ig}^{b,\tor}_{\Sigma,\leq\PP} \ra \Ig^{b}_{\leq \PP}$.
\end{proof}
The proposition implies similar results for individual strata labeled by the $\PP$ conditions as well. More precisely, for each conjugacy class $\PP$, we let $\mathrm{Ig}_{<\PP}^b\subset \Ig_{\leq \PP}^b$ be the closed complement to the stratum $\Ig^b_\PP$, with reduced induced scheme structure. For $? \in \{\diamond,\dagger,\Diamond\}$, we let $\mathrm{Ig}_{<\PP}^{b,?}$ be its associated analytification. This is a closed subsheaf of $\Ig_{\leq \PP}^{b,?}$. We can consider its open complement by defining

\begin{equation}{\label{eqn: ThickenedBoundaryStrata}}
(\Ig^{b,?})_{\PP}:=\Ig^{b,?}_{\leq\PP}\times_{\left|\Ig^{b,?}_{\leq\PP}\right|}\left(\left|\Ig^{b,?}_{\leq\PP}\right|\backslash{\left|\Ig^{b,?}_{<\PP}\right|}\right).
\end{equation}
\begin{remark}[Warning!]\label{remark: differencestrata}
    To the stratum $\Ig^b_{\PP}$, we can also attach a v-sheaf $\mathrm{Ig}^{b,?}_{\PP}$ for $? \in \{\diamond,\dagger,\Diamond\}$,  but this will be very different from $(\Ig^{b,?})_{\PP}$ if $? \in \{\diamond,\dagger\}$. For example, when $\PP=[\mathsf{G}]$, the $v$-sheaf $(\Ig^{b,?})_{\PP}$ is $\Ig^{b,?} \setminus \partial^{?}$, where $\partial \subset \Ig^{b,*}$ is the boundary of the partial minimal compactification. This contains $\Ig^{b,?}$, and their difference is roughly a punctured tubular neighborhood of $\partial^{?}$. Indeed, $\Ig^{b,\dagger}$ will model the fiber of the Hodge-Tate period morphism restricted to the good reduction locus of the Shimura variety, but $(\Ig^{b,\dagger})_{[\mathsf{G}]}$ will model the fiber of the Hodge-Tate period morphism over the entire open Shimura variety (see Corollary \ref{cor: BeyondGoodRed}). 
    
    Relatedly, we note that the constructions $(-)^{\diamond}$ and $(-)^{\dagger}$ in \S \ref{ss: analytificationofperfectschemes} behave well with respect to closed immersions (in the sense of Lemma \ref{lemma: CartesianGivesProper}), but not with respect to open immersions. We deliberately choose the slightly clunky notation $(-)_{\PP}$ to distinguish the two. 
\end{remark}
We now check that these strata are preserved under the group actions.
\begin{corollary}\label{cor: JbPreserveBoundary}
Let $\PP$ be a conjugacy class of admissible parabolics. For $? \in \{\diamond,\dagger\}$, $(\Ig^{b,?})_{\PP}$ is a $\mc{J}_b$-stable sub-v-sheaf of $\Ig^{b,\ast,?}$. For $? \in \{\diamond,\dagger,\Diamond\}$, $(\Ig^{b,?})_{\PP}$ is a $\ul{J_{b}(\bb{Q}_{p})}$-stable sub-$v$-sheaf. 
\end{corollary}
\begin{proof}
    The same argument as Proposition~\ref{prop: JbPreserveBoundary} applies by replacing $\Ig^{b,*,?}$ (resp. $\Ig^{b,*,?}_{\leq \PP}$) with $\Ig^{b,*,?}_{\leq \PP}$ (resp. $(\Ig^{b,*,?})_{\PP}$). 
\end{proof}

We can diagrammatically organize the relationships between these different strata as follows.
\begin{lemma}{\label{lemma: differentstrataofIgs}}
Let $\PP$ be a conjugacy class of admissible parabolics. Then there exists a commutative diagram with cartesian squares 
\[ 
\begin{tikzcd}
& (\Ig^{b,\diamond})_{\PP} \arrow[r] \arrow[d] & \Ig^{b,\diamond}_{\leq \PP} \arrow[d] \arrow[r] & \Ig^{b,*,\diamond} \arrow[d] \\
& (\Ig^{b,\dagger})_{\PP} \arrow[r] \arrow[d] & \Ig^{b,\dagger}_{\leq \PP} \arrow[d] \arrow[r] & \Ig^{b,*,\dagger} \arrow[d] \\
& \Ig^{b,\Diamond}_{\PP} = (\Ig^{b,\Diamond})_{\PP} \arrow[r]  & \Ig^{b,\Diamond}_{\leq \PP} \arrow[r] & \Ig^{b,*,\Diamond}, 
\end{tikzcd}
\]
where the horizontal arrows are the natural inclusions and the rightmost column is defined by evaluating the transformation $(-)^{\diamond} \ra (-)^{\dagger} \ra (-)^{\Diamond}$ on $\Ig^{b,*}$. Similarly, we have a diagram with cartesian squares by replacing each term in the above diagram with their stack quotients by $\ul{J_b(\qp)}$, which is well-defined in light of Proposition~\ref{prop: JbPreserveBoundary} and Corollary \ref{cor: JbPreserveBoundary}. Moreover, if we ignore the bottom row then we have an analogous diagram involving stack quotients by $\mathcal{J}_{b}$.
\end{lemma}
\begin{proof}
If we look at the fpqc subsheaf $\underline{J_{b}(\bb{Q}_{p})} \subset \tilde{J}_{b,\ol{\mathbb{F}}_{p}}$ then we recall that we have an equality  $\underline{J_{b}(\bb{Q}_{p})}^{\diamond} = \underline{J_{b}(\bb{Q}_{p})}^{\dagger} =  \underline{J_{b}(\bb{Q}_{p})}^{\Diamond} = \ul{J_{b}(\bb{Q}_{p})}$ by \cite[Proposition~2.19]{GleasonIvanovMeromorphicVectorBundles}. This reduces the claim on the stacks quotient to the first claim. For the first claim,  the Cartesian squares on the right follow from an immediate application of Lemma \ref{lemma: CartesianGivesProper}. For the cartesian squares on the left, we note, by Lemma \ref{lemma: CartesianGivesProper} applied to the closed immersion $\Ig^{b}_{< \PP} \hookrightarrow \Ig^{b}_{\leq \PP}$, we have a commutative diagram
\[ 
\begin{tikzcd}
& \Ig^{b,\diamond}_{< \PP} \arrow[r]  \arrow[d] &  \Ig^{b,\diamond}_{\leq \PP} \arrow[d] \\
& \Ig^{b,\dagger}_{< \PP} \arrow[r]  \arrow[d] &  \Ig^{b,\dagger}_{\leq \PP} \arrow[d] \\
& \Ig^{b,\Diamond}_{< \PP} \arrow[r] & \Ig^{b,\Diamond}_{\leq \PP}, 
\end{tikzcd}
\]
with cartesian squares. We conclude by passing to the complementary open and using the relationship $\Ig^{b,\Diamond}_{\leq \PP} \setminus \Ig^{b,\Diamond}_{< \PP} = (\Ig^{b}_{\leq \PP} \setminus \Ig^{b}_{< \PP})^{\Diamond} = \Ig^{b,\Diamond}_{\PP}$ from Lemma \ref{lem: ExcisionBigDiamond}.
\end{proof}

\section{Minimally compactified Igusa stacks}\label{sec: IgusaStack}
We can now move our focus to Igusa stacks. We start by defining a Baily--Borel stratification on the minimally compactified Igusa stack and analyzing the geometry of the strata, based on the results from the previous section. This allows us to define intersection complexes attached to \'etale local systems on the open Igusa stack. We show that they recover the corresponding intersection complexes on the Shimura varieties. Along the way, we obtain a dimension formula for the boundary strata on classical partially compactified perfect Igusa varieties.

\subsection{Igusa stacks}\label{sec: IgusaStackNotation}
Let $\gx$ be a Shimura datum satisfying Assumption~\ref{assumption:codimension}, with reflex field $\mathsf{E}$. Fix an isomorphism $\ol{\bb{Q}}_{p} \simeq \bb{C}$, with corresponding $p$-adic place $v$ of $\mathsf{E}$ and denote by $E$ the completion of $\mathsf{E}$ at $v$. Let $\mu$ be the geometric dominant cocharacter of $G := \mathsf{G}_{\bb{Q}_{p}}$ defined by the Hodge cocharacter and this isomorphism. Fix a neat compact open subgroup $K\subset\mathsf{G}(\mathbb{A}_f)$ and let $\mc{S}_K$ be the corresponding Shimura variety at level $K$, viewed as a diamond over $\Spd E$. Recall from \cite[Theorem 9.40]{zhang2023}, \cite[Theorem D]{Kim}\footnote{The original construction in \cite{zhang2023} only gave the good reduction locus of the Igusa stack and \cite{Kim} extends this to the Igusa stack that corresponds to the whole open Shimura variety. The result of the current section gives another way of extending the Igusa stack beyond the good reduction locus. Our result logically does not depend on \cite[Theorem D]{Kim}, see Corollary~\ref{cor: BeyondGoodRed} below.} that we have a corresponding small v-stack $\Igs_{K^p}$ and its minimal compactification $\Igs_{K^p}^\ast$ (both defined over $\Spd \mathbb{F}_p$), together with the following Hecke equivariant cartesian diagrams of $v$-stacks 
\begin{equation}{\label{eqn: CartesiandiagramOpenShimuraVariety}}
 \begin{tikzcd}
\mc{S}_{K^p} \arrow[r,"\pi"] \arrow[d,"\tilde{h}"] & \mathcal{F}\ell_{G,\mu^{-1}} \arrow[d,"h"]\\ 
\Igs_{K^p} \arrow[r,"\overline{\pi}"] & \Bun_{G},
\end{tikzcd}
\end{equation}
and 
\begin{equation}{\label{eqn: CartesiandiagramShimuraVarietyInfiniteLevel}}
 \begin{tikzcd}
\mc{S}^\ast_{K^p} \arrow[r,"\pi_\mathrm{min}"] \arrow[d,"\tilde{h}_\mathrm{min}"] & \mathcal{F}\ell_{G,\mu^{-1}} \arrow[d,"h"] \\ 
\Igs^\ast_{K^p} \arrow[r,"\overline{\pi}_{\mathrm{min}}"] & \Bun_{G},
\end{tikzcd}
\end{equation}
where $\mc{F}\ell_{G,\mu^{-1}}$ is the diamond attached to the partial flag variety for $G$ via the analytification operation $(-)^{\Diamond}$ of \S\ref{ss: AnalytificationinMixedCharacteristic}. It corresponds to the Schubert cell labeled by $\mu^{-1}$ on the $B_{dR}^+$-affine Grassmannian $\operatorname{Gr}_G$. The $v$-stack $\Bun_{G}$ is the moduli stack sending $S \in \Perf$ to the groupoid of $G$-bundles on the relative (adic) Fargues-Fontaine curve $X_{S}$ over $S$. Moreover, the map $\pi$ is the Hodge--Tate period map, and the right vertical map $h$
is the Beauville--Laszlo map. Since the diagrams are equivariant for the $G(\qp)$-actions on the geometric objects, which are trivial on the bottom row, we have similar versions at each level $K_p\subset G(\qp)$, where the top row is replaced by the map between the corresponding stacky quotients. These are related via the commutative diagram
\begin{equation}{\label{eqn: DiagramofIgusaStacks}}
\begin{tikzcd}
\mathcal{S}_K \simeq \left[\mathcal{S}_{K^p}/\ul{K_{p}}\right] \arrow[r,"j_{K_{p}}",hook] \arrow[d,"\tilde{h}_{K_{p}}"] & \left[\mathcal{S}^\ast_{K^p}/\ul{K_{p}}\right] \arrow[d,"\tilde{h}_{K_{p},\min}"] \ar[r, "\pi_{K_p,\min}"]& \left[\mc{F}\ell_{G,\mu^{-1}}/\ul{K_p}\right] \ar[d,"h_{K_p}"]\\ 
\Igs_{K^p} \arrow[r,"\overline{j}",hook] & \Igs^\ast_{K^p} \ar[r,"\overline{\pi}_\mathrm{min}"] & \Bun_G,
\end{tikzcd}
\end{equation}
with cartesian squares.

The bottom rows of the previous diagrams depend only on the fixed compact open subgroup $K^{p} \subset \mathsf{G}(\mathbb{A}_{f}^{p})$ and not on the level at $p$, which justifies the notation. The vertical maps are $\ell$-cohomologically smooth of pure $\ell$-dimension $d := \langle 2\rho_{G}, \mu \rangle$, where $\rho_{G}$ is the half sum of the positive roots. More precisely, we have the following.

\begin{lemma}{\label{lemma: smoothnessofvertmaps}}
For all compact open subgroups $K_{p} \subset G(\mathbb{Q}_{p})$, the map $h_{K_{p}}$  is fdcs $\ell$-cohomologically smooth of pure $\ell$-dimension $d$, and there are natural identifications of functors
\[ h^{!}_{K_{p}}(-) \simeq h^{*}_{K_{p}}(-)[2d](d). \]
In particular, in light of the cartesian diagram (\ref{eqn: DiagramofIgusaStacks}), the maps $\tilde{h}_{K_{p,\min}}$ and $\tilde{h}_{K_{p}}$ are also cohomologically smooth of pure $\ell$-dimension $d$ and surjections on the underlying topological spaces; therefore, both maps are surjections of $v$-sheaves.
\end{lemma}
\begin{proof}
This follows from the fact that the dualizing object of $\Bun_{G}$ is the constant sheaf \cite[Proposition 1.1]{HamannImai}, cf. \cite[Lemma 5.1.8]{DvHKZ2} and the fact that the Beauville--Laszlo map $\left[\mathcal{F}\ell_{G,\mu^{-1}}/\ul{K_{p}}\right]\to \Bun_G$ is $\ell$-cohomologically smooth of $\ell$-dimension $d$ with fibers given by \'etale torsors over the flag variety $\mathcal{F}\ell_{G,\mu}$, cf.~\cite[Proof of Theorem IV.1.19]{FSGeomLLC}. When combined with Proposition \ref{prop: propertiesofthealgebraizationfunctor} (5), this reduces us to the calculation of the dualizing complex on the algebraic flag variety, which, since this is a smooth projective variety of dimension $d$, is given by $\Lambda[2d](d)$. To see that the maps $\tilde{h}_{K_{p,\min}}$ and $\tilde{h}_{K_{p}}$ are surjections of the underlying topological spaces, we note that the image of $\ol{\pi}_{\min}$ and $\ol{\pi}$ is given by the open substack $\Bun_{G,\mu^{-1}} \subset \Bun_{G}$ corresponding to $B(G,\mu^{-1}) \subset B(G)$. Moreover, $h_{K_{p}}$ also has image given by $B(G,\mu^{-1})$, by combining \cite[Proposition~A.9]{RapAppendixtoPadicCohomology} with the interpretation of $h$ as the forgetful map sending a modification of $G$-bundles $\mathcal{E} \dashrightarrow \mathcal{E}_{0}$ of type $\mu$ of the trivial $G$-bundle $\mathcal{E}_{0}$ on the Fargues-Fontaine curve to the $G$-bundle $\mathcal{E}$. This implies the desired claim.
\end{proof}

\subsection{A stratified cartesian diagram}\label{sec: Stratified}

The goal of this section is to show that we have a stratified version of the cartesian diagram exhibited in Equation (\ref{eqn: CartesiandiagramShimuraVarietyInfiniteLevel}), with respect to the Baily--Borel stratification. By Lemma~\ref{lemma: smoothnessofvertmaps} and Lemma~\ref{lemma: stackystratasmooth}, this will lead to the $\ell$-cohomological smoothness of the boundary strata of the Igusa stack, which will play an important role in our construction of the intersection complex on the Igusa stack. From now on,  we fix an algebraically closed field $k/\mathbb{F}_p$ (with discrete topology) and base-change $\Igs^\ast\coloneq \Igs^\ast_{K^p}\gx$ to $\Spd k$. Similarly, we base-change $\Bun_{G}$ and the Igusa varieties $\Ig^{b}$ and its compactifications introduced in \S \ref{sec: pc Igusa varieties} to $k$.

For simplicity of notation, throughout this and the next subsection, we fix the sufficiently small prime-to-$p$ level $K^p$ and write $\Igs^\ast$, $\mathcal{S}^\ast$, $\mathcal{S}_\PP$, $\mathcal{S}_{\leq\PP}$, $\mc{F}\ell$ for $\Igs^\ast_{K^p}$, $\mathcal{S}^\ast_{K^p}$ $\mathcal{S}_{K^p,\PP}$, $\mathcal{S}_{K^p,\leq \PP}$, and $\Fl$. For any of these objects, the superscript $(-)^{b}$ will denote a further base change along $\mathrm{Bun}^b_G\hookrightarrow \mathrm{Bun}_G$, where $i_{b}: \Bun_{G}^{b} \hookrightarrow \Bun_{G}$ is the locally closed Harder-Narasimhan stratum corresponding to $b \in B(G)$.

 \subsubsection{Baily--Borel stratification on the Igusa stack}
\begin{definition}
    For a rational conjugacy class of admissible parabolics $\PP$, we define the substacks
    \[\Igs_{[\mathsf{P}]}\subset \Igs_{\leq [\mathsf{P}]} \subset  \Igs^\ast\]
    to be the $v$-sheaf-theoretic images of $\mathcal{S}_{\PP}$, resp. $\mathcal{S}_{\leq \PP}$ along $\tilde{h}_{\min}: \mathcal{S}^\ast\to \Igs^\ast$. Namely, these are substacks of $\Igs^\ast$ defined by the condition that a map $S\to \Igs^\ast$ factors through $\Igs_{\PP}$, resp. $\Igs_{\leq \PP}$ if and only if there exists a v-cover $\Tilde{S}\to S$ and a map $\Tilde{S}\to \mathcal{S}_{\PP}$, resp. $\Tilde{S}\to \mathcal{S}_{\leq \PP}$ lifting the given map.
\end{definition}

\begin{remark}\label{Rem:FactorThroughP}
By definition, the natural map $\mathcal{S}_{\PP}\to \Igs^\ast$ factors through $\Igs_{\PP}$, similarly if we replace $\PP$ by $\leq \PP$.    
\end{remark}

The following is the main result of this section.
\begin{theorem}\label{Thm:StratifiedCartesian}
    For each rational conjugacy class of admissible parabolics $\PP$, the cartesian diagram (\ref{eqn: CartesiandiagramShimuraVarietyInfiniteLevel}) restricts to a cartesian diagram
    \begin{equation}\label{Eq:StratifiedCartesian}
    \begin{tikzcd}
        \mathcal{S}_{\leq\PP} \ar[r] \ar[d,"\tilde{h}_{\leq [\mathsf{P}]}"] & \mc{F}\ell\ar[d]\\
        \Igs_{\leq\PP} \ar[r] & \Bun_G,
    \end{tikzcd}
    \end{equation}
    where the left vertical map arises as in Remark~\ref{Rem:FactorThroughP}. Similar statements hold if we replace $\mathcal{S}_{\leq\PP}$ and $\Igs_{\leq\PP}$ by $\mathcal{S}_{\PP}$ and $\Igs_{\PP}$, respectively.
\end{theorem}

\begin{remark}{\label{rem: IgusaStackReadyforPeversity}}
    This means we have a set-theoretic locally-closed stratification of the underlying topological space
    \[|\Igs^\ast|=\coprod_{\PP}|\Igs_{\PP}|,\]
    which we will refer to as the \textit{Baily--Borel stratification} on $\Igs^\ast$, since this is true on the $v$-cover $\tilde{h}: \mathcal{S}^{*} \ra \Igs^{*}$ using Lemma \ref{lem: ExcisionBigDiamond}. In particular, this will put us in a situation in which we can apply the discussion of \S~\ref{sec: tStructureVstacks} after showing these strata are cohomologically smooth Artin stacks (Corollary~\ref{cor: IgusaStackStrataSmooth} below).
\end{remark}

We will prove Theorem~\ref{Thm:StratifiedCartesian} by establishing a version of~\eqref{Eq:StratifiedCartesian} after pullback to each individual Newton stratum, based on the following lemma. 

\begin{lemma}\label{Lem:FiberwiseCartesian} Assume that, for each rational conjugacy class of admissible parabolics $\PP$ and each $b\in B(G,\mu^{-1})$, Diagram~\eqref{Eq:StratifiedCartesian} becomes cartesian after pullback along $\Bun_G^b\hookrightarrow\Bun_G$. Then Diagram~\eqref{Eq:StratifiedCartesian} is cartesian. The analogous statement holds if we replace $\leq \PP$ by $\PP$.
\end{lemma}

\begin{proof} 
We spell out the proof for the $\leq \PP$ condition. The one for $\PP$ is similar. Write $X_{\leq\PP}$ for the fiber product $\Igs_{\leq\PP}\times_{\Bun_G}\mc{F}\ell$. The commutativity of~\eqref{Eq:StratifiedCartesian} and the injectivity of the map
\[\mc{S}_{\leq \PP}\to \mc{S}^\ast\simeq \Igs^\ast\times_{\Bun_G}\mc{F}\ell\] 
induces an injective map of v-sheaves 
\begin{equation}\label{Eq:ComparisonMap}
\mathcal{S}_{\leq\PP}\hookrightarrow X_{\leq\PP}.
\end{equation}
Conversely, under our assumption the natural inclusion of topological spaces $|X_{\leq\PP}|\hookrightarrow |\mathcal{S}^\ast|$ factors as 
\[
|X_{\leq\PP}|\hookrightarrow |\mathcal{S}_{\leq\PP}|\hookrightarrow |\mathcal{S}^\ast|.
\]
Indeed, as this is a statement on topological spaces, it can be checked Newton stratum-wise, using the decomposition 
 \[
 |\mathcal{S}^\ast|= \coprod_{b\in B(G,\mu^{-1})} |\mathcal{S}^{\ast,b}|
 \]
and the analogous ones for $|X_{\leq\PP}|$ and $|\mathcal{S}_{\leq\PP}|$. 

This means that, under our assumption, for any affinoid perfectoid space $S$ with a map $S\to X_{\leq\PP}$, the composition 
\[
S\to X_{\leq\PP} \hookrightarrow \mathcal{S}^\ast
\]
factors, on the level of topological spaces, as $|S|\to |\mathcal{S}_{\leq\PP}|\hookrightarrow |\mathcal{S}^\ast|$. Now applying the equivalence of $(i)$ and $(iv)$ in \cite[Proposition 5.3]{Ecod} to $\mathcal{S}_{\leq\PP} \hookrightarrow \mathcal{S}^\ast$, we get a factorization of maps of perfectoid spaces
\[S\to \mathcal{S}_{\leq\PP} \hookrightarrow \mathcal{S}^\ast.\]
This in turn implies that we have a factorization of maps of v-sheaves
\[X_{\leq\PP} \hookrightarrow \mathcal{S}_{\leq\PP}\hookrightarrow \mathcal{S}^\ast.\]
The first map in this factorization provides an inverse to the map in Equation~\eqref{Eq:ComparisonMap} and establishes an isomorphism between the two sides as desired.
\end{proof}

\subsubsection{Newton stratification on the Igusa stack}
The advantage of pullback along $\Bun_G^b\hookrightarrow\Bun_G$ is that we have an explicit description of the fibers of the minimally compactified Igusa stack in terms of partially minimally compactified perfect Igusa varieties. This will allow us to describe the boundary strata explicitly as well. 

More precisely, let $b$ be an element in $B(G,\mu^{-1})$. We fix a completely slope divisible $p$-divisible group $\mathbb{X}/k$ in the isogeny class labeled by $b$. As in \S\ref{Sec:Igusa background}, we can define a perfect Igusa variety $\mathrm{Ig}^b$ and its partial minimal compactification $\mathrm{Ig}^{b,*}$. Let $\mc{J}_b := \widetilde{J}_{b,k}^{\diamond}$ be the v-sheaf of groups attached to the formal group scheme $\widetilde{J}_{b,k}$ over $k$ of self-quasi-isogenies of $\mathbb{X}$ (Definition \ref{defn: automorphismsfpqcsheaf}), by applying the analytification operation $(-)^{\diamond}$ of \S \ref{ss: Analytification}. We recall that, by \cite[Corollary~9.46]{zhang2023}, this is isomorphic to the group diamond of automorphisms of the $G$-bundle $\mathcal{E}_{b}$ attached to $b \in B(G,\mu^{-1})$ on the Fargues-Fontaine curve. In particular, we have an isomorphism $\Bun_{G}^{b} \simeq [\Spd k/\mathcal{J}_{b}]$, as in \cite[Proposition~III.5.3]{FSGeomLLC}.

We set $\mathrm{Ig}^{b,*,\diamond}$ be the v-sheaf over $\Spd k$ attached to $\Ig^{b,\ast}$. As explained in \S \ref{sec: JbActionV}, the action of $\widetilde{J}_{b,\ol{\bb{F}}_{p}}$ on $\Ig^{b,\ast}$ induces an action of $\mc{J}_b$ on $\Ig^{b,\ast,\diamond}$ and further extends to the canonical compactification $\ol{\Ig^{b,\ast, \diamond}} = \Ig^{b,\ast,\dagger}$, by partial properness of $\mathcal{J}_{b}$.

By~\cite[Proposition 9.47]{zhang2023}, we have an isomorphism
\begin{equation}{\label{eqn: strataofIgusaStacks}}
\Igs^{\ast,b} \simeq [\Ig^{b,\ast, \dagger}/\mc{J}_b] 
\end{equation}
such that the morphism to $\Bun_{G}^{b} \simeq [\Spd k/\mathcal{J}_{b}]$ identifies with the natural projection map. Since the $\mc{J}_b$-action preserves the Baily--Borel boundary strata of $\Ig^{b,\ast, \dagger}$ by Proposition~\ref{prop: JbPreserveBoundary} and Corollary~\ref{cor: JbPreserveBoundary}, the $v$-stack $[\Ig^{b,\ast,\dagger}/\mc{J}_b]$ has closed (resp. locally closed) sub-v-stacks $[\Ig_{\leq \PP}^{b,\ast,\dagger}/\mc{J}_b]$ (resp. $[(\Ig^{b,\ast,\dagger})_{\PP}/\mc{J}_b]$).

\subsubsection{A stratified cartesian diagram.}
We now prove the version of Theorem~\ref{Thm:StratifiedCartesian} on Newton strata. Under the identification $\Igs^{\ast,b}\simeq [\Ig^{b,\ast,\dagger}/\mc{J}_b]$ of (\ref{eqn: strataofIgusaStacks}), we obtain a map 
\[\tilde{h}^b_{\leq \PP}:\mathcal{S}^{b}_{\leq\PP}\hookrightarrow \mc{S}^{\ast,b}\to [\mathrm{Ig}^{b,*,\dagger}/\mc{J}_b].\] 

\begin{proposition}\label{prop: FiberwiseStratifiedCartesian}
The map $\tilde{h}^b_{\leq \PP}$ factors through the sub-v-stack 
\[\left[\mathrm{Ig}^{b,\dagger}_{\leq\PP}/\mc{J}_b\right]\hookrightarrow [\mathrm{Ig}^{b,*,\dagger}/\mc{J}_b],\] and we have a cartesian diagram
\begin{equation}\label{Eq:FiberwiseStratifiedCartesian}
\begin{tikzcd}
   \mathcal{S}^{b}_{\leq\PP}  \ar[r]\ar[d]& \mathcal{F}\ell^b \ar[d] \\ 
  \left[\mathrm{Ig}^{b,\dagger}_{\leq\PP}/\mc{J}_b\right] \ar[r]& \mathrm{Bun}_G^b.
\end{tikzcd}
\end{equation}
The same result holds true if we replace $\leq \PP$ by $<\PP$ above.
\end{proposition}

\begin{proof} We spell out the proof for the $\leq \PP$ condition. The proof for the $<\PP$ condition is the same. We have an isomorphism of perfectoid spaces 
\begin{equation}\label{Eq:FiberProductIso}
\mathcal{S}^{*,b} \simeq [\mathrm{Ig}^{b,*,\dagger}/\mc{J}_b]\times_{\mathrm{Bun}_G^b}\mc{F}\ell^b,
\end{equation}
obtained by taking the pullback of the isomorphism (\ref{eqn: strataofIgusaStacks}) along $h: \mathcal{F}\ell \ra \Bun_{G}$. Since the v-sheaves we wish to compare, namely $\mathcal{S}^{b}_{\leq\PP}$ and $[\mathrm{Ig}^{b,\dagger}_{\leq\PP}/\mc{J}_b]\times_{\mathrm{Bun}_G^b}\mathcal{F}\ell^b$, are closed sub-v-sheaves on both sides. By \cite[Proposition 11.20]{Ecod}, it suffices to show their underlying topological spaces match under the isomorphism~\eqref{Eq:FiberProductIso}. 

For this, we let $x:\operatorname{Spa}C\to \mathcal{F}\ell^b$ be a rank one geometric point. Taking fibers over $x$, Equation~\eqref{Eq:FiberProductIso} becomes the identification
\begin{equation}\label{eq: HTfiber}
\mathcal{F}^\ast:=(\pi^\ast_{\HT})^{-1}(x)\simeq \Ig^{b,\ast,\dagger}_C    
\end{equation}
from \cite[Theorem 4.5.1]{CS2}, \cite[Theorem 4.3.9]{santos}\footnote{Since our $C$ is of characteristic $p$, the isomorphism is rather the tilt of the isomorphism in \textit{loc. cit.}.}. It thus suffices to show that the isomorphism~\eqref{eq: HTfiber} respects the boundary strata defined by the $\leq \PP$ conditions, i.e. \begin{equation}\label{Eq: CompareFiberStratum}
\Ig^{b,\dagger}_{\leq \PP} \times_{\Spd k}\Spa C\simeq \mc{F}^\ast\times_{\mc{S}^\ast} \mc{S}_{\leq \PP}=:\mc{F}^\ast_{\leq \PP}.
\end{equation}
This follows from the analogous statement for toroidal compactifications. 

Namely, fix a good compatible choice of cone decompositions $\Sigma$ and form the toroidal compactification $\mc{S}^\tor:=\mc{S}^\tor_{K^p,\Sigma}$. We claim that under \cite[Theorem 4.4.1]{CS2}, cf. \cite[Theorem 4.3.10]{santos}\footnote{Where in the last sentence in the statement, the word ``affinoid'' should be deleted.} we have
\begin{equation}\label{eq: CompareFiberStratumTor}
    \Ig^{b,\tor,\dagger}_{\leq \PP} \times_{\Spd k}\Spa C\simeq \mc{F}^\tor\times_{\mc{S}^\tor} \mc{S}^\tor_{\leq \PP}=:\mc{F}^\tor_{\leq \PP},
\end{equation}
where $\Ig^{b,\tor}_{\leq \PP}$ and $\mc{S}^\tor_{\leq \PP}$ are the preimages of $\Ig^{b}_{\leq \PP}$ and $\mc{S}_{\leq \PP}$ along the natural projection maps from toroidal compactifications to minimal compactifications, and $\mc{F}^\tor$ is the fiber of the Hodge-Tate period map on $\mc{S}^\tor$ over $x$. Since both sides of Equation~\eqref{Eq: CompareFiberStratum} are affinoid and partially proper, they are in particular controlled by the global sections of their structure sheaves. Hence, by \cite[Lemma 4.5.2]{CS2} (cf. Proposition \ref{prop: globalsections})  and Lemma~\ref{lemma: StrataSteinFactorization} below, the claim will imply that Equation~(\ref{Eq: CompareFiberStratum}) is an isomorphism.

To prove the claim, using the partial properness of both sides, we just need to check the map~\eqref{eq: CompareFiberStratumTor} induces natural bijections on rank one geometric points. Let $\Ig^{b,\tor}_{\leq \PP, \CO_{C}}$ be the base-change of $\Ig^{b,\tor}_{\leq \PP}$ to $\Spec \CO_C$, $\hat{\Ig}^{b,\tor}_{\leq \PP, \CO_{C}}$ be its $\varpi$-adic formal completion, and let $C'/C$ be a complete algebraically closed field. By Lemma~\ref{lemma: PerfectDeformation}, a $\Spa C'$-point of the left hand side of~\eqref{eq: CompareFiberStratumTor}, is the same as a $\Spa C'$-point of the adic generic fiber of the formal scheme $\hat{\Ig}^{b,\tor}_{\leq \PP, \CO_{C}}$. Now by examining the proof of \cite[Theorem 4.4.1]{CS2}, cf. \cite[Theorem 4.3.10]{santos}, we see from the explicit data the two sides parameterize that the $\leq\PP$-condition is preserved.
\end{proof}

\begin{lemma}[{Cf. \cite[Lemma 4.5.2]{CS2}}]\label{lemma: StrataSteinFactorization}
Let $f: Y\to X$ be a map of perfectoid spaces, such that the natural map $\CO_X^{+a}/\varpi^n\to f_\ast( \CO^{+a}_Y/\varpi^n)$ is an isomorphism for each $n$, where ``$a$'' means almost and $\varpi$ is a pseudo-uniformizer on $X$. Let $i: X_0\hookrightarrow X$ be a Zariski closed immersion, and $f_0$ be the map $Y_0:=X_0\times_X Y\to X_0$. Then $f_{0,\ast}\CO^{+a}_{Y_0}\simeq \CO^{+a}_{X_0}$. In particular, $f_0$ induces isomorphisms on global sections of the structure sheaves. 
\end{lemma}
\begin{proof} 
By $\varpi$-adic completeness of $\CO^{+a}_{X_0}$ and the fact that $f_{0,\ast}\CO^{+a}_{Y_0}\simeq\varprojlim_{n} f_{0,\ast}(\CO^{+a}_{Y_0}/\varpi^n)$, it suffices to show the natural map  $\CO^{+a}_{X_0}/\varpi^n\to f_{0,\ast}(\CO^{+a}_{Y_0}/\varpi^n)$ is an isomorphism, for all $n$. For this, by \cite[Proposition 14.8]{Ecod}, we may view the sheaves $\CO_X^{+a}/\varpi^n$ and $ \CO^{+a}_Y/\varpi^n$ as quasi-pro-\'etale sheaves. As such, they satisfy base-change for pullback along Zariski closed immersions. In particular, we get
\[ \CO^{+a}_{X_0}/\varpi^n\simeq i^\ast(\CO_X^{+a}/\varpi^n)\xrightarrow{\sim}i^\ast f_\ast( \CO^{+a}_Y/\varpi^n)\simeq f_{0,\ast}(\CO^{+a}_{Y_0}/\varpi^n)\]
as desired. The final statement follows by taking inverse limit of both sides, inverting $\varpi$ and then taking the global sections.
\end{proof}

\begin{corollary}\label{cor: IdentificationOfStrata} Under the identification $\mathrm{Igs}^{\ast,b}\simeq [\mathrm{Ig}^{b,\ast,\dagger}/\mc{J}_b]$ of (\ref{eqn: strataofIgusaStacks}), we have an isomorphism of v-stacks
\[
\mathrm{Igs}^{b}_{\leq\PP}\simeq [\mathrm{Ig}^{b,\dagger}_{\leq\PP}/\mc{J}_b].
\]
\end{corollary}

\begin{proof} 
Proposition~\ref{prop: FiberwiseStratifiedCartesian} shows that $[\mathrm{Ig}^{b,\dagger}_{\leq\PP}/\mc{J}_b]$ is precisely the image of $\mc{S}_{\leq \PP}^b=\mc{S}_{\leq \PP}\times_{\mc{S}^\ast} \mc{S}^{\ast,b}$. The result therefore follows from the definition of $\mathrm{Igs}^{b}_{\leq\PP}$.
\end{proof}

We deduce the result for individual strata as well.
\begin{proposition}\label{prop: FiberwiseStratifiedCartesianEq}
Under the identification $\Igs^{\ast,b}\simeq [\Ig^{b,\ast,\dagger}/\mc{J}_b]$ of (\ref{eqn: strataofIgusaStacks}), we have an isomorphism of $v$-stacks
\[
\mathrm{Igs}^{b}_{\PP}\simeq [(\mathrm{Ig}^{b,\dagger})_{\PP}/\mc{J}_b].
\]
Moreover, Diagram~\eqref{Eq:FiberwiseStratifiedCartesian} restricts to a cartesian diagram
\begin{equation}
\begin{tikzcd}
   \mathcal{S}^{b}_{\PP}  \ar[r]\ar[d]& \mathcal{F}\ell^b \ar[d] \\ 
  \mathrm{Igs}^{b}_{\PP} \ar[r]& \mathrm{Bun}_G^b.
\end{tikzcd}
\end{equation}
\end{proposition}
\begin{proof}
Using Proposition~\ref{prop: FiberwiseStratifiedCartesian}, it suffices to show that the isomorphism
\begin{equation}\label{Eq: CompareNewtonStratumOfBoundaries}
\mathcal{S}^{b}_{\leq\PP}  \simeq {[\mathrm{Ig}^{b,\dagger}_{\leq\PP}/\mc{J}_b]}\times_{\mathrm{Bun}_G^b}\mc{F}\ell^b
\end{equation} 
restricts to an isomorphism of the open sub-v-sheaves
\[ \mathcal{S}^{b}_{\PP}  \simeq {[(\mathrm{Ig}^{b,\dagger})_{\PP}/\mc{J}_b]}\times_{\mathrm{Bun}_G^b}\mc{F}\ell^b.\]
This can again be checked on the level of topological spaces, by \cite[Proposition 11.15]{Ecod}. Using partial properness of both sides, we reduce to check that the above map is a bijection on rank one geometric points. 
This follows from the $<\PP$ version of Proposition~\ref{prop: FiberwiseStratifiedCartesian} and the fact that set theoretically we have 
\[|\mathcal{S}^{b}_{\leq\PP}|=|\mathcal{S}^{b}_{\PP}|\coprod |\mathcal{S}^{b}_{<\PP}|.\]
\end{proof}

Now we can finish the proof of Theorem~\ref{Thm:StratifiedCartesian}.
\begin{proof}[Proof of Theorem~\ref{Thm:StratifiedCartesian}]
    The result for the $\leq \PP$-condition follows by combining Lemma~\ref{Lem:FiberwiseCartesian}, Proposition~\ref{prop: FiberwiseStratifiedCartesian}, and Corollary~\ref{cor: IdentificationOfStrata}. The one for the $\PP$-condition follows by combining Lemma~\ref{Lem:FiberwiseCartesian} and Proposition~\ref{prop: FiberwiseStratifiedCartesianEq}.
\end{proof}

We end this subsection by stating two consequences of Theorem~\ref{Thm:StratifiedCartesian}.  

The first is a construction of the Igusa stack beyond the good reduction locus. The original construction of Igusa stacks in \cite[Theorem 1.3]{zhang2023} only uniformizes the good reduction locus of the Shimura varieties and Kim \cite{Kim} extended this result to the whole Shimura variety. In fact, applying Theorem~\ref{Thm:StratifiedCartesian} to $\PP=[\mathsf{G}]$, we see that taking the image of $\mc{S}_{K^p}$ in $\Igs^\ast$ is an alternative way of constructing the Igusa stack beyond the good reduction locus, logically independent of \cite{Kim}. We summarize this as the corollary below.
\begin{corollary}[Igusa stack beyond the good reduction locus]\label{cor: BeyondGoodRed}
     Under Assumption~\ref{assumption:codimension}, the $v$-stack $\Igs$ and the cartesian diagram in \cite[Conjecture 1.1(1)]{zhang2023} exist, extending the ones on the good reduction locus from \cite[Theorem 1.3]{zhang2023}. For all $b \in B(G,\mu^{-1})$, we have an isomorphism $\Igs^{b} \simeq [(\Ig^{b,\dagger})_{[\mathsf{G}]}/\mathcal{J}_{b}]$, where we note that $(\Ig^{b,\dagger})_{[\mathsf{G}]}$ is precisely the perfectoid space of \cite[Corollary~3.4]{HamannLeeTorsion} when base-changed along $\Spa(C) \ra \Spd k$ for $C$ an algebraically closed perfectoid field, in light of Lemmas \ref{lem: DaggerisCC} and \ref{lemma: PerfectDeformation}.
\end{corollary}

The second consequence is that we have a version of Mantovan's product formula for the Baily--Borel boundaries. Namely, we can formulate Proposition~\ref{prop: FiberwiseStratifiedCartesian} and Proposition~\ref{prop: FiberwiseStratifiedCartesianEq} slightly differently. We take $k$ to be $\ol{\mathbb{F}}_p$. Let $\operatorname{Sht}(G,b,\mu)_{\infty}$ be the fiber of the map $h: \mc{F}\ell\to \Bun_G$ over a $\Spd k$-point corresponding to $b$. This is the local Shimura variety $\Sht(G,b,\mu)_{\infty}$ defined in \cite[Section 23]{SW}, which has a structure map to 
\[\Spd \Breve{E}=\Spd W_{\CO_E}(k)[\tfrac{1}{p}]= \Spd E\times_{\Spd \mathbb{F}_p}\Spd k.\] It is equipped with a natural $\underline{G(\bb{Q}_{p})}$-action (see \S \ref{ss: Variation'sonMantovan'sFiltration}). Then we have the following corollary.  
\begin{corollary}[Mantovan's formula for boundary strata]\label{cor: GeoMantovanBoundary}
    There is a $G(\qp)$-equivariant isomorphism of diamonds over $\Spd \Breve{E}$
    \[\mc{S}^b_{\leq \PP} \simeq \Ig^{b,\dagger}_{\leq \PP} \times^{\mathcal{J}_{b}} \Sht(G,b,\mu)_{\infty}.\]
    Here the $G(\qp)$-action on the right-hand side is through the $\Sht(G,b,\mu)_{\infty}$ factor. The contracted product on the right-hand side is given by the quotient of the product $\Sht(G,b,\mu)_{\infty} \times_{\Spd k} \Ig^{b,\dagger}_{\leq \PP}$ by the diagonal $\mc{J}_b$-action. 

    The analogous statement for the $\PP$-condition is also true. 
\end{corollary}
\begin{remark}
    Note that the base-changes of both $\Ig^{b,\dagger}_{\leq \PP}$ and $\mc{J}_b$ to $\Spd \Breve{E}$ are representable by pre-perfectoid spaces which describe the rigid generic fibers of the Witt vector lifts of the perfect Igusa varieties, and at finite levels the shtuka spaces are also representable by smooth rigid spaces, see \cite[Lecture 24.1, 24.2]{SW}. Hence, this result indeed parallels the classical $p$-adic uniformization.
\end{remark}

\subsection{Dimension formula}
As immediate corollaries to Theorem~\ref{Thm:StratifiedCartesian} and Proposition~\ref{prop: FiberwiseStratifiedCartesianEq}, we obtain dimension formulae for the boundary strata of both the Igusa stack and individual Igusa varieties.
\begin{corollary}{\label{cor: IgusaStackStrataSmooth}}
    Let $d_{[\mathsf{P}]}$ be the dimension of the stratum $\mathcal{S}_{K_p, \PP}$. The stratum $\Igs_{[\mathsf{P}]}$ is an Artin v-stack (in the sense of \cite[Definition~IV.1.1]{FSGeomLLC}) and the structure map $\Igs_{[\mathsf{P}]}\to \Spd k$ is $\ell$-cohomologically smooth of $\ell$-dimension $d_{[\mathsf{P}],\Igs}:= d_{[\mathsf{P}]}-d$.
\end{corollary}
\begin{proof}
    The $v$-stack $\Igs^{\ast}$ is Artin by \cite[Proposition~IV.1.8 (iii)]{FSGeomLLC}, since it admits a suitably representable map, namely $\ol{\pi}_{\min}: \Igs^{\ast} \ra \Bun_{G}$, to the Artin $v$-stack $\Bun_{G}$ (\cite[Theorem IV.1.19.]{FSGeomLLC}). This implies that the stratum $\Igs_{[\mathsf{P}]}$ is also an Artin v-stack, since by Theorem~\ref{Thm:StratifiedCartesian}, it is (locally) closed in the Artin v-stack $\Igs^\ast$. By \cite[Example IV.1.9(iii)]{FSGeomLLC}, the result follows (Similarly, $\Igs_{\leq [\mathsf{P}]}$ is an Artin v-stack, though it is not necessarily cohomologically smooth). 
    
    For the $\ell$-cohomological smoothness, again by Theorem \ref{Thm:StratifiedCartesian} for all open compact subgroups $K_{p}$, we have a cartesian diagram of the form 
    \[
    \begin{tikzcd}
    \left[\mathcal{S}_{K^p,[\mathsf{P}]}/\ul{K_{p}}\right] \arrow[r] \arrow[d] & \left[\Fl/\ul{K_{p}}\right] \ar[d]\\
    \Igs_{\PP} \ar[r] & \Bun_G.
    \end{tikzcd}
    \]
In particular, since the right vertical arrow is separated, fdcs, $\ell$-cohomologically smooth of $\ell$-dimension $d$ (Lemma \ref{lemma: smoothnessofvertmaps}), the same holds true for the left vertical arrow (which is additionally a surjection of v-stacks). Moreover, by \cite[Proposition 23.15]{Ecod}, $\ell$-cohomological smoothness can be checked after passing to a v-cover, so we choose an algebraically closed perfectoid field $C/k$ and base-change everything to it. Now we conclude from Proposition \ref{lemma: stackystratasmooth}, since from the definition \cite[Definition IV.1.11]{FSGeomLLC} it is clear that a map of Artin v-stacks $f: Y\to X$ is cohomologically smooth if there is a separated, fdcs, cohomologically smooth surjection $g: V\to Y$ from an Artin v-stack, such that $f\circ g$ is cohomologically smooth. This also implies the result on $\ell$-dimension, cf. \cite[Definition~IV.1.17]{FSGeomLLC} and \cite[Lemma~3.8]{HamJacCrit}.
\end{proof}

Fix an element $b\in B(G,\mu^{-1})$ and a principal prime-to-$p$ level $K^p=K(N)$, $p\nmid N$. We also deduce an interesting numerical result about the dimension of strata in the Igusa variety $\Ig^b=\Ig^b_{K^p}$ attached to $b$. 

\begin{lemma}\label{lemma: IgusaVarietyStackyQuotient}
    For a principal level $K_b(p^m)\subset J_b(\qp)$, we have a decomposition of the quotient stack 
    \[[\mathrm{Ig}^{b}_{\PP}/\underline{K_b(p^m)}]= \coprod_{[\widetilde{Z}]}[(\mathrm{Ig}^{b}_{[\widetilde{Z}]})^{\perf}/\underline{H_b}],\] where $[\widetilde{Z}]$ runs through Igusa cusp labels at level $K_b(p^m)K^p$ of type $\PP$. If $[\widetilde{Z}]$ lies over a cusp label $[Z]$ at level $K^p$ then $\Ig^{b}_{[\widetilde{Z}]}\simeq \mathrm{Ig}^{b}_{[Z], m}$ is a smaller dimensional Igusa variety attached to $[Z]$ at principal $p^m$-level, cf. \cite[Theorem 3.3.15]{CS2}. Moreover, $H_b$ is some closed subgroup of $K_b(p^m)$ that acts trivially on $\mathrm{Ig}^{b}_{[\widetilde{Z}]}$.
\end{lemma}
\begin{proof}
Fix an Igusa cusp label $\widetilde{Z}$ of type $\PP$ with underlying cusp label $Z$, see Definition~\ref{defn: IgusaCuspLabel}. We denote by $[Z]$ the $K^p$-orbit of $Z$, or equivalently, a cusp label at level $K^p$ that $Z$ projects to. Recall from \cite[Theorem 2.5.8]{CS2} that there is a corresponding smaller dimensional Shimura variety $\mathscr{S}_{[Z]}$ appearing in the boundary of $\mathscr{S}^\ast$, where following Section 4, we use $\mathscr{S}^\ast$ to denote the special fiber of the minimally compactified Shimura variety with hyperspecial level at $p$ and prime-to-$p$ level $K^p$. Denote by $\Ig^b_{[Z]}$ the corresponding smaller dimensional Igusa variety (with infinite level at $p$) above the central leaf $\mathscr{C}^b_{[Z]}\subset \mathscr{S}_{[Z]}$ as in \cite[Theorem 3.3.15]{CS2}. It is attached to an element $b_{\PP}\in B(\mathsf{M}_{h,\qp})$, if $\mathsf{M}_h$ is the reductive group in the Shimura datum for $\mathscr{S}_{[Z]}$; and $b$ is the image of $b_\PP$ under the natural map $B(\mathsf{M}_{h,\qp})\to B(G)$.\footnote{We choose the notation $\mathsf{M}_h$ since $\mathsf{P}$ has a decomposition of the shape $\mathsf{P}=(\mathsf{M}_\ell\times \mathsf{M}_h)\ltimes \mathsf{U}$, where $\mathsf{U}$ is the unipotent radical, $\mathsf{M}_\ell$ is the ``linear'' part of the Levi and $\mathsf{M}_h$ is the ``hermitian'' part, cf. the paragraph after \cite[Definition 3.3.8]{LanStroh}. The part $\mathsf{M}_h$ has a symmetric space that is of hermitian type. They form the Shimura datum of a smaller dimensional Shimura variety appearing in the boundary of $\mathscr{S}^\ast$.}

Let $P_b(\qp)\subset J_b(\qp)$ be the parabolic subgroup stabilizing the filtration on $\mathbb{X}$ given by $Z_b$. It has a decomposition
\[P_b(\qp)=(M_{b,\ell}(\qp)\times M_{b,h}(\qp))\ltimes U_b(\qp),\]
where $M_{b,\ell}$ and $M_{b,h}$ denote the linear and hermitian part of the Levi subgroup $M_b$ of $P_b$, and $U_b$ is a unipotent subgroup of $P_b$.
As in \cite[Section 6.2.1]{CS2}, there is a natural $J_b(\qp)$-equivariant map $\mathrm{Ig}^{b}_{\PP}\to \underline{J_b(\qp)/P_b(\qp)}$. The fiber of this map at the identity, which we denote by $\mathrm{Ig}^{b}_{\PP,0}$, is an infinite  union of Igusa varieties $\Ig^{b}_{[Z]}$, parameterized by a profinite set $\mc{Z}_{\PP, K^p}$ of equivalence classes of Igusa cusp labels of type $\PP$, whose underlying usual cusp label is at level $K^p$. In other words, we have a map\footnote{On the source the disjoint union symbol is an abuse of notation, since there is a profinite topology on the cusp labels.}
\[\mathrm{Ig}^{b}_{\PP,0}=\coprod_{[\widetilde{Z}]\in \mc{Z}_{\PP,K^p}} \Ig^{b}_{[Z]}\to \ul{\mc{Z}_{\PP, K^p}}.\]
This map is $\ul{P_b(\qp)}$-equivariant, where the $P_b(\qp)$-action on $\mc{Z}_{\PP,K^p}$ factors through its projection to $M_{b,\ell}(\qp)\cong \operatorname{Aut}_{\CO_B}(Y\otimes \qp/\mathbb{Z}_p)$, and $Y$ is as in the discussion preceding Equation~\eqref{eq:IgusaCuspLabelDecomposition}. The latter group acts on the Igusa cusp labels by composing $\psi$ with the given automorphism of $Y\otimes \qp/\mathbb{Z}_p$. This action is transitive, and the stabilizer of each equivalence class of a cusp label is $\Gamma_{[Z]}:=\operatorname{Aut}_{B}(Y_\Q)\cap K^p$.

Hence, there remains an action of a group $H_{P_b}$, which sits in an exact sequence
\[1\to M_{b,h}(\qp)\ltimes U_b(\qp)\to H_{P_b}\to \Gamma_{[Z]}\to 1,\] 
on the Igusa variety $\Ig^{b}_{[Z]}$, where the $U_b(\qp)$-action is trivial and the $M_{b,h}(\qp)=J_{b_\PP}(\qp)$-action is the usual one. 

Using this, we can describe the quotient $[\Ig^{b}_{\PP}/\underline{K_b(p^m)}]$ explicitly. Namely, consider the intersection $K_b(p^m)\cap P_b(\qp)$. It has a corresponding decomposition into a linear, hermitian and unipotent part. The linear part $K_b(p^m)\cap M_{b,\ell}(\qp)$ acts on the cusp labels, with stabilizers $\Gamma_{[\widetilde{Z}]}:=K_b(p^m)K^p\cap \operatorname{Aut}_B(Y_{\Q})$. The hermitian part $K_{b_\PP}:=K_b(p^m)\cap M_{b,h}(\qp)$ determines the level of the smaller Igusa variety, giving rise to an equality 
\[[\Ig^{b}_{[Z]}/\underline{K_{b_\PP}}]=(\Ig^{b}_{[Z],m})^{\perf},\]
while the unipotent part $K_{b,U}:=K_b(p^m)\cap U_b(\qp)$ fixes $\Ig^{b}_{[Z],m}$. Hence, define $H_b$ to be the closure of $\Gamma_{[\widetilde{Z}]}\ltimes (K_{b,U})$ in $K_b(p^m)$, and we obtain the desired description of $[\mathrm{Ig}^{b}_{\PP}/\underline{K_b(p^m)}]$.
\end{proof}
We now let $K^{p}$ be a neat level away from $p$. We set $d_{b} := \langle 2\rho_{G},\nu_{b} \rangle$, where $\nu_{b}$ is the slope homomorphism of an element $b \in B(G)$ and $\rho_{G}$ is the half sum of all the positive roots of $G$. Let $d_{b,[\mathsf{P}]}$ denote the dimension of the stratum $\Ig^b_{\PP}$ in the perfect Igusa variety $\Ig^{b,\ast}$. We use the above result to deduce the following lemma.  
\begin{lemma}\label{lemma: StrataFiberwiseSmooth}
     The quotient stacks $[\mathrm{Ig}^{b,\Diamond}_{\PP}/\underline{J_{b}(\bb{Q}_{p})}]$ and $[(\mathrm{Ig}^{b,\diamond})_{\PP}/\underline{J_{b}(\bb{Q}_{p})}]$ are $\ell$-cohomologically smooth over $\Spd k$ of $\ell$-dimension equal to $d_{b,[\mathsf{P}]}$. Similarly, $[(\Ig^{b,\diamond})_{\PP}/\mathcal{J}_{b}]$ is $\ell$-cohomologically smooth of dimension $d_{b,\PP} - d_{b}$.
\end{lemma}
\begin{proof}
    Let $K_b(p^m)\subset J_b(\qp)$ be a pro-$p$ principal level subgroup. Then the natural map
    \[[\mathrm{Ig}^{b,\Diamond}_{\PP}/\underline{K_b(p^m)}]\to [\mathrm{Ig}^{b,\Diamond}_{\PP}/\ul{J_{b}(\bb{Q}_{p})}] \]
    is \'etale; in particular, it is an $\ell$-cohomologically smooth cover of $\ell$-dimension $0$. Hence, it suffices to show $[\mathrm{Ig}^{b,\Diamond}_{\PP}/\underline{K_b(p^m)}]$ is $\ell$-cohomologically smooth of $\ell$-dimension $d_{b,[\mathsf{P}]}$. Using Lemma~\ref{lemma: IgusaVarietyStackyQuotient} (and the fact that the Igusa varieties in question for general level away from $p$ differ from the Igusa varieties at principal level away from $p$ by a finite \'etale morphism), it suffices to show this for each $[\mathrm{Ig}^{b,\Diamond}_{[Z], m}/\underline{H_b}]$ (Note that the functor $(-)^{\Diamond}$ factors through perfection by construction). But this follows from the same argument as in Proposition~\ref{lemma: stackystratasmooth}. Namely, since $H_b\subset K_b(p^m)$ is pro-$p$ and the natural projection
    \[[\mathrm{Ig}^{b,\Diamond}_{[Z], m}/\underline{H_b}]\to \mathrm{Ig}^{b,\Diamond}_{[Z], m}\] is the base-change of $[\Spd k/\underline{H_b}]\to \Spd k$, the map is $\ell$-cohomologically smooth of $\ell$-dimension $0$. The question thus reduces to the $\ell$-cohomological smoothness of $\mathrm{Ig}^{b,\Diamond}_{[Z], m}$. We now conclude, using that we have an equality $\mathrm{Ig}^{b,\Diamond}_{[Z], m} = \mathrm{Ig}^{b,\Diamond}_{[Z],m,\Mant}$, by Remark \ref{rem: comparison with Mantovan} and the fact that $(-)^{\Diamond}$ factors through perfection by definition. As in Remark \ref{rem: comparison with Mantovan}, the subscript $\Mant$ is meant to indicate that these are the finite level Igusa varieties considered in \cite{Mant,CS17,CS2}.  The claim now follows from Corollary \ref{prop: Compatabilitieswithupper!s}. Indeed, $\mathrm{Ig}^{b}_{[Z], m, \Mant}$ is a finite \'etale cover of the central leaf $\mathscr{C}^b_{[Z]}\subset \mathscr{S}_{[Z]}$ and this is smooth by \cite[Proposition~1]{mantovan-PEL}. 
    
    It remains to show the claim for $[(\mathrm{Ig}^{b,\diamond})_{\PP}/\underline{J_{b}(\bb{Q}_{p})}]$. This follows from Lemma \ref{lemma: differentstrataofIgs} and the fact that the map $[\Ig^{b,*,\diamond}/\ul{J_{b}(\bb{Q}_{p})}] \ra [\Ig^{b,*,\Diamond}/\ul{J_{b}(\bb{Q}_{p})}]$ is an open immersion, see Lemma \ref{lemma: IgCCqcqs}. The last claim follows from the first, and the fact that $[(\mathrm{Ig}^{b,\diamond})_{\PP}/\ul{J_{b}(\bb{Q}_{p})}] \ra [(\mathrm{Ig}^{b,\diamond})_{\PP}/\mc{J}_b]$ is $\ell$-cohomologically smooth of dimension $d_{b}$.
\end{proof}
We now have the following.
\begin{proposition}\label{prop: dimensionStrataIg}
 Assuming that the stratum $\mathrm{Ig}^b_{[\mathsf{P}]}$ is non-empty, we have 
 \[
 d_{b,[\mathsf{P}]} - d_b= d_{[\mathsf{P}]} -d.
 \]
\end{proposition}
\begin{proof}
    For each compact open subgroup $K_p\subset G(\qp)$, the map 
    \[\left[\mathcal{F}\ell^b_{G,\mu^{-1}}/\underline{K_p}\right]\to \mathrm{Bun}_G^b\]
    is $\ell$-cohomologically smooth of $\ell$-dimension $d$, by Lemma \ref{lemma: smoothnessofvertmaps}. It then follows from Lemma~\ref{lemma: StrataFiberwiseSmooth} and Proposition~\ref{prop: FiberwiseStratifiedCartesianEq} that 
    \[\left[(\mathrm{Ig}^{b,\diamond})_{\PP}/\mc{J}_b\right]\times_{\mathrm{Bun}_G^b}\left[\mathcal{F}\ell^b_{G,\mu^{-1}}/\underline{K_p}\right]\subset \left[\mathcal{S}^{b}_{K^p,\PP}/\underline{K_p}\right]\]
    is an $\ell$-cohomologically open sub-$v$-sheaf of $\ell$-dimension $d_{b,[\mathsf{P}]}-d_b+d$. By Proposition~\ref{prop: openstrata} (see also~\eqref{eqn: PullingBackStrata}), the above open sub-$v$-sheaf identifies with $[\mathcal{S}_{K^p,\PP}^{b,\cap}/\underline{K}_p]$, where we set $\mathcal{S}_{K^p,\PP}^{b,\cap}:=\mathcal{S}_{K^p}^{*,b,\cap}\times_{\mathcal{S}^*_{K^p}}\mathcal{S}_{K^p,\PP}$. By Lemma~\ref{lem: CapStrataOpen}, $[\mathcal{S}_{K^p,\PP}^{b,\cap}/\underline{K}_p]$ is open in $\mathcal{S}_{K^p,\PP} / \underline{K}_p$.
     However, the $v$-stack $\left[\mathcal{S}_{K^p, \PP}/\underline{K_p}\right]$ is $\ell$-cohomologically smooth of dimension $d_{\PP}$ by Proposition~\ref{lemma: stackystratasmooth}. Combining these observations, this implies that $d_{b,[\mathsf{P}]} - d_b + d = d_{[\mathsf{P}]}$, as desired.
\end{proof}

\begin{remark} Proposition~\ref{prop: dimensionStrataIg} can also be proven in a more elementary way, using the dimension formulae  from~\cite[\S 4.2-4.3]{CS17} and the explicit description of partial compactifications of Igusa varieties. This result can be rephrased as saying that the codimension $d_{b}-d_{b,[\mathsf{P}]}$ of a boundary stratum in the partial minimal compactification of an Igusa variety is independent of the choice of Newton stratum that the Igusa variety lies over. This therefore expresses some sort of transversality between the Newton and boundary stratifications. 

We have chosen to prove the transversality result this way, as this argument will more naturally complement the approach we take to compute the stalks of our relative intersection cohomology sheaf in the next section. This fact underlies the proof given in \S \ref{sec: perversity} that the relative intersection cohomology of the Igusa stack is a perverse sheaf on $\mathrm{Bun}_G$. It also explains that \cite[Proposition 9.23]{zhang2023} is a special case of a more general phenomenon.
\end{remark}

\subsection{Local Systems} \label{sec: EtaleLocalSys} We let $F/\bb{Q}_{\ell}$ be a finite algebraic extension with ring of integers $\mathcal{O}_{F}$. Suppose $\xi$ is a finite dimensional algebraic representation of $\mathsf{G}_{\ol{\bb{Q}}_{\ell}}$ defined over $F$. We suppose that $\Lambda$ is a coefficient system as in Setup~\ref{assumption: coefficientsystemsingeneral}, that is also an $\mathcal{O}_{F}$-algebra. In this situation, we say $\xi$ is defined over $\Lambda$. Using the construction from \S\ref{sec: EtaleSheavesRep}, we construct some $\Lambda$-local systems on the open Igusa stack $\Igs_{K^p}$ (over $k=\ol{k}$ as before), to which we will attach intersection complexes in the next subsection. 

In particular, $\xi$ defines the action of $\mathsf{G}(F)$ on a finite-dimensional $F$-vector space $V_{\xi}$. Suppose $\Lambda$ is as in Setup~\ref{assumption: coefficientsystemsingeneral} (1). Given a compact open subgroup $K_{\ell} \subset \mathsf{G}(\bb{Q}_{\ell})$, consider the natural induced map 
\[ K_\ell\subset \mathsf{G}(\bb{Q}_{\ell}) \ra \GL(V_{\xi})(F).\]
Since $K_{\ell}$ is compact, it fixes some $\mathcal{O}_{F}$-lattice $L_{\xi}\subset V_{\xi}$. We can then consider the $K_{\ell}$-representation $L_{\xi} \otimes_{\mathcal{O}_{F}} \Lambda$. This defines an object in $\mathrm{D}(K_{\ell},\Lambda)$ and in turn an \'etale $\Lambda$-sheaf on the classifying stack $[\ast/\ul{K_\ell}]$, and hence on $\Igs_{K^p}$ by pullback along the natural maps
\begin{equation}{\label{eqn: TorsoronIgusaStack}}
\Igs_{K^p}\to [\ast/\ul{K^p}]\to [\ast/\ul{K_\ell}],
\end{equation}
where $\ast=\Spd k$. We denote this $\Lambda$-local system by $\mc{L}_\xi$. Similarly, if $\Lambda$ is as in Setup~\ref{assumption: coefficientsystemsingeneral} (2) then we can define an object in $\Detale([\ast/\underline{K_{\ell}}],\Lambda)$  by taking inverse limits, and by pulling back we get a lisse $\Lambda$-adic sheaf on $\Igs_{K^p}$.

Similarly, one obtains a $\Lambda$-local system $\tilde{\mathcal{L}}_\xi$ on the Shimura variety by further pullback along $\mc{S}_K\to \Igs_{K^p}$, for any $K=K^pK_p$. Applying Lemma \ref{lemma: ComparisonbetweenPinksFunctorandClassifyingStack}, we see that this agrees with the analytification of Pink's construction applied to $L_{\xi}$, hence also with the usual $\Lambda$-local system attached to $\xi$ upon inverting $\ell$, see \cite[Section~5.1]{PinkHigherDirectImages}. 

\subsection{Intersection complex}\label{sec: IConIgs}
Let $\Lambda/\mathbb{Z}_\ell$ be a coefficient system as in Setup~\ref{assumption: coefficientsystemsingeneral} and $\xi$ be an algebraic representation of $\mathsf{G}_{\ol{\mathbb{Q}}_{\ell}}$ defined over $\Lambda$ (i.e of the form described in \S~\ref{sec: EtaleLocalSys}). Let $\mc{L}:=\mc{L}_\xi$ be the attached \'etale local system on $\Igs_{K^p}$ as in \S\ref{sec: EtaleLocalSys} above. We will consider the following sheaves.
\begin{definition}{\label{defn: ICIgs*}}
We let $\mc{L}:=\mc{L}_\xi$ be the attached \'etale local system on $\Igs_{K^p}$, as in \S\ref{sec: EtaleLocalSys} above. If $\Lambda$ is as in Setup~\ref{assumption: coefficientsystemsingeneral} (1) then, by Remark \ref{rem: IgusaStackReadyforPeversity}, Corollary~\ref{cor: IgusaStackStrataSmooth}, and Definition~\ref{defn: IntersectionComplex}, we can define its intersection complex $\IC_{\Igs^\ast}(\mc{L})$, for the Baily--Borel stratification on $\Igs^\ast$ by taking intermediate extensions. If $\Lambda$ is as Setup \ref{assumption: coefficientsystemsingeneral} (2) then we define its intersection complex to be given by the inverse limit 
\[ \lim_{n \geq 1} \IC_{\Igs^{*}}(\mathcal{L}/\ell^{n}) \in \Detale(\Igs^{*},\Lambda) \]
over the mod $\ell^{n}$-reduction maps. Similarly, denote the pullback of $\mc{L}$ to $\mc{S}_K$ by $\tilde{\mc{L}}$. We can define $\IC_{[\mc{S}^\ast_{K^p}/\ul{K_p}]}(\tilde{\mc{L}}[d](\tfrac{d}{2}))$, using the Baily--Borel stratification and Proposition~\ref{lemma: stackystratasmooth} with torsion coefficients and in general by taking inverse limits. We drop the stratification from the notation for simplicity. This is justified by Corollary~\ref{cor: CompareIC} below, since, at least when passing to $\mc{S}_{K}^\ast$, the complexes agree with the analytification of the intrinsic intersection complex on the algebraic Shimura variety.
\end{definition}
\begin{remark}{\label{rem: TakingInverseLimitisPresumablyNotNecessary}}
We note that when $\Lambda$ is as in Setup \ref{assumption: coefficients} (2), we have given a somewhat unorthodox definition for $\IC_{\Igs^{\ast}}(\mathcal{L})$. Namely, instead of directly defining it to be an intermediate extension with respect to the perverse $t$-structure on $\Detale(\Igs^{*},\Lambda)$ defined by the Baily-Borel stratification, we have instead done this only with torsion coefficients and taken an inverse limit. We are unsure whether these two definitions agree. However, this inverse limit will still have the correct formal properties for our purposes in light of the fact that, in the context of schemes of finite type over a field, these two definitions agree (see Theorem \ref{thm: inverselimit}). 
\end{remark}
We now have the following.
\begin{proposition}\label{prop: PullbackofIC}
    Following the notation in \S\ref{sec: IgusaStackNotation}, for each compact open subgroup $K_p\subset G(\qp)$, there is an isomorphism in $\Detale([\mc{S}_{K^p}^\ast/\ul{K_p}],\Lambda)$
    \[ \tilde{h}_{K_p,\min}^{\ast}\IC_{\Igs^\ast}(\mc{L})[d](\tfrac{d}{2}) \simeq \IC_{[\mc{S}^\ast_{K^p}/\ul{K_p}]}(\tilde{\mc{L}}[d](\frac{d}{2})).\]
\end{proposition}
\begin{proof}
    Since $\tilde{h}_{K_p,\min}$ is $\ell$-cohomologically smooth of $\ell$-dimension $d$ and the stratifications on its source and target are compatible by Theorem \ref{Thm:StratifiedCartesian}, the functor $\tilde{h}_{K_p,\min}^\ast[d](\tfrac{d}{2})$ is perverse $t$-exact with respect to the $t$-structure defined by the Baily--Borel stratification, by smooth and proper base change, see \cite[Proposition 22.19, Proposition 23.16(i)]{Ecod}. Moreover, the functor $\tilde{h}_{K_{p},\min}^{*}$ commutes with inverse limits.
\end{proof}

Let $\mc{L}^\mathrm{alg}$ be the usual \'etale $\Lambda$-local system attached to $\xi$ on the algebraic Shimura variety $\Sh_K:=\Sh\gx_K$ and write $(-)^\Diamond$ for the analytification functor $c_{\Shstar_K}^\ast$. We have the following corollary.
\begin{corollary}\label{cor: CompareIC}
    Let $K_p \subset G(\qp)$ be a compact open subgroup whose pro-order is coprime to $\ell$ and $q_{K_p}: [\mc{S}^\ast_{K^p}/\ul{K_p}]\to \mc{S}_K^\ast$ be the natural map. For $K =K_pK^p$ neat, we have isomorphisms
    \[ q_{K_p,\ast}\tilde{h}_{K_p,\min}^{\ast}\IC_{\Igs^\ast}(\mc{L})[d](\tfrac{d}{2}) \simeq \IC_{\mc{S}^\ast_K}(\tilde{\mc{L}}[d](\frac{d}{2}))\simeq \IC_{\Shstar_K}(\mc{L}^\mathrm{alg}[d](\frac{d}{2}))^\Diamond\]
    and 
    \[ q_{K_p,\ast}\tilde{h}_{K_p,\min}^{\ast}\bb{D}_{\Igs_{K^{p}}^{*}}\IC_{\Igs^\ast}(\mc{L})[d](\tfrac{d}{2}) \simeq \bb{D}_{\mc{S}^\ast_K}\IC_{\mc{S}^\ast_K}(\tilde{\mc{L}}[d](\frac{d}{2}))\simeq (\bb{D}_{\Shstar_{K}}\IC_{\Shstar_K}(\mc{L}^\mathrm{alg}[d](\frac{d}{2})))^\Diamond\]
\end{corollary}
\begin{proof}
The second part of the claim follows from the first by applying Verdier duality, using Lemma \ref{lemma: smoothnessofvertmaps}, Proposition \ref{prop: qKpProper}, and Proposition \ref{prop: propertiesofthealgebraizationfunctor} (5). Hence it suffices to prove the maps in the first row are isomorphisms.

Using Theorem \ref{thm: inverselimit} and the fact that both $q_{K_p,\ast}$ and $\tilde{h}^\ast_{K_p,\min}$ commute with inverse limits (for the latter, use cohomological smoothness as in Lemma \ref{lemma: smoothnessofvertmaps}), we can reduce to the case that $\Lambda$ is torsion. 

Now the first isomorphism in the first row follows from Proposition~\ref{prop: PullbackofIC} and the fact that $q_{K_p}$ is $\ell$-cohomologically proper of $\ell$-dimension $0$ and respects the stratification by definition of the strata on the stacky quotient, see Proposition~\ref{prop: qKpProper}. Indeed, $q_{K_p,\ast}$ is exact, so if we present both intersection complexes using Deligne's formula (Corollary~\ref{cor: DeligneFormula}), it is easy to see that $q_{K_p,\ast}$ preserves the IC complexes. The second isomorphism follows from Pink's formula \cite[Theorem 4.2.1]{PinkHigherDirectImages} and Corollary~\ref{cor: diamonddescriptionoftheICsheaf}.

\end{proof}

\subsubsection{Universal local acyclicity}
We denote by $\mathbb{D}=\mathbb{D}_{\Igs_{K^p}^\ast}$ the Verdier duality functor on $\Igs_{K^p}^\ast$. 
\begin{proposition}\label{prop: ICStratifiedULA}
    Suppose $\Lambda$ is regular. Both $\IC_{\Igs^\ast}(\mc{L})$ and its Verdier dual $\mathbb{D}(\IC_{\Igs^\ast}(\mathcal{L}))$ are stratified ULA with respect to $\Igs_{K^p}^\ast\to\Spd k$ and the Baily--Borel stratification $\{\Igs_{K^p,\PP}\}_{\PP}$, in the sense of Definition~\ref{defn: StratifiedULA}. 
\end{proposition}
\begin{proof}
We give the proof for $\IC_{\Igs^{*}}(\mathcal{L})$ with the one for $\mathbb{D}(\IC_{\Igs^\ast}(\mathcal{L}))$ being virtually identical. We first show that $\IC_{\Igs^\ast}(\mc{L})$ is ULA over $\Spd k$. By \cite[Proposition IV.2.13]{FSGeomLLC}, ULAness is cohomologically smooth local on the source. Hence, it suffices to prove the ULAness of $\tilde{h}_{K_p,\min}^\ast\IC_{\Igs^\ast}(\mc{L})[d](\tfrac{d}{2})$ with respect to $[\mc{S}^\ast_{K^p}/K_p]\to \Spd k$. By Proposition~\ref{prop: PullbackofIC}, this complex identifies with $\IC_{[\mc{S}_{K^p}^\ast/\ul{K_p}]}(\tilde{\mc{L}}[d](\frac{d}{2}))$.
    Consider the factorization 
    \begin{equation*}
    q_{K_p}: [\mc{S}_{K^p}^{*}/\ul{K_{p}}] \xrightarrow{q^{1}_{K_{p}}} [\mc{S}^\ast_{K}/\ul{K_{p}}] \xrightarrow{q_{K_{p}}^{2}} \mc{S}^\ast_{K} 
    \end{equation*}
    from (\ref{eqn: factorizationofqKp}). Since $q^1_{K_p}$ is quasi-pro-\'etale and proper by Lemma~\ref{lemma: q1isquasiproetaleandproper}, it suffices to show the ULAness of $A:=q^1_{K_p,\ast}\IC_{[\mc{S}_{K^p}^\ast/\ul{K_p}]}(\tilde{\mc{L}}[d](\frac{d}{2}))$, see \cite[Proposition IV.2.28]{FSGeomLLC}. For this, we use Proposition~\ref{prop: ULA-Criterion}. Namely, for each open subgroup $K_p'\subset K_p$, consider the diagram
    \[[\mc{S}^\ast_{K}/\ul{K_{p}}] \xleftarrow{a_{K_p'}}[\mc{S}^\ast_{K}/\ul{K'_{p}}]\xrightarrow{b_{K_{p}'}} \mc{S}^\ast_{K}.\]
    Then, by Proposition~\ref{prop: ULA-Criterion}, it suffices to check $b_{K_p',\ast} a_{K_p'}^\ast A$ is ULA over $\Spd k$ for a cofinal system of $K_p'$'s. By proper base change along
    \[\begin{tikzcd}
        \left[\mc{S}^\ast_{K^p}/\ul{K_p'}\right]\ar[d]\ar[r,"q^1_{K_p,K_p'}"] & \left[\mc{S}^\ast_{K}/\ul{K_p'}\right]\ar[d,"a_{K_p'}"]\\
        \left[\mc{S}^\ast_{K^p}/\ul{K_p}\right] \ar[r,"q^1_{K_p}"] & \left[\mc{S}^\ast_{K}/\ul{K_p}\right]
    \end{tikzcd}\]
    one has that $b_{K_p',\ast}a_{K_p'}^\ast A$ is isomorphic to 
    \[  b_{K_p',\ast}q^1_{K_p,K_p',\ast}\IC_{[\mc{S}_{K^p}^\ast/\ul{K_p'}]}(\tilde{\mc{L}}[d](\frac{d}{2}))\simeq \pi_{K_p',K_p,\ast}q_{K_p',\ast}\IC_{[\mc{S}_{K^p}^\ast/\ul{K_p'}]}(\tilde{\mc{L}}[d](\frac{d}{2}))\simeq \pi_{K_p',K_p,\ast} \IC_{\Shstar_{K'}}(\mc{L}^\mathrm{alg}[d](\frac{d}{2}))^\Diamond,\]
    where $\pi_{K_p',K_p}: \mc{S}_{K^{p}K_{p}'}^\ast\to \mc{S}_{K^{p}K_{p}}^\ast$ is the natural map and the last isomorphism follows from Corollary~\ref{cor: CompareIC}. But the last term is ULA over $\Spd k$, since $\IC_{\Shstar_{K'}}(\mc{L}^\mathrm{alg}[d](\frac{d}{2}))^\Diamond$ is so and $\pi_{K_p',K_p}$ is proper, see Proposition~\ref{prop: constructiblealgebrizestoULA}(2) and \cite[Proposition IV.2.11]{FSGeomLLC}. This concludes the proof of the universal local acyclicity of $\IC_{\Igs^\ast}(\mc{L})$. 

    For ULAness of the restriction to strata, for each $\PP$, consider the cartesian diagram 
     \[\begin{tikzcd}
        \left[\mc{S}_{K^p,\PP}/\ul{K_p}\right]\ar[d,"\tilde{h}_{K_p,\PP}"]\ar[r,"\tilde{i}_{\PP}"] & \left[\mc{S}^\ast_{K^p}/\ul{K_p}\right]\ar[d,"\tilde{h}_{K_p,\min}"]\\
        \Igs_{K^p,\PP} \ar[r,"i_{\PP}"] & \Igs^\ast_{K^p},
    \end{tikzcd}\]
    where the vertical arrows are cohomologically smooth surjections by Lemma \ref{lemma: smoothnessofvertmaps} and \cite[Proposition IV.2.13]{FSGeomLLC}. Note that we have
    \[\tilde{h}^{*}_{K_p,\PP}i_{\PP}^\ast \IC_{\Igs^\ast}(\mc{L})[d](\frac{d}{2})\simeq \tilde{i}_{\PP}^\ast\tilde{h}^{*}_{K_p,\min} \IC_{\Igs^\ast}(\mc{L})[d](\frac{d}{2}).\]
    It therefore suffices to show ULAness of 
    \[\tilde{i}_{\PP}^\ast\tilde{h}^{*}_{K_p,\min} \IC_{\Igs^\ast}(\mc{L})\simeq \tilde{i}_{\PP}^\ast \IC_{[\mc{S}_{K^p}^\ast/\ul{K_p}]}(\tilde{\mc{L}}[d](\frac{d}{2})),\]
    for each $\PP$. Therefore, the desired statement follows from the stratified ULAness of $\IC_{\Shstar_{K'}}(\mc{L}^\mathrm{alg}[d](\frac{d}{2}))^\Diamond$  which follows from Corollary \ref{cor: diamonddescriptionoftheICsheaf} (2), and the same argument as in the first paragraph, with $\mc{S}^\ast_{K^p}$ (resp. $\mc{S}^\ast_{K}$) replaced by $\mc{S}_{K^p,\PP}$ (resp. $\mc{S}_{K,\PP}$).
\end{proof}

In particular, when $\Lambda$ is as in Setup~\ref{assumption: coefficientsystemsingeneral} and regular, the complex $\IC_{\mathrm{Igs}^\ast}(\mathcal{L})$ is ULA over a point by Proposition~\ref{prop: ICStratifiedULA} and Proposition~\ref{prop: constructiblealgebrizestoULA}(2). We note that this in particular implies that $\IC_{\mathrm{Igs}^{*}}(\mathcal{L})$ is also overconvergent as $\Igs^\ast$ is Artin (see the proof of Corollary \ref{cor: IgusaStackStrataSmooth}). Indeed, overconvergence may be checked after pullback along a $v$-cover (see \cite[Proposition~4.4.5 (1)]{GHILZIsocComparison}), e.g. a cohomologically smooth surjection from an atlas of an Artin $v$-stack given by a locally spatial diamond. On such an atlas of $\Igs^\ast$, being ULA over a point implies overconvergence by definition (\cite[Remark~IV.2.2]{FSGeomLLC}). 

\subsubsection{Overconvergence and reflexivity}
In fact, the ULA condition in the analytic context contains two parts: overconvergence and a certain finiteness condition (perfect constructibility of its sections on quasicompact \'etale neighborhoods), see \cite[Definition IV.2.1]{FSGeomLLC}. The finiteness condition for $\IC_{\Igs^\ast}(\mc{L})$ will fail if we drop the regularity assumption on our coefficient ring $\Lambda$. Indeed, for modules over rings such as $\bb{Z}/\ell^n\bb{Z}$, the truncation functors for the standard $t$-structure can produce $\ell$-torsion, which will break perfect constructibility. However, the overconvergence condition will still hold without the regularity assumption on $\Lambda$. Going forward, the overconvergence property for the sheaves $\IC_{\Igs^{*}}(\mathcal{L})$ will play a delicate, but important technical role. Therefore, we document the following weaker version of Proposition~\ref{prop: ICStratifiedULA} for more general coefficients. 
\begin{proposition}{\label{prop: ICIgsstaroverconvergent}}
Let $\Lambda$ be as in Setup~\ref{assumption: coefficientsystemsingeneral}. The sheaves $\IC_{\Igs^{*}}(\mathcal{L})$ and $\bb{D}(\IC_{\Igs^{*}}(\mathcal{L}))$ are overconvergent.
\end{proposition}
\begin{proof}
We explain the claim for $\IC_{\Igs^{*}}(\mathcal{L})$ with the proof for $\bb{D}_{\Igs^{*}}(\IC_{\Igs^{*}}(\mathcal{L}))$ being virtually the same. Fix a pro-$p$ compact open $K_{p} \subset G(\bb{Q}_{p})$. We recall that the map $\tilde{h}_{K_{p},\min}: [\mathcal{S}_{K_{p}}^{*}/\ul{K_{p}}] \ra \Igs^{*}_{K^{p}}$ is a $v$-cover by Lemma \ref{lemma: smoothnessofvertmaps} for $K_{p}$ a pro-$p$ group. In particular, in light of \cite[Proposition~4.4.5]{GHILZIsocComparison}, we may reduce to checking the claim for $\IC_{[\mathcal{S}^{*}_{K^{p}}/\ul{K_{p}}]}(\tilde{\mathcal{L}})$. Moreover, note that we know the following. 
\begin{lemma}{\label{obs: ProofofOverconvergenceObservations}}
\begin{enumerate}
\item The pushforward $q_{K_{p}*}\IC_{[\mathcal{S}^{*}_{K^{p}}/\ul{K_{p}}]}(\tilde{\mathcal{L}})$ is overconvergent.
\item The pullback followed by pushforward of $\IC_{[\mathcal{S}^{*}_{K^{p}}/\ul{K_{p}}]}$ along the diagram 
\begin{equation}{\label{eqn: pullpushdiagram}}
\begin{tikzcd}
\left[\mathcal{S}^{*}_{K^{p}}/\ul{K_{p}'}\right] \arrow[r,"q_{K_{p}}'"] \arrow[d] & \mathcal{S}^{*}_{K^{p}K_p'} \\
\left[\mathcal{S}^{*}_{K^{p}}/\ul{K_{p}}\right] &
\end{tikzcd}
\end{equation}
is also overconvergent.
\end{enumerate}
\end{lemma}
\begin{proof}
Indeed, this follows from Corollary \ref{cor: CompareIC} and Proposition \ref{prop: propertiesofthealgebraizationfunctor} (6).
\end{proof}
Now observe that the inclusions of the Baily--Borel strata $[\mathcal{S}^{*}_{K^{p},\PP}/\ul{K_{p}}] \hookrightarrow [\mathcal{S}^{*}_{K^{p}}/\ul{K_{p}}]$ are partially proper, since they are pulled back from the corresponding strata on $\mc{S}_K^\ast$, where here the strata are partially proper, as they arise by applying $(-)^{\Diamond}$ to an algebraic stratification. Therefore, a $\Spa(C,C^{+})$-point $x: \Spa(C,C^{+}) \ra [\mathcal{S}^{*}_{K^{p}}/\ul{K_{p}}]$ must factor through some Baily--Borel stratum. We replace $[\mathcal{S}^{*}_{K^{p}}/\ul{K_{p}}]$ with the relevant stratum $[\mathcal{S}^{*}_{K^{p},\PP}/\ul{K_{p}}]$. Using the description of the space and the map $q_{K_{p},\PP}: [\mathcal{S}^{*}_{K^{p},\PP}/\ul{K_{p}}] \ra \mathcal{S}^{*}_{K^{p}K_{p},\PP}$ given in the proof of Lemma \ref{lemma: stackystratasmooth}, we see that the map $\Spa(C,C^{+}) \ra [\mathcal{S}^{*}_{K^{p}, \PP}/\ul{K_{p}}]$ factors through the quotient map to the classifying stack $[\Spa(C,C^{+})/\hat{H}_{p}]$ for some closed subgroup $\hat{H}_{p} \subset K_{p}$ acting trivially. We write $A$ for the pullback of $\IC_{[\mathcal{S}^{*}_{K^{p}}/\ul{K_{p}}]}(\tilde{\mathcal{L}})$ to this classifying stack. By definition of overconvergence and qcqs base-change, we are reduced to checking that the adjunction map $\mathrm{id}\to j_{*}j^{*}$ evaluated on $A$ is an isomorphism. Here
\[ j: [\Spa(C,\mathcal{O}_{C})/\hat{H}_{p}] \ra [\Spa(C,C^{+})/\hat{H}_{p}] \]
is the natural inclusion. Hence Diagram~\eqref{eqn: pullpushdiagram} reduces to 
\[ 
\begin{tikzcd}
\left[\Spa(C,C^{+})/\hat{H}'_{p}\right] \arrow[r] \arrow[d] & \Spa(C,C^{+}) \\
\left[\Spa(C,C^{+})/\hat{H}_{p}\right] &
\end{tikzcd}
\]
where $\hat{H}'_{p} \subset K'_{p}$ is a closed subgroup. By Lemma~\ref{obs: ProofofOverconvergenceObservations}, we know that if we pull and push $A$ along this diagram, we get an overconvergent sheaf. On the other hand, if we apply the pull push to the map $\mathrm{id}\to j_{*}j^{*}$, then we get the corresponding adjunction map on $\Spa(C,C^{+})$ (noting that the vertical arrow is smooth and the right vertical arrow is proper), which will then be an isomorphism. The claim now follows from the conservativity established in Lemma \ref{lemma: conservativitylemma}. 
\end{proof}

Another implication of the ULA condition is reflexivity, see Remark \ref{remark: ULAReflexive}. As shown in Proposition~\ref{prop: VerdierDualIC}, stratified reflexivity (see Definition~\ref{defn: stratifiedreflexive}) for intersection complexes implies Verdier self-duality of the intermediate extension functor, which helps us analyze the complexes $\IC_{\Igs^{*}}(\mathcal{L})$. Again, dropping the regularity condition on $\Lambda$, the stronger property of being stratified ULA no longer holds for $\IC_{\Igs^{*}}(\mathcal{L})$. Yet a similar analysis as our proof of overconvergence shows stratified reflexivity still holds.

Let $\Lambda$ be as in Setup~\ref{assumption: coefficientsystemsingeneral}. We start with some straightforward lemmas.
\begin{lemma}{\label{lemma: reflexivitystabilityproperties}}
Let $f: X \ra Y$ be a map of Artin $v$-stacks. The following is true. 
\begin{enumerate}
\item If $f$ is $\ell$-cohomologically proper and  $A \in \Detale(X,\Lambda)$ then the pushforward of the biduality map of (\ref{eqn: DoubleDual})
\[ A \ra \bb{D}_{X}\bb{D}_{X}(A). \]
is the biduality map of $f_{*}A$
\[ f_{*}A \ra \bb{D}_{Y}\bb{D}_{Y}(f_{*}A).\]
\item If $f$ is $\ell$-cohomologically smooth and $A \in \Detale(Y,\Lambda)$ then the pullback of the biduality map 
\[ A \ra \bb{D}_{Y}\bb{D}_{Y}(A) \]
is the biduality map of $f^{*}A$
\[ f^{*}A \ra \bb{D}_{X}\bb{D}_{X}f^{*}A.\]
\end{enumerate}
\end{lemma}
\begin{proof}
We always have identities 
\[ \bb{D}_{Y}f_{!} \simeq f_{*}\bb{D}_{X}\ \mathrm{and}\ 
\bb{D}_{X}f^{*} \simeq f^{!}\bb{D}_{Y}, \]
by projection formula and Yoneda. If $f$ is proper then this implies that we also have 
\[ \bb{D}_{Y}f_{*} \simeq f_{!}\bb{D}_{X} \]
as $f_{!} \simeq f_{*}$, which implies part (1). 
Similarly, if $f$ is $\ell$-cohomologically smooth then we have 
\[ \bb{D}_{X}f^{!} \simeq \bb{D}_{X}(f^{*}(-) \otimes f^{!}\Lambda) \simeq \bb{D}_{X}(f^{*}(-))\otimes (f^!\Lambda)^{-1} \simeq f^{!}\bb{D}_{Y}(-)\otimes(f^!\Lambda)^{-1} \simeq f^{*}\bb{D}_{Y}(-) \]
using the $\ell$-cohomological smoothness of $f$ and invertibility (in particular dualizability) of the dualizing sheaf. This implies part (2).
\end{proof}

\begin{lemma}{\label{lemma: appliedconservativitycriterion}}
The following statements are true. 
\begin{enumerate}
\item The family of functors 
\[ q_{K_{p}*}\tilde{h}_{\min,K_{p}}^{*}: \Detale(\Igs^{*}_{K^{p}},\Lambda) \ra \Detale(\mathcal{S}^{*}_{K^{p}K_{p}},\Lambda) \]
is conservative for $K_{p} \subset G(\bb{Q}_{p})$ ranging over pro-$p$ compact open subgroups.
\item For $\PP$ a conjugacy class of admissible parabolics, the family of functors 
\[ q_{K_{p},\PP*}\tilde{h}_{K_{p},\PP}^{*}: \Detale(\Igs^{*}_{K^{p},\PP},\Lambda) \ra \Detale(\mathcal{S}^{*}_{K^{p}K_{p},\PP},\Lambda) \]
is conservative for $K_{p} \subset G(\bb{Q}_{p})$ ranging over pro-$p$ compact open subgroups.
\end{enumerate}
\end{lemma}
\begin{proof}
We only prove the first claim with the proof of the second claim being strictly easier. Since reduction mod $\ell$ is conservative on $\Detale(\Igs_{K^p}^\ast,\Lambda)$, we reduce to the case $\Lambda$ is finite by using properness of $q_{K_{p}}$ (Proposition \ref{prop: qKpProper}) and projection formula. Note that $\tilde{h}_{\min,K_{p}}$ is a $v$-cover by Lemma \ref{lemma: smoothnessofvertmaps}, and hence $\tilde{h}_{\min,K_{p}}^\ast$ is conservative, using that $\mc{D}_{\et}(-,\Lambda)$ is a $v$-sheaf. Similarly, factoring $q_{K_{p}} := q^{2}_{K_{p}} \circ q^{1}_{K_{p}}$ as in Lemma \ref{lemma: q1isquasiproetaleandproper}, where $q^{1}_{K_{p}}$ is a quasi-pro-\'etale proper map, then the functor $q^{1}_{K_{p}*}$ is also conservative. Indeed, this may be checked on stalks by proper base-change, where it reduces to the fact that $R\Gamma(T,-)$ for $T$ a profinite set is a conservative functor (cf. the proof of \cite[Proposition~IV.2.28]{FSGeomLLC}). This reduces us to Lemma \ref{lemma: conservativitylemma} (cf. the proof of Proposition \ref{prop: ICStratifiedULA}).
\end{proof}

Below let us again abbreviate $\bb{D}_{\Igs^{*}}$ as $\bb{D}$ to lighten notation.
\begin{proposition}{\label{prop: stratifiedreflexivity}}
Let $\Lambda$ be as in Setup~\ref{assumption: coefficientsystemsingeneral}. Then $\IC_{\Igs^{*}}(\mathcal{L})$ and $\bb{D}(\IC_{\Igs^{*}}(\mathcal{L}))$ are stratified reflexive.
\end{proposition}
\begin{proof}
We explain the proof for $\IC_{\Igs^{*}}(\mathcal{L})$ with the proof for $\bb{D}(\IC_{\Igs^{*}}(\mathcal{L}))$ being identical. Consider the biduality map 
\[ \IC_{\Igs^{*}}(\mathcal{L}) \ra \bb{D}\bb{D}(\IC_{\Igs^{*}}(\mathcal{L})). \]
By Lemma \ref{lemma: appliedconservativitycriterion} (1), it suffices to show that this map is an isomorphism after applying $q_{K_{p}*}\tilde{h}_{\min,K_{p}}^{*}(-)[d](\frac{d}{2})$ for $K_{p} \subset G(\bb{Q}_{p})$ running through pro-$p$ groups. However, by Lemma \ref{lemma: reflexivitystabilityproperties}, this identifies with the biduality map applied to $q_{K_{p}*}\tilde{h}_{\min,K_{p}}^{*}(\IC_{\Igs^{*}}(\mathcal{L}))[d](\frac{d}{2})$, which, in light of Corollary \ref{cor: CompareIC}, identifies with the biduality map of $\IC_{\Shstar_K}(\mc{L}^\mathrm{alg}[d](\frac{d}{2}))^\Diamond$. The latter is an isomorphism by Corollary~\ref{cor: diamonddescriptionoftheICsheaf} (2). The same applies to its pullback to the strata by replacing Lemma~\ref{lemma: appliedconservativitycriterion} (1) with Lemma~\ref{lemma: appliedconservativitycriterion} (2).
\end{proof}
We now use this to deduce the following
\begin{corollary}\label{cor: VerdierSelfDualityofICIgs}
    Let $\Lambda$ be as in Setup~\ref{assumption: coefficientsystemsingeneral} and self-injective. Let $\xi$ be an algebraic representation of $\mathsf{G}_{\ol{\bb{Q}}_\ell}$ defined over $\Lambda$. Then there is a natural identification 
    \[\mathbb{D}(\IC_{\Igs^\ast}(\mc{L})) \simeq \IC_{\Igs^\ast}(\mc{L}^\vee).\]
    In particular, for $\mc{L}=\Lambda$, the intersection complex $\IC_{\Igs^\ast}\coloneq\IC_{\Igs^\ast}(\Lambda)$ is Verdier self-dual.
\end{corollary}
\begin{proof}
    This follows from Proposition~\ref{prop: stratifiedreflexivity} above and Proposition~\ref{prop: VerdierDualIC} (recalling that the dualizing sheaf on $\Igs$ is the constant sheaf). 
\end{proof}

\subsubsection{Relation to the intersection complex on the Shimura variety}
Consider the natural map $\tilde{f}_{K_p}:\mc{S}^\ast_{K^p}\to \mc{S}_K^\ast$ and factor it as before
\begin{equation}\label{Eq: MapfromInfiniteLevel}
\mc{S}^\ast_{K^p}\xrightarrow{f_{K_p}} [\mc{S}^\ast_{K^p}/\ul{K_p}]\xrightarrow{q_{K_p}}\mc{S}_K^\ast.
\end{equation}
For any $\mathcal{G}$ in $\Detale(\mathrm{Igs}^*_{K^p},\Lambda)$, we define $\mathcal{G}_{K_p}:=q_{K_p,*}\tilde{h}_{K_p,\min}^*\mathcal{G}$, a complex on the coarse quotient $\mathcal{S}^*_{K^pK_p}$. The adjunction map $q_{K_p}^\ast q_{K_p,\ast}\to \operatorname{id}$ for each $K_p$ induces a morphism
\[
\alpha_{\mathcal{G}}:\varinjlim_{K_p} \tilde{f}^*_{K_p}\mathcal{G}_{K_p} \to
\varinjlim_{K_p} f^*_{K_p}\tilde{h}^*_{K_p,\min}\mathcal{G}
\simeq \tilde{h}_{\min}^*\mathcal{G}
\]
of complexes on $\mathcal{S}^*_{K^p}$. The pullback of $\mathcal{G}$ along $[\mathcal{S}_{K^p}^*/\underline{G(\mathbb{Q}_p)}] \ra \Igs^{*}$ induces a $G(\mathbb{Q}_p)$-equivariant structure on the colimit $\varinjlim_{K_p} \tilde{f}^*_{K_p}\mathcal{G}_{K_p}$. The map $\alpha_{\mathcal{G}}$ is $G(\qp)$-equivariant for this and the natural $G(\qp)$-equivariant structure on $\tilde{h}^\ast_{\min} \mc{G}$.  

\begin{lemma}\label{lem: InfiniteLevelisColimit}
    Let $\Lambda$ be as in Setup~\ref{assumption: coefficientsystemsingeneral} (1). Assume that the sheaves $\mathcal{G}_{K_p}$ are overconvergent, for all $K_p$. Then the morphism $\alpha_{\mathcal{G}}$ is a $G(\mathbb{Q}_p)$-equivariant isomorphism. 
\end{lemma}
\begin{proof} 
We may restrict the colimit to those $K_p$ that are pro-$p$. To check that $\alpha_{\mathcal{G}}$ is an isomorphism, it suffices to do this after pulling back to geometric points $x: \Spa(C',C'^{+}) \ra \mathcal{S}^\ast_{K^p}$. By assumption, the sheaves $\mathcal{G}_{K_p}$ are all overconvergent. Then the same proof as in Proposition~\ref{prop: ICIgsstaroverconvergent} shows that $\mc{G}$ is also overconvergent. Therefore, we can reduce to $x: \Spa C' \ra \mathcal{S}^\ast_{K^p}$ being a rank one point. Assume it factors through some stratum $\mathcal{S}_{K^p,[\mathsf{P}]}$. Pulling back to this point, and using that the map $q_{K_p}$ is $\ell$-cohomologically proper (Proposition~\ref{prop: qKpProper}), we see that the diagram (\ref{Eq: MapfromInfiniteLevel}) degenerates to 
\[\tilde{f}_{K_{p}} \circ x: 
\Spa C' \xrightarrow{f_{K_{p}} \circ x} [\Spa C'/\hat{H}_{\mathsf{U},p}] \xrightarrow{q_{K_p}}\Spa C',
\]
where we have used the cartesian diagram \eqref{eqn: ClassifyingStackDescriptionofBoundary} and the notation there. Now as $K_{p} \ra \{1\}$, we have that $\hat{H}_{\mathsf{U},p} \ra \{1\}$, since $\hat{H}_{\mathsf{U},p} \subset K_{p}$ is a closed subgroup of $K_{p}$. The claim hence follows from the fact that $\tilde{h}^*_{K_p,\min}\mathcal{G}$, when pulled back to $ [\Spa C'/\hat{H}_{\mathsf{U},p}]$, defines a smooth representation of $\hat{H}_{\mathsf{U},p}$, which will be the colimit of its fixed vectors under open subgroups of $\hat{H}_{\mathsf{U},p}$. 
\end{proof}

\begin{corollary}\label{cor: InfiniteLevelisColimitIC}
Let $\Lambda$ be as in Setup~\ref{assumption: coefficientsystemsingeneral} (1). Let $\xi$ be an algebraic representation of $\mathsf{G}_{\ol{\bb{Q}}_\ell}$ defined over $\Lambda$, and $\mc{L}$, $\tilde{\mc{L}}$, $\mc{L}^\mathrm{alg}$ be the corresponding local systems on $\Igs$, $[\mathcal{S}^\ast_{K^p}/\ul{K_{p}}]$ and $\Shstar_{K}$ as before. The natural map
\[ \alpha_{\mathrm{IC}}: \varinjlim_{K_p} \tilde{f}_{K_{p}}^{*}\IC_{\Shstar_{K}}(\mc{L}^\mathrm{alg}[d](\frac{d}{2}))^{\Diamond} \ra \varinjlim_{K_p} f_{K_{p}}^{*}\IC_{[\mathcal{S}^\ast_{K^p}/\ul{K_{p}}]}(\tilde{\mc{L}}[d](\frac{d}{2})) \xrightarrow{\simeq} \tilde{h}_{\min}^\ast(\IC_{\Igs^{*}}(\mathcal{L})[d](\frac{d}{2})) \]
is a $G(\bb{Q}_{p})$-equivariant isomorphism, where the second map is obtained from Proposition~\ref{prop: PullbackofIC}.
\end{corollary}

\begin{proof}
This follows from Lemma~\ref{lem: InfiniteLevelisColimit}  applied to $\mathcal{G}=\mathrm{IC}_{\mathrm{Igs}^*}(\mathcal{L})[d](\frac{d}{2})$. Indeed, Corollary~\ref{cor: CompareIC} implies that we can identify $\mathcal{G}_{K_p}$ with $\IC_{\Shstar_{K}}(\mathcal{L}^\mathrm{alg}[d](\frac{d}{2}))^\Diamond$ in this case and these (and hence also $\mathcal{G}$) are overconvergent as discussed in the proof of Proposition \ref{prop: ICIgsstaroverconvergent}. 
\end{proof}
We are now ready to define the \emph{relative intersection cohomology of the Igusa stack}. 

\begin{definition}\label{defn:FIC}
    Let $\Lambda/\mathbb{Z}_{\ell}$ be a coefficient system as in Setup~\ref{assumption: coefficientsystemsingeneral} and let $\xi$ be an algebraic representation of $\mathsf{G}_{\ol{\mathbb{Q}}_{\ell}}$  defined over $\Lambda$. These determine an \'etale local system $\mathcal{L}_{\xi}$ on $\Igs:=\mathrm{Igs}_{K^p}$ for an appropriate choice of tame level $K^p$, as in \S \ref{sec: EtaleLocalSys}. We define the sheaf 
    \[
    \hat{\mathcal{F}}_{\xi,\IC, \Lambda}:=R\ol{\pi}_{\min,*}\IC_{\Igs^*}(\mathcal{L}_{\xi}).
    \]
    in $\Detale(\Bun_G, \Lambda)$, where $\IC_{\Igs^{*}}(\mathcal{L}_{\xi})$ is as defined in Definition \ref{defn: ICIgs*}. When $\xi$ is trivial, we write $\hat{\mathcal{F}}_{\IC,\Lambda}$ for $\hat{\mathcal{F}}_{\xi,\IC,\Lambda}$. 
\end{definition}

\section{Stalks of \texorpdfstring{$\mathcal{F}_\mathrm{IC}$ and perversity}{}}\label{sec: Stalks}

We fix a PEL type AC Shimura datum $(\mathsf{G},\mathsf{X})$ as in Assumption~\ref{assumption:codimension} before and consider the cartesian diagrams in Equation~\eqref{eqn: CartesiandiagramShimuraVarietyInfiniteLevel}. The goal of this section is to determine the stalks of the sheaf $\hat{\mc{F}}_{\xi,\IC,\Lambda}$ and compare them with the intersection cohomology of Igusa varieties. As a direct consequence of this computation, we prove the (semi-)perversity of $\hat{\mc{F}}_{\xi,\IC,\Lambda}$ with respect to the natural $t$-structure on $\mathrm{Bun}_G$ induced by the Harder--Narasimhan stratification. Let $\Lambda/\bb{Z}_\ell$ be a coefficient system as in Setup~\ref{assumption: coefficientsystemsingeneral}. We fix $k=\ol{\bb{F}}_p$ and base change Diagram~\eqref{eqn: CartesiandiagramShimuraVarietyInfiniteLevel} to $k$. 

\subsection{Intersection cohomology on perfect Igusa varieties}{\label{ss: IntCohonPerfectIgusaVarieties}}
We take $b\in B(G,\mu^{-1})$ and choose a $p$-divisible group $\bb{X}=\bb{X}_b$ over $k$ with $G$-structures representing the corresponding isogeny class. Fix a compact open subgroup $K^p\subset \mathsf{G}(\bb{A}_f^p)$, and consider the corresponding perfect Igusa variety $\Ig^b=\Ig^b_{K^p}$ over $k$,  as in \S \ref{ss: IgusaVarieties}. It has a partial minimal compactification $\Ig^{b,\ast}$, which admits a stratification $\Ig^{b,\ast}=\coprod_{\PP}\Ig^b_\PP$ labeled by conjugacy classes of admissible rational parabolic subgroups of $\mathsf{G}$, pulled back from the analogous stratification on the central leaf $\mathscr{C}^{\bb{X},\ast}$. We will refer to this stratification as the Baily--Borel stratification for $\Ig^{b,\ast}$. There are finite-level variants as well: as introduced in \S \ref{ss: IgusaVarieties}, we have a tower of Igusa varieties with their partial minimal compactifications $g_{b,m}: \Ig^{b}_{m} \hookrightarrow \Ig^{b,\ast}_{m}$, for varying $m \geq 1$. Each of these acquires a stratification $\Ig^{b,\ast}_{m}=\coprod_{\PP} \Ig^{b}_{m,\PP}$ pulled back from $\mathscr{C}^{\bb{X},\ast}$. 

Let $K_{\ell} \subset K^{p}$ denote the level at $\ell$. For each $m$, $\Ig^b_m$ supports a pro-\'etale $K_{\ell}$-torsor by passing to infinite level at $\ell$. This gives rise to a map of groupoid-valued fpqc sheaves
\[ \eta_{b,m}: \Ig^{b}_{m} \ra [\Spec(\ol{\mathbb{F}}_{p})/\underline{K_{\ell}}].\]
As in \S \ref{ss: analytificationofperfectschemes}, we have an analytification functor 
\[c^{*}: \D_{\et}([\Spec(\ol{\mathbb{F}}_{p})/\underline{K_{\ell}}],\Lambda) \ra \D_{\et}([\Spd(\ol{\mathbb{F}}_{p})/\underline{K_{\ell}}],\Lambda).\] 
This is an equivalence of categories, as can be shown by reducing first to the case $\Lambda$ is finite (by $\ell$-adic completeness) and then arguing as in the proof of \cite[Lemma 7.1.7]{GHILZIsocComparison}. In particular, this identifies $\D_{\et}([\Spec(\ol{\mathbb{F}}_{p})/\underline{K_{\ell}}],\Lambda)$ with the $\ell$-completed derived category of smooth representations  $\widehat{\D}(K_{\ell},\Lambda)$. 

In particular, as in \S \ref{sec: EtaleLocalSys}, for a finite dimensional algebraic representation $\xi$ of $\mathsf{G}_{\ol{\mathbb{Q}}_{\ell}}$ defined over $\Lambda$, we can construct a $\Lambda$-local system on $\Ig^{b}_{m}$, which we denote by $\mathcal{L}_{\xi,b,m}$. 

On the other hand, since the schemes $\Ig^{b,*}_{m}$ are of finite type over $\ol{\bb{F}}_{p}$, we can define two perverse $t$-structures
on $\D^{\mathrm{b}}_{\mathrm{c}}(\Ig^{b,\ast}_m,\Lambda)$ as in \S \ref{ss: ComparisonWithAlgebraicIC}. One is the intrinsic perverse $t$-structure $(\phantom{}^{p}\D_{\mathrm{c}}^{\mathrm{b},\leq 0}(\Ig^{b,\ast}_m,\Lambda),\phantom{}^{p}\D_{\mathrm{c}}^{\mathrm{b},\geq 0}(\Ig^{b,\ast}_{m},\Lambda))$ and the other is $(\phantom{}^{p}\D_{\mathrm{c}}^{\mathrm{b},\BB,\leq 0}(\Ig^{b,\ast}_m,\Lambda),\phantom{}^{p}\D_{\mathrm{c}}^{\mathrm{b},\BB,\geq 0}(\Ig^{b,\ast}_{m},\Lambda))$, the one defined with respect to the Baily--Borel stratification. This allows us to form a priori distinct intersection cohomology complexes
\[ \IC_{\Ig^{b,*}_{m}}^{\BB}(\mathcal{L}_{\xi,b,m}[d_{b}]),\, \IC_{\Ig^{b,*}_{m}}(\mathcal{L}_{\xi,b,m}[d_{b}]) \]
by taking the intermediate extension along $g_{b,m}$ of $\mathcal{L}_{\xi,b,m}[d_{b}]$ with respect to these two $t$-structures, where $d_b=\langle2\rho,\nu_b\rangle$ is the dimension of $\Ig^b$. 

\begin{proposition}{\label{prop: TwoDefinitionsofIntCohonperfectIgusavarieties}}
For all $m \geq 1$, we have an isomorphism in $\Dbc(\Ig^{b,*}_m,\Lambda)$
\[ \IC_{\Ig^{b,*}_{m}}(\mathcal{L}_{\xi,b,m}[d_{b}]) \simeq \IC_{\Ig^{b,*}_{m}}^{\BB}(\mathcal{L}_{\xi,b,m}[d_{b}]).\]
\end{proposition}
\begin{proof}
Using Lemma~\ref{lemma: NaturalandStratifiedIC}, it suffices to show that $g_{b,m*}(\mathcal{L}_{\xi,b,m})$ is constructible with respect to the Baily--Borel stratification. 

First, we reduce to the case of principal level away from $p$ by noting that in general $\Ig^{b,*}_{m,K^{p}}$ differs from an Igusa variety with principal level away from $p$ by a morphism which respects the Baily--Borel stratification and which, moreover, is finite \'etale after restricting to each boundary stratum (cf. Lemma \ref{lemma: proetaletorsorboundary}, which carries over to Lan's integral models of the Shimura varieties over $\mathcal{O}_E$, and therefore to the corresponding Igusa varieties as well). In particular, this implies that (\'etale locally) $g_{b,m*}(\mathcal{L}_{\xi,b,m})$ is a direct summand of the $*$-pushforward of the analogous pushforward at principal level reducing us to the principal case. 

Now, as in the proof of Corollary \ref{cor: CompareIC}, this will follow from the analogue of Pink's formula,~\cite[Theorem 4.2.1]{PinkHigherDirectImages}, for the partial minimal compactifications $\mathrm{Ig}^{b,*}_m$. The axiomatic properties of the partially compactified Igusa tower needed to establish Pink's formula, as the level at $\ell\not =p$ increases, are guaranteed by the structure of the toroidal boundary charts of Theorem~\ref{thm:IgusaStratumDiagram}. See also~\cite[Lemma 4.3.2]{LanStroh}, where these axiomatic properties are stated and~\cite[Theorem 4.3.10]{LanStroh}, where they are invoked to establish Pink's formula in various settings. The key point is that, as the level at $\ell$ increases, the transition morphisms in the boundary chart should be given, \'etale locally, by multiplication by powers of $\ell$ on an abelian scheme, respectively on a torus. 
\end{proof}

In light of this comparison, we simply denote these sheaves by $\IC_{\Ig^{b,*}_{m}}(\mathcal{L}_{\xi,b,m}[d_b])$ in what follows. We note that, since the transition morphisms $f_{mn}: \Ig^{b,*}_{m} \ra \Ig^{b,*}_{n}$ are finite and only ramified along the boundary, we have natural transition morphisms 
\[ f_{mn}^{*}\IC_{\Ig^{b,*}_{n}}(\mathcal{L}_{\xi,b,n}[d_{b}]) \ra \IC_{\Ig^{b,*}_{m}}(\mathcal{L}_{\xi,b,m}[d_{b}]) \]
obtained by adjunction, which allows us to form the following smooth $J_{b}(\bb{Q}_{p})$-representations with $\Lambda$-coefficients
\begin{equation}{\label{eqn: defofIC}}
V_{\xi,\IC,\Lambda,b} := \varinjlim_{m \geq 1} R\Gamma(\Ig^{b,*}_{m},\IC_{\Ig^{b,*}_{m}}(\mathcal{L}_{\xi,b,m}[d_{b}])).   
\end{equation}
and\footnote{We note that $V_{\xi,\IC,\Lambda,b}^{D}$ is genuinely different from $V_{\xi,\IC,\Lambda,b}$ if $\Lambda = \mathcal{O}_{F}$, as $\mathcal{O}_{F}$ is not self-injective.} 
\begin{equation}{\label{eqn: defofICdual}}
V_{\xi,\IC,\Lambda,b}^{D} := \varinjlim_{m \geq 1} R\Gamma(\Ig^{b,*}_{m},\bb{D}_{\Ig^{b,*}_{m}}\IC_{\Ig^{b,*}_{m}}(\mathcal{L}_{\xi,b,m}[d_{b}])).   
\end{equation}

We also define the following $\ell$-complete objects
\begin{equation}{\label{eqn: defofICcomplete}}
\hat{V}_{\xi,\IC,\Lambda,b} := \varprojlim_{n \geq 1} \varinjlim_{m \geq 1} R\Gamma(\Ig^{b,*}_{m},\IC_{\Ig^{b,*}_{m}}(\mathcal{L}_{\xi,b,m} \otimes_{\Lambda} \Lambda/\ell^{n}[d_{b}])).  
\end{equation}
and
\begin{equation}{\label{eqn: defofICcompletedual}}
\hat{V}_{\xi,\IC,\Lambda,b}^{D} := \varprojlim_{n \geq 1} \varinjlim_{m \geq 1} R\Gamma(\Ig^{b,*}_{m},\bb{D}_{\Ig^{b,*}_{m}}\IC_{\Ig^{b,*}_{m}}(\mathcal{L}_{\xi,b,m} \otimes_{\Lambda} \Lambda/\ell^{n}[d_{b}])).  
\end{equation}

If $\Lambda$ is torsion then we have $\hat{V}_{\xi,\IC,\Lambda,b} = V_{\xi,\IC,\Lambda,b}$ and $\hat{V}^{D}_{\xi,\IC,\Lambda,b} = V^{D}_{\xi,\IC,\Lambda,b}$. 
\begin{remark}
By Remark \ref{rem: comparison with Mantovan} and the topological invariance of the \'etale site, we can also describe $V_{\xi,\IC,\Lambda,b}$ as
\[ \varinjlim_{m \geq 1} R\Gamma(\Ig^{b,*}_{m,\Mant},\IC_{\Ig^{b,*}_{m,\Mant}}(\mathcal{L}_{\xi,b,m}[d_{b}])), \]
where $\Ig^{b,*}_{m,\Mant}$ are the normalizations with respect to $\mathscr{C}^{b} \hookrightarrow \mathscr{C}^{b,*}$, the partial minimal compactification of the Igusa varieties $\Ig^{b}_{m,\Mant}$ considered in \cite{CS17} and \cite{mantovan-PEL}, and $\IC_{\Ig^{b,*}_{m,\Mant}}(\mathcal{L}_{\xi,b,m}[d_{b}])$ is the analogous intersection complex on these algebraic varieties. A similar result holds for the $\ell$-completed and Verdier dual versions.
\end{remark}

Let $i_b: \Bun_G^b\hookrightarrow \Bun_G$ be the inclusion of the stratum labeled by $b$. The main result of this section is the following.
\begin{theorem}\label{thm: stalks}
Let $\Lambda$ be as in Setup~\ref{assumption: coefficientsystemsingeneral} and $\xi$ be an algebraic representation of $\mathsf{G}_{\ol{\mathbb{Q}}_{\ell}}$ defined over $\Lambda$. Then, for each $b \in B(G,\mu^{-1})$, under the isomorphism $\Detale(\Bun_G^b,\Lambda)\simeq \widehat{\D}(J_b(\qp),\Lambda)$, we have an identification $i_b^\ast \hat{\mc{F}}_{\xi,\IC,\Lambda} \simeq \hat{V}_{\xi,\IC,\Lambda,b}[-d_b]$ and $i_b^\ast \bb{D}_{\Bun_{G}}\hat{\mc{F}}_{\xi,\IC,\Lambda} \simeq \hat{V}^{D}_{\xi,\IC,\Lambda,b}[-d_b]$ of $\ell$-completed smooth $J_{b}(\bb{Q}_{p})$-representations.
\end{theorem}
We make some preparations for this computation.

\subsection{Generic and special Newton polygon stratifications}\label{ss: GenericSpecialNP}
Recall from \S\ref{sec: IgusaStack} that, using the cartesian diagram \eqref{eqn: CartesiandiagramShimuraVarietyInfiniteLevel}, we may define a Newton polygon stratification of the minimally compactified Shimura variety by pulling back the Harder--Narasimhan stratification on $\Bun_G$ by the substacks $i_{b}: \Bun_{G}^{b} \hookrightarrow \Bun_{G}$. We will call this the \textit{generic Newton polygon stratification} (or simply Newton stratification), and denote the strata by $\mc{S}^{\ast,b}_{K^p}$, or sometimes by $\mc{S}^{\ast,b, \eta}_{K^p}$ to distinguish it from the \textit{special Newton polygon stratification} that we will now discuss. Throughout this section, we endow $B(G)$ with its order topology, with the partial order on $B(G)$ given as in the Notation and Convention section. Namely, $b\in\ol{\{b'\}}$ if and only if $b\geq b'$. Note that this is opposite to the specialization relation for $G$-isocrystals, but compatible with the topology on $\Bun_G$.

Pick a hyperspecial level subgroup $K_p\subset G(\qp)$. The Shimura variety at level $K=K_pK^p$ has a canonical integral model $\mathscr{S}_K$ over $\Spf \mc{O}_E$, which possesses a minimal compactification $\mathscr{S}_K^\ast$. For simplicity, below we base-change it to $\mc{O}_{\bb{C}_p}$, and maintain the same notation. The (geometric) special fiber $\mathscr{S}_{K,\ol{\bb{F}}_p}$ has a Newton polygon stratification according to the isogeny class of the universal $p$-divisible group with $G$-structures. Using the theory of well-positioned subschemes, see \cite[\S 3.3]{LanStroh}, this stratification extends over $\mathscr{S}_{K,\ol{\bb{F}}_p}^\ast$. Consider the specialization map from the underlying topological space of the adic generic fiber $\mc{S}^\ast_K=\mc{S}^\ast_{K,\bb{C}_p}$ to that of $\mathscr{S}_{K,\ol{\bb{F}}_p}^\ast$
\[\operatorname{sp}:|\mc{S}^\ast_{K}|\to |\mathscr{S}_{K,\ol{\bb{F}}_p}^\ast|.\]
It sends a point $x=\Spa(K,K^+)\to \mc{S}^\ast_{K}$ (where $K$ is a non-archimedean field over $\bb{C}_p$) to the image of the unique closed point, under the corresponding map $\Spf K^+ \to \mathscr{S}^\ast_{K}$.
\begin{definition}
    For each $b\in B(G,\mu^{-1})$, we define the special Newton polygon stratum of $\mc{S}^\ast_{K}$ attached to $b$ to be the locally closed sub-v-sheaf 
    \[\mc{S}^{\ast,b, s}_{K}:=\mc{S}^\ast_{K}\times_{|\mc{S}^\ast_{K}|,\,\operatorname{sp}} |\mathscr{S}_{K,\ol{\bb{F}}_p}^{\ast,b}|.\]
    For any $K_p'\subset K_p$, we define the special Newton polygon strata $\mc{S}^{\ast,b, s}_{K_p'K^p}$ by taking the preimages of $\mc{S}^{\ast,b, s}_{K}\subset \mc{S}_{K}^*$.
\end{definition}
\begin{remark}\label{rem: GenericSpecial}
    Given $b\in B(G,\mu^{-1})$, the generic and special Newton polygon strata labeled by $b$ have the same rank one points. For a higher rank point $x=\Spa(K,K^+)\to \mc{S}_K^\ast$, it lies in $\mc{S}^{\ast,b,\eta}_K$ if and only if, under the corresponding map $\Spf K^+ \to \mathscr{S}^\ast_{K}$, the image of the generic point lies in $|\mathscr{S}_{K,\ol{\bb{F}}_p}^{\ast,b}|$; while it lies in $\mc{S}^{\ast,b,s}_K$ if and only if under $\Spf K^+ \to \mathscr{S}^\ast_{K}$, the image of the closed point lies in $|\mathscr{S}_{K,\ol{\bb{F}}_p}^{\ast,b}|$. This explains our choice of terminology. From this description, we also see that the generic Newton polygon strata are partially proper, cf. \cite[Definition 18.4]{Ecod}.
\end{remark}
\begin{remark}
    Because of the difference on higher rank points, the specialization relations among the generic and special polygon strata are reversed relative to each other. Namely, for a generic Newton polygon stratum $\mc{S}^{\ast,b',\eta}_K$ to lie in the closure of $\mc{S}^{\ast,b,\eta}_K$, one must have that $\Bun_G^{b'}$ lies in the closure of the Harder--Narasimhan stratum $\Bun_{G}^{b} \subset \Bun_{G}$. This closure is contained in (and in fact equal to, by \cite[Theorem~1.1]{Vi}) the locus of $G$-bundles in $\Bun_{G}$ with Kottwitz element $b' \geq b$. On the other hand, for a special Newton polygon stratum $\mc{S}^{\ast,b',s}_K$ to lie in the closure of $\mc{S}^{\ast,b,s}_K$, one must have $b' \leq b$ due to the reversal of these closure relationships on the special fiber. (The inequalities are with respect to our convention for the partial order on $B(G)$ that is fixed in \S\ref{Notation}.)
\end{remark}
We can study the interaction between the two stratifications. For each $b$, we consider the v-sheaf 
\[\mc{S}^{\ast,b,\cap}_K:=\mc{S}^{\ast,b,\eta}_K\times_{\mc{S}^{\ast}_K} \mc{S}^{\ast,b,s}_K.\] 
One result that will be useful to us is the following.
\begin{lemma}\label{lem: CapStrataOpen}
    For each $b$, the natural map of v-sheaves $\mc{S}^{\ast,b,\cap}_K\hookrightarrow \mc{S}^{\ast}_K$
    is an open immersion. 
\end{lemma}
\begin{proof}
    This follows since we can rewrite $\mc{S}^{\ast,b,\cap}_K$ as the fiber product of two open subfunctors
    \[\mc{S}^{\ast,b,\cap}_K\simeq \left(\bigcup_{b'\leq b}
    \mc{S}^{\ast,b',\eta}_K\right)\times_{\mc{S}^{\ast}_K} \left(\bigcup_{b''\geq b}\mc{S}^{\ast,b'',s}_K\right).\]
 Indeed, this identification can be checked on the level of underlying topological spaces. Clearly, $|\mc{S}^{\ast,b, \cap}_K|$ is contained in the topological space on the right-hand side. But, given any $x=\Spa(C,C^+)$-point on the right-hand side, where $C/\bb{C}_p$ is algebraically closed, we obtain a corresponding $\Spf C^+$-point of $\mathscr{S}^\ast_K$. According to Remark~\ref{rem: GenericSpecial}, its generic point lies in $|\mathscr{S}^{\ast,b'}_{\ol{\bb{F}}_p}|$ for some $b'\leq b$. On the other hand, its closed point lies in $|\mathscr{S}^{\ast,b''}_{\ol{\bb{F}}_p}|$ for some $b''\geq b$. Together, these two inequalities imply $b'\leq b''$. However, the closed point is a specialization of the generic point. Since the specialization relation for the universal family of $G$-isocrystals on $\mathscr{S}^\ast_{\ol{\bb{F}}_p}$ is inverse to the specialization relation of $G$-bundles on the Fargues--Fontaine curve, we must also have $b'\geq b''$. Combining all the inequalities, we see that $b'=b=b''$, which means $x$ has to lie in $\mc{S}^{\ast,b,\cap}_K$ as desired.
\end{proof}

\subsection{Open substacks in \texorpdfstring{$\Igs^\ast$}{}}\label{sec: OpenStratainIgs}
Lemma~\ref{lem: CapStrataOpen} leads to an interesting consequence on the geometry of $\Igs^\ast$. Namely, for each $b\in B(G,\mu^{-1})$, the substack $[\Ig^{b,\ast,\diamond}/\mc{J}_b]$ of $\Igs^{\ast,b}$ is in fact open in $\Igs^{\ast}$.

\begin{proposition}\label{prop: openstrata}
    The composition of the natural maps 
    \[\alpha_b: [\Ig^{b,\ast,\diamond}/\mc{J}_b]\hookrightarrow [\Ig^{b,\ast,\dagger}/\mc{J}_b]\simeq \Igs^{\ast,b}\hookrightarrow\Igs^\ast\]
    is an open immersion. Moreover, we have an identification 
    \begin{equation}\label{eq: CapStrata}
        \mc{S}_{K^p}^{\ast,b,\cap}\simeq [\Ig^{b,\ast,\diamond}/\mc{J}_b]\times_{\Igs^\ast}\mathcal{S}^*_{K^p}.
    \end{equation}
\end{proposition}
\begin{proof}
     By \cite[Proposition 10.11]{Ecod}, we can check that $\alpha_b$ is an open immersion after pullback to a v-cover. By combining with Lemma~\ref{lem: CapStrataOpen}, we see that the second part of the proposition implies the first part. For the second part, note that we have 
     \[\mc{S}_{K^p}^{\ast,b}\simeq [\Ig^{b,\ast,\dagger}/\mc{J}_b]\times_{\Igs^\ast}\mathcal{S}^*_{K^p}.\]
     This is a locally spatial diamond, since it is locally closed in the spatial diamond $\mc{S}_{K^p}^\ast$. \footnote{Indeed, factoring the injection $\mc{S}_{K^p}^{\ast,b}\hookrightarrow \mc{S}_{K^p}^{\ast}$ into an open immersion followed by a closed immersion, we see this by combining \cite[Corollary 11.26, Proposition 11.19]{Ecod}.} By \cite[Proposition 11.20]{Ecod} and Lemma~\ref{lemma: IgCCqcqs}(1), we just need to identify the underlying topological spaces of the two sides in Equation~\eqref{eq: CapStrata}. We may do this after restricting to each geometric fiber of the Hodge--Tate period map. 
     
     Fix a geometric point $x=\Spa C$ of $\Fl$ and identify the fiber $\pi^{-1}(x)$ of the Hodge--Tate period map over $x$ with $\Ig^{b,\ast,\dagger}_C$. Let $k$ denote the residue field of $C$. We need to show that, under this identification, we have an isomorphism
     \begin{equation}\label{eq: CapFibers}
         \mc{S}_{K^p}^{\ast,b,\cap}\times_{\mc{S}_{K^p}^\ast} \pi^{-1}(x)\simeq \Ig^{b,\ast,\diamond}_C.
     \end{equation}

     For this, let us observe that one can characterize the subfunctor $\Ig^{b,\ast,\diamond}_C\hookrightarrow \Ig^{b,\ast,\dagger}_C$ as follows: Take a perfectoid space $S=\Spa(R,R^+)$ over $\Spa C$. By the partial properness of $\Ig^{b,\ast,\dagger}_C$, giving a map $f:S\to \Ig^{b,\ast,\dagger}_C$ is the same as giving a map $\Spec R^\circ \to \Ig^{b,\ast}_k$. Then $f$ factors through $\Ig^{b,\ast,\diamond}_C$ if and only if this map extends over $\Spec R^+$, i.e. there is a dotted arrow filling in the commutative diagram 
     \[\begin{tikzcd}
         \Spec R^\circ\ar[r,"f"]\ar[d, hook]& \Ig^{b,\ast}_k\\
         \Spec R^+\ar[ur,dashed,"g", swap]& .
     \end{tikzcd}\]

     The left hand side of Equation~\eqref{eq: CapFibers} has a similar description. Namely, let $K\subset \mathsf{G}(\bb{A}_f)$ be hyperspecial at $p$, consider the natural map $q: \Ig_k^{b,\ast}\to \mathscr{S}^{\ast,b}_{K,k}$, which factors through the perfection $\mathscr{S}^{\ast,b,\mathrm{perf}}_{K,k}$ of the latter. Then, for $S=\Spa(R,R^+)$ over $\Spa C$ as above, a map from $S$ to the LHS is the same as a map $f: \Spec R^\circ \to \Ig^{b,\ast}$, such that there is a dotted arrow filling in the commutative diagram
     \[\begin{tikzcd}
         \Spec R^\circ\ar[r,"f"]\ar[d, hook]& \Ig_k^{b,\ast}\ar[d,"q"]\\
         \Spec R^+ \ar[r,dashed,"g", swap]& \mathscr{S}^{\ast,b,\mathrm{perf}}_{K,k}.
     \end{tikzcd}\]
     
     For our purpose, it suffices to compare the above two subfunctors of $\Ig^{b,\ast,\dagger}_C$, when taking test objects to be geometric points. It then follows from the valuative criterion of properness for the map $q$ that these two subfunctors agree. Note that the valuative criterion holds for $q$, since it factors as a pro-proper map $\Ig_k^{b,\ast}\to \mathscr{C}^{b,\ast,\mathrm{perf}}_{k}$ to the central leaf, and a closed immersion $\mathscr{C}^{b,\ast,\mathrm{perf}}_k\hookrightarrow \mathscr{S}^{\ast,b,\mathrm{perf}}_{K,k}$. (The fact that $\mathscr{C}^{b,\ast}\hookrightarrow \mathscr{S}^{\ast,b}_{K}$ is a closed immersion follows from \cite[Proposition 3.4.1(2)]{LanStroh} and the construction of partial minimal compactifications of well-positioned subsets.)
\end{proof}

\subsection{Proof of Theorem~\ref{thm: stalks}}
We are now ready to prove Theorem~\ref{thm: stalks}. For notational simplicity, we will prove that claim for $\hat{\mathcal{F}}_{\xi,\IC,\Lambda}$ with the claim for $\bb{D}_{\Bun_{G}}\hat{\mathcal{F}}_{\xi,\IC,\Lambda}$ being completely analogous (by replacing $\IC$ with $\bb{D}\IC$ everywhere and noting that Corollary \ref{cor: CompareIC} and Proposition \ref{prop: ICIgsstaroverconvergent} are both true for the Verdier dual). This will take a few steps.
\subsubsection{Step 1} We first consider the restriction of $\IC_{\Igs^\ast}(\mc{L}_\xi)$ to the open substacks discussed in \S \ref{sec: OpenStratainIgs}. We may pull back the Baily--Borel stratum $\Igs_\PP$ of $\Igs^{*}$ along the map $\alpha_{b}: [\Ig^{b,*,\diamond}/\mathcal{J}_{b}] \ra \Igs^{*}$ of Proposition~\ref{prop: openstrata}. Combining Proposition~\ref{prop: FiberwiseStratifiedCartesianEq} and Lemma~\ref{lemma: differentstrataofIgs}, this identifies with the locally closed stratum $[(\Ig^{b,\diamond})_{\PP}/\mathcal{J}_{b}]$, as defined in \eqref{eqn: ThickenedBoundaryStrata}. In other words, we have an isomorphism 
\begin{equation}{\label{eqn: PullingBackStrata}}
 [(\Ig^{b,\diamond})_{\PP}/\mathcal{J}_{b}] \simeq [\Ig^{b,\ast,\diamond}/\mathcal{J}_{b}] \times_{\Igs^\ast} \Igs_{\PP}.  
\end{equation}
By Lemma \ref{lemma: StrataFiberwiseSmooth}, each $[(\Ig^{b,\diamond})_\PP/\mc{J}_b]$ is an Artin v-stack that is $\ell$-cohomologically smooth over $\Spd \ol{\bb{F}}_p$ of pure $\ell$-dimension $d_{b,\PP} - d_{b} = d_{\PP} -d$, where $d$ and $d_\PP$ are the dimensions of the Shimura variety and its $\PP$-stratum, respectively, and $d_{b,\PP}$ is the dimension of $\Ig^{b}_{\PP}$ and the equality follows from Proposition~\ref{prop: dimensionStrataIg}.

In particular, we can apply the constructions in \S\ref{subsub: PerverseSetup} and define a perverse $t$-structure on $\Detale([\Ig^{b,\ast,\diamond}/\mc{J}_b],\Lambda)$.  

Denote by $\mathcal{L}_{\xi, b}$ the restriction of $\mc{L}_\xi$ to
\[[(\Ig^{b,\diamond})_{[\mathsf{G}]}/\mathcal{J}_{b}]\simeq [\Ig^{b,\ast,\diamond}/\mathcal{J}_{b}] \times_{\Igs^\ast}\Igs.\]
When $\Lambda$ is as in Setup \ref{assumption: coefficients} (1), we can consider the intersection complex for $\mathcal{L}_{\xi, b}$ along the open immersion 
\[[(\Ig^{b,\diamond})_{[\mathsf{G}]}/\mc{J}_b]\hookrightarrow[\Ig^{b,\ast,\diamond}/\mc{J}_b]\]
and denote it by $\IC_{[\Ig^{b,\ast,\diamond}/\mc{J}_b]}(\mathcal{L}_{\xi, b})$. When $\Lambda$ is as in Setup \ref{assumption: coefficients} (2), we define it via taking inverse limits, as in Definition \ref{defn: ICIgs*}. Here, let us again caution the reader that $(\Ig^{b,\diamond})_{[\mathsf{G}]}$ is not the same as $\Ig^{b,\diamond}$; the difference between these objects is explained in Remark~\ref{remark: differencestrata}. 
\begin{lemma}\label{lem: restrictionICtocapstrata}
    There is a natural identification
    \[\alpha_b^\ast\IC_{\Igs^\ast}(\mathcal{L}_{\xi}) \simeq \IC_{[\Ig^{b,\ast,\diamond}/\mc{J}_b]}(\mathcal{L}_{\xi, b}) \]
    in $\Detale([\Ig^{b,*.\diamond}/\mathcal{J}_{b}],\Lambda)$. 
\end{lemma}
\begin{proof}
    Both $\IC_{\Igs^\ast}(\mathcal{L}_{\xi})$ and $\IC_{[\Ig^{b,\ast,\diamond}/\mc{J}_b]}(\mathcal{L}_{\xi,b})$ have a presentation using Deligne's formula with torsion coefficients, see Corollary~\ref{cor: DeligneFormula}. Since $\alpha_b$ is an open immersion by Proposition~\ref{prop: openstrata}, pullback along $\alpha_b$ commutes with push-forwards (also with truncation since $\alpha_b^\ast$ is exact). Hence, we obtain the desired result with torsion coefficients, by combining (\ref{eqn: PullingBackStrata}) and the equality $d_{b,\PP} - d_{b} = d_{\PP} - d$, of the $\ell$-codimensions appearing in the definition of $\IC_{\Igs^\ast}(\mathcal{L}_{\xi})$ and $\IC_{[\Ig^{b,\ast,\diamond}/\mc{J}_b]}(\mathcal{L}_{\xi, b})$. With integral coefficients, we note that, since $\alpha_{b}$ is an open immersion, $*$-pullback commutes with inverse limits.
\end{proof}

\subsubsection{Step 2} We need to relate the cohomology of $\IC_{[\Ig^{b,\ast,\diamond}/\mc{J}_b]}(\mathcal{L}_{\xi,b})$ with the colimit of the cohomology of the algebraic intersection cohomology groups $\IC_{\Ig^{b,\ast}_{m}}(\mc{L}_{\xi,b,m}^\mathrm{alg})$, defined in \S \ref{ss: IntCohonPerfectIgusaVarieties}. Here we added a superscript ${(-)}^\mathrm{alg}$ for the local systems on the algebraic varieties to distinguish them from the ones on the corresponding v-sheaves. For this purpose, we reinterpret $\IC_{[\Ig^{b,\ast,\diamond}/\mc{J}_b]}(\mathcal{L}_{\xi,b})$ in terms of the functor $t^\ast c^\ast$ introduced in \S \ref{ss: analytificationofperfectschemes} 

More precisely, we consider the finite level Igusa varieties obtained by taking the coarse moduli spaces of the quotient $[\Ig^{b,\ast}/\underline{K_b(p^m})]$. These coarse moduli spaces are given by the perfections $(\Ig^{b,\ast}_{m})^{\perf}$ of the partially compactified Igusa varieties $\Ig^{b,\ast}_{m}$ defined in \S \ref{ss: IgusaVarieties}, for compact open subgroups $K_b(p^m)\subset J_b(\qp)$ at principal levels (cf. Lemma \ref{lemma: IgusaVarietyStackyQuotient}). 

Similar to $\Ig^{b,\diamond}$, the v-sheaf $(\mathrm{Ig}^{b,\ast}_m)^{\perf,\diamond}= \mathrm{Ig}^{b,\ast,\diamond}_m$ has a stratification by the locally closed sub-v-sheaves
\[(\Ig^{b,\diamond}_{m})_{\PP}:=\Ig^{b,\ast,\diamond}_{m} \times_{\left|{\Ig^{b,\ast,\diamond}_{m}}\right|}\left(\left|{\Ig^{b,\diamond}_{m,\leq\PP}}\right|\backslash{\left|{\Ig^{b,\diamond}_{m,<\PP}}\right|}\right), \]
as in (\ref{eqn: ThickenedBoundaryStrata}). Similar to Lemma~\ref{lemma: differentstrataofIgs}, for each $\PP$, we have that 
\begin{equation}{\label{eqn: finitelevelthickenedboundary}}
(\Ig^{b,\diamond}_{m})_{\PP}\simeq \Ig^{b,\Diamond}_{m,\PP}\times_{{\Ig^{b,\ast,\diamond}_{m}}}\Ig^{b,\ast,\diamond}_m, 
\end{equation}
where we used the open immersion $t_{\Ig^{b,*,\diamond}_{m}}: \Ig^{b,*,\diamond}_{m} \ra \Ig^{b,*,\Diamond}_{m}$ in forming the above fiber product, cf. Lemma \ref{lemma: IgCCqcqs}. It follows that $(\Ig^{b,\diamond}_{m})_{\PP}$ is $\ell$-cohomologically smooth of dimension $d_{b,\PP}$, since this is true for $\Ig^{b,\Diamond}_{m,\PP}$, cf. Lemma \ref{lemma: StrataFiberwiseSmooth}.

Therefore, we can consider the local system $\mc{L}_{\xi,b,m}$ on $(\Ig^{b,\diamond}_{m})_{[\mathsf{G}]}$ obtained by restricting $c_{\Ig^b_m}^\ast \mc{L}^{\alg}_{\xi,b,m}$ along $(\Ig^{b,\diamond}_{m})_{[\mathsf{G}]}\hookrightarrow \Ig_{m}^{b,\ast, \Diamond}$ and taking its intermediate extension, denoted $\IC_{\Ig_{m}^{b,\ast,\diamond}}(\mc{L}_{\xi,b,m})$, with respect to this stratification. We have the following.
\begin{lemma}\label{lem: diamondICviaAnalytification}
With notation as above, we have a natural isomorphism
    \[\IC_{\Ig_{m}^{b,\ast,\diamond}}(\mc{L}_{\xi,b,m}[d_{b}])\simeq t_{\Ig_{m}^{b,\ast}}^\ast c_{\Ig_{m}^{b,\ast}}^\ast \IC_{\Ig^{b,\ast}_{m}}(\mc{L}_{\xi,b,m}^\mathrm{alg}[d_{b}]).\]
\end{lemma}
\begin{proof}
    By Lemma~\ref{lemma: intermediateextensioncompatiblewithanalytification}, we have  
    \begin{equation}\label{eq: CompareICDiamonds}
        c_{\Ig_{m}^{b,\ast}}^\ast \IC_{\Ig^{b,\ast}_{m}}(\mc{L}_{\xi,b,m}^\mathrm{alg}[d_{b}])\simeq  \IC_{\Ig^{b,\ast,\Diamond}_{m}}(c_{\Ig_{m}^{b,\ast}}^\ast\mc{L}_{\xi,b,m}^\mathrm{alg}[d_{b}]),
    \end{equation}
    where, on the right-hand side, the intersection complex is formed with respect to the stratification on $\Ig^{b,\ast,\Diamond}_{m}$ given by the $\Ig^{b,\ast,\Diamond}_{m,\PP}$s. Since $t_{\Ig_m^{b,\ast}}$ is an open immersion, we have by smooth base-change and (\ref{eqn: finitelevelthickenedboundary}) an isomorphism
    \[\IC_{\Ig_{m}^{b,\ast,\diamond}}(\mc{L}_{\xi,b,m}[d_{b}])\simeq t_{\Ig_m^{b,\ast}}^\ast\IC_{\Ig^{b,\ast,\Diamond}_m}(c_{\Ig_{m}^{b,\ast}}^\ast\mc{L}_{\xi,b,m}^\mathrm{alg}[d_{b}]).\]
    Combining this with Equation~\eqref{eq: CompareICDiamonds}, we get the desired claim.
\end{proof}

\subsubsection{Step 3}
Now we can identify the cohomology of $\IC_{[\Ig^{b,\ast,\diamond}/\mc{J}_b]}(\mc{L}_{\xi,b})$. Let $\beta_b: [\Ig^{b,\ast,\diamond}/\mc{J}_b]\to [\ast/\mc{J}_b]$ be the natural map, and let $V_{\xi,\IC,\Lambda,b}$ be the representation of $J_b(\qp)$ defined in (\ref{eqn: defofIC}).
\begin{proposition}\label{prop: CohoSmallDiamond}
    Under the identification $\Detale([\ast/\mc{J}_b],\Lambda)\simeq \widehat{\D}(J_b(\qp),\Lambda)$, we have an isomorphism
    \[R\beta_{b,\ast}\IC_{[\Ig^{b,\ast,\diamond}/\mc{J}_b]}(\mathcal{L}_{\xi, b}) \simeq \hat{V}_{\xi,\IC,\Lambda,b}[-d_b].\]
\end{proposition}
\begin{proof}
Since $\IC_{[\Ig^{b,\ast,\diamond}/\mc{J}_b]}(\mathcal{L}_{\xi, b})$ is an inverse limit over the mod $\ell^{n}$ IC-sheaves by Definition, we reduce to the case of torsion coefficients using Theorem \ref{thm: inverselimit} to rewrite the RHS of the desired isomorphism as an inverse limit over the mod $\ell^{n}$-reductions of the IC-sheaf. Since the identification $\Detale([\ast/\mc{J}_b],\Lambda)\simeq \D(J_b(\qp),\Lambda)$ is made via pullback along 
\begin{equation}\label{eq: restricttoJbqp}
    u_b: [\ast/\ul{J_b(\qp)}]\to [\ast/\mc{J}_b],
\end{equation} 
which is an $\ell$-cohomologically smooth map of $\ell$-dimension $d_b=\langle2\rho,\nu_b\rangle$, we just need to compute the pushforward of 
\[\IC_{[\Ig^{b,\ast,\diamond}/\mc{J}_b]}(\mathcal{L}_{\xi, b})|_{[\Ig^{b,\ast,\diamond}/\ul{J_b(\qp)}]}\simeq \IC_{[\Ig^{b,\ast,\diamond}/\ul{J_b(\qp)}]}(\mathcal{L}_{\xi, b}[d_{b}])[-d_b]\]
along the map
\[\tilde{\beta}_b: [\Ig^{b,\ast,\diamond}/\ul{J_b(\qp)}]\to [\ast/\ul{J_b(\qp)}], \]
where $\IC_{[\Ig^{b,\ast,\diamond}/\ul{J_b(\qp)}]}(\mathcal{L}_{\xi, b}[d_{b}])$ is defined completely analogously to $\IC_{[\Ig^{b,\ast,\diamond}/\mc{J}_b]}(\mathcal{L}_{\xi, b})$\footnote{Note in particular that the map $[\ast/\underline{J_{b}(\qp)}\ra [\ast/\mc{J}_{b}]]$ is $\ell$-cohomologically smooth, so we may use smooth base-change to compare Deligne's formula on both sides in the case of torsion coefficients. Moreover, $*$-pullback commutes with inverse limits, so we may reduce to the torsion case.} using that $\ul{J_{b}(\bb{Q}_{p})}$ preserves the boundary strata of $\Ig^{b,\ast,\diamond}$ by Proposition \ref{prop: JbPreserveBoundary}.
For each $m \geq 1$, we have a cofinal system of compact open subgroups $K_{b}(p^{m}) \subset J_{b}(\bb{Q}_{p})$, as in Lemma \ref{lemma: IgusaVarietyStackyQuotient}. We then have a diagram 
\begin{equation}{\label{eqn: KeyDiagram}}
\begin{tikzcd}
\left[ \Ig^{b,\ast,\diamond}/\ul{K_{b}(p^{m})} \right] \arrow[r,"q_{m}"] \arrow[d,"h_{m}"] & \Ig^{b,\ast,\diamond}_{m}\\
\left[ \Ig^{b,\ast,\diamond}/\Jbul \right] &. 
\end{tikzcd}
\end{equation}
Here $h_{m}$ is the natural map on stack quotients, which is \'etale. The map $q_{m}$ is the map from the stack quotient to the coarse quotient, and it will be cohomologically proper of cohomological amplitude $0$ if we choose $m$ sufficiently large so that $K_{b}(p^{m})$ is pro-$p$, cf. Proposition \ref{prop: qKpProper}. 

In particular, we obtain 
\begin{align*}
    R\tilde{\beta}_{b,\ast}\IC_{[\Ig^{b,\ast,\diamond}/\ul{J_b(\qp)}]}(\mathcal{L}_{\xi, b}[d_{b}])[-d_b] &\simeq \varinjlim_{m} R\Gamma([ \Ig^{b,\ast,\diamond}/\ul{K_{b}(p^{m})}], h_{m}^{*}\IC_{[ \Ig^{b,\ast,\diamond}/\ul{J_{b}(\qp)}]}(\mathcal{L}_{\xi, b})[d_b])[-d_{b}]\\
    &\simeq \varinjlim_{m} R\Gamma(\Ig^{b,\ast,\diamond}_{m}, \IC_{[ \Ig^{b,\ast,\diamond}/\ul{K_{b}(p^m)}]}(\mathcal{L}_{\xi, b,m}[d_b]))[-d_b] \\
    &= \varinjlim_{m} R\Gamma(\Ig^{b,\ast,\diamond}_{m}, \IC_{\Ig^{b,\ast,\diamond}_{m}}(\mathcal{L}_{\xi,b,m}[d_{b}]))[-d_b].
\end{align*}
By Lemma~\ref{lem: diamondICviaAnalytification} and Proposition~\ref{prop: SmallDiamondvsPerfectScheme}, this identifies with the algebraic intersection cohomology
\[\varinjlim_{m} R\Gamma(\Ig^{b,\ast}_{m},\IC_{\Ig^{b,\ast}_{m}}(\mathcal{L}_{\xi,b,m}^\mathrm{alg}[d_b]))[-d_{b}],\]
which is precisely $\hat{V}_{\xi,\IC,\Lambda,b}[-d_b]$, as desired.
\end{proof}

\subsubsection{Step 4}
To conclude, we need to control the difference between the spaces $\Ig^{b,\ast,\diamond}$ and $\Ig^{b,\ast,\dagger}$ (see Lemma~\ref{lem: DaggerisCC}), since it is the latter that appear as fibers of $\ol{\pi}_{\min}$. This difference does not matter for overconvergent sheaves, and in particular, the IC sheaves.

Recall that if $X$ is a $v$-stack then we have the full subcategory of overconvergent sheaves $\Detale^{\oc}(X,\Lambda)  \subset \Detale(X,\Lambda)$. It is easy to see that the condition of being overconvergent is stable under pullback. We now have the following.

\begin{proposition}{\label{prop: overconvergentreplacement}}
Consider the natural map \[a_{b}: [\Ig^{b,*,\diamond}/\Jbul] \ra [\Ig^{b,*,\dagger}/\Jbul],\] 
then the functors $a_{b}^{*}$ and $a_{b*}$ define inverse equivalences
\[ \Detale^{\oc}([\Ig^{b,*,\diamond}/\Jbul],\Lambda) \simeq \Detale^{\oc}([\Ig^{b,\dagger}/\Jbul],\Lambda).\]
\end{proposition}
\begin{proof}
By Lemma~\ref{lemma: IgCCqcqs}, the map $a_{b}: [\Ig^{b,*,\diamond}/\Jbul] \ra [\Ig^{b,*,\dagger}/\Jbul]$ is a qcqs open immersion of $v$-stacks. Hence, the claim follows from \cite[Lemma~4.4.7]{GHILZIsocComparison}.
\end{proof}

We obtain Theorem~\ref{thm: stalks} as a corollary of this.
\begin{corollary}[{Theorem~\ref{thm: stalks}}]
For all $b \in B(G,\mu^{-1})$, under the equivalence 
\[\Detale(\Bun_{G}^{b},\Lambda) \simeq \widehat{\D}(J_{b}(\bb{Q}_{p}),\Lambda),\] 
the complex $i_b^\ast \mc{F}_{\xi,\IC,\Lambda}$ is identified with the $\ell$-complete smooth representation $\hat{V}_{\xi,\IC,\Lambda,b}[-d_b]$.
\end{corollary}
\begin{proof}
Consider the natural maps 
\[\begin{tikzcd}
    \left[\ast/\Jbul\right] &\left[\Ig^{b,\ast,\dagger}/\Jbul\right] \ar[l, "\tilde{\beta}_b^\dagger", swap]\ar[r,"k_b"]& \Igs^\ast.
\end{tikzcd}\]
We know that $\IC_{\Igs^\ast}(\mathcal{L}_{\xi})$ is overconvergent by Proposition \ref{prop: ICIgsstaroverconvergent}. Since overconvergence is stable under pullback, $k_b^\ast\IC_{\Igs^\ast}(\mathcal{L}_{\xi})$ is also overconvergent. Hence, we have
\begin{align*}
u_b^\ast i_b^\ast \mc{F}_{\xi,\IC,\Lambda} &:=u_b^\ast i_b^\ast R\ol{\pi}_{\min,\ast}(\IC_{\Igs^\ast}(\mathcal{L}_{\xi})) \simeq R\tilde{\beta_b}^\dagger_\ast k_b^\ast \IC_{\Igs^\ast}(\mathcal{L}_{\xi}) \\
&\simeq R\tilde{\beta_b}_\ast\IC_{[\Ig^{b,\ast,\diamond}/\Jbul]}(\mathcal{L}_{\xi,b})[-d_b] \\
&\simeq \hat{V}_{\xi,\IC,\Lambda,b}[-d_b],
\end{align*}
where the first isomorphism follows from proper and smooth base-changes; the second isomorphism follows from Proposition~\ref{prop: overconvergentreplacement} and Lemma~\ref{lem: restrictionICtocapstrata}. The last isomorphism follows from Proposition~\ref{prop: CohoSmallDiamond}.
\end{proof}
This completes our computation of the stalks of $\mc{F}_{\xi,\IC,\Lambda}$. 
\subsubsection{Aside on completed versus smooth cohomology}
We let $\Lambda$ be as in setup \ref{assumption: coefficientsystemsingeneral} (2). Let $G/\bb{Q}_{p}$ be a connected reductive group and $\ast = \Spd \ol{\bb{F}}_p$. We recall from \cite[Section~2.2]{HansenSupercuspidalCohomology} that we have an exact symmetric monoidal functor 
\begin{equation}{\label{eqn: LCompletionFunctor}}
\gamma: \D(G,\Lambda) \ra \widehat{\D}(G,\Lambda) 
\end{equation}
\[ A \mapsto \lim_{n \geq 1} A \otimes_{\Lambda} \Lambda/\ell^{n},\]
which admits a right adjoint functor 
\begin{equation}{\label{eqn: SmoothVectorsofLcomplete}}
\delta: \widehat{\D}(G,\Lambda) \ra \D(G,\Lambda) 
\end{equation}
\[ A \mapsto \colim_{K_{p} \ra \{1\}} R\hat{\Gamma}(K_{p},A), \]
where $K_{p}$ runs over compact open pro-$p$ subgroups. Here $R\hat{\Gamma}(K_{p},A)$ for an object $A \in \widehat{\D}(G,\Lambda)$ represented by an inverse system $(A_{n})$ of smooth representations under the mod $\ell^{n}$-reduction maps is given by $\lim_{n \ra \infty} R\Gamma(K_{p},A_{n})$, the $\ell$-completed version of continuous cohomology (we recall that the limit here is derived). Alternatively, under the identification $\widehat{\D}(G(\bb{Q}_{p}),\Lambda) \simeq \Detale([\ast/G(\bb{Q}_{p})],\Lambda)$ (\cite[Proposition~2.6 (1)]{HansenSupercuspidalCohomology}), we have that 
\[ \delta(A) := \colim_{K_{p} \ra 1} a_{K_{p}*}b_{K_{p}}^{*}(A), \]
where $a_{K_{p}}$ and $b_{K_{p}}$ are as in the diagram 
\[[\ast/\underline{G(\bb{Q}_{p}})] \xleftarrow{b_{K_{p}}} [\ast/\underline{K_{p}}] \xrightarrow{a_{K_{p}}} \ast, \] 
and the functors are computed in the category of \'etale $\Lambda$-adic sheaves. If we write $\D(G,\Lambda)_{\Adm} \subset \D(G,\Lambda)$ for the full subcategory of sheaves such that its $K_{p}$-invariants for all compact open pro-$p$ subgroups have cohomology sheaves given by finitely generated $\Lambda$-modules then it follows from \cite[Proposition~2.6 (3)]{HansenSupercuspidalCohomology} that $\gamma$ defines a fully faithful embedding $\D(G,\Lambda)_{\Adm} \hookrightarrow \widehat{\D}(G,\Lambda)$. Moreover, its essential image is easily verified to be the full subcategory of $\ell$-complete representations, which are represented by an inverse system $(A_{n})$ over the mod $\ell^{n}$-reduction maps such that $A_{n}$ has $K_{p}$-invariants given by an admissible representation for all $n \geq 1$. We denote this subcategory by $\widehat{\D}(G,\Lambda)_{\Adm} \subset \widehat{\D}(G,\Lambda)$. In other words, this adjoint pair induces an exact equivalence
\begin{equation}{\label{eqn: AdmissibleComplexesLCompletetoSmooth}}
\gamma: \D(G,\Lambda)_{\Adm} \xrightarrow{\simeq} \widehat{\D}(G,\Lambda)_{\Adm}. 
\end{equation}
Now we specialize this discussion to the case that $G = J_{b}$ and the representations appearing in the stalks of $\mathcal{F}_{\xi,\IC,\Lambda}$.
\begin{lemma}{\label{lemma: smoothvectorsinthestalk}}
For $\Lambda$ as in Setup \ref{assumption: coefficientsystemsingeneral} (2), we have 
\[ \delta(\hat{V}_{\xi,\IC,\Lambda,b}) = V_{\xi,\IC,\Lambda,b} \]
and 
\[ \delta(\hat{V}^{D}_{\xi,\IC,\Lambda,b}) = V^{D}_{\xi,\IC,\Lambda,b}. \]
\end{lemma}
\begin{proof}
We explain the proof for $\hat{V}_{\xi,\IC,\Lambda,b}$ with the proof for $V^{D}_{\xi,\IC,\Lambda,b}$ being completely analogous. We look at the diagram
\[ 
\begin{tikzcd}
& \Ig^{b,\ast,\diamond}_{m} \arrow[dr] & & \\
\left[ \Ig^{b,\ast,\diamond}/\ul{K_{b}(p^{m})} \right]  \arrow[r] \arrow[ur,"q_{m}"] \arrow[d,"h_{m}"] & \left[\ast/\ul{K_{b}(p^{m})}\right] \arrow[d,"b_{m}"] \arrow[r,"a_{m}"] & \ast \\
\left[ \Ig^{b,\ast,\diamond}/\Jbul \right] \arrow[r,"\beta_{b}"] & \left[\ast/\ul{J_{b}(\bb{Q}_{p})}\right],& 
\end{tikzcd}
\]  
which extends (\ref{eqn: KeyDiagram}) with a cartesian square. In particular, by Proposition \ref{prop: CohoSmallDiamond}, the sheaf $\delta(\hat{V}_{\xi,\IC,\Lambda,b})$ identifies with
\[ \colim_{K_{p} \ra 1} a_{m*}b_{m}^{*}\beta_{b*}\IC_{[\Ig^{b,*}_{m}/\ul{J_{b}(\bb{Q}_{p})}]}(\mathcal{L}_{\xi,b,m})[d_b]. \]
However, by applying base change to the above diagram and arguing as in the proof of Proposition \ref{prop: CohoSmallDiamond} (We use Theorem \ref{thm: inverselimit} and that in Setup~\ref{assumption: coefficientsystemsingeneral}, $\Lambda$ is finite over $\mathbb{Z}_\ell$.) recalling that $\IC_{[\Ig^{b,*}_{m}/\ul{J_{b}(\bb{Q}_{p})}]}(\mathcal{L}_{\xi,b,m})$ is defined as an inverse limit of the mod $\ell^{n}$-IC-sheaves with integral coefficients, we deduce that this is isomorphic to 
\[ \colim_{K_{p} \ra 1} R\Gamma(\Ig^{b,*}_{m},\IC_{\Ig^{b,*}_{m}}(\mathcal{L}_{\xi,b,m}^{\alg}))[d_b], \]
which is precisely the smooth representation $V_{\xi,\IC,\Lambda,b}$, as desired. 
\end{proof}

\subsection{(Semi-)perversity}\label{sec: perversity} Let $\Lambda$ be a coefficient system as in Setup \ref{assumption: coefficientsystemsingeneral}. Fix an algebraically closed field $k/\mathbb{F}_p$ and consider $\Bun_G$ over $\ast:=\Spd k$. We recall that $\Detale(\Bun_{G},\Lambda)$ is equipped with a natural perverse $t$-structure, where, for $b \in B(G)$, with associated Harder--Narasimhan stratum $i_{b}: \Bun_{G}^{b} \hookrightarrow \Bun_{G}$, we say an object $A \in \phantom{}^{p}\Detale^{\leq 0}(\Bun_{G},\Lambda)$ (resp. $\phantom{}^{p}\Detale^{\geq 0}(\Bun_{G},\Lambda)$) if $i_{b}^{*}(A) \in \widehat{\D}^{\leq d_{b}}(J_{b}(\bb{Q}_{p}),\Lambda)$ (resp. $i_b^!A\in \widehat{\D}^{\geq d_{b}}(J_{b}(\bb{Q}_{p}),\Lambda)$) under the isomorphism $\Detale(\Bun_{G}^{b},\Lambda) \simeq \widehat{\D}(J_{b}(\bb{Q}_{p}),\Lambda)$. Here $d_{b} := \langle 2\rho_{G}, \nu_{b} \rangle$, where $\nu_{b}$ denotes the slope homomorphism of $b$. That this defines a $t$-structure follows from the formalism in \S\ref{sec: tStructureVstacks}, or see \cite[Proposition 8.1.5]{DvHKZ} for a proof of this specific case. We also recall that, for $b \in B(G,\mu^{-1})$, we have an equality $d_{b}= \dim(\Ig_k^{b,*})$. 

Let $\xi$ be an algebraic representation of $\mathsf{G}_{\ol{\bb{Q}}_{\ell}}$ defined over $\Lambda$,
$\mc{L}_\xi$ be the associated $\Lambda$-local system on $\Igs_{K^p,\et}$, and $\IC_{\Igs^\ast}(\mc{L}_\xi)$ be its intersection complex with respect to the Baily--Borel stratification, as introduced in \S\ref{sec: IConIgs}. We consider $\hat{\mc{F}}_{\xi,\IC,\Lambda} := R\overline{\pi}_{\min,\ast}\IC_{\Igs^\ast}(\mc{L}_\xi)\in \Detale(\Bun_G,\Lambda)$ and, when $\xi$ is the trivial representation, we simply denote this by $\hat{\mc{F}}_{\IC,\Lambda}$ as before. Since the image $\overline{\pi}_{\min,\ast}$ is given by the open substack $j_{\leq \mu}: \Bun_{G,\mu^{-1}} \hookrightarrow \Bun_{G}$, $\hat{\mc{F}}_{\xi,\IC,\Lambda}$ is supported on this substack. We will abuse notation and also write $\hat{\mc{F}}_{\xi,\IC,\Lambda}$ for its restriction.

\begin{theorem}{\label{thm: perversityofICIgs}}
For $\Lambda$ as in Setup \ref{assumption: coefficientsystemsingeneral}, we have that $\hat{\mc{F}}_{\xi,\IC,\Lambda}$ lies in $\phantom{}^{p, \leq 0}\Detale(\Bun_{G,\mu^{-1}},\Lambda)$ for the $t$-structure on $\Detale(\Bun_{G,\mu^{-1}},\Lambda)$ coming from the Harder--Narasimhan stratification. Moreover, the following is true. 
\begin{enumerate}
\item If $\Lambda$ is self-injective, then $\mathbb{D}_{\Bun_{G,\mu^{-1}}}(\hat{\mc{F}}_{\xi,\IC,\Lambda})\simeq \hat{\mc{F}}_{\xi^{\vee},\IC,\Lambda}$. In particular, if $\xi$ is the trivial representation then $\hat{\mc{F}}_{\IC,\Lambda} = \hat{\mathcal{F}}_{\IC,\xi,\Lambda}$ is Verdier self-dual, and $\hat{\mc{F}}_{\xi,\IC,\Lambda}$ is perverse, since $\bb{D}(\phantom{}^{p, \leq 0}\Detale(\Bun_{G,\mu^{-1}},\Lambda)) \subset \phantom{}^{p, \geq 0}\Detale(\Bun_{G,\mu^{-1}},\Lambda)$ by Remark \ref{rem: dualityinterchangingtstructures}.
\item If $\Lambda$ is regular, then $\hat{\mathcal{F}}_{\xi,\IC,\Lambda}$ is ULA with respect to $\mathrm{Bun}_{G,\mu^{-1}} \to *$.
\end{enumerate}
\end{theorem}
\begin{proof}
Assume first that $\Lambda$ is self-injective and hence finite in our situation. By Corollary \ref{cor: VerdierSelfDualityofICIgs} and the properness of the morphism $\overline{\pi}_{\min}$ over $\Bun_{G,\mu^{-1}}$, we know that $\hat{\mc{F}}_{\xi,\IC,\Lambda}$ is Verdier dual to $\hat{\mc{F}}_{\xi^{\vee},\IC,\Lambda}$, where $\xi^{\vee}$ denotes the representation dual to $\xi$.  It therefore suffices to show, for all $\xi$, that $\hat{\mc{F}}_{\xi,\IC,\Lambda}$ lies in ${}^p\Detale^{\leq 0}(\Bun_G,\Lambda)$. This can be checked on $*$-stalks. By Theorem~\ref{thm: stalks}, it is enough to show that $\hat{V}_{\xi,\IC,b}\in \widehat{\D}^{\leq 0}(J_b(\Q_p),\Lambda)$. Recall that we have
\[
\hat{V}_{\IC,\xi,\Lambda,b} = V_{\IC,\xi,\Lambda,b} = \mathrm{colim}_{m\geq 1} R\Gamma(\mathrm{Ig}^{b,*}_m, \IC_{\mathrm{Ig}^{b,*}_m}(\mathcal{L}_{\xi,b,m}^{\alg}))[d_b].
\]
Hence the desired result follows from combining the affineness of each $\mathrm{Ig}^{b,\ast}_m$, the perversity of the intermediate extension $\IC_{\mathrm{Ig}^{b,\ast}_m}(\mathcal{L}_{\xi,b,m}^{\alg})$ (for the usual perverse $t$-structure on finite type $\ol{\mathbb{F}}_p$-schemes), and the fact that pushforward along affine morphisms is perverse right $t$-exact, cf.~\cite[Theorem 4.1.1]{BBD} (Here, we have implicitly invoked the comparison of intermediate extensions on $\mathrm{Ig}^{b,*}_m$ established in Proposition \ref{prop: TwoDefinitionsofIntCohonperfectIgusavarieties}).

When $\Lambda$ is regular, then part (2) follows from Proposition \ref{prop: ICStratifiedULA} and the fact that ULA sheaves are preserved under proper pushforward \cite[Proposition~IV.2.11]{FSGeomLLC}. This in particular implies that the stalks $i_{b}^{*}\hat{\mathcal{F}}_{\xi,\IC,\Lambda} = \hat{V}_{\xi,\IC,\Lambda,b}[-d_{b}] \in \widehat{\D}(J_{b}(\bb{Q}_{p}),\Lambda)$ lie in the full subcategory $\widehat{\D}_{\Adm}(J_{b}(\bb{Q}_{p}),\Lambda)$, using the classification of ULA sheaves on $\Bun_{G}$ in the torsion case \cite[Theorem~V.7.1.]{FSGeomLLC}. Therefore, by applying Lemma \ref{lemma: smoothvectorsinthestalk} and the equivalence (\ref{eqn: AdmissibleComplexesLCompletetoSmooth}), we deduce that $\hat{V}_{\xi,\IC,\Lambda,b} = \gamma\delta(\hat{V}_{\xi,\IC,\Lambda,b}) = \gamma(V_{\xi,\IC,\Lambda,b})$. It is easy to see that $V_{\xi,\IC,\Lambda,b}[-d_{b}]$ is concentrated in degrees $\leq d_{b}$ by the same reasoning as above, and therefore we conclude the same for $\hat{V}_{\xi,\IC,\Lambda,b}$ by the exactness of $\gamma$. This gives us the semi-perversity statement in this case.
\end{proof}
\begin{remark}
One should not expect $\hat{\mathcal{F}}_{\xi,\IC,\Lambda}$ to be ULA without the assumption that $\Lambda$ is cohomologically regular. Indeed, the invariants under various compact open subgroups could fail to be perfect complexes if the intersection cohomology groups of Igusa varieties at finite level have torsion and $\Lambda$ is a ring like $\bb{Z}/\ell^{n}\bb{Z}$ for $n \geq 2$. As a result, the finiteness part of the ULA condition on $\hat{\mathcal{F}}_{\xi,\IC,\Lambda}$ could fail. 
\end{remark}

\section{Applications}\label{sec: applications}

\subsection{Intersection cohomology of Shimura varieties}

In this section, we give a formula for the intersection cohomology of Shimura varieties in terms of the Hecke action on the complex $\hat{\mc{F}}_{\xi, \IC,\Lambda}$ and various sheaves obtained from it. To establish the most general possible results, including working with both rational and integral coefficients, we need to use the lisse-\'etale sheaves introduced in \cite[Chapter VII]{FSGeomLLC}. We first review some of the basic structure of $\Dlis(\Bun_{G},\Lambda)$ in this level of generality. We fix $k=\ol{\bb{F}}_p$ and consider $\Bun_G$ over $\Spd k$.

\subsubsection{Recollection on lisse-\'etale sheaves}
We assume for the rest of this section that $\Lambda$ is a $\bb{Z}_{\ell}[\sqrt{p}]$-algebra for $\ell \neq p$ unless otherwise stated, considered as a relatively (over $\bb{Z}_\ell$) discrete condensed ring via $\Lambda:=\bb{Z}_{\ell}\otimes^{\bs}_{\bb{Z}_{\ell,\mathrm{disc}}}\Lambda_{\mathrm{disc}}$, see \cite[Chapter VII.6, 2nd paragraph]{FSGeomLLC}, where $\otimes^{\bs}$ denotes the solid tensor product. We recall that, for such a $\Lambda$, we have the category $\Dlis(\Bun_{G},\Lambda)$, the derived category of lisse-\'etale $\Lambda$-sheaves on the Artin $v$-stack $\Bun_{G}$ constructed in \cite[Definition~VII.6.1]{FSGeomLLC}. This is a full subcategory $\Dlis(\Bun_{G},\Lambda) \hookrightarrow \mathrm{D}_{\bs}(\Bun_{G},\Lambda)$ of solid $\Lambda$-sheaves on $\Bun_{G}$, as described in \cite[Definition~VII.1.10]{FSGeomLLC}. If $\Lambda$ is a torsion ring then we have an identification $\Dlis(\Bun_{G},\Lambda) \simeq \Detale(\Bun_{G},\Lambda)$ with the usual left-completed category of \'etale $\Lambda$-modules on $\Bun_{G}$ by \cite[Proposition~VII.6.6]{FSGeomLLC} combined with \cite[Proposition~V.3.5]{FSGeomLLC} and \cite[Theorem~V.3.7]{FSGeomLLC}.

\subsubsection{Hecke operators on \texorpdfstring{$\Dlis(\Bun_{G},\Lambda)$}{}}{\label{s: HeckeOperatorsandDlisse}}
We consider the $L$-group $\phantom{}^{L}G := \widehat{G} \rtimes W_{\bb{Q}_{p}}/\Lambda$ and write $\Rep(\phantom{}^{L}G)$ for the category of algebraic representations of the group $\phantom{}^{L}G := \widehat{G} \rtimes Q$, where $Q$ is the finite quotient of $W_{\bb{Q}_{p}}$ through which the action of $W_{\bb{Q}_{p}}$ on $\widehat{G}$ factors.

We recall that, for $V \in \Rep(\phantom{}^{L}G)$ an algebraic representation of the $L$-group over $\Lambda$, we have the Hecke operator 
\begin{equation}{\label{eqn: HeckeOperator}} T_{V}: \Dlis(\Bun_{G},\Lambda) \ra \Dlis(\Bun_{G},\Lambda)^{BW_{\mathbb{Q}_{p}}} \end{equation}
\[ A \mapsto h_{G\natural}^{\ra}(h_{G}^{\la*}(-) \otimes \mathrm{IC}'_{V}), \]
as constructed in \cite[Section~IX.2]{FSGeomLLC}. Here the morphisms are given by
\begin{equation}{\label{eqn: FullHeckeDiagram}}
 \Bun_{G} \xleftarrow{h^{\la}_{G}} \mathrm{Hck}_{G} \xrightarrow{h^{\ra}_{G}} \Bun_{G} \times \Div^{1}, 
\end{equation}
where $\Hck_{G}$ is the $v$-stack parameterizing, for $S \in \Perf$, degree 1 closed Cartier divisors on the relative Fargues-Fontaine curve $X_S$, and modifications of $G$-bundles $\mathcal{E}_{1} \dashrightarrow \mathcal{E}_{2}$ on $X_{S}$, which is meromorphic along the Cartier divisor determined by a map $S \ra \Div^{1}$ (Recall that $\Div^{1}$ is the moduli space of degree $1$ closed Cartier divisors on the Fargues-Fontaine curve). The maps $h_{G}^{\la}$ and $h_{G}^{\ra}$ remember the bundles $\mathcal{E}_{1}$ and $\mathcal{E}_{2}$, respectively. While $\IC'_{V} \in \D_{\bs}(\Hck_{G},\Lambda)$ is the solid sheaf attached to $V$ under the geometric Satake equivalence, which is pulled back from the local Hecke stack to the global one (see \cite[Section~IX.2]{FSGeomLLC}) and $h_{G\natural}^{\ra}$ is the natural pushforward in the category of solid sheaves (the left adjoint to $*$-pullback), as defined in \cite[Section~VII.3]{FSGeomLLC}. 

More precisely, given a small Artin $v$-stack $X$ over a diamond $S$ (in the sense of \cite[Section~IV.1]{FSGeomLLC}), we have a natural fully faithful embedding 
\begin{equation}{\label{eqn: dualembedding}}
 \Detale^{\ULA}(X/S,\bb{Z}_{\ell}[\sqrt{p}]) \hookrightarrow \mathrm{D}_{\bs}(X,\bb{Z}_{\ell}[\sqrt{p}]),  \\
 A \mapsto \bb{D}_{X}(A)^{\vee}, 
\end{equation}
where $(-)^{\vee} := R\mathcal{H}\mathrm{om}_{\mathcal{D}_{\bs}(X,\Lambda)}(-,\Lambda)$ is the dual embedding, as defined in \cite[Section~VII.5]{FSGeomLLC}. We recall that the category appearing on the RHS is the homotopy category of the inverse limit (as $\infty$-categories) $\cDetale^{\ULA}(X/S,\bb{Z}_{\ell}[\sqrt{p}]) := \lim_{n} \cDetale^{\ULA}(X/S,\bb{Z}_{\ell}[\sqrt{p}]/\ell^{n})$, where $\cDetale^{\ULA}(X/S,\bb{Z}_{\ell}[\sqrt{p}]/\ell^{n})$ is the $\infty$-category of \'etale $\bb{Z}_{\ell}[\sqrt{p}]/\ell^{n}$-sheaves which are ULA over $S$. We write $\mc{H}\mathrm{ck}_{G}$ for the local Hecke stack of $G$, which is equipped with a natural map $\mathrm{Hck}_{G} \ra \mc{H}\mathrm{ck}_{G}$ from the global Hecke stack, given by formally completing the modification of vector bundles along the divisor. For $V \in \Rep(\phantom{}^{L}G)$, we may first look at the sheaf in $\Detale^{\ULA}(\mc{H}ck_{G}/[\Div^{1}/L^{+}G],\bb{Z}_{\ell}[\sqrt{p}])$ attached to $V$ under geometric Satake and then apply the functor $A\mapsto \bb{D}(A)^\vee$ to obtain a sheaf in $\D_{\bs}(\mc{H}ck_{G},\bb{Z}_{\ell}[\sqrt{p}])$, where the Verdier duality functor is with respect to $\mc{H}ck_{G}\to[\Div^{1}/L^{+}G]$. We write $\IC_{V} \in \Detale(\Hck_{G},\bb{Z}_{\ell}[\sqrt{p}])$ and $\IC'_{V} \in \D_{\bs}(\Hck_{G},\bb{Z}_{\ell}[\sqrt{p}])$ for the sheaves produced by this procedure and pulling back along $\mathrm{Hck}_{G} \ra \mc{H}\mathrm{ck}_{G}$ (see the beginning of \cite[Section~IX.2]{FSGeomLLC}).

We note it follows by \cite[Proposition~IX.2.1]{FSGeomLLC} that the Hecke operator factors through a natural map $\Dlis(\Bun_{G},\Lambda)^{BW_{\bb{Q}_{p}}} \ra \mathrm{D}_{\bs}(\Bun_{G} \times \Div^{1})$ and this gives rise to the operator described in (\ref{eqn: HeckeOperator}).

We assume that $G$ is quasi-split (e.g it is unramified, as in Assumption \ref{assumption:codimension}) for simplicity. We fix a maximal torus and Borel $T \subset B \subset G$, and let $\bb{X}_{*}(T_{\ol{\mathbb{Q}}_{p}})^{+}$ be the set of geometric dominant cocharacters. For $\mu \in \bb{X}_{*}(T_{\ol{\bb{Q}}_{p}})^{+}$ a geometric dominant cocharacter of $G$, we write $\mu^{\Gamma} \in \bb{X}_{*}(T_{\ol{\mathbb{Q}}_{p}})^{+}/\Gamma$ for the corresponding Galois orbit.  Attached to $\mu$, we can define an associated highest weight tilting representation of $\widehat{G}$, denoted $\mathcal{T}_{\mu}$, and a representation $\mathcal{T}_{\mu^{\Gamma}}$ of $\phantom{}^{L}G$ attached to the Galois orbit, given by extending $\mathcal{T}_{\mu}$ to a representation of $\widehat{G} \rtimes W_{E_{\mu}}$ as in \cite[Lemma~2.1.2]{KottwitzTOO}, and then inducing it along the embedding $W_{E_{\mu}} \subset W_{\mathbb{Q}_{p}}$ (see \cite[Section~9.1]{HamGeomES}). We then have a natural Hecke operator
\[ T_{\mu^{\Gamma}}: \Dlis(\Bun_{G},\Lambda) \ra \Dlis(\Bun_{G},\Lambda)^{BW_{\mathbb{Q}_{p}}} \]
attached to the cocharacter $\mu$. 

We will be interested in a slight variant of the operation $T_{\mu^{\Gamma}}$; in particular, if $E_{\mu}$ denotes the reflex field of $\bb{Q}_{p}$ then the operation $T_{\mu^{\Gamma}}$ factors over the natural functor 
\[ \Dlis(\Bun_{G},\Lambda)^{BW_{E_{\mu}}} \ra \Dlis(\Bun_{G},\Lambda)^{BW_{\mathbb{Q}_{p}}} \] 
induced by induction from $W_{E_{\mu}}$ to $W_{\bb{Q}_{p}}$ (See the discussion on \cite[Page~322]{FSGeomLLC}). We write
\[ T_{\mu}: \Dlis(\Bun_{G},\Lambda) \ra \Dlis(\Bun_{G},\Lambda)^{BW_{E_{\mu}}} \]
for the  resulting natural operator attached to $\mathcal{T}_{\mu}$. It has a similar formula as \eqref{eqn: HeckeOperator}, with $\IC_\mu'$, the lisse version of the intermediate extension for the Schubert cell labeled by $\mu$ on $\mc{H}ck_G$ (restricted to $\Hck_G$), being the kernel of the correspondence. When we want to emphasize the coefficients, we denote this by $T_{\mu}^{\Lambda}$.

Similarly, if $\Lambda$ is as in Setup \ref{assumption: coefficientsystemsingeneral} (2) then we will also consider the functor
\begin{equation}{\label{eqn: LCompleteHeckeOperator}}
\hat{T}_{\mu}: \Detale(\Bun_{G},\Lambda) \ra \Detale(\Bun_{G},\Lambda)^{BW_{E_{\mu}}} 
\end{equation}
\[ A \mapsto h_{G!}^{\ra}(h_{G}^{\la*}(-) \otimes \IC_\mu) \]
defined by the Hecke operator on the $\ell$-complete category.

For $C$ an algebraically closed perfectoid field in characteristic $0$, under the equivalence $\Dlis(\Bun_{G},\Lambda) \simeq \Dlis(\Bun_{G,C},\Lambda)$ given by \cite[Corollary~V.2.3]{FSGeomLLC}, the composite of the functor $T_{\mu}$ with the forgetful functor 
\[ \Dlis(\Bun_{G},\Lambda)^{BW_{E_{\mu}}} \ra \Dlis(\Bun_{G},\Lambda) \]
identifies with the functor
\[ \phantom{}_{x}h^{\ra}_{G\natural}(\phantom{}_{x}h^{\la*}_{G}(-) \otimes \IC'_{\mu}) \] 
given by the analogous pull and push along the diagram
\begin{equation}{\label{eqn: basechangedHeckeDiagram}}
\Bun_{G} \xleftarrow{\phantom{}_{x}h^{\la}_{G}} \phantom{}_{x}\mathrm{Hck}_{G} \xrightarrow{\phantom{}_{x}h^{\ra}_{G}} \Bun_{G,C}, 
\end{equation}
obtained by base-changing $h^{\ra}_{G}$ along the natural map $x: \Spd C \ra \Div^{1}$ corresponding to the degree $1$ Cartier divisor in $X_{C^{\flat}}$ defined by the characteristic $0$ untilt of $C^{\flat}$ given by $C$. Here the sheaf $\IC'_{\mu}$ is the analogous solid sheaf attached to $\mathcal{T}_{\mu} \in \Rep(\widehat{G})$ described above by applying geometric Satake to get a sheaf $\IC_{\mu}$ and then applying the embedding (\ref{eqn: dualembedding}).

As we will be interested in studying these operators in the context of Shimura varieties, we now assume for simplicity that $\mu$ is minuscule, and consider the pullback of the Hecke correspondence (\ref{eqn: basechangedHeckeDiagram}) along the map $h_{G}^{\ra}$ to the open substack $i_{1}: \Bun_{G}^{1} \hookrightarrow \Bun_{G}$ corresponding to the trivial $G$-bundle. The Hecke correspondence pulled back to this locus is easily identified with the $v$-stack quotient $[\Fl/\underline{G(\bb{Q}_{p})}]$ of the flag variety by the action of the group $v$-sheaf $\ul{G(\bb{Q}_{p})}$. In particular, after rigidifying the $G(\bb{Q}_{p})$-action, this will parameterize minuscule modifications $\mathcal{E} \dashrightarrow \mathcal{E}_{0}$ to the trivial $G$-bundle, via the Bialynicki-Birula isomorphism \cite[Proposition~3.4.3,Theorem~3.4.5]{CS17}. We therefore have the cartesian diagram
\begin{equation}{\label{eqn: CartesianDiagramForFlagVariety}}
\begin{tikzcd}
\left[\mathcal{F}\ell_{G,\mu^{-1},C}/\ul{G(\mathbb{Q}_{p})}\right] \arrow[r,"h_1"] \arrow[d] & \left[\Spd C/\ul{G(\bb{Q}_{p})}\right] \arrow[d,"i_{1}"] & \\
\phantom{}_{x}\Hck_{G, \leq \mu} \arrow[r,"h_{G}^{\ra}"] & \Bun_{G,C}. & 
\end{tikzcd}
\end{equation}
If $\mu$ is minuscule then one has an isomorphism $V^{\mu} = \mathcal{T}_{\mu} = V_{\mu}$ between the standard, tilting, and costandard highest weight representations. The sheaf $\IC_{\mu}$ is simply identified with the constant sheaf $\Lambda[d](\frac{d}{2})$ on $[\Fl_{,C}/\underline{G(\bb{Q}_{p})}]$, where $d = \langle 2\rho_{G},\mu \rangle$ for $\rho_{G}$ the half sum of the positive roots of $G$. If we write $h_{2}: [\mathcal{F}\ell_{G,\mu^{-1},C}/\ul{G(\bb{Q}_{p}})] \ra \Bun_{G}$ for the map induced by $h^{\la}_{G}$ and $h_{1}: [\mathcal{F}\ell_{G,\mu^{-1},C}/\ul{G(\bb{Q}_{p}})] \ra [\Spd C/\ul{G(\bb{Q}_{p})}]$ for the map induced by $h^{\ra}_{G}$ then this implies the following.
\begin{lemma}{\label{lemma: computingminusculeHeckecorrespondences}}
For $\Lambda/\bb{Z}_{\ell}[\sqrt{p}]$, we have an isomorphism
\[ i_{1}^{*}T_{\mu}(-) \simeq h_{1\natural}h^{*}_{2}(-)[-d](-\frac{d}{2}),  \]
in $\D(G(\bb{Q}_{p}),\Lambda)^{BW_{E_{\mu}}}$
where the Tate twist on the RHS means we twist the Weil group action by $|\cdot|^{-d/2}$ for $|\cdot|$ the norm character of $W_{E_{\mu}}$, using the chosen square root of $p$. Moreover, if $\Lambda/\bb{Z}_{\ell}[\sqrt{p}]$ is a torsion algebra then under the equivalence $\mathrm{D}_{\mathrm{lis}}(\Bun_{G},\Lambda) \simeq \Detale(\Bun_{G},\Lambda)$ given by the naive embedding, we have an identification 
\[ i_{1}^{*}T_{\mu}(-) \simeq  h_{1!}h_{2}^{*}(-)[d](\frac{d}{2}) \simeq h_{1*}h_{2}^{*}(-)[d](\frac{d}{2}) \]
\end{lemma}
\begin{proof}
The first part follows by applying the base-change result \cite[Proposition~VII.3.1 (iii)]{FSGeomLLC} to the cartesian diagram (\ref{eqn: CartesianDiagramForFlagVariety}) and noting that the sheaf $\IC_{V_{\mu}}$ identifies under the functor $A\mapsto \bb{D}(A)^\vee$ with $\Lambda[-d](-\frac{d}{2})$.

It remains to justify why one can write $h_{1*}$ instead of $h_{1\natural}$ in the torsion case. As the map $h_1$ is $\ell$-cohomologically smooth of dimension $d$, the left adjoint $h_{1\natural}$  to $*$-pullback satisfies 
\[h_{1\natural}\simeq h_{1!}[2d](d).\] Finally, we use that $h_1$ is proper (its fibers are given by $\mathcal{F}\ell_{G,\mu}$) to write $h_{1*}$ instead of $h_{1!}$. 
\end{proof}

\noindent This lemma will allow us to write the cohomology of Shimura varieties and its boundary strata in terms of Hecke operators applied to certain sheaves on $\Bun_{G}$. It will also be important for us to study the perversity properties of these sheaves. We describe the relevant piece of structure now. 

\subsubsection{The semi-orthogonal decomposition and the perverse $t$-structure}{\label{s: SemiOrthogonalDecompositionPerverseTStruture}}
By applying excision with respect to the locally closed stratification $i_{b}: \Bun_{G}^{b} \hookrightarrow \Bun_{G}$, we recall that $\Dlis(\Bun_{G},\Lambda)$ has a semi-orthogonal decomposition with graded pieces given by $\Dlis(\Bun_{G}^{b},\Lambda)$, the category of lisse \'etale sheaves on the strata $\Bun_{G}^{b}$. By \cite[Proposition~VII.7.1]{FSGeomLLC}, this identifies with $\mathrm{D}(J_{b}(\bb{Q}_{p}),\Lambda)$, the derived category of smooth $J_b(\qp)$-representations on discrete $\Lambda$-modules. If $\Lambda$ is torsion then this follows from the identification $\Dlis(\Bun_{G},\Lambda) \simeq \Detale(\Bun_{G},\Lambda)$  with usual \'etale sheaves, which is equipped with excision triangles coming from this sheaf theory. In general however, $!$-pushforward is not defined in the context of lisse-\'etale sheaves; however, in the particular case of $\Dlis(\Bun_{G},\Lambda)$, where only locally closed immersions are relevant, one can construct the functors $i_{b!}: \Dlis(\Bun_{G}^{b},\Lambda) \ra \Dlis(\Bun_{G},\Lambda)$ and show that there is excision (See the discussion on \cite[Page~11-12]{ImaConv}, cf. \cite[Proposition VII.6.7]{FSGeomLLC}). In particular, with respect to inclusions of locally closed strata of $\Bun_{G}$, we have the full six operations satisfying the usual compatibilities.

Recall that we also have the full subcategory $\Dlis^{\ULA}(\Bun_{G},\Lambda) \subset \Dlis(\Bun_{G},\Lambda)$ of objects which are ULA over $\Spd k$, as defined in \cite[Definition~VII~7.8]{FSGeomLLC}. By~\cite[Proposition~VII.7.9]{FSGeomLLC}, this can be described explicitly using the semi-orthogonal decomposition as the subcategory of complexes $A \in \Dlis(\Bun_{G},\Lambda)$ such that, for each $b\in B(G)$, the pullback $i_{b}^{*}A \in \Dlis(\Bun_{G}^{b},\Lambda) \simeq \mathrm{D}(J_{b}(\bb{Q}_{p}),\Lambda)$ lies in the full subcategory $\D_{\PerfAdm}(J_{b}(\bb{Q}_{p}),\Lambda)$ of admissible objects. More precisely, $M \in \mathrm{D}(J_{b}(\bb{Q}_{p}),\Lambda)$ lies in $\D_{\PerfAdm}(J_{b}(\bb{Q}_{p}),\Lambda)$ if $M^K$ is a perfect complex of $\Lambda$-modules for any compact open pro-$p$ subgroup $K\subset J_b(\bb{Q}_p)$. Here we use ``$\mathrm{Pf}$'' in the subscript for ``perfect'', to distinguish it from the condition that all cohomological groups of $M^K$ are finitely generated $\Lambda$-modules.

Moreover, $\Dlis(\Bun_{G},\Lambda)$ is equipped with a Verdier duality operation $\bb{D}_{\Bun_{G}}$ given by $R\mathcal{H}\mathrm{om}_{\Dlis(\Bun_{G},\Lambda)}(-,\Lambda)$, since the dualizing sheaf on $\Bun_{G}$ is the constant sheaf (see \cite[Proposition~3.1]{HamannImai}). We will denote this by $\bb{D}_{\Bun_{G},\Lambda}$ when we want to emphasize the coefficients.

The existence of this semi-orthogonal decomposition allows us to define a perverse $t$-structure, by arguing just as in Definition/Proposition \ref{defn: tstructuresonIgusaStacks}.  
\begin{definition/proposition}{\label{defn: existenceofperversetstructure}}
There exists a $t$-structure $(\phantom{}^{p}\Dlis^{\leq 0}(\Bun_{G},\Lambda),\phantom{}^{p}\Dlis^{\geq 0}(\Bun_{G},\Lambda))$ on $\Dlis(\Bun_G,\Lambda)$ characterized by the following property:
\begin{itemize}
\item $A \in \phantom{}^{p}\Dlis^{\leq 0}(\Bun_{G},\Lambda)$ (resp. $A\in \phantom{}^{p}\Dlis^{\geq 0}(\Bun_{G},\Lambda))$) if and only if for all $b \in B(G)$ we have that $i_{b}^{*}A \in \D^{\leq d_{b}}(J_{b}(\bb{Q}_{p}),\Lambda)$ (resp. $i_{b}^{!}A \in \D^{\geq d_{b}}(J_{b}(\bb{Q}_{p}),\Lambda)$) where $d_{b} := \langle 2\rho_{G},\nu_{b} \rangle$ for $\nu_{b}$ the slope homomorphism of $b$ and $\rho_{G}$ the half sum of the positive roots.
\end{itemize}
We denote by $\Perv_{\mathrm{lis}}(\Bun_{G},\Lambda)$ the heart of this $t$-structure. 
\end{definition/proposition}
We endow $B(G)$ with the order topology as in \S~\ref{ss: GenericSpecialNP}. For a locally closed subset $V \subset B(G)$, we define $\Bun_{G}^V:= \Bun_{G}\times_{\ul{B(G)}}\ul{V}$. It is easy to see that the semi-orthogonal decomposition, perverse $t$-structure, and characterization of ULA sheaves pass analogously to the subcategory $j_{V!}: \Dlis(\Bun_{G}^{V},\Lambda) \hookrightarrow \Dlis(\Bun_{G},\Lambda)$. 

\subsubsection{Intersection cohomology via Hecke operators}
We start with the case of torsion coefficients, so we first let $\Lambda$ be a $\mathbb{Z}_{\ell}[\sqrt{p}]$-algebra as in Setup~\ref{assumption: coefficientsystemsingeneral} (1) and $\xi$ an algebraic representation of $\mathsf{G}_{\ol{\mathbb{Q}}_{\ell}}$ defined over $\Lambda$. We let $\mathcal{L}_{\xi}^{\mathrm{alg}}$ be the associated \'etale $\Lambda$-local system on the algebraic Shimura variety with coefficients in $\Lambda$, as defined in \cite{PinkHigherDirectImages}. We define the complex 
\begin{equation}{\label{eqn: InfiniteLevelIntersectionCohomology}}
\IC(\mathsf{G},\mathsf{X},\xi)_{K^p, \Lambda}:=\varinjlim_{K_{p} \ra \{1\}} R\Gamma\left(\Shstar\gx_{K^pK_p,\ol{\mathsf{E}}},\IC_{\Shstar_{K_{p}}}(\mathcal{L}_{\xi}^{\mathrm{alg}}[d](\frac{d}{2}))\right)
\end{equation}
of smooth $G(\mathbb{Q}_p)$-representations with admissible cohomology. It is equipped with a continuous action of $W_E$ and an action of an away-from-$p$ Hecke algebra. We omit $\xi$ from the notation if it is the trivial representation.

\begin{theorem}\label{thm: Hecke operator torsion}
For $\Lambda$ as in Setup~\ref{assumption: coefficientsystemsingeneral} (1), we have an identification
\[ 
i_{1}^{*}T_{\mu}\hat{\mathcal{F}}_{\xi,\IC,\Lambda} \simeq \IC(\mathsf{G},\mathsf{X},\xi)_{K^p, \Lambda}
\]
of $G(\bb{Q}_{p}) \times W_{E}$-representations, compatibly with the away-from-$p$ Hecke action. 
\end{theorem}

\begin{proof} We combine Lemma~\ref{lemma: computingminusculeHeckecorrespondences} with proper base change applied to the diagram~\eqref{eqn: CartesiandiagramShimuraVarietyInfiniteLevel}. This gives a Hecke- and $G(\mathbb{Q}_p)$-equivariant isomorphism
\[
i_1^*T_{\mu}\hat{\mathcal{F}}_{\xi,\IC,\Lambda} \simeq R\Gamma(\mathcal{S}^*_{K^p}, \tilde{h}_{\min}^*\IC_{\mathrm{Igs}^*}(\mathcal{L}_{\xi})). 
\]
We then apply Corollary~\ref{cor: InfiniteLevelisColimitIC} and commute the colimit with taking cohomology using the fact that the transition morphisms in the tower of Shimura varieties are qcqs, cf. Stacks Project \cite[Tag~073E]{stacks-project}. For the $W_{E}$-action, we argue as in \cite[Theorem 8.4.10]{DvHKZ}.
\end{proof}
We now extend our result to the coefficient systems of Setup~\ref{assumption: coefficients}. Case (1) is already done. For Cases (2)-(5), we need some preparations. First take $\Lambda = \mathcal{O}_{F}$, as in case (2) of \ref{assumption: coefficients}. We note that there is a natural $\ell$-adic completion functor
\[ 
\mc{D}_{\mathrm{lis}}(\Bun_{G},\mathcal{O}_{F}) \ra \varprojlim_{n} \mc{D}_{\mathrm{lis}}(\Bun_{G},\mathcal{O}_{F}/\ell^{n}) \simeq \varprojlim_{n} \cDetale(\Bun_{G},\mathcal{O}_{F}/\ell^{n}) =: \cDetale(\Bun_{G},\mathcal{O}_{F})
\] 
\[ A \mapsto \varprojlim_{n} A \otimes_{\mathcal{O}_{F}} \mathcal{O}_{F}/\ell^{n},\]
where we identify the lisse and \'etale categories for torsion coefficients. We write 
\begin{equation}{\label{eqn: lcompletinglissesheaves}}
\gamma_{\Bun_{G}}: \mathrm{D}_{\mathrm{lis}}(\Bun_{G},\mathcal{O}_{F}) \ra \Detale(\Bun_{G},\mathcal{O}_{F})
\end{equation}
for the resulting functor on the homotopy category. For a locally closed subset $V \subset B(G)$, we have the locally closed substack $\Bun_{G}^{V}$ as before. We have an analogously defined functor
\[ \gamma_{V}: \mathrm{D}_{\mathrm{lis}}(\Bun_{G}^{V},\mathcal{O}_{F}) \ra \Detale(\Bun_{G}^{V},\mathcal{O}_{F}).\]
For an inclusion $V_{1} \subset V_{2}$, we write $i_{V_{1}V_{2}}: \Bun_{G}^{V_{1}} \ra \Bun_{G}^{V_{2}}$ for the induced inclusion. We will compare the semi-orthogonal decomposition on the lisse and the $\ell$-complete categories under these functors.

Write $\Convex_{B(G)}$ for the category of convex subsets (in the sense that they can be written as the difference of two closed subsets) of $B(G)$ with morphisms given by inclusion. Then the semi-orthogonal decompositions on $\Detale(\Bun_{G},\mathcal{O}_{F})$ and $\Dlis(\Bun_{G},\mathcal{O}_{F})$ can be viewed as functors
\begin{equation}{\label{eqn: DetaleSemiOrthogonalDecomp}}
\bb{S}^{\mathcal{O}_{F}}_{B(G),\et}: \Corr(\Convex_{B(G)},\mathrm{all}) \ra \LinCat_{\CO_F}^{\mathrm{sm}}
\end{equation}
\[ \{W_{1} \la W_{2} \ra W_{3}\} \ra \{ i_{W_{2}W_{3}!}i_{W_{2}W_{1}}^{*}: \cDetale(\Bun_{G}^{W_{1}},\mathcal{O}_{F}) \ra \cDetale(\Bun_{G}^{W_{3}},\mathcal{O}_{F})\} \]
and 
\begin{equation}{\label{eqn: DlisSemiOrthogonalDecomp}}
 \bb{S}^{\mathcal{O}_{F}}_{B(G),\lis}: \Corr(\Convex_{B(G)},\mathrm{all}) \ra \LinCat^{\mathrm{sm}}_{\CO_F}
\end{equation}
\[ \{W_{1} \la W_{2} \ra W_{3}\} \ra \{ i_{W_{2}W_{3}\lis!}i_{W_{2}W_{1}\lis}^{*}: \cDlis(\Bun_{G}^{W_{1}},\mathcal{O}_{F}) \ra \cDlis(\Bun_{G}^{W_{3}},\mathcal{O}_{F})\} \]
satisfying certain properties, see \cite[Definition~6.2.1]{GHILZIsocComparison}, with notation as in \emph{loc.cit}. Similarly, for an open subset $U \subset B(G)$, we write $\mathbb{S}^{\mathcal{O}_{F}}_{U,\et}$ and $\mathbb{S}^{\mathcal{O}_{F}}_{U,\lis}$ for the analogous semi-orthogonal decomposition on $\Detale(\Bun_{G}^{U},\mathcal{O}_{F})$ and $\Dlis(\Bun_{G}^U,\mathcal{O}_{F})$, respectively. In fact, using the characterization of ULA objects in $\cDetale(\Bun_{G}^{U},\Lambda)$ and $\cDlis(\Bun_{G}^{U},\Lambda)$ in terms of $*$-pullbacks to HN-strata, it is clear that these even restrict to semi-orthogonal decompositions on the ULA subcategories $\cDlis^{\ULA}(\Bun_{G}^{U},\Lambda)$ and $\cDetale^{\ULA}(\Bun_{G}^{U},\Lambda)$, which we denote by $\bb{S}_{U,\lis}^{\mathcal{O}_{F},\ULA}$ and $\bb{S}_{U,\et}^{\mathcal{O}_{F},\ULA}$.\footnote{In the current version of \cite{GHILZIsocComparison}, the definition of semi-orthogonal decomposition only allows for values in the presentable categories (which notably excludes the subcategory of ULA sheaves). The point being that the $\infty$-categorical adjoint functor theorem then automatically guarantees the existence of the right adjoint functors. However, one can also encode these adjoints as additional data (as in \cite[Definition~3.2.1]{HeyerMann}) and allow for values in stable $\infty$-categories. Much of the discussion still goes through.} 

Since $\Bun_{G}$ is an increasing union of $\Bun_{G}^{U}$ for $U \subset B(G)$ quasi-compact open, by excision we have that $*$-pullback induces equivalences
\begin{equation}{\label{eqn: lissesubstackequivalence}}
\mc{D}_{\mathrm{lis}}(\Bun_{G},\mathcal{O}_{F}) \xrightarrow{\simeq} \lim_{U} \mc{D}_{\mathrm{lis}}(\Bun_{G}^{U},\mathcal{O}_{F}) 
\end{equation}
and 
\begin{equation}{\label{eqn: etalesubstackequivalence}}
 \cDetale(\Bun_{G},\mathcal{O}_{F}) \xrightarrow{\simeq} \lim_{U} \cDetale(\Bun_{G}^{U},\mathcal{O}_{F}). 
\end{equation}
and similarly on the ULA subcategories. It is easy to check that, under the identifications 
\[\widehat{\D}(J_{b}(\bb{Q}_{p}),\mathcal{O}_{F}) \simeq \Detale(\Bun_{G}^{b},\mathcal{O}_{F}),\,\,\D(J_{b}(\bb{Q}_{p}),\mathcal{O}_{F}) \simeq \Dlis(\Bun_{G}^{b},\mathcal{O}_{F}),\] the functor  $\gamma_{b}$ is given by the natural map
\begin{equation}{\label{eqn: lcompletioninJb}}
 \D(J_{b}(\bb{Q}_{p}),\mathcal{O}_{F}) \ra \widehat{\D}(J_{b}(\bb{Q}_{p}),\mathcal{O}_{F}) 
\end{equation}
from \eqref{eqn: LCompletionFunctor} applied to $G = J_{b}$. We will identify $\gamma_{b}$ with this functor in what follows.

Let $\D_{\PerfAdm}(J_{b}(\bb{Q}_{p}),\mathcal{O}_{F}) \subset \D(J_{b}(\bb{Q}_{p}),\mathcal{O}_{F})$  (resp. $\widehat{\D}_{\PerfAdm}(J_{b}(\bb{Q}_{p}),\mathcal{O}_{F}) \subset \widehat{\D}(J_{b}(\bb{Q}_{p}),\mathcal{O}_{F})$) be the subcategory of (perfect) admissible objects (resp. the inverse limit of perfect admissible objects) as before. We note that the functor (\ref{eqn: lcompletioninJb}) restricts to an equivalence
\begin{equation}{\label{eqn: lcompletioninJbperfadmissible}}
 \gamma_{b}^{\ULA}: \D_{\PerfAdm}(J_{b}(\bb{Q}_{p}),\mathcal{O}_{F}) \xrightarrow{\simeq} \widehat{\D}_{\PerfAdm}(J_{b}(\bb{Q}_{p}),\mathcal{O}_{F}),
\end{equation}
which agrees with the restriction of (\ref{eqn: AdmissibleComplexesLCompletetoSmooth}) to the perfect admissible objects. We note that $\D_{\PerfAdm}(J_{b}(\bb{Q}_{p}),\mathcal{O}_{F})$ (resp. $\widehat{\D}_{\PerfAdm}(J_{b}(\bb{Q}_{p}),\mathcal{O}_{F})$) are precisely the graded pieces of the semi-orthogonal decompositions on $\Dlis^{\ULA}(\Bun_{G},\mathcal{O}_{F})$ (resp. $\Detale^{\ULA}(\Bun_{G},\mathcal{O}_{F})$). In particular, we can verify the following.
\begin{lemma}{\label{lemma: FromLCompletetoSmooth}}
The following is true. 
\begin{enumerate}
\item For $V_{1} \subset V_{2}$ an inclusion of locally closed subsets of $B(G)$, we have natural isomorphisms $\gamma_{V_{2}} i_{V_{1}V_{2}\lis!} \simeq i_{V_{1}V_{2}!}\gamma_{V_{1}}$ and $\gamma_{V_{1}} i_{V_{1}V_{2}\lis}^{*} \simeq i_{V_{1}V_{2}}^{*}\gamma_{V_{2}}$. In particular, under the equivalences (\ref{eqn: lissesubstackequivalence}) and (\ref{eqn: etalesubstackequivalence}), the functor $\gamma$ identifies with $\lim_{U} \gamma_{U}$, where $U \subset B(G)$ ranges over quasi-compact opens, and the functor $\gamma_{\Bun_{G}}$ upgrades to a natural transformation of semi-orthogonal decompositions $\bb{S}_{B(G),\lis}^{\mathcal{O}_{F}} \ra \bb{S}_{B(G),\et}^{\mathcal{O}_{F}}$. 
\item  The functor $\gamma_{\Bun_{G}}$ restricts to an equivalence
\begin{equation}{\label{eqn: EquivalenceonULASheaves}} 
\gamma_{\Bun_{G}}^{\ULA}: \Dlis^{\ULA}(\Bun_{G},\mathcal{O}_{F}) \xrightarrow{\simeq} \Detale^{\ULA}(\Bun_{G},\mathcal{O}_{F})
\end{equation}
on the full subcategories of ULA objects. Moreover, it intertwines Verdier duality on the category of lisse-\'etale ULA sheaves with Verdier duality on the category of $\ell$-complete \'etale ULA sheaves.
\item The equivalence $\gamma^{\ULA}_{\Bun_{G}}$ is Hecke-equivariant. In other words, for all $\mu \in \bb{X}_{*}(T_{\ol{\bb{Q}}_{p}})^{+}$ we have a natural equivalence 
\[ \hat{T}_{\mu}\gamma^{\ULA}_{\Bun_{G}} \simeq \gamma^{\ULA}_{\Bun_{G}}T_{\mu} \]
of functors $\Dlis^{\ULA}(\Bun_{G},\mathcal{O}_{F}) \ra \Detale^{\ULA}(\Bun_{G},\mathcal{O}_{F})^{BW_{E_{\mu}}}$.
\end{enumerate}
\end{lemma}
\begin{proof}
Point (1) for $*$-pullback follows from the fact that $*$-pullback is a symmetric monoidal and that $i_{V_{1}V_{2}}^{*}$ commutes with limits (since it admits an exceptional left adjoint, see \cite[Proposition~VII.7.2]{FSGeomLLC} and \cite[Proposition~9.5]{GHILZIsocComparison}). For $!$-pushforward, this follows from projection formula and the fact that $i_{V_{1}V_{2}!}$ commutes with limits (since it admits an exceptional left adjoint (See \cite[Proposition~9.5]{GHILZIsocComparison})\footnote{Strictly speaking, the proof is only for torsion coefficients. However, the proof with rational coefficients is exactly the same.}). 

Combining (1), excision, and the fact that $\gamma_{b}$ sends perfect admissible objects to perfect admissible complexes, we see that $\gamma_{U}$ induces a functor 
\[ \gamma_{U}^{\ULA}: \cDlis^{\ULA}(\Bun_{G}^{U},\mathcal{O}_{F}) \ra \cDetale^{\ULA}(\Bun_{G}^{U},\mathcal{O}_{F}) \]
on the ULA subcategories. We may now appeal to \cite[Proposition~6.4.4]{GHILZIsocComparison} and \cite[Proposition~6.4.5]{GHILZIsocComparison} with respect to the semi-orthogonal decompositions $\bb{S}^{\mathcal{O}_{F}}_{U,\et}$ and $\bb{S}^{\mathcal{O}_{F}}_{U,\lis}$, for varying quasi-compact open $U \subset B(G)$, to reduce part (2) to the claim the functor $\gamma_{b}^{\ULA}$ is an equivalence, which follows from the equivalence \eqref{eqn: lcompletioninJbperfadmissible}. The claim on Verdier duality follows from the fact that it sends the constant sheaf to the constant sheaf.

Point (3) follows from \cite[Proposition VII.5.2]{FSGeomLLC} by taking the sheaf $A$ in \emph{loc.cit} to be $\IC_{\mu} \in \Detale(\Hck_{G},\bb{Z}_{\ell}[\sqrt{p}])$, where we recall that this sheaf gives rise to $\IC'_{\mu} \in \D_{\bs}(\Hck_{G},\bb{Z}_{\ell}[\sqrt{p}])$ under the embedding $(-)\mapsto \bb{D}(-)^\vee$.
\end{proof}
To deal with the case of rational coefficients, we note that, since the condensed ring $F:=\CO_F{\otimes^{\bs}_{\mathcal{O}_{F,\mathrm{disc}}} F_{\mathrm{disc}}}$ agrees with $\colim_{n}\CO_F\cdot \tfrac{1}{\ell^n}$ as a $\CO_F$-module and the lisse category is stable under colimits, it is itself a lisse $\CO_F$-sheaf. We thus have an embedding
\begin{equation}\label{eq: FasOFmodule}
    \Dlis(\Bun_{G},{F}) \hookrightarrow \Dlis(\Bun_{G},\CO_F),
\end{equation}
by forgetting the $F$-structure. It has a left adjoint given by the solid tensor product
\begin{equation}{\label{eqn: RationalCoefficients}}
- \otimes^{\bs}_{\mathcal{O}_{F}} F: \Dlis(\Bun_{G},\mathcal{O}_{F}) \ra \Dlis(\Bun_{G},F)
\end{equation}
\[ A \mapsto A \otimes^{\bs}_{\mathcal{O}_{F}} F,\]
where this is well-defined by \cite[Proposition~VII.6.2]{FSGeomLLC}. We write $\bb{S}_{\lis}^{F}$ for the semi-orthogonal decomposition on $\mc{D}_{\lis}(\Bun_{G},F)$ as above. Below we omit $\bs$ from notation when the context is clear.
\begin{lemma}{\label{lemma: RationalCoefficients}}
The functor (\ref{eqn: RationalCoefficients}) satisfies the following.
\begin{enumerate}
\item It upgrades to a natural transformation $\bb{S}_{\lis}^{\mathcal{O}_{F}} \ra \bb{S}_{\lis}^{F}$ of semi-orthogonal decompositions.
\item It restricts to a functor 
\[ - \otimes_{\mathcal{O}_{F}} F: \Dlis^{\ULA}(\Bun_{G},\mathcal{O}_{F}) \ra \Dlis^{\ULA}(\Bun_{G},F), \]
such that the natural map below is an isomorphism.
\[  \bb{D}_{\Bun_{G},\mathcal{O}_{F}}(-) \otimes_{\mathcal{O}_{F}} F\to \bb{D}_{\Bun_{G},F}(- \otimes_{\mathcal{O}_{F}} F)\]

\item It is Hecke equivariant. In other words, for all $\mu \in \bb{X}_{*}(T_{\ol{\bb{Q}}_{p}})^{+}$ we have a natural equivalence
\[ T_{\mu}^{\mathcal{O}_{F}}(-) \otimes_{\mathcal{O}_{F}} F \simeq T_{\mu}^{F}(- \otimes_{\mathcal{O}_{F}} F). \]
\end{enumerate}
\end{lemma}
\begin{proof}
Part (1) follows from the fact that both $*$-pullback and $!$-pushforward along $i_{V_{1}V_{2}}: \Bun_{G}^{V_{1}} \ra \Bun_{G}^{V_{2}}$ are compatible with $- \otimes_{\mathcal{O}_{F}} F$. This is because $*$-pullback commutes with colimits, while $i_{V_{1}V_{2}!}$ satisfies the projection formula.

For part (2), we argue as in the proof of \cite[Proposition IV.2.19]{FSGeomLLC}. Namely, denote the structure map $\Bun_G\to \ast$ by $\pi$. Let $A$ be an object in $\Dlis^\mathrm{ULA}(\Bun_G, \CO_F)$, and write $A_{F}$ for the image under (\ref{eqn: RationalCoefficients}). We consider the natural transformation 
\begin{equation}\label{eq: DualityRationalization}
    \bb{D}_{\Bun_{G},\mathcal{O}_{F}}(A) \otimes_{\mathcal{O}_{F}} \pi^\ast(-)\to \RHomint_{\Dlis(\Bun_{G},\CO_F)}(A, \pi^\ast(-)).
\end{equation} 
We only need to show that, when applied to the constant sheaf $F$, this becomes an isomorphism, noting that by adjunction, we have (under \eqref{eq: FasOFmodule})
\[ \RHomint_{\Dlis(\Bun_{G},F)}(A, F)\simeq \RHomint_{\Dlis(\Bun_{G},F)}(A_F, F)\simeq \bb{D}_{\Bun_{G},F}(A_F).\]
For this purpose, it suffices to show that both sides commute with filtered colimits, because $F=\colim_n \CO_F\cdot \tfrac{1}{\ell^n}$ and \eqref{eq: DualityRationalization} evaluated on $\CO_F$ is clearly an isomorphism. The left hand side commutes with filtered colimits, since both $\pi^\ast$ and tensor product do. For the right-hand side, this follows from the adjunction
\[ \pi_\natural(A\otimes-) \dashv R\mc{H}om_{\Dlis(\Bun_{G},\CO_F)}(A, \pi^\ast(-)),\]
and the fact that the left adjoint $\pi_\natural(A\otimes-)$ preserves compact objects (See \cite[Lemma~IV.2.20]{FSGeomLLC}). Indeed, we can consider the compact generators $A_K^b$ defined in \cite[Proposition VII.7.4]{FSGeomLLC} for varying $b \in B(G)$ and $K\subset J_b(\qp)$ an open pro-$p$ subgroup. It follows from \cite[Proposition VII.7.6]{FSGeomLLC} that  we have an isomorphism 
\[ \pi_{\natural}(A \otimes A_K^b) \simeq \RHom_{\Dlis(\Bun_{G},\CO_F)}(\bb{D}_{\mathrm{BZ}}A_K^b, A), \]
where $\bb{D}_{\mathrm{BZ}}$ is the Bernstein-Zelevinsky duality functor. However, the Bernstein--Zelevinsky dual of $A_K^b$ is up to shift and twist $i_{b,!}\operatorname{c-Ind}_K^{J_b(\qp)}(\Lambda)$ (cf. proof of \cite[Proposition 7.6]{FSGeomLLC}). Hence the right-hand side is isomorphic to (up to shift and twist) $(i_b^!A)^K$, which is a perfect complex of $\mathcal{O}_{F}$-modules by ULAness of $A$ (\cite[Proposition VII.7.9]{FSGeomLLC}), which are in particular compact objects of $\Dlis(\ast,\mathcal{O}_{F}) \simeq \D(\mathcal{O}_{F})$, as desired.

Point (3) is an application of projection formula for $h^{\ra}_{\natural}$ in the solid formalism (\cite[Proposition VII.3.1 (i)]{FSGeomLLC}) and the fact that $*$-pullback is symmetric monoidal. 
\end{proof}
We now restrict the coefficients to Setup~\ref{assumption: coefficients} and make the following construction in $\Dlis(\Bun_G,\Lambda)$, similar to Definition~\ref{defn:FIC}. In the case of rational coefficients, it is important for us to work with lisse-\'etale sheaves, rather than the isogeny category of $\Detale(\Bun_G,\CO_F)$, due to potential unbounded $\ell$-torsion in $\hat{\mc{F}}_{\IC, \xi, \CO_F}$. Namely, for an admissible smooth representation $M$ with $\CO_F$-coefficients, we would like the operation of passing to $F$ coefficients to kill all $\ell$-torsion in $M$. But it can happen that $M$ has unbounded torsion, which builds up a torsion-free part in the $\ell$-adic completion. This part will survive in the isogeny category, though undesired (cf. \cite[Remark~4.10]{BhattHanSixFunctors}). Therefore, we need to work with a category which contains all (discrete) $\bb{Z}_\ell$-modules, instead of just $\ell$-complete ones, where we can form tensor product with $\bb{Q}_\ell$ as $\bb{Z}_\ell$-modules. The formalism of $\Dlis$ is suitable for this purpose, see the comments on Page 24 of \cite{FSGeomLLC}.

\begin{construction}{\label{cons: ConstructionofSheaveswithQellCoefficients}}
If $\Lambda$ is as in \ref{assumption: coefficients} (1) and $\xi$ is an algebraic representation of $\mathsf{G}_{\ol{\mathbb{Q}}_{\ell}}$ defined over $\Lambda$, we define 
\[\mathcal{F}_{\xi,\IC,\Lambda} = \hat{\mathcal{F}}_{\xi,\IC,\Lambda} \in \Detale(\Bun_{G},\Lambda) \simeq \Dlis(\Bun_{G},\Lambda).\] 
If $\Lambda = \mathcal{O}_{F}$ and $\xi$ is defined over $\CO_F$, we define $\mathcal{F}_{\xi,\IC,\mathcal{O}_{F}} \in \Dlis^{\ULA}(\Bun_{G},\mathcal{O}_{F})$ to be the sheaf corresponding to
$\hat{\mathcal{F}}_{\xi,\IC,\mathcal{O}_{F}} \in \Detale^{\ULA}(\Bun_{G},\mathcal{O}_{F})$  (where we recall that this lies in the ULA subcategory by Theorem \ref{thm: perversityofICIgs} and the fact that $\mathcal{O}_{F}$ is regular) under the equivalence (\ref{eqn: EquivalenceonULASheaves}). Similarly, if $\Lambda = F$ is the fraction field, then we let $\mathcal{F}_{\xi,\IC,F} := \mathcal{F}_{\xi,\IC,\mathcal{O}_{F}} \otimes^{\bs}_{\mathcal{O}_{F,,\disc}} F_{\disc} \in \Dlis^{\ULA}(\Bun_{G},F)$ denote the rationalization as in \eqref{eqn: RationalCoefficients}. If $\Lambda = \mathcal{O}_{F},F$ for $F/\bb{Q}_{\ell}$ a (possibly infinite) algebraic extension, we define $\mathcal{F}_{\xi,\IC,\Lambda} \in \Dlis(\Bun_{G},\Lambda)$ to be the colimit of the intersection cohomology of the sheaves attached to the finite subextensions of $F/\bb{Q}_{\ell}$ (equivalently, since $\otimes^{\blacksquare}$ commutes with colimits it is given by $\mathcal{F}_{\xi,\IC,L} \otimes^{\blacksquare}_{L} F$, where $L/\bb{Q}_{\ell}$ is a finite subextension contained in $F$ and over which $\xi$ is defined). Similarly, we define it in the case where $\Lambda$ is an infinite extension of $\bb{F}_{\ell}$ by taking colimits over the torsion case. As before, these sheaves are supported on the open substack $\Bun_{G,\mu^{-1}} \hookrightarrow \Bun_{G}$ and we will occasionally abuse notation and also write $\mathcal{F}_{\xi,\IC,\Lambda} \in \Dlis(\Bun_{G,\mu^{-1}},\Lambda)$ for their restrictions. 
\end{construction}

We now assume that $\Lambda = \mathcal{O}_F$ or $\Lambda = F$ for $F/\bb{Q}_{\ell}$ an algebraic extension and $\xi$ is an algebraic representation of $\mathsf{G}_{\ol{\mathbb{Q}}_{\ell}}$ defined over $F$ or that $\Lambda$ is an algebraic extension of $\bb{F}_{\ell}$ and $\xi$ is defined over a finite subextension of $\bb{F}_{\ell}$. Let $(\mathsf{G},\mathsf{X})$ be a Shimura datum as in Assumption~\ref{assumption:codimension} with reflex field $\mathsf{E}$ and Hodge cocharacter $\mu$. We let $\mathcal{L}^{\mathrm{alg}}_{\xi}$ be the lisse sheaf with coefficients in $\Lambda$ determined by $\xi$ on the inverse system of Shimura varieties $\{\Shstar\gx_{K.\ol{\mathsf{E}}}\}_K$, for $K \subset \mathsf{G}(\bb{A}_{f})$ sufficiently small compact opens, as defined in \cite[Section~5.1]{PinkHigherDirectImages}. We set $\IC(\gxno,\xi)_{K,\Lambda} := R\Gamma(\Shstar\gx_{K,\ol{\mathsf{E}}},\IC_{\Shstar_{K}}(\mathcal{L}_{\xi}^{\mathrm{alg}}))$ to be the intersection cohomology at level $K$, and, as in the torsion case, consider the smooth admissible $G(\bb{Q}_{p})$-representation with coefficients in $\Lambda$
\begin{equation}{\label{eqn: IntersectionComplexAtInfiniteLevel}}
\IC(\gxno,\xi)_{K^{p},\Lambda}:= \varinjlim_{K_{p} \ra \{1\}} \IC(\gxno,\xi)_{K^{p}K_{p},\Lambda}.
\end{equation}
It is also equipped with a continuous $\operatorname{Gal}_{\mathsf{E}}$-action and an away-from-$p$ Hecke action. 

With Construction~\ref{cons: ConstructionofSheaveswithQellCoefficients} in hand, we can now carry over the results of Theorems~\ref{thm: stalks} and~\ref{thm: Hecke operator torsion} to the case of the coefficient systems considered in assumption \ref{assumption: coefficients}. As before, we fix an isomorphism $\ol{\bb{Q}}_{p} \simeq \bb{C}\supset \mathsf{E}$, and let $E$ be the completion of $\mathsf{E}$ at the induced $p$-adic place and $\mu$ the cocharacter corresponding to the Hodge cocharacter of the Shimura datum under this isomorphism.
\begin{theorem}{\label{thm: Hecke operator}}
We let $\Lambda/\bb{Z}_{\ell}[\sqrt{p}]$ be as in Setup~\ref{assumption: coefficients} and $\xi$ be an algebraic representation of $\mathsf{G}_{\ol{\mathbb{Q}}_{\ell}}$ defined over $\Lambda$. The following is true. 
\begin{enumerate}
\item (Stalks) For all $b \in B(G,\mu^{-1})$, we have natural, Hecke-equivariant identifications
\[
V_{\xi,\IC,\Lambda,b}[-d_b]\simeq i_b^*\mc{F}_{\xi, \IC, \Lambda}. 
\]
as smooth $J_{b}(\bb{Q}_{p})$-representations, where the LHS is defined exactly as in~\eqref{eqn: defofIC} (without requiring $\Lambda$ to be as in Setup \ref{assumption: coefficientsystemsingeneral}). 
\item (Verdier self-duality) If $\Lambda$ is self-injective, then we have $\bb{D}_{\Bun_{G,\mu^{-1}}}(\mathcal{F}_{\xi,\IC,\Lambda}) \simeq \mathcal{F}_{\xi^{\vee},\IC,\Lambda}$.
\item (Perversity) If $\Lambda$ is self-injective, then we have $\mathcal{F}_{\xi,\IC,\Lambda} \in \Perv_{\mathrm{lis}}(\Bun_{G,\mu^{-1}},\Lambda)$. 
\item (ULAness) If $\Lambda$ is regular, then we have  $\mathcal{F}_{\xi,\IC,\Lambda} \in \D^{\ULA}_{\mathrm{lis}}(\Bun_{G},\Lambda)$.
\item (Hecke operator) We have an away-from-$p$ Hecke- and $G(\Q_p)\times W_E$-equivariant isomorphism 
\[ 
i_{1}^{*}T_{\mu}\mc{F}_{\xi, \IC, \Lambda} \simeq \IC(\mathsf{G},\mathsf{X}, \xi)_{K^p,\Lambda}.
\] 
\end{enumerate}
\end{theorem}
\begin{proof}
For $\Lambda$ finite, part (1) is already established in Theorem~\ref{thm: stalks}, (2)-(4) in Theorem \ref{thm: perversityofICIgs} and part (5) is established in Theorem \ref{thm: Hecke operator torsion}. 

Let $F/\bb{Q}_{\ell}$ be a finite extension. For $\Lambda=\CO_F$, since it is not self-injective, only part (1), (4), (5) are relevant. We recall that the inverse functor to the equivalence $\gamma_{b}^{\ULA}$ of \eqref{eqn: lcompletioninJbperfadmissible} is given by the functor of taking smooth vectors in an $\ell$-adically complete representation, as in \eqref{eqn: SmoothVectorsofLcomplete}. Hence part (1) follows from Theorem \ref{thm: stalks}, Lemma \ref{lemma: smoothvectorsinthestalk}, Lemmas \ref{lemma: FromLCompletetoSmooth} (1)-(2). Part (4) follows from Theorem~\ref{thm: perversityofICIgs} and Lemma~\ref{lemma: FromLCompletetoSmooth}(2). Part (5) again follows from the $\ell$-complete case (obtained by taking an inverse limit of Theorem~\ref{thm: Hecke operator torsion} and using Theorem \ref{thm: inverselimit}, recalling our definition of $\IC_{\Igs^{*}}(\mathcal{L})$ in this case (Definition \ref{defn: ICIgs*})), using Lemma~\ref{lemma: FromLCompletetoSmooth} (3), and the isomorphism ~\eqref{eqn: InvariantsinCompletedIntersectionCohomology} below.

For $\Lambda=F$, part (1) follows from the case of $\Lambda=\CO_F$ and Lemma \ref{lemma: RationalCoefficients} (1). Part (2) needs some justification, we note that, for all $n \geq 1$, we have a natural map
\[ \IC_{\Igs^{*}}(\mathcal{L}_{\xi}) \otimes_{\mathcal{O}_{F}} \mathcal{O}_{F}/\ell^{n} \ra \IC_{\Igs^{*}}(\mathcal{L}_{\xi}/\ell^{n}) \]
following from the definition of $\IC_{\Igs^{*}}(\mathcal{L}_{\xi})$ as an inverse limit (Definition \ref{defn: ICIgs*}).Now, by applying Verdier duality to this map and precomposing with the morphism from Construction \ref{cons: DualityMap} applied to $\IC_{\Igs^{*}}(\mathcal{L}_{\xi}/\ell^{n}) \in \Detale(\Igs^{*},\mathcal{O}_{F}/\ell^{n})$, we obtain a natural morphism 
\[ \IC_{\Igs^{*}}(\mathcal{L}_{\xi}/\ell^{n}) \ra \bb{D}_{\Igs^{*}}(\IC_{\Igs^{*}}(\mathcal{L}_{\xi}) \otimes_{\mathcal{O}_{F}} \mathcal{O}_{F}/\ell^{n})  \]
passing to inverse limits over $n \geq 1$, and then applying the proper pushforward $\ol{\pi}_{\min*}$, we obtain a natural map 
\begin{equation}{\label{eqn: maptodualwithFCoefficients}}
 \hat{\mathcal{F}}_{\xi,\IC,\mc{O}_{F}} \ra \bb{D}_{\Bun_{G,\mu^{-1}}}\hat{\mathcal{F}}_{\xi^{\vee},\IC,\mc{O}_{F}}, 
\end{equation}
which, by applying the equivalence $\gamma^{\ULA}_{\Bun_{G,\mu^{-1}}}$ given by Lemma \ref{lemma: FromLCompletetoSmooth} (1) and (2), gives us a natural map 
\[ \mathcal{F}_{\xi,\IC,\mc{O}_{F}} \ra \bb{D}_{\Bun_{G,\mu^{-1}}}\mathcal{F}_{\xi^{\vee},\IC,\mc{O}_{F}} \]
in $\Dlis^{\ULA}(\Bun_{G,\mu^{-1}},\mathcal{O}_{F})$. Applying (\ref{eqn: RationalCoefficients}) and using Lemma \ref{lemma: RationalCoefficients} (2), we obtain a natural map 
\[ \mathcal{F}_{\xi,\IC,F} \ra \bb{D}_{\Bun_{G,\mu^{-1}}}\mathcal{F}_{\xi^{\vee},\IC,F}. \]
As discussed above, we have computed the $*$-stalks of the source of the map as
\[ \colim_{m \geq 1} R\Gamma(\Ig^{b,*}_{m},\IC_{\Ig^{b,*}_{m}}(\mathcal{L}_{\xi,b,m}^{\alg}))[-d_{b}]. \]
in terms of a colimit of the usual algebraic intersection cohomology with $F$-coefficients. Similarly, we can compute the $*$-stalks of the target using Theorem \ref{thm: stalks} and Lemma \ref{lemma: smoothvectorsinthestalk} as 
\[ \colim_{m \geq 1} R\Gamma(\Ig^{b,*}_{m},\bb{D}_{\Ig^{b,*}_{m}}(\IC_{\Ig^{b,*}_{m}}(\mathcal{L}_{\xi^\vee,b,m}^{\alg}))[-d_{b}] \]
Moreover, through a lengthy but straightforward calculation (using that Construction \ref{cons: DualityMap} respects perverse $t$-exact proper pushforwards and smooth pullbacks) the natural map on $*$-stalks identifies with the natural map 
\[ \colim_{m \geq 1} R\Gamma(\Ig^{b,*}_{m},\IC_{\Ig^{b,*}_{m}}(\mathcal{L}_{\xi,b,m}^{\alg}))[-d_{b}] \ra \colim_{m \geq 1} R\Gamma(\Ig^{b,*}_{m},\bb{D}_{\Ig^{b,*}_{m}}(\IC_{\Ig^{b,*}_{m}}(\mathcal{L}_{\xi^\vee,b,m}^{\alg}))[-d_{b}] \]
on algebraic intersection cohomology with coefficients in $F$ induced by Construction \ref{cons: DualityMap}, which is an isomorphism by classical results on the intersection cohomology of algebraic varieties. By varying over all $b \in B(G,\mu^{-1})$, this in particular implies that (\ref{eqn: maptodualwithFCoefficients}) is also an isomorphism. 

Part (3) follows from (1)-(2) and Artin vanishing as in the proof of \ref{thm: perversityofICIgs}. Part (4) follows from the case $\Lambda=\CO_F$, combined with Lemma~\ref{lemma: RationalCoefficients}(2). Part (5) reduces to the case of $\Lambda = \mathcal{O}_{F}$ using Lemma \ref{lemma: RationalCoefficients} (3). 

Parts (1)-(5) when $\Lambda = \mathcal{O}_{F},F$ for $F/\bb{Q}_{\ell}$ an infinite algebraic extension is easily checked to follow formally from the case of a finite extension by taking colimits, where, for Part (2), we argue as in the proof of Lemma \ref{lemma: RationalCoefficients} (2). Similarly, the case where $\Lambda$ is an infinite extension of $\bb{F}_{\ell}$ follows from the case of a finite torsion algebra by taking colimits. 
\end{proof}

\subsection{The Mantovan product formula for intersection cohomology}

In this section, we establish a version of the Mantovan product formula for the intersection cohomology of Shimura varieties. Let $\gx$ be a Shimura datum satisfying Assumption~\ref{assumption:codimension} as before.

We follow the approach of~\cite[Corollary~3.22]{HamannLeeTorsion} and~\cite[Theorem 8.5.7]{DvHKZ}. We already showed that the intersection cohomology groups of $\Sh\gx_K$ is obtained by applying a Hecke operator to the complex $\mathcal{F}_{\xi,\IC,\Lambda}$. The desired formula follows from using excision with respect to the Newton stratification on $\Bun_G$, and rewriting the graded pieces in terms the (compact support) cohomology of local Shimura varieties and the intersection cohomology of partial minimal compactification of Igusa varieties. The last step uses our computation of stalks of $\mathcal{F}_{\xi,\IC,\Lambda}$ obtained in Theorem~\ref{thm: Hecke operator}.

Again, we let $k=\ol{\bb{F}}_p$ and consider Diagram~\eqref{eqn: CartesiandiagramforClosedShimuraVariety} over $\ast=\Spd k$. Let $\mu$, $\mathsf{E}$, $\bb{C}\simeq \ol{\bb{Q}}_p$, $E$, $d=\langle2\rho, \mu\rangle$ be as before, and $\Lambda$ be as in Setup~\ref{assumption: coefficients}.
\subsubsection{}{\label{ss: Variation'sonMantovan'sFiltration}}
 Recall from \S~\ref{cor: GeoMantovanBoundary} that for $b \in B(G,\mu^{-1})$, we have the moduli space of local shtukas $\Sht(G,b,\mu)_{\infty}$, which is a $v$-sheaf over $\Breve{E}$, the compositum of $E$ with the maximal unramified extension of $\bb{Q}_p$. This represents the functor sending $S \in \Perf$ to the set of all pairs $(S^\#,\alpha)$ where $S^\#$ is an untilt of $S$ over $\Breve{E}$, and $\alpha$ is a modification $\mathcal{E}_{0,S} \dashrightarrow \mathcal{E}_{b,S}$ from the trivial $G$-bundle with meromorphy along $S^\#$ and bounded by $\mu$. Here $\mathcal{E}_{b,S}$ (resp. $\mathcal{E}_{0,S}$) is the pullback of the $G$-bundle $\mathcal{E}_{b}$ attached to $b$ on $X_{\Spd k}$ along the natural map $X_{S} \ra X_{\Spd k}$. The space $\Sht(G,b,\mu)_{\infty}$ is equipped with a natural action of $\underline{G(\bb{Q}_{p})}$ (resp. $\underline{J_{b}(\bb{Q}_{p})}$) via automorphisms of $\mathcal{E}_{0}$ (resp. $\mathcal{E}_{b}$).  

We define the complex
\begin{equation}{\label{eqn: CohofLocalShimuraVariety}}
R\Gamma_{c}(G,b,\mu,\Lambda(d_{b})) := \varinjlim_{K_{p} \ra \{1\}}R\Gamma_{c}\left(\Sht(G,b,\mu)_{\infty,C}/\ul{K_{p}},\Lambda(d_{b})\right) 
\end{equation}
of $G(\mathbb{Q}_{p}) \times J_{b}(\mathbb{Q}_{p}) \times W_{E}$-modules, where $K_p$ runs over compact open subgroups of $G(\Q_p)$ and where $\Lambda(d_{b})$ is the sheaf with $J_{b}(\mathbb{Q}_{p})$-action given as in \cite[Lemma~7.4]{KoTorsionPaper}. 
\begin{corollary}\label{cor: product formula}
    (Mantovan product formula) For $\Lambda$ a coefficient system as in Setup~\ref{assumption: coefficients}, the complex $\IC(\mathsf{G},\mathsf{X}, \xi)_{K^p,\Lambda}$ is equipped with a Hecke- and $G(\mathbb{Q}_p)\times W_E$-equivariant filtration indexed by $b\in B(G,\mu^{-1})$, with graded pieces isomorphic to 
\[
(R\Gamma_{c}(G,b,\mu,\Lambda(d_{b}))[d] \otimes |\cdot|^{d/2}) \otimes_{\mathcal{H}(J_{b})} V_{\xi,\IC,\Lambda,b}[d_{b}]. 
\]
Here $|\cdot|: W_{E} \ra \Lambda^{*}$ denotes the norm character and $\mc{H}(J_b)$ is the Hecke algebra for $J_b(\qp)$ with $\Lambda$-coefficients.
\end{corollary}
\begin{proof}
This follows from Theorem~\ref{thm: Hecke operator} and excision applied to the sheaf $\mathcal{F}_{\xi,\IC,\Lambda}$ with respect to the semi-orthogonal decomposition. To identify the graded pieces, we combine the computation in Theorem~\ref{thm: Hecke operator} part (1) and (5), with \cite[Lemma~10.1]{HamGeomES}, just as in the proof of \cite[Corollary~3.22]{HamannLeeTorsion}. 
\end{proof}

\subsubsection{Variants of the Mantovan product formula}\label{sec: variants}
For $\PP$ a conjugacy class of admissible parabolic subgroups, we recall the cartesian diagram 
\[ 
    \begin{tikzcd}
        \mathcal{S}_{\leq\PP} \ar[r,"\pi_{\leq [\mathbf{P}]}"] \ar[d,"\tilde{h}_{\leq [\mathsf{P}]}"] & \Fl \ar[d]\\
        \Igs_{\leq\PP} \ar[r,"\ol{\pi}_{\leq [\mathsf{P}]}"] & \Bun_G,
    \end{tikzcd}
\]
of Theorem~\ref{Thm:StratifiedCartesian}, and the analogous diagram for the $\PP$-condition in place of $\leq \PP$. We can use these diagrams to obtain variants of the Mantovan product formula for the cohomology of the boundary. Indeed, we pass to the quotient by $G(\bb{Q}_{p})$ to obtain a commutative diagram
\begin{equation}{\label{eqn: CartesianDiagramGQpqquotient}}
    \begin{tikzcd}
        \left[\mathcal{S}_{\leq\PP}/\underline{G(\bb{Q}_{p})}\right] \ar[r,"\pi_{\leq \PP}"] \ar[d,"\tilde{h}_{\leq [\mathsf{P}]}"] & \left[\Fl/\underline{G(\bb{Q}_{p})}\right] \arrow[r,"h_{2}"]  \ar[d,"h_{1}"] & \left[\ast/\underline{G(\bb{Q}_{p})}\right]\\
        \Igs_{\leq\PP} \ar[r,"\ol{\pi}_{\leq [\mathsf{P}]}"] & \Bun_G &,
    \end{tikzcd}
\end{equation}
where $h_{1}$ and $h_{2}$ are as in Lemma  \ref{lemma: computingminusculeHeckecorrespondences} and the left square is cartesian. 

\begin{definition}{\label{defn: SheavesonBunG}}Fix a conjugacy class $\PP$ of admissible parabolic subgroups and a $\mathbb{Z}_\ell[\sqrt{p}]$-algebra $\Lambda$ as in Setup~\ref{assumption: coefficients}.
\begin{enumerate}
    \item If $\Lambda$ is finite, we let $\mathcal{F}_{\leq \PP,\Lambda} := R\ol{\pi} _{\leq \PP*}\Lambda$ and $\mathcal{F}_{!,\PP,\Lambda} := R\ol{\pi}_{\PP!}\Lambda$ in $\Detale(\Bun_G,\Lambda)$. 
    \item More generally, if $\Lambda$ is as in cases (2)-(5) of Assumption~\ref{assumption: coefficients}, we define $\mathcal{F}_{\leq \PP,\Lambda} \in \Dlis^{\ULA}(\Bun_{G},\Lambda)$ and $\mathcal{F}_{!,\PP,\Lambda} \in \Dlis^{\ULA}(\Bun_{G},\Lambda)$ by the same recipe as in Construction~\ref{cons: ConstructionofSheaveswithQellCoefficients}, where the required ULAness follows from the same proof as Theorem \ref{thm: Hecke operator}, part (4).
\end{enumerate}
\end{definition}

\noindent For these sheaves, we obtain the following analog of Theorem~\ref{thm: Hecke operator} and Corollary~\ref{cor: product formula}, which is proved in essentially the same way, but is strictly easier (cf. the proof of \cite[Proposition~3.6]{HamannLeeTorsion}). We leave the details to the interested reader.  

\begin{theorem}\label{thm: product formula boundary}\leavevmode
Let $\Lambda$ be as in Setup \ref{assumption: coefficients}.
\begin{enumerate}

\item (Stalks) For all $b \in B(G,\mu^{-1})$, write $g_{b,\PP}: \Ig^{b}_{\PP} \hookrightarrow \Ig^{b}_{\leq \PP}$ for the natural open immersion of $\Ig^{b}_{\PP}$ into its closure. We have natural, prime-to-$p$ Hecke- and $J_{b}(\bb{Q}_{p})$-equivariant identifications
\[ V_{b,\Lambda,\leq\PP}:= R\Gamma(\Ig^{b}_{\leq \PP},\Lambda)^{\mathrm{sm}} \simeq i_{b}^{*}\mathcal{F}_{\leq \PP,\Lambda} 
\]
and
\[
V_{b,\Lambda, \PP}:=R\Gamma(\Ig^{b}_{\leq\PP},g_{b, \PP!}\Lambda)^{\mathrm{sm}} \simeq i_{b}^{*}\mathcal{F}_{!,\PP,\Lambda}, \]
where the LHS denotes smooth vectors in the usual \'etale cohomology with $\Lambda$-coefficients, i.e. the colimit of the cohomology of the strata on finite level Igusa varieties.
\item (Semi-perversity) We have $\mathcal{F}_{!,\PP,\Lambda}$, $\mathcal{F}_{\leq \PP,\Lambda} \in \phantom{}^{p}\D_{\mathrm{lis}}^{\leq 0}(\Bun_{G,\mu^{-1}},\Lambda)$. 
\item (ULAness) The sheaves $\mathcal{F}_{!,\PP,\Lambda}$, $\mathcal{F}_{\leq \PP,\Lambda}$ are ULA if $\Lambda$ is regular.
\item (Hecke operator) We have prime-to-$p$ Hecke and $G(\mathbb{Q}_p)\times W_E$-equivariant isomorphisms
\[ 
i_{1}^{*}T_{\mu}\mathcal{F}_{\leq \PP,\Lambda} \simeq \varinjlim_{K_{p} \ra \{1\}} R\Gamma(\Shstar\gx_{K^{p}K_{p},\leq \PP,\ol{\mathsf{E}}},\Lambda)[d](\frac{d}{2}), 
\]
and 
\[ 
i_{1}^{*}T_{\mu}\mathcal{F}_{!,\PP,\Lambda} \simeq \varinjlim_{K_{p} \ra \{1\}} R\Gamma_{c}(\Shstar\gx_{K^{p}K_{p},\PP,\ol{\mathsf{E}}},\Lambda)[d](\frac{d}{2}). 
\]

\item (Mantovan product formula) The colimit $\varinjlim_{K_{p} \ra \{1\}} R\Gamma(\Shstar\gx_{K^{p}K_{p},\leq \PP,\ol{\mathsf{E}}},\Lambda)[d](\frac{d}{2})$ of cohomology groups is equipped with a Hecke- and $G(\mathbb{Q}_p)\times W_E$-equivariant filtration indexed by $b\in B(G,\mu^{-1})$, with graded pieces isomorphic to 
\[
 (R\Gamma_{c}(G,b,\mu, \Lambda(d_{b}))  \otimes_{\mathcal{H}(J_{b})} V_{b,\leq \PP,\Lambda}[2d_{b}])
\]
The analogous result holds for $\varinjlim_{K_{p} \ra \{1\}} R\Gamma_{c}(\Shstar\gx_{K^{p}K_{p},\PP,\ol{\mathsf{E}}},\Lambda)$ and $V_{b,\PP,\Lambda}$. 
\end{enumerate}
\end{theorem}

\subsection{Vanishing Results}
In this section, we will use our main results on the construction of $\mathcal{F}_{\xi,\IC, \Lambda}$ and its properties, to prove that under a certain genericity assumption, the intersection cohomology of the minimal compactification of Shimura varieties is concentrated in degree $0$. This is consistent with the conjectures of Koshikawa--Shin~\cite{KoshikawaShin}, as explained in \S \ref{ss: comparisonwithKoshikawaShin} below. 

Throughout this section, we fix an algebraically closed perfectoid field $C/\qp$ and base-change Diagram~\eqref{eqn: CartesiandiagramforClosedShimuraVariety} to $\ast=\Spd C$. We moreover impose the following assumption.
\begin{assumption}{\label{assumption: coefficientsfortexactness}}
The coefficient ring $\Lambda/\bb{Z}_{\ell}[\sqrt{p}]$ is an algebraically closed field as in Setup~\ref{assumption: coefficients}. If $\Lambda$ is torsion then we also assume that $\ell \nmid \pi_{0}(Z(G))$.
\end{assumption}

We let $\Phi^{\mathrm{ss}}(G)$ denote the set of $\widehat{G}$-conjugacy classes of continuous semisimple homomorphisms $\phi: W_{\mathbb{Q}_{p}} \ra \widehat{G}(\Lambda)$, see \cite[Definition VIII.3.1]{FSGeomLLC}. We recall \cite[Appendix~A]{HamannLeeTorsion} that, for each $\phi \in \Phi^{\mathrm{ss}}(G)$, we can form the full subcategory $\Dlis(\Bun_{G},\Lambda)_{\phi} \hookrightarrow \Dlis(\Bun_{G},\Lambda)$ of objects for which the endomorphism defined by $f \in \mathcal{Z}^{\mathrm{spec}}(G) \setminus \mf{m}_{\phi}$ is an isomorphism. Here $\mathcal{Z}^{\mathrm{spec}}(G) := \mathcal{O}_{\mathfrak{X}_{\widehat{G}}}(\mathfrak{X}_{\widehat{G}})$ is the ring of global functions of $\mathfrak{X}_{\widehat{G}}/\Spec \Lambda$, the stack of Langlands parameters constructed in \cite[Section~VII.1]{FSGeomLLC} base-changed along $\bb{Z}_\ell\to\Lambda$, $\mf{m}_{\phi} \subset\mathcal{Z}^\mathrm{spec}(G)$ is the maximal ideal corresponding to the semisimple parameter $\phi$ via \cite[Proposition~VIII.3.2]{FSGeomLLC}, and $\mathcal{Z}^\mathrm{spec}(G)$ acts on $\Dlis(\Bun_{G},\Lambda)$ by the excursion action. The localized category is a full subcategory of $\Dlis(\Bun_{G},\Lambda)$ and the natural inclusion $\Dlis(\Bun_{G},\Lambda)_{\phi} \hookrightarrow \Dlis(\Bun_{G},\Lambda)$ admits a left adjoint given by an idempotent localization map
\[ (-)_{\phi}:  \Dlis(\Bun_{G},\Lambda) \ra \Dlis(\Bun_{G},\Lambda)_{\phi}, \]
as described in \cite[Appendix~A]{HamannLeeTorsion}. We recall, by \cite[Lemma~4.2 (1)]{HamannLeeTorsion}, that if $\Lambda$ is a field, the category $\Dlis(\Bun_{G},\Lambda)_{\phi}$ has (and is essentially characterized by) the property that any Schur-irreducible object $A \in \Dlis(\Bun_{G},\Lambda)_{\phi}$  has Fargues--Scholze parameter $\phi$. By \cite[Lemma~4.2]{HamannLeeTorsion}, the Hecke action $T_{V}: \Dlis(\Bun_{G},\Lambda) \ra \Dlis(\Bun_{G},\Lambda)^{BW_{\mathbb{Q}_{p}}}$ for $V \in \Rep_{\Lambda}(\phantom{}^{L}G)$ described in \S \ref{s: HeckeOperatorsandDlisse}, induces a localized Hecke operator
\begin{equation}{\label{eqn: LocalizedHeckeOperator}}
 T_{V,\phi}: \Dlis(\Bun_{G},\Lambda)_{\phi} \ra \Dlis(\Bun_{G},\Lambda)^{BW_{\mathbb{Q}_{p}}}_{\phi}, 
\end{equation}
which sits in an obvious commutative square involving $T_{V}$ and the localization map $(-)_{\phi}$.

For $\mu \in \bb{X}_{*}(T_{\ol{\mathbb{Q}}_{p}})^{+}$ a geometric dominant cocharacter of $G$ with reflex field $E_{\mu}$, we recall the Hecke operator attached to the tilting module $\mathcal{T}_{\mu} \in \Rep(\widehat{G})$
\[ T_{\mu}: \Dlis(\Bun_{G},\Lambda) \ra \Dlis(\Bun_{G},\Lambda)^{BW_{E_{\mu}}}.\]
This gives rise to a localized Hecke operator $T_{\mu,\phi}$.

Consider the inclusion of the open locus 
\[i_{\mu^{-1}}: \Bun_{G,\mu^{-1}} \hookrightarrow \Bun_{G}\]
corresponding to $B(G,\mu^{-1}) \subset B(G)$ under the bijection $B(G) \simeq |\Bun_{G}|$. By \cite[Proposition~A.9]{RapAppendixtoPadicCohomology}, the restriction of the Hecke operator $i_{1}^\ast T_{\mu}$ factors through the idempotent operation
\[\Dlis(\Bun_{G},\Lambda)\xrightarrow{i_{\mu^{-1}}^\ast}  \Dlis(\Bun_{G,\mu^{-1}},\Lambda) \xrightarrow{i_{\mu^{-1},!}}\Dlis(\Bun_{G},\Lambda).\]
We abusively still write $i_{1}^{*}T_{\mu}$ for the corresponding operator 
\[ i_{1}^{*}T_{\mu}i_{\mu^{-1},!}: \Dlis(\Bun_{G,\mu^{-1}},\Lambda) \ra \Dlis(\Bun_G^1,\Lambda)\simeq\D(G(\mathbb{Q}_{p}),\Lambda).\]

The localization map $(-)_{\phi}$ respects the semi-orthogonal decomposition of $\Dlis(\Bun_{G},\Lambda)$ described in \S~\ref{s: SemiOrthogonalDecompositionPerverseTStruture}. In particular, if we restrict to the full subcategory $i_{1!}: \D(G(\mathbb{Q}_{p}),\Lambda) \hookrightarrow \D(\Bun_{G},\Lambda)$ corresponding to the locus of $\Bun_{G}$ defined by the trivial $G$-bundle then we obtain a localization map 
\[ (-)_{\phi}: \D(G(\mathbb{Q}_{p}),\Lambda) \ra \D(G(\mathbb{Q}_{p}),\Lambda)_{\phi}, \]
which we may apply to the complex (\ref{eqn: IntersectionComplexAtInfiniteLevel}) of $G(\mathbb{Q}_{p})$-representations. 

By the perfect admissibility\footnote{This is guaranteed by the assumption that $\Lambda$ is a field.} of the complex in \eqref{eqn: IntersectionComplexAtInfiniteLevel}, these localization maps induce a direct sum decomposition
\begin{equation}{\label{eqn: directsumdecompositionoverparameters}}
 \IC(\gxno,\xi)_{K^{p},\Lambda} \simeq \bigoplus_{\phi \in \Phi^{\mathrm{ss}}(G)} \IC(\gxno,\xi)_{K^{p},\Lambda,\phi}, 
\end{equation}
by \cite[Proposition~A.5.]{HamannLeeTorsion}. We have the following variant of Theorem \ref{thm: Hecke operator}(5) linking these direct summands to the localized Hecke operators (\ref{eqn: LocalizedHeckeOperator}) which follows from combining the above discussion and Theorem \ref{thm: Hecke operator}~(5).

\begin{corollary}{\label{cor: phiisotyicpartcalculation}}
Let $\Lambda$ be a coefficient system as in Assumption~\ref{assumption: coefficientsfortexactness}, and $\xi$ be an algebraic representation of $\mathsf{G}$ defined over $\Lambda$. We let $\mathcal{F}_{\xi,\IC,\Lambda}$ be obtained via Construction~\ref{cons: ConstructionofSheaveswithQellCoefficients}.  
We have an isomorphism
\[ i_{1}^{*}T_{\mu,\phi}(\mathcal{F}_{\xi,\IC,\Lambda})_{\phi} \simeq \IC(\gxno,\xi)_{K^{p},\Lambda,\phi}, \]
of objects in $\D(G(\mathbb{Q}_{p}),\Lambda)^{BW_{E}}$, i.e complexes of smooth $G(\mathbb{Q}_{p})$ representations with a commuting continuous action of $W_{E}$.
\end{corollary}

The following assumption is key to bounding the cohomological amplitude of $\IC(G,X,\mathcal{L})_{K^p,\phi}$. 
\begin{assumption}{\label{assumption: perversetexactness}}
The parameter $\phi \in \Phi^{\mathrm{ss}}(G)$ has the property that the localized operator
\[ i_{1}^{*}T_{\mu,\phi}: \Dlis^{\ULA}(\Bun_{G,\mu^{-1}},\Lambda)_{\phi} \ra \D_{\PerfAdm}(G(\mathbb{Q}_{p}),\Lambda)_{\phi} \]
is $t$-exact with respect to the perverse $t$-structure on $\Dlis^{\ULA}(\Bun_{G,\mu^{-1}},\Lambda)_{\phi}$ and the standard $t$-structure on $ \D_{\PerfAdm}(G(\mathbb{Q}_{p}),\Lambda)_{\phi} \subset \D(G(\bb{Q}_{p}),\Lambda)_{\phi}$, where we note that the perverse $t$-structure is well-defined on the ULA subcategory by the assumption that $\Lambda$ is a field \cite[Proposition 1.2.1]{HansenBeijingNotes}.
\end{assumption}
\begin{remark}
We have restricted to the ULA subcategory (which is preserved under Hecke operators by \cite[Theorem~I.7.2]{FSGeomLLC}), since in several cases where we can verify this assumption, we only know on this subcategory that the assumption is satisfied, see for example \cite[Corollary 4.27]{HamannLeeTorsion}. We do not know whether we should expect such a thing to hold on the whole $\Dlis(\Bun_{G},\Lambda)_\phi$, cf. Conjecture \cite[Conjecture~6.4]{HamannLeeTorsion}. Indeed, \cite[Theorem~3.38]{YangZhuTorsion} establishes the claim on the full category of sheaves defined by the formal completion around the parameter $\phi$ in the stack of parameters. However, one can only compare the localization around $\phi$ with the formal completion around $\phi$ after restricting ULA objects \cite[Remark~3.28]{YangZhuTorsion}. It seems essential to work with the localization as opposed to the formal completion in general, as the former will be readily related to the localization around maximal ideals in the spherical Hecke algebra.
\end{remark}

As a consequence, we deduce the following theorem.
\begin{theorem}{\label{thm: axiomatizedvanishingstatement}} Let $\Lambda$ be a coefficient system as in Assumption~\ref{assumption: coefficientsfortexactness}, and assume Assumption~\ref{assumption: perversetexactness} holds. The complex $\IC(\gxno,\xi)_{K^{p},\Lambda,\phi}$ is concentrated in degree $0$. This holds, in particular, if $\Lambda = \ol{\bb{F}}_{\ell}$ and $G$, $\phi$, $\ell$, and $p$ are of the form described in \cite[Corollary~4.29]{HamannLeeTorsion} or \cite[Theorem~8.2.1]{peng2025farguesscholzeparameterstorsionvanishing} by results in \textit{loc.cit}.
\end{theorem} 
\begin{proof}
This follows from combining Corollary \ref{cor: phiisotyicpartcalculation} with Part (3) and (4) of Theorem \ref{thm: Hecke operator}.
\end{proof}

We now want to extract an analogous statement at finite level for compact open subgroups $K_{p} \subset G(\bb{Q}_{p})$ by applying a Hochschild--Serre type argument. For this, we will impose an additional banality assumption on the prime $\ell$ when $\Lambda$ is torsion, because passing to a finite level $K_{p}$ can be done fairly routinely when $K_p$ is pro-$p$. We denote by $\mathcal{H}_{K^p}$ the Hecke algebra of bi-$K_p$-invariant, compactly supported, $\Lambda$-valued smooth functions on $G(\mathbb{Q}_p)$. 

\begin{lemma}{\label{lemma: HochschildSerreforIntCohomology}}
For $\Lambda$ a coefficient system as in Assumption \ref{assumption: coefficients}, $\xi$ an algebraic representation of $\mathsf{G}_{\ol{\mathbb{Q}}_{\ell}}$ defined over $\Lambda$, and $K_{p} \subset G(\mathbb{Q}_{p})$ a compact open subgroup whose pro-order is coprime to $\ell$, we have an $\mathcal{H}_{K_p}$-equivariant identification 
\[
R\Gamma(K_{p},\IC(\gxno, \xi)_{K^{p},\Lambda}) \simeq \IC(\gxno,\xi)_{K^{p}K_{p},\Lambda}, 
\]
where the LHS denotes continuous group cohomology, where we equip the $\Lambda$-module with the discrete topology.
\end{lemma}
\begin{proof}
We set $\mathcal{L} = \mathcal{L}_{\xi}$ to be the \'etale $\Lambda$-local system on $\Igs^{*}$ defined by $\xi$. We first explain the proof in the case of finite coefficients. We have an identification of complexes of smooth $G(\mathbb{Q}_p)$-representations
\begin{equation}{\label{eqn: InvariantCohomologyatinfinitelevelproofI}}
R\Gamma(\mathcal{S}_{K^{p}}^{*},\tilde{h}_{\min}^*(\IC_{\Igs^{*}}(\mathcal{L}))) \simeq \IC(\gxno,\xi)_{K^{p},\Lambda}  \end{equation}
by applying $R\Gamma(\mathcal{S}^*_{K^p}, -)$ to the isomorphism $\alpha_{\IC}$ of $G(\mathbb{Q}_p)$-equivariant sheaves of Corollary \ref{cor: InfiniteLevelisColimitIC}. 

The LHS of~\eqref{eqn: InvariantCohomologyatinfinitelevelproofI}, a priori just a sheaf on $\ast$, in fact comes from a sheaf on $[\ast/\underline{K_p}]$, via base-change along the diagram
\[
\begin{tikzcd}
    \mathcal{S}_{K^{p}}^{*} \ar[r]\ar[d]& {[\mathcal{S}_{K^{p}}^{*}/\ul{K_{p}}]} \ar[d]\\
    \ast \ar[r]& {[\ast/\underline{K_{p}}]}.
\end{tikzcd}
\]

This justifies the fact that we have an isomorphism 
\begin{equation}{\label{eqn: InvariantCohomologyatinfinitelevelproofII}} R\Gamma(K_{p},R\Gamma(\mathcal{S}_{K^{p}}^{*},\tilde{h}_{\min}^{*}(\IC_{\Igs^{*}}(\mathcal{L}))) \simeq R\Gamma([\mathcal{S}_{K^{p}}^{*}/\underline{K_{p}}],\tilde{h}_{K_{p},\min}^{*}(\IC_{\Igs^{*}}(\mathcal{L}))), 
\end{equation}
By combining (\ref{eqn: InvariantCohomologyatinfinitelevelproofI}) and (\ref{eqn: InvariantCohomologyatinfinitelevelproofII}), we obtain an isomorphism 
\begin{equation}{\label{eqn: InvariantCohomologyatinfinitelevelproofIII}}
 R\Gamma(K_{p},\IC(\gxno,\xi)_{K^{p},\Lambda}) \simeq R\Gamma([\mathcal{S}_{K^{p}}^{*}/\underline{K_{p}}],\tilde{h}_{K_{p}}^{\min*}(\IC_{\Igs^{*}}(\mathcal{L}))). \end{equation}

Now, by factoring the morphism $[\mathcal{S}^{*}_{K^{p}}/\underline{K_{p}}] \ra \ast$ as $[\mathcal{S}^{*}_{K^{p}}/\underline{K_{p}}] \xrightarrow{q_{K_{p}}} \mathcal{S}^{*}_{K^{p}K_{p}} \ra \ast$, we can rewrite the RHS of \eqref{eqn: InvariantCohomologyatinfinitelevelproofIII} via the isomorphism
\[ R\Gamma([\mathcal{S}_{K^{p}}^{*}/\underline{K_{p}}],\tilde{h}_{K_{p}}^{*}(\IC_{\Igs^{*}}(\mathcal{L}))) \simeq R\Gamma(\mathcal{S}_{K^{p}K_{p}}^{*},q_{K_{p}*}\tilde{h}_{K_{p}}^{*}(\IC_{\Igs^{*}}(\mathcal{L}))), \]
and this identifies with $\IC(\gxno,\mathcal{L})_{K^{p}K_{p}}$ by Corollary \ref{cor: CompareIC} (where the condition on the pro-order of $K_p$ is used) and Proposition \ref{prop: propertiesofthealgebraizationfunctor} (3), as desired. Let $F/\bb{Q}_{\ell}$ be a finite extension. The case of rational coefficients formally follows from the case of $\Lambda = \mathcal{O}_{F}$ by tensoring with $F$. For the case of $\Lambda = \mathcal{O}_{F}$, we may deduce from the torsion case and Theorem \ref{thm: inverselimit} that we have an isomorphism
\begin{equation}{\label{eqn: InvariantsinCompletedIntersectionCohomology}}
R\hat{\Gamma}(K_{p},\widehat{\IC}(\gxno,\xi)_{K^{p},\mathcal{O}_{F}}) = \IC(\gxno,\xi)_{K^{p}K_{p},\mathcal{O}_{F}} 
\end{equation}
for any compact open $K_{p} \subset G(\bb{Q}_{p})$ of pro-order prime to $\ell$. Here 
\[ \widehat{\IC}(\gxno,\xi)_{K^{p},\mathcal{O}_{F}} := \lim_{n \ra \infty} \colim_{K_{p} \ra \{1\}} \IC(\gxno,\xi)_{K^{p}K_{p},\mathcal{O}_{F}/\ell^{n}}  \] 
and $R\hat{\Gamma}(K_{p},-)$ is as defined in (\ref{eqn: SmoothVectorsofLcomplete}).

However, as in the proof of Lemma \ref{lemma: smoothvectorsinthestalk}, (\ref{eqn: InvariantsinCompletedIntersectionCohomology}) implies that 
\[ \delta(\widehat{\IC}(\gxno,\xi)_{K^{p},\mathcal{O}_{F}}) \in \D(G(\bb{Q}_{p}),\mathcal{O}_{F}) \]
is precisely $\IC(\gxno,\xi)_{K^{p},\mathcal{O}_{F}}$. Then it is a straightforward calculation to show that 
\[ R\Gamma(K_{p},\delta(\widehat{\IC}(\gxno,\xi)_{K^{p},\mathcal{O}_{F}})) = R\hat{\Gamma}(K_{p},\widehat{\IC}(\gxno,\xi)_{K^{p},\mathcal{O}_{F}}), \]
see the proof of \cite[Proposition~2.6]{HansenSupercuspidalCohomology}. This gives the desired claim when combined with (\ref{eqn: InvariantsinCompletedIntersectionCohomology}). When $\Lambda$ is as in Setup \ref{assumption: coefficients} (4)-(5), the claim follows by taking colimits of the coefficient systems for $F$ running over finite extensions of $\bb{Q}_\ell$ or $\bb{F}_{\ell}$.
\end{proof}

 We have the following lemma, which gives us slightly stronger results when working with rational coefficients. 

\begin{lemma}{\label{lemma: descentlemma}}
For $K_{p} \subset G(\qp)$ a compact open subgroup such that the pro-order of $K_{p}$ is invertible in $\Lambda$ as in Setup \ref{assumption: coefficients}, we have a natural $\mathcal{H}_{K_{p}}$-equivariant
\[ \IC(\gxno,\xi)_{K^{p}K_{p},\Lambda} \ra R\Gamma(K_{p},\IC(\gxno,\xi)_{K^{p},\Lambda}) \]
map, which induces an injection on each cohomology group. 
\end{lemma}
\begin{proof}
Let $K_{p}' \subset K_{p}$ denote an open and normal pro-$p$ subgroup and let $d := [K_{p}:K'_p]$ denote the index. We set $\mathcal{L} := \mathcal{L}_{\xi}^{\alg}$ to be resulting local systems on the tower $\Shstar(\gxno)_{K,\ol{\mathsf{E}}}$ defined by $\xi$. We consider the finite \'etale Galois covering 
\[ f: \Sh\gx_{K^{p}K_{p}',\ol{\mathsf{E}}} \ra \Sh\gx_{K^{p}K_{p},\ol{\mathsf{E}}} \]
of degree $d$ of finite type schemes over $\ol{\mathsf{E}}$ and write $\overline{f}$ for the corresponding finite morphism on the minimal compactification. Recall that we have the trace and restriction morphisms 
\[ \mathcal{L} \xrightarrow{\mathrm{res}_{f}} f_{*}f^{*}\mathcal{L} \xrightarrow{\mathrm{tr}_{f}} \mathcal{L}, \]
whose composite is multiplication by $d$. We now apply the intermediate extension operation along the open immersion $\Sh\gx_{K^{p}K_{p}',\ol{\mathsf{E}}} \ra \Shstar\gx_{K^{p}K_{p},\ol{\mathsf{E}}}$, which induces a diagram
\begin{equation}{\label{eqn: methodelatraceforIntCoh}}
\IC_{\Shstar_{K^{p}K_{p}}}(\mathcal{L}) \ra \overline{f}_{*}(\IC_{\Shstar_{K^{p}K_p'}}(\mathcal{L})) \ra \IC_{\Shstar_{K^{p}K_{p}}}(\mathcal{L}),  
\end{equation}
where the composite arrow is also given by multiplication by $d$ by $\Lambda$-linearity of the intermediate extension. Here we have used the perverse $t$-exactness and properness of $\overline{f}_{*}$ (since $\ol{f}$ is a finite morphism) to identify the middle term with $\overline{f}_{*}(\IC_{\Shstar_{K^{p}K_p'}}(\mathcal{L}))$. 

Applying $R\Gamma(\Shstar\gx_{K_{p}^{\mathrm{hs}}K^{p},\ol{\mathsf{E}}},-)$ to (\ref{eqn: methodelatraceforIntCoh}), gives us a diagram 
\begin{equation}{\label{eqn: methodelatraceforIntCoh2}}
\IC(\gxno,\xi)_{K^{p}K_{p},\Lambda} \ra \IC(\gxno,\xi)_{K^{p}K_{p}',\Lambda} \ra \IC(\gxno,\xi)_{K^{p}K_{p},\Lambda}, 
\end{equation}
where the composite is multiplication by $d$. Here the first map is precisely the transition morphism in the colimit (\ref{eqn: IntersectionComplexAtInfiniteLevel}). In particular, the first map factors over the invariants $\IC(\gxno,\xi)_{K^{p}K_{p},\Lambda}^{K_{p}/K_{p}'}$ and the resulting map 
\[\IC(\gxno,\xi)_{K^{p}K_{p},\Lambda} \ra \IC(\gxno,\xi)_{K^{p}K_{p}',\Lambda}^{K_{p}/K_{p}'}\] 
is $\mathcal{H}_{K_{p}}$-equivariant and is injective on cohomology groups, because $d$ is invertible in $\Lambda$ and hence the composite map in \eqref{eqn: methodelatraceforIntCoh2} is an isomorphism. We may rewrite the middle term as 
\[\IC(\gxno,\xi)_{K^{p}K_{p}'} \simeq R\Gamma(K'_{p},\IC(\gxno,\xi)_{K^{p},\Lambda})\] 
since $K'_{p}$ is pro-$p$, using Lemma \ref{lemma: HochschildSerreforIntCohomology}. Moreover, we have an isomorphism 
\[(R\Gamma(K_{p}',\IC(\gxno,\xi)_{K^{p},\Lambda}))^{K_{p}/K_{p}'} \simeq R\Gamma(K_{p},\IC(\gxno,\xi)_{K^{p},\Lambda}),\] 
since 
\[R\Gamma(K_{p}/K_{p}',R\Gamma(K_{p},\IC(\gxno,\xi)_{K^{p},\Lambda})) \simeq R\Gamma(K_{p},\IC(\gxno,\xi)_{K^{p},\Lambda})\] 
and the higher cohomology of $R\Gamma(K_{p}^{'}/K_{p},-)$ is all trivial (using that $d$ is invertible in $\Lambda$ by our assumption on $K_{p}$). This gives us the desired map $\IC(\gxno,\xi)_{K^{p}K_{p},\Lambda} \ra R\Gamma(K_{p},\IC(\gxno,\xi)_{K^{p},\Lambda})$. 

To see that it is injective on each cohomology group, we note that by construction it fits into the commutative diagram (\ref{eqn: methodelatraceforIntCoh2}) as follows 
\begin{equation}{\label{eqn: methoddelatraceforIntCoh3}}
 \begin{tikzcd}
&  R\Gamma(K_{p},\IC(\gxno,\xi)_{K^{p},\Lambda})  \arrow[dr]  &  \\
 \IC(\gxno,\xi)_{K^{p}K_{p},\Lambda} \arrow[rr] \arrow[ur] & & R\Gamma(K_{p}',\IC(\gxno,\xi)_{K^{p},\Lambda}).
\end{tikzcd} 
\end{equation}
This forces the desired injectivity by injectivity of the bottom map on cohomology groups.
\end{proof}
\begin{remark}
One should also be able to deduce Lemma \ref{lemma: descentlemma} more directly using a version of Lemma \ref{lemma: HochschildSerreforIntCohomology} that works with rational coefficients $\mathbb{Q}_{\ell}$. However, arguing as in the proof of this lemma would require a six functor formalism with $\bb{Q}_{\ell}$-coefficients in which $[\ast/\underline{K}] \ra \ast$ is $\ell$-cohomologically proper, without assuming that $K$ is pro-$p$, which would take us too far afield. For this reason, we have opted for the above argument. Although it establishes a weaker claim, it suffices for our purposes.
\end{remark}
We now use this to deduce the following Corollary.
\begin{corollary}{\label{cor: HeckeAlgebraResult}}
Let $K_p^{\mathrm{hs}} \subset G(\qp)$ be a hyperspecial subgroup (whose existence is guaranteed by our running Assumption \ref{assumption:codimension}). Let $\Lambda$ be as in Assumption~\ref{assumption: coefficientsfortexactness} and assume that the pro-order of $K_{p}$ is invertible in $\Lambda$. For a maximal ideal $\mf{m} \subset \mathcal{H}_{K_{p}^{\hs}}$ in the spherical Hecke algebra, let $\phi_{\mf{m}} \in \Phi^{\mathrm{ss}}(G)$ denote the corresponding semisimple L-parameter. Assume Assumption \ref{assumption: perversetexactness} holds for $\phi_{\mf{m}}$. We then have that 
\[ \IC(\gxno,\xi)_{K^{p}K_{p}^{\mathrm{hs}},\Lambda,\mf{m}} \]
is concentrated in degree $0$. In particular, if $\Lambda = \ol{\bb{F}}_{\ell}$, $\ell$ is banal with respect to $G$, $p > 2$, $\mf{m}$ is generic in the sense of \cite[Definition~1.1]{HamannLeeTorsion}, and $G$ is a product of groups of the form in \cite[Table~(1)]{HamannLeeTorsion} and \cite[Theorem~8.2.1]{peng2025farguesscholzeparameterstorsionvanishing} then by \textit{loc.cit} the complex $\IC(\gxno,\xi)_{K^{p}K_{p}^{\mathrm{hs}},\Lambda,\mf{m}}$ is concentrated in degree $0$.
\end{corollary}
\begin{proof}
We look at the decomposition  
\[ R\Gamma(K_{p}^{\hs},\IC(\gxno,\xi)_{K^{p},\Lambda}) \simeq \bigoplus_{\phi \in \Phi^{\mathrm{ss}}(G)} R\Gamma(K_{p}^{\hs},\IC(\gxno,\xi)_{K^{p},\Lambda,\phi}) \]
induced by the direct sum decomposition (\ref{eqn: directsumdecompositionoverparameters}).

It follows by \cite[Lemma~4.2 (1)]{HamannLeeTorsion} and the fact that the $L$-parameter given by usual Satake and the Fargues--Scholze correspondence agree for unramified irreducible representations that $R\Gamma(K_{p}^{\hs},\IC(\gxno,\xi)_{K^{p},\Lambda,\phi})_{\mf{m}} = 0$ unless $\phi = \phi_{\mf{m}}$ for $\phi_{\mf{m}}$ the semisimple $L$-parameter attached to a maximal ideal $\mf{m} \subset H_{K_{p}^{\hs}}$ in the spherical Hecke algebra. Moreover, again by \cite[Lemma~4.2 (3)]{HamannLeeTorsion}, for such a maximal ideal $\mf{m} \subset \mathcal{H}_{K_{p}^{\hs}}$ one has an isomorphism $R\Gamma(K_{p}^{\hs},\IC(\gxno,\xi)_{K^{p},\Lambda,\phi_{\mf{m}}}) \simeq R\Gamma(K_{p}^{\hs},\IC(\gxno,\xi)_{K^{p},\Lambda})_{\mf{m}}$ of $\mathcal{H}_{K_{p}^{\hs}}$-modules, where the RHS is the usual localization under the spherical Hecke algebra. Now the claim follows immediately from Lemma \ref{lemma: HochschildSerreforIntCohomology} and Lemma \ref{lemma: descentlemma}, using the assumption on $K_{p}$. 
\end{proof}
\subsubsection{Relationship to the work of Koshikawa--Shin}{\label{ss: comparisonwithKoshikawaShin}}
For $K = K^{p}K_{p} \subset \mathbf{G}(\bb{A}_{f})$ a sufficiently small open compact subgroup, we consider $R\Gamma_{(2)}(\Sh_{K}\gx(\bb{C})^{\mathrm{an}},\bb{C})$, the $L^{2}$-cohomology with coefficients in $\bb{C}$ of the complex analytic variety $\Sh_{K}\gx(\bb{C})^{\mathrm{an}}$ attached to the Shimura variety of level $K$ together with its structure as a real manifold. We then pass to infinite level at $p$
\begin{equation}{\label{eqn: L2Cohomology}} 
R\Gamma_{(2)}(\gxno)_{K^{p}} := \varinjlim_{K_{p} \ra \{1\}} R\Gamma_{(2)}(\Sh_{K}\gx(\bb{C})^{\mathrm{an}},\bb{C}). 
\end{equation}
It is well known that tempered representations contribute only to the middle degree of this cohomology complex. The recent work of Koshikawa--Shin \cite{KoshikawaShin} studies the question how representations whose Arthur parameters have nontrivial Arthur $\mathrm{SL}_2$ could spread out in different degrees of cohomology. This is connected to the intersection cohomology via results of Looijenga~\cite{Looijenga88} and Saper--Stern~\cite{SaperStern}. Namely, they prove, after fixing an isomorphism $i: \bb{C} \simeq \ol{\mathbb{Q}}_{\ell}$, that one has an identification
\begin{equation}{\label{eqn: IntCohtoL2Coh}}
 \IC(\mathsf{G},X)_{K^{p},\ol{\mathbb{Q}}_{\ell}} \simeq^{\iota}  R\Gamma_{(2)}(\gxno)_{K^{p}}[d]. 
\end{equation}
In light of this, one may wonder how our vanishing results (Theorem \ref{thm: axiomatizedvanishingstatement}) compare with the conjecture of Koshikawa--Shin (\cite[Conjecture~1.2.1]{KoshikawaShin}) on $L^{2}$-cohomology. 

To understand this, we recall that we expect the optimal condition under which Assumption~\ref{assumption: perversetexactness} is expected to hold is when $\phi$ is \emph{of weakly Langlands--Shahidi type} in the sense of \cite{HamannLeeTorsion} (see \cite[Conjecture~6.4]{HamannLeeTorsion}). More precisely, assume that $\phi$ is induced (up to $\widehat{G}$-conjugation) from a supercuspidal parameter $\phi_{M}: W_{\bb{Q}_{p}} \ra \phantom{}^{L}M(\ol{\bb{Q}}_{\ell})$ for a Levi subgroup $M \subset G$, under the natural (up to $\widehat{G}$-conjugacy) embedding $\phantom{}^{L}M(\ol{\bb{Q}}_{\ell}) \ra \phantom{}^{L}G(\ol{\bb{Q}}_{\ell})$. The condition that $\phi$ is of weakly Langlands--Shahidi type says that 
\[ H^{2}(R\Gamma(W_{\bb{Q}_{p}},r_{\mathrm{ad},P} \circ \phi_{M})) \]
vanishes for all proper parabolics $P = MU$ with Levi factor $M$, where $r_{\mathrm{ad},P}$ is the natural representation of $\phantom{}^{L}M(\ol{\bb{Q}}_{\ell})$ induced by the action of $\widehat{M}$ on the Lie algebra of $\widehat{U}$ by the adjoint action.

On the other hand, \cite[Conjecture~1.2.1]{KoshikawaShin} bounds the degrees of cohomology in which a smooth representation $\pi_{p}$ of $G(\mathbb{Q}_{p})$ can contribute to $\IC(\gxno)_{K^{p}}$. This bound (see the discussion before \cite[Conjecture~1.2.1]{KoshikawaShin}) is computed in terms of the size of the Weil--Deligne $\SL_{2}$ of the $L$-parameter of $\hat{\pi}_{p}$, the Aubert dual of $\pi_{p}$, under a local Langlands correspondence satisfying certain desiderata.\footnote{In this discussion, we do not distinguish the $\ell$-adic and Weil-Deligne forms of $L$-parameters, see \cite[Remark 7.1.7]{DvHKZ2}.} By the closure ordering conjecture \cite[Section~3.1]{KoshikawaShin}, this should be an optimal bound (in the sense of the natural ordering of nilpotent orbits) for the size of the Arthur $\SL_{2}$ of any local $A$-parameter $\psi$ such that $\pi_{p} \in \Pi_{\psi}(G)$ lies in a local A-packet for $G$. In particular, if the Arthur $\SL_{2}$ of all such $\psi$ is trivial, then the conjecture of Koshikawa--Shin predicts that $\pi_{p}$ can only contribute to middle degree.

Now if $\pi_{p}$ contributes to the direct summand $\IC(\gxno)_{K^{p},\ol{\mathbb{Q}}_{\ell},\phi}$ then, as discussed above, this forces that $\phi^{\mathrm{FS}}_{\pi_{p}} = \phi$. If $\phi$ is of weakly Langlands--Shahidi type, then, under the belief that the Fargues--Scholze correspondence is the semi-simplification of any reasonable local Langlands correspondence with coefficients in $\bb{C}$ compared under the isomorphism $i$, we claim that it is impossible for $\pi_{p}$ to lie in a local $A$-packet $\Pi_{\psi}(G)$ attached to a local $A$-parameter $\psi$ with non-trivial Arthur $\SL_{2}$. Indeed, suppose such a $\psi$ existed, then, by semi-simplifying the Weil--Deligne $\SL_{2}$ using the natural map 
\[W_E\hookrightarrow W_E\times\mathrm{SL}_2(\bb{C}): g \mapsto \begin{pmatrix} |g|^{\frac{1}{2}} & 0 \\ 0 & |g|^{-\frac{1}{2}}  \end{pmatrix},\] 
we would obtain a non-Frobenius semi-simple deformation of the full semi-simplification of $\psi$, which we would expect to be precisely $\phi^{\mathrm{FS}}_{\pi_{p}}$ up to $\widehat{G}$-conjugacy (i.e that the Fargues-Scholze correspondence is compatible with our assignment of A-parameters under the fixed choice of isomorphism $i$, where we normalize the Fargues-Scholze correspondence with respect to the choice of square root $i(\sqrt{p})$). In other words, there would be a non-trivial lift
\[ 
\begin{tikzcd}
& \phantom{}^{L}P(\ol{\bb{Q}}_{\ell}) \arrow[d] & \\
W_{\bb{Q}_{p}} \arrow[r] \arrow[ur,dotted,"\exists"] & \phantom{}^{L}M(\ol{\bb{Q}}_{\ell}),  &
\end{tikzcd}
\]
which is non-Frobenius semisimple. However, this would give rise to a class in $H^{2}(W_{\bb{Q}_{p}},r_{\mathrm{ad},P} \circ \phi_{M})$ (see \cite[Remark~6.5]{HamannLeeTorsion}), which would contradict the assumption that $\phi$ is of weakly Langlands-Shahidi type. 

Analogously, one can connect our results with the conjectures of Koshikawa--Shin with torsion coefficients (\cite[Conjecture~6.2.1 (iii)]{KoshikawaShin}), by using that the Langlands--Shahidi condition mod $\ell$ also prevents the existence of non Frobenius-semisimple deformations for any lift of the parameter to characteristic $0$ (see \cite[Remark~5.1.1.]{KoshikawaShin}). 

\subsection{Eichler--Shimura Relation}\label{sec:Eichler-Shimura}

Following the strategy of \cite{XiaoZhu}, \cite{KoshikawaES}, \cite[Section 9]{DvHKZ} and \cite{vdHove}, we use Theorem~\ref{thm: Hecke operator} to establish the Eichler--Shimura congruence relation on the intersection cohomology of Shimura varieties at Iwahoric levels. 

Let $\Lambda$ over $\mathbb{Z}_\ell[\sqrt{p}]$ be a coefficient ring as in Setup~\ref{assumption: coefficients} such that $\ell \nmid \pi_{0}(Z(G))$ if $\Lambda$ is integral or torsion. Let $\gx$ be a Shimura datum as in Assumption~\ref{assumption:codimension}, and we have $\mu$, $\mathsf{E}$, $E$ as before. Let $q$ be the cardinality of the residue field of $E$. Fix a hyperspecial subgroup $K_p^\mathrm{hs}\subset G(\qp)$, as well as an Iwahori subgroup $I_p\subset K_p^\mathrm{hs}$. Recall that $\mathcal{H}_{K_p^\mathrm{hs}}$ and $\mathcal{H}_{I_p}$ denote the corresponding Hecke algebras with $\Lambda$-coefficients and let $\mathcal{Z}_{I_p}$ denote the center of the Iwahori Hecke algebra $\mathcal{H}_{I_p}$. It can be identified with the spherical Hecke algebra $\mc{H}_{K_p^\mathrm{hs}}$ via the Bernstein isomorphism, see \cite[Theorem 1.2]{Vigneras} for the statement in this generality. For $K^p \subset G(\mathbb{A}_f^p)$ a neat compact open subgroup, we set $K:=I_pK^p$. Let $\xi$ be an algebraic representation of $\mathsf{G}_{\ol{\mathbb{Q}}_{\ell}}$ defined over $\Lambda$. The intersection cohomology complex $\mathrm{IC}(\mathsf{G},\mathsf{X},\xi)_{K,\Lambda}$ has an action of $\mathcal{Z}_{I_p}\times W_E$. 

Consider the stack of $L$-parameters
\[\mathfrak{X}_{\widehat{G}} :=[Z^1(W_{\qp},\widehat{G})_\Lambda/\widehat{G}]\] 
as constructed in \cite[Section VIII]{FSGeomLLC}. This is an algebraic stack over $\Lambda$. We write $\phi_E^\mathrm{univ}$ for the universal $L$-parameter on $Z^1(W_{\qp},\widehat{G})$, restricted to $W_E\subset W_{\qp}$. Let $r_\mu$ be the representation of $\widehat{G}\rtimes W_E$, whose restriction to $\widehat{G}$ is the highest weight representation\footnote{We note that this is well-defined and agrees with the tilting module $\mathcal{T}_{\mu}$ with $\Lambda$-coefficients because $\mu$ is minuscule.} labeled by $\mu$, with $W_E$ acting trivially on the highest weight vector. It gives rise to a vector bundle $\mathcal{V}_\mu$ on $\ParG$, whose spectral action (see \cite[Section IX.2, X]{FSGeomLLC}) on $\mc{D}_{\mathrm{lis}}(\Bun_G,\Lambda)$ agrees with the Hecke operator $T_\mu$. The vector bundle $\mc{V}_\mu$ carries an action by $W_E$, because the pullback of $\mc{V}_\mu$ to $Z^1(W_{\qp},\widehat{G})$ is the trivial bundle $V_\mu\times Z^1(W_{\qp},\widehat{G})$. It carries a $W_E$-action through $r_\mu\circ \phi_E^\mathrm{univ}$ and this $W_E$-action commutes with the diagonal action of $\widehat{G}$ via 
\[(v,\phi)\mapsto (r_\mu(g)v, g\phi g^{-1}).\] 
Hence, it descends to an action on $\mc{V}_\mu$.

Now let $\mathcal{Z}^\mathrm{spec}(G):=\mathcal{O}_{\mathfrak{X}_{\widehat{G}}}(\mathfrak{X}_{\widehat{G}})$ be the spectral Bernstein center, as before. Choose a lift $\sigma_E\in W_E$ of the (arithmetic) $q$-Frobenius.  As in \cite[Definition 9.3.2]{DvHKZ}, we can define the spectral Hecke polynomial attached to $\sigma_E$ as the characteristic polynomial of its action on $\mc{V}_\mu$
\begin{align}
 H_{\mu,\sigma_E}^\mathrm{spec}(X):= \operatorname{det}\left(X-r_{\mu} \circ  \phi_E^{\operatorname{univ}}(\sigma_E)\right) \in \mathcal{Z}^\mathrm{spec}(G)[X].
\end{align}
Theorem~\ref{thm: Hecke operator} implies that 
there is a $\Lambda$-algebra homomorphism
\[\mathrm{End}_{\ParG}(\mathcal{V}_\mu)\to \mathrm{End}_{D(G(\qp),\Lambda)}(\operatorname{IC}(\mathsf{G},\mathsf{X},\xi)_{K^p,\Lambda}),\]
cf. \cite[Theorem 9.1.4]{DvHKZ}, which restricts to a $\Lambda$-algebra map
\[\mathcal{Z}^\mathrm{spec}(G)[X]/(H_{\mu,\sigma_E}^\mathrm{spec}(X))\to \mathrm{End}_{D(G(\qp),\Lambda)}(\operatorname{IC}(\mathsf{G},\mathsf{X},\xi)_{K^p,\Lambda}),\]
sending $X$ to $\sigma_E$. The latter acts on $\operatorname{IC}(\mathsf{G},\mathsf{X},\xi)_{K^p,\Lambda}$ through its action on $\mc{V}_\mu$ and the formula
\[i_1^\ast\mc{V}_\mu\ast (\mc{F}_{\xi, \IC, \Lambda})=i_{1}^{*}T_{\mu}\mc{F}_{\xi, \IC, \Lambda} \simeq \IC(\mathsf{G},\mathsf{X}, \xi)_{K^p,\Lambda}.\]
from Theorem~\ref{thm: Hecke operator} (5). But the spectral Weil group action agrees with the usual geometric action coming from the structure map of the Shimura variety to $\Spa E$, see \cite[Theorem 9.1.4]{DvHKZ}. Hence we deduce that $H_{\mu,\sigma_E}^\mathrm{spec}(X)$ annihilates $\sigma_E$ for its usual geometric action on $\IC(\mathsf{G},\mathsf{X}, \xi)_{K^p,\Lambda}$.

In addition to $\ell \nmid \pi_{0}(Z(G))$, suppose that $I_{p}$ has pro-order coprime to $\ell$ if $\Lambda$ is integral or torsion (e.g $\ell \nmid q - 1$ if $G$ is split). In that case, we can pass to Iwahori level as in~\cite[Theorem 5.4]{vdHove}, using the Hochschild--Serre spectral sequence for intersection cohomology in the case of torsion or integral coefficients, cf. Lemma~\ref{lemma: HochschildSerreforIntCohomology} and Lemma~\ref{lemma: descentlemma}.
We deduce that the following diagram is commutative
\[
\begin{tikzcd}
    \mathcal{Z}^\mathrm{spec}(G)[X]/(H_{\mu,\sigma_E}^\mathrm{spec}(X))\ar[r]\ar[d]& \mathrm{End}_{D(G(\qp),\Lambda)}(\operatorname{IC}(\mathsf{G},\mathsf{X},\xi)_{K^p,\Lambda})\ar[d]\\
    \mathcal{Z}_{I_p}[X]/(H_{\mu,\sigma_E}^\mathrm{spec}(X))\ar[r]& \mathrm{End}_{D(\Lambda)}(\operatorname{IC}(\mathsf{G},\mathsf{X},\xi)_{K,\Lambda}),
\end{tikzcd}
\]
where the left vertical map exists through \cite[Proposition 4.4]{vdHove}, cf. \cite[Lemma 9.4.4, Proposition 9.4.8]{DvHKZ} and all other maps are the natural ones. Moreover, under the Bernstein isomorphism $\mathcal{Z}_{I_p}\simeq\mathcal{H}_{K_p^\mathrm{hs}}$, the image of $H_{\mu,\sigma_E}^\mathrm{spec}(X)$ in $\mathcal{Z}_{I_p}[X]$ agrees with the Hecke polynomial \[H_\mu(X)=\mathrm{det}(X-r_\mu(g,\sigma_E))\in \mathcal{H}_{K_p^\mathrm{hs}}[X],\] 
which no longer depends on the choice of the lift $\sigma_E$ of the $q$-Frobenius. So all in all, we have reached the following theorem:

\begin{theorem}\label{thm: ES} Suppose $\Lambda$ is as in Setup~\ref{assumption: coefficients}, and that $I_p$ has pro-order coprime to $\ell$ and $\ell \nmid |\pi_0(Z(G))|$ if $\Lambda$ is integral or torsion (e.g $\ell \nmid q - 1$ if $G$ is split). 
For $K=I_pK^p$, the inertia subgroup $I_E \subset W_E$ acts unipotently on $\operatorname{IC}(\mathsf{G},\mathsf{X},\xi)_{K,\Lambda}$. Moreover, the action of any lift of the (arithmetic) $q$-Frobenius $\sigma_E \in W_E$ on $\operatorname{IC}(\mathsf{G},\mathsf{X},\xi)_{K,\Lambda}$ satisfies $H_{\mu}(\sigma_E)=0$.
\end{theorem}

\subsection{A conjectural connection between Arthur's conjectures and categorical local Langlands}{\label{ss: EigensheafConjecture}
In view of the existing results on the compactly supported and usual cohomology of Shimura varieties, the applications we have discussed so far are not completely surprising once we established our main result, Theorem~\ref{thm: main}. However, we envision further applications of our construction of the perverse sheaf $\mathcal{F}_{\mathrm{IC},\Lambda}$. 

One of the most interesting facts about the intersection cohomology of Shimura varieties, unlike compactly supported or usual cohomology, is that it admits a very clean and simple description without passing to any localization by a Hecke algebra, via the relationship to $L^{2}$-cohomology given in (\ref{eqn: IntCohtoL2Coh}).  In particular, for $K^{p} \subset \mathsf{G}(\bb{A}_{f}^{p})$ a neat compact open subgroup, we recall, by the work of Borel--Casselman \cite{BorelCasselman}, that we have the following $\mathcal{H}_{K^{p}} \times G(\mathbb{Q}_{p})$-equivariant isomorphism 
\begin{equation}{\label{eqn: BorelCasselman}}
 R\Gamma_{(2)}(\gxno)_{K^{p}} \simeq \bigoplus_{\pi} m(\pi) \otimes \pi_{p} \otimes (\pi^{p})^{K^{p}} \otimes R\Gamma((\mathfrak{g},K_{\infty}),\pi_{\infty}), 
\end{equation}
where 
\begin{enumerate}
\item $\pi := \pi_{\infty} \otimes \pi_{p} \otimes \pi^{p}$ runs over discrete automorphic representations of $\mathsf{G}(\mathbb{A})$, with multiplicity $m(\pi)$ in the discrete automorphic spectrum; we may restrict to those automorphic representations that are cohomological (regular algebraic) of trivial weight at $\infty$.
\item $R\Gamma((\mathfrak{g},K_{\infty}),\pi_{\infty})$ is the complex of $\mathbb{C}$-vector spaces attached to the  $(\mathfrak{g},K_{\infty})$-cohomology of the archimedean component $\pi_{\infty}$ of $\pi$, where $K_{\infty}$ is the centralizer in $\mathbf{G}(\bb{R})$ of an element of $X$.
\item $\mathcal{H}_{K^{p}}$ denotes the away from $p$ Hecke algebra with coefficients in $\mathbb{C}$; it acts on the space of invariants $(\pi^p)^{K^p}$. 
\end{enumerate}
Now, fixing an isomorphism  $\iota: \ol{\mathbb{Q}}_{\ell} \simeq \mathbb{C}$ as before, by (\ref{eqn: IntCohtoL2Coh}) and Theorem \ref{thm: Hecke operator}, we obtain a chain of $\mathcal{H}_{K^p} \times G(\bb{Q}_{p})$-equivariant isomorphisms
\[ 
 i_{1}^{*}T_{\mu}\mathcal{F}_{\mathrm{IC},\ol{\mathbb{Q}}_{\ell}} \simeq \IC(\mathsf{G},X)_{K^{p},\ol{\mathbb{Q}}_{\ell}} \simeq^{\iota}  R\Gamma_{(2)}(\gxno)_{K^{p}}[d], 
\]
where $\mathcal{F}_{\mathrm{IC},\ol{\mathbb{Q}}_{\ell}}$ denotes the perverse ULA sheaf in $\mathrm{D}_{\mathrm{lis}}(\Bun_{G},\ol{\mathbb{Q}}_{\ell})$ constructed in \ref{cons: ConstructionofSheaveswithQellCoefficients}. By combining with (\ref{eqn: BorelCasselman}), this gives us an identification
\begin{equation}{\label{eqn: GKCohomologyFormulaHeckeOperator}}
  i_{1}^{*}T_{\mu}\mathcal{F}_{\mathrm{IC},\ol{\mathbb{Q}}_{\ell}}\simeq^{\iota} \bigoplus_{\pi} m(\pi) \otimes \pi_{p} \otimes (\pi^{p})^{K^{p}} \otimes R\Gamma((\mathfrak{g},K_{\infty}),\pi_{\infty})[d] 
\end{equation}
under the fixed isomorphism $\iota$. If we compare both sides, this leads to the following natural question.
\begin{question}
Does the sheaf $\mathcal{F}_{\mathrm{IC},\ol{\mathbb{Q}}_{\ell}}$ decompose as a direct sum of sheaves in $\mathrm{D}_{\mathrm{lis}}(\Bun_{G},\ol{\mathbb{Q}}_{\ell})$ in such a way that, after applying $i_{1}^{*}T_{\mu}$, one recovers the direct sum decomposition appearing in (\ref{eqn: GKCohomologyFormulaHeckeOperator})?
\end{question}

We strongly suspect that the answer should be indeed yes. To explain further, let us write $A_{\pi} := R\Gamma((\mathfrak{g},K_{\infty}),\pi_{\infty}) \in \mathrm{D}(\mathbb{C})$ for the complex of $\mathbb{C}$-vector spaces attached to the $(\mathfrak{g},K_{\infty})$-cohomology of $\pi_{\infty}$. Recall that the cohomological amplitude of $A_{\pi}$ admits an alternative description in terms of the Arthur $\SL_{2}$ attached to a global representation $\pi$ equipped with an embedding into $L^{2}_{\disc}$, an explicit algebraic representation, denoted $\psi_{\pi}|_{\SL_{2}(\mathbb{C})}: \SL_{2}(\mathbb{C}) \ra \widehat{G}(\mathbb{C})$ attached to the global representation $\pi$ in which $\pi_{\infty}$ occurs under Arthur's conjectures (see \cite[Section~2]{KoshikawaShin} for a discussion of Arthur parameters and the Arthur $\mathrm{SL}_2$). Writing $r_{\mu}: \widehat{G} \ra \GL(V_{\mu})$ for the algebraic representation attached to the Hodge cocharacter of the Shimura datum, we may decompose $V_{\mu} := \bigoplus_{n \in \mathbb{Z}} V_{\mu}^{\psi_{\pi},n}$ in terms of the weight spaces of the diagonal $\mathbb{G}_{m} \hookrightarrow \SL_{2}(\mathbb{C})$ acting on $V_{\mu}$ via $\psi_{\pi}|_{\SL_{2}(\mathbb{C})}$. Now one can show, using the theory of Adams--Johnson packets \cite{AdamsJohnson} and the Vogan--Zuckerman classification of cohomological representations \cite{VoganZuckerman}, that the degrees of cohomology of $R\Gamma((\mathfrak{g},K_{\infty}),\pi_{\infty})$ will be related to the degrees of cohomology of $\bigoplus_{n \in \mathbb{Z}} V_{\mu}^{\psi_{\pi},n}[-d+n]$ (see \cite[Proposition~6.1.4]{KoshikawaShin} and \cite[Proposition~9.1]{ArthurUnipotentConjectures} for a more precise discussion), where as before $d$ is the dimension of the Shimura variety. 

This kind of ``shearing'' by the weights of the Arthur $\SL_{2}$ directly mimics the behavior of certain purely local generalized eigensheaves that are being constructed in forthcoming work of  Bertoloni-Meli and Koshikawa~\cite{CuspidalSheavesBertiKoshikawa} which is strongly inspired by the $p$-adic ABV theory described in \cite{padicABVtheory}. In particular, for each local Arthur parameter $\psi_{p}$ at $p$, assuming categorical local Langlands and the compatibility with geometric Eisenstein series, they construct a perverse ULA Hecke eigensheaf $\mathcal{F}_{\psi_{p}} \in \mathrm{D}_{\mathrm{lis}}(\Bun_{G},\ol{\mathbb{Q}}_{\ell})$, whose stalks at the neutral components are the local $A$-packets (Strictly speaking, these will be the ABV packets and it is of separate interest to show that these ABV packets agree with the usual Arthur packets (See \cite[Section~8.3,Conjecture~1]{padicABVtheory} and \cite{cunningham2022proofvogansconjecturearthur})). Moreover, this is a ``sheared Hecke eigensheaf'' in the sense that one has an isomorphism 
\[ T_{\mu}(\mathcal{F}_{\psi_{p}}) \simeq \mathcal{F}_{\psi_{p}} \otimes \iota^{-1}(\bigoplus_{n \in \mathbb{Z}} V_{\mu}^{\psi_{p},n}[n]) \]
of sheaves on $\Bun_{G}$ for all cocharacters $\mu$ after forgetting the Weil group action, where $V_{\mu}^{\psi_{p},n}[-n]$ is defined as above with respect to the Arthur $\SL_{2}$ of the local parameter $\psi_{p}$. This leads to a natural expectation that $\mathcal{F}_{\mathrm{IC},\ol{\mathbb{Q}}_{\ell}}$ should decompose in terms of certain direct summands of $\mathcal{F}_{\psi_{p}}$ (after restricting to the open substack $\Bun_{G,\mu^{-1}} \hookrightarrow \Bun_{G}$ corresponding to $B(G,\mu^{-1}) \subset B(G)$) in a way that is consistent with the behavior of the $(\mathfrak{g},K_{\infty})$-cohomology described above, as well as with Arthur's description of the discrete spectrum in terms of local A-packets. In future work, we will formulate the precise statement of this conjecture. We believe that this so-called ``eigensheaf conjecture'' is the correct interpretation of Fargues' local-global compatibility conjecture, as discussed in \cite[Section~7]{Fargues-Overview}.

\begin{remark}\label{rem: function field version}
It would be very interesting if, in the context of function fields and moduli of global shtukas, one could give a construction of the analog of the sheaf $\mathcal{F}_{\IC,\ol{\bb{Q}}_{\ell}}$ over all of $\Bun_{G}$, cf. \cite[Theorem~E]{lihuerta2026courbesetfibresvectoriels}. In particular, this analog should have a very similar direct sum decomposition. The spectral counterpart of this sheaf seems to be intimately related with recent conjectures in unramified geometric Langlands, describing the discrete $L^{2}$-spectrum in terms of the global unramified geometric Langlands conjecture (see \cite[Conjecture~3.4.1]{RaskinAnnouncement}).
\end{remark}

\appendix
\section{Intersection complex with adic coefficients as an inverse limit}
In this appendix, we record the following result: in the algebraic setting, the intersection complex with $\ell$-adically complete coefficients is the inverse limit of the ones with torsion coefficients. The proof is due to O. Gabber. We thank him for communicating the ideas to us. Any error below is due to us.

Let $X$ be a separated, finite type scheme over $k$ (or the perfection of such a scheme when $k$ is in characteristic $p$), and $\Lambda$ be the ring of integers in a finite extension of $\bb{Q}_\ell$. We consider the usual perverse t-structure on $\mc{D}^b_{\mathrm{c}}(X, \Lambda)$, the $\infty$-category of bounded constructible sheaves with $\Lambda$-coefficients on $X$, and denote it by $(\phantom{}^{p}\mc{D}_{\mathrm{c}}^{\leq 0}(X,\Lambda), \phantom{}^{p}\mc{D}_{\mathrm{c}}^{\geq 0}(X,\Lambda))$ (or $(\phantom{}^{p}\mc{D}_{\mathrm{c}}^{\leq 0}, \phantom{}^{p}\mc{D}_{\mathrm{c}}^{\geq 0})$ for simplicity), where p is the middle perversity. For a morphism of perverse sheaves $f: K_1\to K_2$ in the heart of the $t$-structure ${}^p\mc{D}^{=0}_c(X,\Lambda)$, we write $\operatorname{ker}(f)$, $\operatorname{coker}(f)$ for the kernel and cokernel taken in the abelian category of perverse sheaves. Similarly, we write $K/\ell$ for ${}^p\mc{H}^0(K\otimes^{\bb{L}}\Lambda/\ell)$ the cokernel of multiplication by $\ell$ and so on.

Let $j: U\hookrightarrow X$ be the inclusion of a dense open subscheme of pure dimension $d_U$. For each perverse sheaf $K$ on $U$, we can form the intermediate extension ${}^pj_{!\ast}K$, as in Definition \ref{defn: intermediateextension}. Note that when $K$ is torsion-free, each $K/\ell^n$ is again a perverse sheaf, by considering the long exact sequence of perverse cohomology induced by the short exact sequences
\[0\to K\xrightarrow{\ell^n}K\to K/\ell^n\to 0.\]

\begin{theorem}{\label{thm: inverselimit}}
    Let $j: U\hookrightarrow X$ be as above. Assume $K \in \mc{D}_{c}^{b}(U,\Lambda)$ is a torsion-free perverse sheaf with $\Lambda$-coefficients. We have an isomorphism
    \[{}^pj_{!\ast}K\simeq \varprojlim_n{}^pj_{!\ast}(K/\ell^n).\]
\end{theorem}
To prepare for the proof, let us first recall from \cite[\S~3.3]{BBD}, \cite[Section 2]{JuteauDecompositionNumbers} that Verdier duality on $X$ (denoted by $\bb{D}_X$) does not preserve the p-perverse t-structure, since $\Lambda$ is not self-injective. Rather, the Verdier duality interchanges the p-perverse t-structure with the p$^+$-perverse t-structure $(\phantom{}^{p^+}\mc{D}_{\mathrm{c}}^{\leq 0}, \phantom{}^{p^+}\mc{D}_{\mathrm{c}}^{\geq 0})$, which can be described explicitly as follows:
\begin{itemize}
    \item $A\in \phantom{}^{p^+}\mc{D}_{\mathrm{c }}^{\leq 0}(X,\Lambda)\iff$ $A \in \phantom{}^p\mc{D}_{\mathrm{c }}^{\leq 1}(X,\Lambda)$ and $\phantom{}^p\mc{H}^1(A)$ is torsion.
    \item $A\in \phantom{}^{p^+}\mc{D}_{\mathrm{c }}^{\geq 0}(X,\Lambda)\iff$ $A \in \phantom{}^p\mc{D}_{\mathrm{c }}^{\geq 0}(X,\Lambda)$ and $\phantom{}^p\mc{H}^0(A)$ is torsion-free.
\end{itemize}
One can form intermediate extensions of p$^+$-perverse sheaves along a dense open immersion $j:U\hookrightarrow X$ as well. Note that when $K\in \phantom{}^p\mc{D}^{=0}(U,\Lambda)$ is torsion-free (as a p-perverse sheaf, meaning that it has no torsion subsheaf that is also p-perverse), then it is also automatically a p$^+$-perverse sheaf.

\begin{lemma}\label{lemma: pTorsionfree}
    Let $j:U\hookrightarrow X$ be a dense open immersion and $K$ be a torsion-free p-perverse sheaf on $U$. Then ${}^pj_{!\ast}K$ is a torsion-free p-perverse sheaf on $X$.
\end{lemma}
\begin{proof}
    Assume there is a nonzero map $\mc{T}\to {}^pj_{!\ast} K$ from a torsion p-perverse sheaf $\mc{T}$. Post-composing with the natural injection ${}^pj_{!\ast} K \ra \phantom{}^{p}j_{\ast} K$ defining ${}^pj_{!\ast} K$, we get a nonzero map $\mc{T}\to {}^pj_{!\ast} K\hookrightarrow {}^pj_\ast K$. By adjunction between ${}^pj_{\ast}$ and ${}^pj^\ast=j^\ast$, this gives rise to a nonzero map $j^\ast\mc{T}\to K$. Since ${}^pj_{!\ast} K$ does not have any p-perverse subsheaf supported on the complement of $U$, $j^*\mc{T}\neq 0$. This contradicts the assumption that $K$ is torsion-free.
\end{proof}
\begin{corollary}\label{cor: pIC+perverse}
    For $j$ and $K$ as above, ${}^pj_{!\ast}K$ is a perverse sheaf for the p$^+$-perversity.
\end{corollary}
\begin{proof}
    This follows from the definition of the p$^+$-t-structure and the torsion-freeness of ${}^pj_{!\ast}K$.
\end{proof}
We now have the following Verdier dual statement.
\begin{corollary}\label{cor: p+ICperverse}
    For $j$ and $K$ as above, ${}^{p^+}j_{!\ast} K$ is a perverse sheaf for the $p$-perversity.   
\end{corollary}
\begin{proof}
    Since ${}^{p^+}j_{!\ast} K$ is p$^+$-perverse, it lies in ${}^p\mc{D}^{[0,1]}(X,\Lambda)$ and we have a distinguished triangle
    \[^p\mc{H}^0({}^{p^+}j_{!\ast} K)[0]\to {}^{p^+}j_{!\ast} K \to {}^p\mc{H}^1({}^{p^+}j_{!\ast} K)[-1].\]
    It therefore suffices to show that the torsion sheaf ${}^p\mc{H}^1({}^{p^+}j_{!\ast} K)$ vanishes. Suppose it is killed by $\ell^m$ for some $m\in \bb{N}$. We need to show that multiplication by $\ell^m$ is epimorphic on ${}^{p^+}j_{!\ast} K$ as a p$^+$-perverse sheaf. This follows from the dual statement. Namely,
    \[\bb{D}_X({}^{p^+}j_{!\ast} K)\simeq {}^{p}j_{!\ast} \bb{D}_U(K),\]
    and it suffices to show that ${}^{p}j_{!\ast} \bb{D}_U(K)$ is torsion-free as a p-perverse sheaf. But since $K\in {}^p\mc{D}^{=0}\cap {}^{p^+}\mc{D}^{=0}$, so is $\bb{D}_UK\in {}^p\mc{D}^{=0}\cap {}^{p^+}\mc{D}^{=0}$. This means $\bb{D}_UK={}^p\mc{H}^0(\bb{D}_UK)$ is torsion-free by the second condition describing the p$^+$-t-structure, and in particular torsion-free as a p-perverse sheaf. Now the desired statement follows by applying Lemma~\ref{lemma: pTorsionfree} to $\bb{D}_U(K)$. 
\end{proof}

\begin{proposition}\label{prop: InjectiveAlpha}
    Let $j:U\hookrightarrow X$ be a dense open immersion. Assume $U$ is of pure dimension $d_U$, and $K$ is a p-perverse sheaf on $U$ that is torsion-free (as a p-perverse sheaf). Then there is a unique, injective map of p-perverse sheaves on $X$
    \begin{equation}{\label{eqn: alphamap}}
    \alpha:{}^pj_{!\ast} K\to {}^{p+}j_{!\ast}K,
    \end{equation}
    extending the identity on $K$.
\end{proposition}
\begin{proof}
    Under the assumption that $K$ is torsion-free, the p-perversity of ${}^{p+}j_{!\ast}K$ is established in Corollary~\ref{cor: p+ICperverse}. Pick a stratification of $X$ for which $K$ is constructible and assume that $U$ is a union of strata. We write $X$ as an increasing union
    \[X=U_n\supset U_{n-1}\supset \cdots \supset U_1\supset U_0,\]
    with inclusion maps $j_i: U_{i-1}\hookrightarrow U_i$. We assume that for each $i\geq 1$, $U_i\backslash U_{i-1}$ is a union of strata of pure codimension $i$. By definition of being constructible, we know that $K$ is (up to shift) a local system when restricted to $U_0$. We can moreover assume that $K|_{U_0}$ is nonzero, since otherwise one can replace $X$ by the closure of the support of $K$, and $U$ by its intersection with the support of $K$.
    
    Then we have Deligne's formula for the intermediate extensions in both the p- and p$^+$-perverse t-structures, cf. Corollary~\ref{cor: DeligneFormula}, \cite[\S 2.3]{JMW}. Namely, we have
    \begin{equation*}
    {}^pj_{!\ast}K\simeq \tau_{\leq -d_U+n-1}j_{n,\ast}\cdots \tau_{\leq -d_U+1} j_{2,\ast}\tau_{\leq -d_U} j_{1,\ast}K
    \end{equation*}
    \begin{equation*}
    {}^{p^+}j_{!\ast}K\simeq \tau^+_{\leq -d_U+n-1}j_{n,\ast}\cdots \tau^+_{\leq -d_U+1} j_{2,\ast}\tau^+_{\leq -d_U} j_{1,\ast}K,
    \end{equation*}
    where $\tau^+$ denotes the truncation functor for the linear dual of the standard t-structure, which we denote by $({}^+\mc{D}^{\leq 0}, {}^+\mc{D}^{\geq 0})$. Explicitly, for $M\in \mc{D}^b_{\mathrm{c}}(U_i,\Lambda)$ for some $i$, $\tau^+_{\leq 0}M\to M$ is an isomorphism in degrees $\leq 0$ and in degree 1, it is an isomorphism onto the torsion part. 
    
    Since, by definition, we have an embedding $\iota: \mc{D}^{\leq 0}\subset {}^{+}\mc{D}^{\leq 0}$, there is a natural transformation
    \[\tau_{\leq 0}\to \tau^+_{\leq 0}.\]
    After successively shifting and composing with the $j_{i\ast}$'s, we obtain the desired map $\alpha$. The uniqueness follows from applying adjunctions inductively.

    To show injectivity of $\alpha$, note that ${}^pj_{!\ast}K$ has no p-perverse sub- or quotient-object that is supported on  $Z:=X\backslash U$. But since $\alpha$ is a map of p-perverse sheaves and is the identity on $K$ when restricted to $U$, $\operatorname{ker}(\alpha)$ is p-perverse and is supported on $Z$. This forces it to be zero. 
\end{proof}
\begin{remark}
    We warn the reader that it does not follow from Corollary~\ref{cor: pIC+perverse} and Corollary~\ref{cor: p+ICperverse} that the map $\alpha$ is an isomorphism. In particular, if we consider the map as a map of p$^+$-perverse sheaves, then it will be an epimorphism, since ${}^{p+}j_{!\ast}K$ does not have p$^+$-perverse quotient that is supported on $Z$. However, as a map of p$^+$-perverse sheaves, $\alpha$ is no longer injective, as ${}^pj_{!\ast}K$ could have a subsheaf supported on $Z$ that is p$^+$-perverse. It genuinely matters in which abelian category we are taking the kernel or cokernel. For example, consider the open immersion $j:U\to X=\bb{A}^1$ over $\bb{C}$, where $U$ is the complement of the origin. Take $K$ to be the rank one \'etale $\bb{Z}_\ell$-local system on $U$ attached to the representation $\hat{\bb{Z}}\to \bb{Z}_\ell^\times: 1\mapsto \ell+1$, shifted by $[1]$. Then one can compute that there is a short exact sequence of p-perverse sheaves
    \[0\to {}^{p}j_{!\ast}K\xrightarrow{\alpha} {}^{p+}j_{!\ast}K \to i_\ast\bb{F}_\ell[0]\to 0, \]
    but $i_\ast\bb{F}_\ell[0]$ is not p$^+$-perverse, while $i_\ast\bb{F}_\ell[-1]$ is. Hence, when viewing $\alpha$ as a map between p$^+$-perverse sheaves, it sits in another short exact sequence
    \[0\to i_\ast\bb{F}_\ell[-1]\to {}^{p}j_{!\ast}K\xrightarrow{\alpha} {}^{p+}j_{!\ast}K\to 0.\]
\end{remark}
\begin{proposition}\label{prop: ImagePerverse}
    Taking reduction of the map $\alpha$ of (\ref{eqn: alphamap}) modulo $\ell^n$ for $n\gg 0$, the image of $({}^{p}j_{!\ast}K)/\ell^n\xrightarrow{\alpha/\ell^n} ({}^{p+}j_{!\ast}K)/\ell^n$ in the category of p-perverse sheaves is isomorphic to ${}^{p}j_{!\ast}(K/\ell^n)$.
\end{proposition}
We need the following lemma.
\begin{lemma}
    Denote the Verdier duality functor on $\mc{D}^b_{\mathrm{c}}(X,\Lambda/\ell^n)$ (resp. $\mc{D}^b_{\mathrm{c}}(X,\Lambda)$) by $\bb{D}_{X,\Lambda/\ell^n}$ (resp. $\bb{D}_{X,\Lambda}$). Let $\iota: \mc{D}^b_{\mathrm{c}}(X,\Lambda/\ell^n)\to \mc{D}^b_{\mathrm{c}}(X,\Lambda)$ be the natural forgetful functor. Then we have a natural isomorphism of functors from $\mc{D}^b_{\mathrm{c}}(X,\Lambda/\ell^n)^\mathrm{op}$ to $\mc{D}^b_{\mathrm{c}}(X,\Lambda)$ 
    \[\iota\circ \bb{D}_{X,\Lambda/\ell^n}(-)\simeq \bb{D}_{X,\Lambda}\circ \iota(-) [1].\]
\end{lemma}
\begin{proof}
    This follows from the definition and adjunction. Indeed, for any $A\in \mc{D}^b_{\mathrm{c}}(X,\Lambda/\ell^n)$, we have the following computation:
    \begin{align*}
        \bb{D}_{X,\Lambda}\circ \iota( A)[1]&\simeq R\ul{\Hom}_{\Lambda}(\iota(A),\omega_X)[1]\\
        &\simeq \iota\, R\ul{\Hom}_{\Lambda/\ell^n}(A,R\ul{\Hom}_{\Lambda}(\Lambda/\ell^n,\omega_X))[1]\\
        &\simeq \iota\, R\ul{\Hom}_{\Lambda/\ell^n}(A,\omega_X\otimes^\bb{L} \Lambda/\ell^n[-1])[1]\\
        &\simeq \iota\, R\ul{\Hom}_{\Lambda/\ell^n}(A,\omega_{X,\Lambda/\ell^n})=(\iota\circ \bb{D}_{X,\Lambda/\ell^n})(A),
    \end{align*}
    where $\omega_{X,\Lambda/\ell^n}\in \mc{D}^b_\mathrm{c}(X,\Lambda/\ell^n)$ and $\omega_{X,\Lambda}\in \mc{D}^b_\mathrm{c}(X,\Lambda)$ denote the dualizing complexes.
\end{proof}
\begin{proof}[Proof of Proposition~\ref{prop: ImagePerverse}]
    The intermediate extension ${}^{p}j_{!\ast}(K/\ell^n)$ is the unique p-perverse extension of $K/\ell^n$ to $X$ that does not have sub- or quotient objects supported on $Z=X\backslash U$. Now any quotient p-perverse sheaf of $\operatorname{Im}(\alpha/\ell^n)$ supported on $Z$ is also a quotient of ${}^pj_{!\ast}K$, so it has to be zero. 
    
    On the other hand, note that we can view $({}^{p+}j_{!\ast}K)/\ell^n$ as a perverse sheaf in $\mc{D}^b_{\mathrm{c}}(X,\Lambda/\ell^n)$. Then by the lemma above, we have (by abuse of notation, we omit the $\iota$'s)
    \begin{align*}
    \bb{D}_{X,\Lambda/\ell^n}(({}^{p+}j_{!\ast}K)/\ell^n)&\simeq {}^p\mc{H}^0(\bb{D}_{X,\Lambda}(({}^{p+}j_{!\ast}K)\otimes^\bb{L} \Lambda/\ell^n)[1])\\
    &\simeq {}^p\mc{H}^0(({}^{p}j_{!\ast}\bb{D}_{U,\Lambda}K)\otimes^\bb{L}\Lambda/\ell^n)\\
    &\simeq ({}^{p}j_{!\ast}\bb{D}_{U,\Lambda}K)/\ell^n,
    \end{align*}
    where the quotients are taken in the category of p-perverse sheaves in $\mc{D}^b_\mathrm{c}(X,\Lambda)$. Hence, if $\operatorname{Im}(\alpha/\ell^n)$ has a p-perverse subsheaf supported on $Z$, then after applying $\bb{D}_{X,\Lambda/\ell^n}$, we would obtain a p-perverse quotient of $({}^{p}j_{!\ast}\bb{D}_{U,\Lambda}K)/\ell^n$. Composing with the surjection 
    \[{}^{p}j_{!\ast}\bb{D}_{U,\Lambda}K\to ({}^{p}j_{!\ast}\bb{D}_{U,\Lambda}K)/\ell^n,\]
    we would obtain a p-perverse quotient of ${}^{p}j_{!\ast}\bb{D}_{U,\Lambda}K$ supported on $Z$, which is necessarily zero. Hence $\operatorname{Im}(\alpha/\ell^n)$ has no p-perverse subsheaf supported on $Z$.
    Combining the above, we conclude that $\operatorname{Im}(\alpha/\ell^n)\simeq {}^pj_{!\ast}(K/\ell^n)$ as desired.
\end{proof}

Now we can finish the proof of Theorem~\ref{thm: inverselimit}.
\begin{proof}[Proof of Theorem~\ref{thm: inverselimit}]
    We take the inverse limit of the inverse system of short exact sequences
    \[0\to\operatorname{ker}(\alpha/\ell^n)\to ({}^pj_{!\ast}K)/\ell^n\to {}^pj_{!\ast}(K/\ell^n)\to 0.\]
    Since $\operatorname{ker}(\alpha/\ell^n)$ is a constructible p-perverse sheaf, the cohomology groups of its stalks are finite abelian groups. Hence, the system $\{\operatorname{ker}(\alpha/\ell^n)\}_n$ satisfies the Mittag-Leffler condition. Therefore, after taking inverse limit, we obtain a short exact sequence
    \[0\to\varprojlim_n\operatorname{ker}(\alpha/\ell^n)\to \varprojlim_n({}^pj_{!\ast}K)/\ell^n\to \varprojlim_n{}^pj_{!\ast}(K/\ell^n)\to 0, \]
    but $\varprojlim_n\operatorname{ker}(\alpha/\ell^n)=\operatorname{ker}(\alpha)=0$ by Proposition~\ref{prop: InjectiveAlpha}. Also, we know ${}^pj_{!\ast}K$ is $\ell$-complete, so 
    \[{}^pj_{!\ast}K\simeq \varprojlim_n ({}^pj_{!\ast}K)/\ell^n.\] 
    Thus, the theorem follows.
\end{proof}

\printbibliography
\end{document}